\documentclass{lmcs}
\pdfoutput=1

\usepackage{lastpage}
\lmcsdoi{17}{4}{8}
\lmcsheading{}{\pageref{LastPage}}{}{}%
{Dec.~31,~2020}{Nov.~16,~2021}{}

\usepackage[utf8]{inputenc}

\usepackage{etex}

\usepackage{xcolor}  % Coloured text etc.

\usepackage{amsmath,amsthm,amssymb,amsfonts,graphicx}
\usepackage{stmaryrd}
\usepackage{cmll}
\usepackage{tikz}
\usepackage{proof}
\usepackage{mathtools}
\usepackage{xypic}
\usepackage{lscape}
\usepackage{multicol}
\usepackage{leftidx}
\usepackage{bbm}
\usepackage{rotating}

\usepackage{array}   % for \newcolumntype macro
\newcolumntype{L}{>{$}l<{$}} % math-mode version of "l" column type

\tikzstyle{strings}=[baseline={([yshift=-.5ex]current bounding box.center)}]

\usetikzlibrary{topaths}

\tikzset{every picture/.append style={scale=.5}, transform shape, strings}
\tikzset{map/.style={draw,fill=white, rectangle}}
\tikzset{simple/.style={}}
\tikzset{nothing/.style={outer sep=-3.4pt}}

\tikzset{%
symbol/.style={%
draw=none,
every to/.append style={%
edge node={node [sloped, allow upside down, auto=false]{$#1$}}}
}
}

\usetikzlibrary{shapes.geometric}
\usetikzlibrary{patterns}
\usetikzlibrary{fit}
\usetikzlibrary{positioning}
\usetikzlibrary{calc}
\usetikzlibrary{arrows}
\usetikzlibrary{decorations.markings}
\usetikzlibrary{decorations.pathreplacing}
\usetikzlibrary{shapes}

\pgfdeclarelayer{nodelayer}
\pgfdeclarelayer{edgelayer}
\pgfsetlayers{edgelayer,nodelayer,main}

\tikzset{none/.style={}}
\tikzset{none/.style={outer sep=-3.4pt}}

\tikzset{dot/.style={thick, fill=black, circle, scale=1, inner sep = .05cm}}

\tikzset{oa/.style={draw, scale=0.9,minimum height=.1cm,circle,append after command={
[shorten >=\pgflinewidth, shorten <=\pgflinewidth,]
(\tikzlastnode.north) edge (\tikzlastnode.south)
(\tikzlastnode.east) edge (\tikzlastnode.west) } } }

\tikzset{ox/.style={draw, scale=0.9,minimum height=.1cm,circle,append after command={
[shorten >=\pgflinewidth, shorten <=\pgflinewidth,]
(\tikzlastnode.north west) edge (\tikzlastnode.south east)
(\tikzlastnode.north east) edge (\tikzlastnode.south west) } } }
        
\tikzset{circ/.style={
shape=circle, inner sep=1pt, draw}}

\tikzstyle{none}=[inner sep=-1pt]
\tikzstyle{circle}=[shape=circle,draw]

\newcommand*{\StrikeThruDistance}{0.15cm}%
\newcommand*{\StrikeThru}{\StrikeThruDistance,\StrikeThruDistance}%

\tikzset{wires/.style={}}

\tikzset{box/.style={inner sep=0pt, thick, draw=black, text height=1.5ex, text depth=.25ex, text centered, minimum height=3em, anchor=center}}

\newcommand{\<}{\langle}
\renewcommand{\>}{\rangle}
\newcommand{\X}{\mathbb{X}}
\newcommand{\A}{\mathbb{A}}
\newcommand{\B}{\mathbb{B}}
\newcommand{\C}{\mathbb{C}}
\newcommand{\D}{\mathbb{D}}
\newcommand{\I}{\mathbb{I}}
\newcommand{\J}{\mathbb{J}}
\newcommand{\N}{\mathbb{N}}
\newcommand{\U}{\mathbb{U}}
\newcommand{\V}{\mathbb{V}}
\newcommand{\Z}{\mathbb{Z}}
\newcommand{\Y}{\mathbb{Y}}
\newcommand{\dsa}{$\dag$-$*$-autonomous}  
\newcommand{\dldc}{$\dag$-LDC}  
\newcommand{\m}{{\sf m}}
 
\newcommand{\nat}{\text{nat. }} 
\newcommand{\id}{\text{id}} 
\newcommand{\CP}{\mathsf{CP}}
\newcommand{\ox}{\otimes}
\newcommand{\pr}{\oplus}
\newcommand{\oa}{\oplus}
\newcommand{\op}{\mathsf{op}}
\newcommand{\rev}{\mathsf{rev}}
\newcommand{\mx}{\mathsf{mx}}
\newcommand{\Chu}{\mathsf{Chu}}
\newcommand{\Core}{\mathsf{Core}}
\newcommand{\Unitary}{\mathsf{Unitary}}
\newcommand{\dual}{\text{\reflectbox{$\Vdash$}}}
\newcommand{\fin}{\mathsf{FinSp}}
\newcommand{\lollipop}{\ensuremath{\!-\!\!\circ}}
\renewcommand{\bar}[1]{\overline{#1}}
\newcommand{\x}{\times}
\newcommand {\poppilol} {\reflectbox{$\multimap$}}
\newcommand{\cal}[1]{\mathcal{#1}}

\newcommand{\dashvv}{\dashv \!\!\!\!\! \dashv}  

\newcommand{\Hilb}{{\sf Hilb}}
\newcommand{\FHilb}{{\sf FHilb}}

\renewcommand{\phi}{\varphi}
\newcommand{\Asp}{\mathsf{Asp}}
\makeatletter

\newdimen\w@dth

\def\setw@dth#1#2{\setbox\z@\hbox{\scriptsize $#1$}\w@dth=\wd\z@
\setbox\@ne\hbox{\scriptsize $#2$}\ifnum\w@dth<\wd\@ne \w@dth=\wd\@ne \fi
\advance\w@dth by 1.2em}

\def\t@^#1_#2{\allowbreak\def\n@one{#1}\def\n@two{#2}\mathrel
{\setw@dth{#1}{#2}
\mathop{\hbox to \w@dth{\rightarrowfill}}\limits
\ifx\n@one\empty\else ^{\box\z@}\fi
\ifx\n@two\empty\else _{\box\@ne}\fi}}
\def\t@@^#1{\@ifnextchar_ {\t@^{#1}}{\t@^{#1}_{}}}

\def\t@left^#1_#2{\def\n@one{#1}\def\n@two{#2}\mathrel{\setw@dth{#1}{#2}
\mathop{\hbox to \w@dth{\leftarrowfill}}\limits
\ifx\n@one\empty\else ^{\box\z@}\fi
\ifx\n@two\empty\else _{\box\@ne}\fi}}
\def\t@@left^#1{\@ifnextchar_ {\t@left^{#1}}{\t@left^{#1}_{}}}

\def\two@^#1_#2{\def\n@one{#1}\def\n@two{#2}\mathrel{\setw@dth{#1}{#2}
\mathop{\vcenter{\hbox to \w@dth{\rightarrowfill}\kern-1.7ex
                 \hbox to \w@dth{\rightarrowfill}}%
       }\limits
\ifx\n@one\empty\else ^{\box\z@}\fi
\ifx\n@two\empty\else _{\box\@ne}\fi}}
\def\tw@@^#1{\@ifnextchar_ {\two@^{#1}}{\two@^{#1}_{}}}

\def\tofr@^#1_#2{\def\n@one{#1}\def\n@two{#2}\mathrel{\setw@dth{#1}{#2}
\mathop{\vcenter{\hbox to \w@dth{\rightarrowfill}\kern-1.7ex
                 \hbox to \w@dth{\leftarrowfill}}%
       }\limits
\ifx\n@one\empty\else ^{\box\z@}\fi
\ifx\n@two\empty\else _{\box\@ne}\fi}}
\def\t@fr@^#1{\@ifnextchar_ {\tofr@^{#1}}{\tofr@^{#1}_{}}}

\newdimen\W@dth
\def\setW@dth#1#2{\setbox\z@\hbox{$#1$}\W@dth=\wd\z@
\setbox\@ne\hbox{$#2$}\ifnum\W@dth<\wd\@ne \W@dth=\wd\@ne \fi
\advance\W@dth by 1.2em}

\def\T@^#1_#2{\allowbreak\def\N@one{#1}\def\N@two{#2}\mathrel
{\setW@dth{#1}{#2}
\mathop{\hbox to \W@dth{\rightarrowfill}}\limits
\ifx\N@one\empty\else ^{\box\z@}\fi
\ifx\N@two\empty\else _{\box\@ne}\fi}}
\def\T@@^#1{\@ifnextchar_ {\T@^{#1}}{\T@^{#1}_{}}}

\def\T@left^#1_#2{\def\N@one{#1}\def\N@two{#2}\mathrel{\setW@dth{#1}{#2}
\mathop{\hbox to \W@dth{\leftarrowfill}}\limits
\ifx\N@one\empty\else ^{\box\z@}\fi
\ifx\N@two\empty\else _{\box\@ne}\fi}}
\def\T@@left^#1{\@ifnextchar_ {\T@left^{#1}}{\T@left^{#1}_{}}}

\def\Tofr@^#1_#2{\def\N@one{#1}\def\N@two{#2}\mathrel{\setW@dth{#1}{#2}
\mathop{\vcenter{\hbox to \W@dth{\rightarrowfill}\kern-1.7ex
                 \hbox to \W@dth{\leftarrowfill}}%
       }\limits
\ifx\N@one\empty\else ^{\box\z@}\fi
\ifx\N@two\empty\else _{\box\@ne}\fi}}
\def\T@fr@^#1{\@ifnextchar_ {\Tofr@^{#1}}{\Tofr@^{#1}_{}}}

\def\Two@^#1_#2{\def\N@one{#1}\def\N@two{#2}\mathrel{\setW@dth{#1}{#2}
\mathop{\vcenter{\hbox to \W@dth{\rightarrowfill}\kern-1.7ex
                 \hbox to \W@dth{\rightarrowfill}}%
       }\limits
\ifx\N@one\empty\else ^{\box\z@}\fi
\ifx\N@two\empty\else _{\box\@ne}\fi}}
\def\Tw@@^#1{\@ifnextchar_ {\Two@^{#1}}{\Two@^{#1}_{}}}

\def\to{\@ifnextchar^ {\t@@}{\t@@^{}}}
\def\from{\@ifnextchar^ {\t@@left}{\t@@left^{}}}
\def\tofro{\@ifnextchar^ {\t@fr@}{\t@fr@^{}}}
\def\To{\@ifnextchar^ {\T@@}{\T@@^{}}}
\def\From{\@ifnextchar^ {\T@@left}{\T@@left^{}}}
\def\Two{\@ifnextchar^ {\Tw@@}{\Tw@@^{}}}
\def\Tofro{\@ifnextchar^ {\T@fr@}{\T@fr@^{}}}

\makeatother

\newcommand{\pullbackcorner}[1][ul]{\save*!/#1+1.2pc/#1:(1,-1)@^{|-}\restore}
\newcommand{\pushoutcorner}[1][dr]{\save*!/#1-1.2pc/#1:(-1,1)@^{|-}\restore}

\newcommand{\mulmap}[2]{
\begin{tikzpicture}[scale=#1, rotate=180]
	\begin{pgfonlayer}{nodelayer}
		\node [style=circle,scale=0.4, fill=#2] (0) at (-1, 3) {};
		\node [style=none] (1) at (-1.25, 2.75) {};
		\node [style=none] (2) at (-0.75, 2.75) {};
		\node [style=none] (3) at (-1, 3.25) {};
	\end{pgfonlayer}
	\begin{pgfonlayer}{edgelayer}
		\draw [style=none] (3.center) to (0);
		\draw [style=none, in=90, out=180, looseness=1.50] (0) to (1.center);
		\draw [style=none, in=90, out=0, looseness=1.50] (0) to (2.center);
	\end{pgfonlayer}
\end{tikzpicture}
}

\newcommand{\leftaction}[2]{
\begin{tikzpicture}[scale=#1]
	\begin{pgfonlayer}{nodelayer}
		\node [style=none] (0) at (0.5, 0.25) {};
		\node [style=none] (1) at (0.5, 0.75) {};
		\node [style=none] (2) at (0, 0.75) {};
		\node [style=none] (3) at (0.5, -0.75) {};
		\node [style=none] (4) at (0.5, 2) {};
		\node [style=none] (5) at (-0.75, 2) {};
	\end{pgfonlayer}
	\begin{pgfonlayer}{edgelayer}
		\draw[fill=#2] (0.center) -- (1.center) -- (2.center) -- (0.center);
		\draw (3.center) to (0.center);
		\draw (1.center) to (4.center);
		\draw [bend left, looseness=1.00] (2.center) to (5.center);
	\end{pgfonlayer}
\end{tikzpicture}
}

\newcommand{\unitmap}[2]{
\begin{tikzpicture}[scale=#1]
	\begin{pgfonlayer}{nodelayer}
		\node [style=circle, scale=0.4, fill=#2] (0) at (0, 0) {};
		\node [style=none] (1) at (0, -0.5) {};
	\end{pgfonlayer}
	\begin{pgfonlayer}{edgelayer}
		\draw [style=none] (0) to (1.center);
	\end{pgfonlayer}
\end{tikzpicture}
}

\newcommand{\counitmap}[2]{
\begin{tikzpicture}[scale=#1, rotate=180]
	\begin{pgfonlayer}{nodelayer}
		\node [style=circle, scale=0.4, fill=#2] (0) at (0, 0) {};
		\node [style=none] (1) at (0, -0.5) {};
	\end{pgfonlayer}
	\begin{pgfonlayer}{edgelayer}
		\draw [style=none] (0) to (1.center);
	\end{pgfonlayer}
\end{tikzpicture}
}

\newcommand{\bialgunitmap}[1]{
\begin{tikzpicture}[scale=#1]
	\begin{pgfonlayer}{nodelayer}
		\node [style=none] (0) at (-1, 2) {};
		\node [style=none] (1) at (-0.5, 2) {};
		\node [style=none] (2) at (-0.75, 1.75) {};
		\node [style=none] (3) at (-0.75, 1.25) {};
	\end{pgfonlayer}
	\begin{pgfonlayer}{edgelayer}
		\draw (0.center) to (1.center);
		\draw (1.center) to (2.center);
		\draw (2.center) to (0.center);
		\draw (2.center) to (3.center);
		\draw[fill=white] (0) -- (1) -- (2) -- (0);
	\end{pgfonlayer}
\end{tikzpicture}
}

\newcommand{\bialgcounitmap}[1]{
\begin{tikzpicture}[scale=#1]
	\begin{pgfonlayer}{nodelayer}
		\node [style=none] (0) at (-1, 1) {};
		\node [style=none] (1) at (-0.5, 1) {};
		\node [style=none] (2) at (-0.75, 1.25) {};
		\node [style=none] (3) at (-0.75, 1.75) {};
	\end{pgfonlayer}
	\begin{pgfonlayer}{edgelayer}
		\draw (0.center) to (1.center);
		\draw (1.center) to (2.center);
		\draw (2.center) to (0.center);
		\draw (3.center) to (2.center);
		\draw[fill=white] (0) -- (1) -- (2) -- (0);
	\end{pgfonlayer}
\end{tikzpicture}
}

\newcommand{\twistedmulmap}[2]{
\begin{tikzpicture}[scale=#1]
	\begin{pgfonlayer}{nodelayer}
		\node [style=none] (0) at (1.75, -0.5) {};
		\node [style=circle, scale=0.4, fill=#2] (1) at (1.5, -1.25) {};
		\node [style=none] (2) at (1.25, -0.5) {};
		\node [style=none] (3) at (1.5, -1.5) {};
	\end{pgfonlayer}
	\begin{pgfonlayer}{edgelayer}
		\draw [style=none, in=150, out=-90, looseness=2.00] (0.center) to (1);
		\draw [style=none, in=-90, out=30, looseness=1.75] (1) to (2.center);
		\draw [style=none] (1) to (3.center);
	\end{pgfonlayer}
\end{tikzpicture}
}

\newcommand{\productunitmap}[3]{
\begin{tikzpicture}[scale=#1]
	\begin{pgfonlayer}{nodelayer}
		\node [style=circle, scale=0.4, fill=#2] (0) at (0, 0) {};
		\node [style=none] (1) at (0, -0.5) {};
		\node [style=none] (2) at (0.25, -0.5) {};
		\node [style=circle, scale=0.4, fill=#3] (3) at (0.25, 0) {};
	\end{pgfonlayer}
	\begin{pgfonlayer}{edgelayer}
		\draw [style=none] (0) to (1.center);
		\draw [style=none] (3) to (2.center);
	\end{pgfonlayer}
\end{tikzpicture}
}

\newcommand{\productmulmap}[3]{
\begin{tikzpicture}[scale=#1]
	\begin{pgfonlayer}{nodelayer}
		\node [style=circle, scale=0.4, fill=#2] (0) at (0, -0.5) {};
		\node [style=circle, scale=0.4, fill=#3] (1) at (0.5, -0.5) {};
		\node [style=none] (2) at (-0.5, 0) {};
		\node [style=none] (3) at (0, 0) {};
		\node [style=none] (4) at (0.5, 0) {};
		\node [style=none] (5) at (1, 0) {};
		\node [style=none] (6) at (0, -1) {};
		\node [style=none] (7) at (0.5, -1) {};
	\end{pgfonlayer}
	\begin{pgfonlayer}{edgelayer}
		\draw [style=none, bend right, looseness=1.00] (2.center) to (0);
		\draw [style=none, bend right, looseness=1.00] (0) to (4.center);
		\draw [style=none, bend right=15, looseness=1.25] (3.center) to (1);
		\draw [style=none, bend right, looseness=1.00] (1) to (5.center);
		\draw [style=none] (0) to (6.center);
		\draw [style=none] (1) to (7.center);
	\end{pgfonlayer}
\end{tikzpicture}
}

\newcommand{\comulmap}[2]{
\begin{tikzpicture}[scale=#1]
	\begin{pgfonlayer}{nodelayer}
		\node [style=circle,scale=0.4, fill=#2] (0) at (-1, 3) {};
		\node [style=none] (1) at (-1.25, 2.75) {};
		\node [style=none] (2) at (-0.75, 2.75) {};
		\node [style=none] (3) at (-1, 3.25) {};
	\end{pgfonlayer}
	\begin{pgfonlayer}{edgelayer}
		\draw [style=none] (3.center) to (0);
		\draw [style=none, in=90, out=180, looseness=1.50] (0) to (1.center);
		\draw [style=none, in=90, out=0, looseness=1.50] (0) to (2.center);
	\end{pgfonlayer}
\end{tikzpicture}
}

\newcommand{\conjugateaction}[2]{
\begin{tikzpicture}[scale=#1]
	\begin{pgfonlayer}{nodelayer}
		\node [style=none] (0) at (0, 1) {};
		\node [style=none] (1) at (0, 0.75) {};
		\node [style=none] (2) at (-0.25, 1) {};
		\node [style=circle] (3) at (0, 1.5) {};
		\node [style=none] (4) at (0, 2) {};
		\node [style=none] (5) at (-0.75, 2) {};
		\node [style=none] (6) at (0, 0.5) {};
		\node [style=none] (7) at (0, 1.5) {$s$};
	\end{pgfonlayer}
	\begin{pgfonlayer}{edgelayer}
		\draw (4.center) to (3);
		\draw (3) to (0.center);
		\draw (2.center) to (0.center);
		\draw (0.center) to (1.center);
		\draw (1.center) to (2.center);
		\draw [bend right, looseness=1.00] (5.center) to (2.center);
		\draw (1.center) to (6.center);
	\end{pgfonlayer}
	\draw[fill=#2] (0.center) -- (1.center) -- (2.center) -- (0.center);
\end{tikzpicture}
}

\newcommand{\tinycomulmap}{
{\tiny \comulmap{0.5}{white}}
}
\begin{document}

\title[Dagger linear logic for categorical quantum mechanics]{Dagger linear logic\texorpdfstring{\\}{} for categorical quantum mechanics}

\author[R.~Cockett]{Robin Cockett\rsuper{a}}	%required

\author[C.~Comfort]{Cole Comfort\rsuper{b}}

\author[P.~V.~Srinivasan]{Priyaa V. Srinivasan\rsuper{a}}	%optional

\address{Department of  Computer Science, University of Calgary}	%required
\email{robin@ucalgary.ca, priyaavarshinee.srin@ucalgary.ca}  %optional
\address{Department of  Computer Science, University of Oxford}
\email{cole.comfort@new.ox.ac.uk}

\thanks{Partially supported by NSERC, Canada}	%optional

\begin{abstract}
Categorical quantum mechanics exploits the dagger compact closed structure of finite dimensional Hilbert spaces, and 
uses the graphical calculus of string diagrams to facilitate reasoning about finite dimensional processes.  A significant portion of 
quantum physics, however, involves reasoning about infinite dimensional processes, and it is well-known that the category of all 
Hilbert spaces is not compact closed.  Thus, a limitation of using dagger compact closed categories is that one cannot directly 
accommodate reasoning about infinite dimensional processes.  

A natural categorical generalization of compact closed categories, in which infinite dimensional spaces
 can be modelled, is $*$-autonomous categories and, more generally, linearly distributive categories.  
 This article starts the development of this direction of generalizing categorical quantum mechanics.  
 An important first step is to establish the behaviour of the dagger in these more general settings.  
 Thus, these notes simultaneously develop the categorical semantics of multiplicative dagger linear logic. 

The notes end with the definition of a mixed unitary category.  It is this structure which is subsequently 
used to extend the key features of categorical quantum mechanics. 
\end{abstract}

\maketitle

%%%%%%%%%%%%%%%%%%%%%%%%%%%%%%%%%%%%%%%%%%%%%%%%
%%%%%%%%%%%%%%%%%%%%%%%%%%%%%%%%%%%%%%%%%%%%%%%%

\section{Introduction}

Categorical quantum mechanics (CQM), as described in \cite{CoK17,HeV19}, employs a graphical calculus for  \dag-compact closed categories 
(\dag-KCCs) to study quantum processes within the \dag-KCC of finite dimensional Hilbert spaces ($\FHilb$).  From a logical perspective, the graphical calculus is 
the proof theory of a compact fragment of multiplicative \dag-linear logic \cite{Dun06}.  This programme of CQM was initiated by Abramsky and Coecke's seminal paper 
\cite{AC04} and it has allowed much of the structure of $\FHilb$ to be abstracted away and absorbed into the graphical calculus.

 A well-known limitation of compact closed categories, is that, while they model finite dimensional Hilbert spaces,  they do not model infinite dimensional 
 spaces \cite{Heunen16}.  A categorical generalization of compact closed categories, in which infinite dimensions spaces can be modelled, 
 however, is *-autonomous categories. Thus, from a categorical perspective, an inviting direction to accommodate infinite dimensional 
 quantum processes is to utilize these *-autonomous settings as extensions of the compact closed settings of CQM. This can be taken a 
 step further by generalizing to linearly distributive categories (LDC) in which the existence of dual objects is not assumed. These linear 
 settings come with a proof theory -- a graphical calculus -- so that we need not abandon this key aspect of CQM.  They also provide a 
 natural setting for Frobenius structures (see Egger \cite{Egg10}), thus another important innovation of CQM can potentially be extended. 
 The downside is that these linear settings are more complicated as the single tensor product of CQM splits into two tensor products 
 (the tensor, $\ox$, and the ``par'', $\oa$\footnote{In the linear logic community the par is often written $\parr$.}).
 
 A first step towards extending CQM from compact closed categories to these linear settings is to understand how the dagger is expressed in 
 $*$-autonomous and linearly distributive categories.  This issue is the main focus of this paper.  In CQM, the dagger structure directly 
 determines the important notion of a unitary isomorphism.  As is discussed further below, the expression of unitary structure in linear 
 settings, is more complicated and leads ultimately to the introduction of mixed unitary categories (MUCs).  A mixed unitary category may be 
 viewed as an essentially traditional CQM setting extended by a dagger linear setting in which (potentially) infinite dimensional objects can reside. 
 In our subsequent papers, \cite{CS19, CoS20} we show how some fundamental features of CQM can be realized in mixed unitary categories.  
 In \cite{CS19} we show, following the ideas in \cite{CoH16}, that completely positive maps can be expressed in mixed unitary categories and, 
 in \cite{CoS20}, we explore Frobenius structure and complementarity in the mixed unitary setting.
   
The problem of extending CQM to include infinite dimensional processes has been explored in a number of different ways.  
A fundamental feature of CQM is the replacement of the notion of an orthonormal basis by the algebraic structure of a special commutative 
$\dagger$-Frobenius algebra.  In the category of finite dimensional Hilbert spaces the correspondence between orthonormal bases 
and special commutative $\dagger$-Frobenius algebras is precise \cite{CPV12}.  Taking this correspondence as fundamental,  
Abramsky and Heunen, in \cite{AbH12}, showed how Ambrose's $H^*$-algebras \cite{Amb45} could be used, in much the 
same way, to replace orthonormal bases in infinite dimensional Hilbert spaces.   However, there was a cost: they had to move 
to semi-Frobenius algebras (that is they had to drop the units from the Frobenius structure) and to require a special property 
(which they appropriately called $H^*$) on the maps from the unit.  Gogoiso and Genovese \cite{GG17} proposed an 
interesting approach to reinstating the units using techniques from non-standard analysis \cite{Rob96}.  
They considered $~^\star${\bf Hilb} the category of non-standard separable Hilbert Spaces and linear maps.  This they claimed is a 
dagger compact closed category, in which, among other things, the semi-Frobenius algebras of Abramsky and Heunen can be modelled.  
Furthermore, the units can be reinstated because formal infinite sums are permitted.  

In \cite{CoH16}, Coecke and Heunen, in order to include infinite dimensional  quantum processes, take the simple step of dropping the 
requirement of being compact closed and work in dagger symmetric monoidal categories ($\dagger$-SMCs).   The category of Hilbert 
spaces is the prototypical example of a $\dagger$-SMC.  In particular, they show how to build, for an arbitrary $\dagger$-SMC, the 
category of completely positive maps.  Of course, it is also possible to consider special commutative $\dagger$-Frobenius algebras in an 
arbitrary $\dagger$-SMC: although these, for Hilbert spaces are, of course, just the finite dimensional objects again.  In \cite{HeR18}, 
Heunen and Reyes considered a different $\dagger$-SMC, namely, the category of Hilbert $\C^*$-modules.  Its objects can be equivalently viewed as 
bundles of Hilbert spaces over a locally compact Hausdorff space.  They characterized the special commutative $\dagger$-Frobenius algebras in this
 category as bundles of finite dimensional Hilbert spaces (with dimensions are uniformly bounded).  These objects, while being far from finite, do retain a 
 (uniform) locally finite nature.  This example,  by using vector bundles and ideas from differential geometry, enters the domain of traditional 
 theoretical physics, and serves as a reminder that $\dagger$-Frobenius algebras are not only of interest for Hilbert spaces. 

The approach taken in this article to the problem of allowing infinite dimensional objects is different again.  We accept CQM and its finite dimensionality as a 
feature not a bug.  However, we also accept that it is useful to have access to infinite dimensional structures. Rather than insisting that these infinite 
dimensional structures are concretely related to Hilbert spaces, we allow that they may be a system of formal types which extend an essentially traditional 
CQM core.  A prototypical example of this is provided by the extension of finite dimensional complex matrices to infinite dimensional ``finiteness matrices'' 
(see Example \ref{Sec: Finiteness matrices}).  

The system of types extending the core are minimally expected to organize themselves into a linear setting with a dagger, that is, into a dagger 
linearly distributive category (a $\dagger$-LDC).   There are various reasons for this expectation. To start with the CQM core with its dagger involution is 
already an example of a $\dagger$-LDC (albeit a ``compact'' one, that is one in which the two tensor products, tensor and par, coincide).  Formally 
extending such a CQM core with (infinite) limits and colimits freely -- using the bicompletion of Joyal \cite{Joy95} -- causes the tensor product of the 
core to separate into two linearly related tensor structures.  These form the tensor, $\ox$, and the par, $\oa$, of linear logic and of a $\dagger$-LDC.   
Extending a setting to accommodate infinite objects may, dually, be viewed as a process of extracting a core of ``finite objects'' from that larger setting.  
Taking this latter view, suggests that we should start with a linear setting with a dagger, namely a $\dagger$-LDC, and extract from it a traditional CQM core. 

At the very outset of such a program, there are some immediate -- and perhaps paradigm breaking -- questions to be faced.  The most immediate one is the 
question of what a  $\dagger$-LDC might be.  To answer this question is also to answer the question of what (multiplicative) ``dagger linear logic'' is.  
This article essentially focuses on this question.  It lays down the categorical groundwork for dagger linear logic and, thereby, begins to build a bridge to CQM.

Recall, following the lead of \cite{AC04,Sel07}, that it is now standard in CQM to interpret the dagger functor as a stationary on objects $(A = A^\dagger)$ involution.  
However, in the setting of linear logic and LDCs, one expects an involution to flip the tensor and par structure so that $A^\dagger \ox B^\dagger = (A \oa B)^\dagger$.  This has 
the -- perhaps uncomfortable -- effect of implying that,  in these more general settings, the involution can no longer be viewed as being stationary on objects.  Of course, having started 
down this road, it seems prudent also to replace the equality above by a coherent isomorphism $\lambda_\ox: A^\dagger \ox B^\dagger \to (A \oa B)^\dagger$ and indeed the involution by an isomorphism $\iota_A: A \to (A^\dagger)^\dagger$.  

At this juncture it is perhaps appropriate to acknowledge that these are not new ideas.  Models for quantum mechanics in $*$-autonomous categories are often described 
as ``toy models'' \cite{Abr12} and were, in particular discussed by Pavlovic \cite{Pav11} where some very similar directions were advocated.   Indeed, Egger \cite{Egg11}, in initiating the 
development of ``involutive'' categories, was also implicitly suggesting that dagger functors should not necessarily be regarded as being stationary on objects in these more 
general settings.  Section~\ref{daggers-duals-conjugation} is essentially a realization of Egger's ideas: we have, however, changed his terminology preferring to talk of ``conjugation'' 
rather than ``involution'' as we think of the contravariant dagger, $(\_)^\dagger$, as an involution.  Nonetheless, conjugation and involution are closely related (see section \ref{daggers-duals-conjugation}) as in the presence of dualization having conjugation is equivalent to having an involution. 

There are some significant complications attendant on allowing a non-stationary dagger.  The first and foremost amongst these is that 
 one gets inundated by coherence issues.  This article does provide a path through these coherence issues, and, hopefully shows -- once one has assimilated all 
 the structure -- that these issues are not nearly as terrible as might be expected.  However, we are forced to concede that they are non-trivial.  The next problem, 
 once the coherences are under control, is that one would like to be able to say what a unitary isomorphism is with respect to a non-stationary involution.  How 
 this may be accomplished seems, at first glance, less than obvious.  

The fact that the dagger  functor is an involution with a coherent isomorphism $\iota_A: A \to (A^\dagger)^\dagger$ makes it natural to view a 
{\bf unitary object} as an object with an isomorphism $\varphi_A: A \to A^\dagger$, such that $\iota_A = \varphi_A (\varphi_A^{-1})^\dagger$: we refer to $\varphi_A$ as 
the {\bf unitary structure} of $A$.  Considering this, with the expected coherent behaviour of this unitary structure, leads one to realize that, for unitary objects, we 
must also have $A \ox B \simeq A \oa B$.  This, in turn, leads one to ask how this can happen in an LDC.  Fortunately, there is a theory which has been developed 
for this situation: namely that of LDCs with mix \cite{CS97,Pav11}.  An LDC with mix has a coherent map, from the par unit to the tensor unit ${\sf m}: \bot \to \top$, called the {\bf mix map}, which, in turn, induces a map ${\sf mx}: A \ox B \to A \oa B$, called the {\bf mixor}.  In a mix category we say an object $A$ is in the {\bf core} \cite{BCS00} in case the mixor for that object with any other object, ${\sf mx}_{A,X}: A \ox X \to A \oa X$, is an isomorphism.  This allows our earlier expectation to be expressed as the requirement that unitary objects be in the core.  As $\dagger$-LDCs with mix play an important role in this development we refer to them as a $\dagger$-mix category.  If, further, we want our unitary objects to form a {\bf compact} LDC -- that is one in which the tensor and par structures are equivalent -- then the tensor and par structures must agree at their units, which means we must ask that the mix map, ${\sf m}: \bot \to \top$, be an isomorphism.   The first milestone of the paper is therefore to collect all this structure into what we call a {\bf $\dagger$-isomix category}.   

Amidst the introduction of all this structure, the astute reader may notice that we have still failed to elucidate how the unitary isomorphisms arise.  Let us hasten to correct 
this oversight: a unitary isomorphism $f: A \to B$ is an isomorphism between unitary objects which is (twisted) natural with respect to the unitary structures, $\varphi_A$ and $\varphi_B$\footnote{This formulation of unitary isomorphisms is not completely original as a lively discussion of whether $\dagger$-categories were ``evil''  led 
Peter Lumsdaine to suggest in the math overflow forum \cite{Lums15} how they might be made a little less evil.  These ideas never took off, perhaps because it was pointed 
out by Peter Selinger that, when one regarded something as evil if it was not preserved by equivalence, then it was impossible for dagger not to be evil! Here we quite explicitly 
have ``unitary structure'' which is of course is not only not necessarily unique but also will not necessarily be preserved by an equivalence.  Thus, like all structure, it is thoroughly evil!},
in the sense that  the following diagram is rendered commutative:
\[ \xymatrix{ A \ar[d]_{\varphi_A} \ar[rr]^f & & B \ar[d]^{\varphi_B} \\ A^\dagger  && B^\dagger  \ar[ll]^{f^\dagger} } \]
Note that, when the  unitary structure is the identity map, we recover the usual notion of unitary isomorphism.
The coherence requirements on unitary structure then have the pleasing effect of forcing the coherence isomorphisms, for unitary objects, 
to be unitary maps.  

A {\bf unitary category} is a $\dagger$-isomix category in which {\em all\/} objects have a unitary structure.   Unitary categories are necessarily compact LDCs and so are rather special. 
In fact, we show that they are $\dagger$-linearly equivalent to the more standard CQM notion of a $\dagger$-SMC -- and, furthermore, that a closed 
unitary category is linearly equivalent to a dagger compact  closed category, $\dagger$-KCC.   One may -- with some justification -- feel that one has, at this stage, come full circle as the standard structures from categorical quantum mechanics seem to be emerging.   However, notice that it is our more complex notion of unitary that allows the extraction of a unitary core from a larger $\dagger$-isomix category: the larger $\dagger$-isomix category, with its possibly infinite dimensional objects, is still available.

A mixed unitary category (MUC) is a strong $\dagger$-isomix functor $M: \U \to \C$,  which has domain a ``small'' unitary category (the unitary core) and codomain of a ``large'' $\dagger$-isomix category.  
A mixed unitary category can be represented schematically as follows:
\[  \begin{tikzpicture}[scale=1.5] %MUC.tikz
	\begin{pgfonlayer}{nodelayer}
		\node [style=circle, scale=14] (0) at (4.5, -0) {};
		\node [style=none] (1) at (-3, 2) {};
		\node [style=none] (2) at (-5, -0) {};
		\node [style=none] (3) at (-3, -2) {};
		\node [style=none] (4) at (-1, -0) {};
		\node [style=none] (5) at (-3, -0.75) {$A \to^{\varphi_A}_{\simeq} A^\dagger$};
		\node [style=none] (6) at (-3, 0.5) {Unitary};
		\node [style=none] (7) at (-3, -0) {category};
		\node [style=none] (8) at (-0.75, -0) {};
		\node [style=none] (9) at (2.25, -0) {};
		\node [style=none] (10) at (0.25, 0.25) {$\dagger$-isomix};
		\node [style=none] (11) at (0.25, -0.25) {functor};
		\node [style=none] (12) at (4.5, 2) {$\dagger$-isomix};
		\node [style=none] (13) at (4.5, 1.5) {category};
		\node [style=none] (14) at (4.5, -1.5) {$B$};
		\node [style=none] (15) at (5.5, -2) {$B^\dagger$};
		\node [style=none] (16) at (3, 0.75) {};
		\node [style=none] (17) at (2.25, -0.5) {};
		\node [style=none] (18) at (3.25, -1.5) {};
		\node [style=none] (19) at (4, -0.25) {};
	\end{pgfonlayer}
	\begin{pgfonlayer}{edgelayer}
		\draw (1.center) to (2.center);
		\draw (2.center) to (3.center);
		\draw (3.center) to (4.center);
		\draw (4.center) to (1.center);
		\draw [->] (8.center) to (9.center);
		\draw [dotted] (16.center) to (17.center);
		\draw [dotted] (17.center) to (18.center);
		\draw [dotted] (18.center) to (19.center);
		\draw [dotted] (19.center) to (16.center);
	\end{pgfonlayer}
\end{tikzpicture} \]

One may think of the unitary category as acting on the larger category, much as a field $K$ acts on a $K$-algebra as scalars.   
Expressing these categories in this manner allows an obvious notion of functor as a square of $\dagger$-Frobenius functors, whose component on the unitary categories preserves unitary structure, and which commutes up to a linear natural isomorphism.  
For any $\dagger$-isomix category one can build a unitary category by collecting the ``pre-unitary'' objects which are in the core. This give a way of generating a mixed unitary category from a $\dagger$-isomix category, which we call the unitary construction.

We have provided examples throughout the text.  An important example, closely related to traditional CQM, is the $*$-autonomous category of ``finiteness matrices'', ${\sf FMat}(\C)$, over the complex numbers \cite{Ehr05} (see Example \ref{Sec: Finiteness matrices}).  Here the maps are infinite dimensional matrices whose support is carefully controlled by the finiteness structure.  The dagger on the category is given by simultaneously transposing and conjugating the matrices: on objects it is given by taking the dual finiteness space.  ${\sf FMat}(\C)$ forms a $\dagger$-isomix category which, furthermore, is $\dagger$-$*$-autonomous.  An object is in the core if and only if it's underlying finiteness space is finite and these objects are also exactly the unitary objects.  The unitary structure of a finite object in this case is the identity map (so the unitary structure is ``trivial'') -- which means that the coherence requirements are immediately satisfied.

Another source of examples (see Example \ref{Section: Chu}) is from the Chu construction, \cite{Barr06}, where the dualizing object is set to the tensor unit.  Considering the Chu construction over complex vector spaces there is an obvious notion of conjugation which means that this category forms a $\dagger$-isomix category.  From there one can obtain a non-trivial MUC by extracting the pre-unitary objects, or, more directly, by using the fact that the category of Hilbert spaces embeds into this category.  To obtain a MUC one must then restrict this last embedding to the finite dimensional Hilbert spaces.

\paragraph*{Notation:} Throughout the paper we use diagrammatic order for composition of maps so $fg: A \to C = A \to^{f} B \to^{g} C$.   Functors, however, are usually written using applicative notation so that $GF(A) = F(G(A))$.   We write par, the dual linear logic tensor, as $\oa$ throughout.  Circuit diagrams should be read from top to bottom (in the same direction as the pull of  gravity).  

%%%%%%%%%%%%%%%%%%%%%%%%%%%%%%%%%%%%%%%%%%%%%%%%

\section{Linearly distributive categories}

%%%%%%%%%%%%%%%%%%%%%%%%%%%%%%%%%%%%%%%%%%%%%%%%

This section recalls some background concepts from the theory of linearly distributive categories. The definition of linearly distributive categories is available in \cite{CS97,BCST96}.  Here we briefly recall the definition of linear functors and their transformations \cite{CS99}, the notion of a linear adjoint \cite{CKS00} -- which we shall refer to as a linear dual -- and the notion of a mix category and its core \cite{CS97a,BCS00}.

\subsection{Linearly distributive categories, functors, and transformations}

A {\bf linearly distributive category (LDC)} is a category, $\X$, with two monoidal structures 
$$(\ox, \top, a_\ox, u_\ox^L, u_\ox^R)~~~\mbox{ and } ~~~ (\oa, \bot, a_\oa, u_\oa^L, u_\oa^R)$$ 
linked by natural transformations called the linear distributors:
\begin{align*}
&\partial^L: A \ox (B \oa C) \rightarrow  (A \ox B) \oa C \\
& \partial^R: (B \oa C) \ox A \rightarrow B \oa (C \ox A)
\end{align*} 
such that the monoidal natural isomorphisms - associators and unitors - interact coherently with the linear distributors, see \cite{BCST96, CS97} for more details.  A symmetric LDC is an LDC in which both monoidal structures are symmetric,  with symmetry maps $c_\ox$ and $c_\oa$, such that
$\partial^R = c_\ox (1 \ox c_\oa) \partial^L (c_\ox \oa 1) c_\oa$. 

LDCs provide categorical semantics for the multiplicative linear logic (MLL).  LDCs come equipped with a graphical calculus \cite{BCST96} that contains the calculus for monoidal categories. 

 In this section, we review the fundamentals of the graphical calculus for LDCs. For detailed exposition, see \cite{BCST96, CS97}.  
The following are the generators of LDC circuits: wires represent objects and circles represent maps. The input wires of a map are tensored (with $\ox$), and the output wires are ``par''ed (with $\oa$). The following diagram represents a map $f: A \ox B \to C \oa D$. 
 \[ \begin{tikzpicture}
	\begin{pgfonlayer}{nodelayer}
		\node [style=circle, scale=2] (0) at (-0.25, -0.5) {};
		\node [style=none] (1) at (-0.25, -0.5) {$f$};
		\node [style=none] (2) at (-1, 0.5) {$A$};
		\node [style=none] (3) at (0.5, 0.5) {$B$};
		\node [style=none] (4) at (-1, -1.5) {$C$};
		\node [style=none] (5) at (0.5, -1.5) {$D$};
		\node [style=ox] (6) at (-0.25, 1) {};
		\node [style=oa] (7) at (-0.25, -2) {};
		\node [style=none] (8) at (-0.25, -2.75) {};
		\node [style=none] (9) at (-0.25, 1.75) {};
		\node [style=none] (10) at (-0.5, -1.65) {};
		\node [style=none] (11) at (0, -1.65) {};
		\node [style=none] (12) at (0, 0.65) {};
		\node [style=none] (13) at (-0.5, 0.65) {};
		\node [style=none] (14) at (-0.25, -3) {$f: A \oa B \to C \ox D$};
	\end{pgfonlayer}
	\begin{pgfonlayer}{edgelayer}
		\draw [bend right=45, looseness=1.25] (6) to (0);
		\draw [bend left=45, looseness=1.25] (6) to (0);
		\draw [bend left=45, looseness=1.25] (7) to (0);
		\draw [bend right=45, looseness=1.25] (7) to (0);
		\draw (9.center) to (6);
		\draw (7) to (8.center);
	\end{pgfonlayer}
\end{tikzpicture}  \]

The $\ox$-associator, the $\oa$-associator, the left linear distributor, and the right linear distributors are, respectively, drawn as follows:
\[(a) ~~~ \begin{tikzpicture}
	\begin{pgfonlayer}{nodelayer}
		\node [style=ox] (0) at (1, -2.5) {};
		\node [style=ox] (1) at (0.5, -0) {};
		\node [style=ox] (2) at (1.75, -1.5) {};
		\node [style=ox] (3) at (1.5, 1.25) {};
		\node [style=none] (4) at (1.5, 2) {};
		\node [style=none] (5) at (1, -3.25) {};
		\node [style=none] (6) at (0.75, -0.5) {};
		\node [style=none] (7) at (0.25, -0.5) {};
		\node [style=none] (8) at (1.25, 1) {};
		\node [style=none] (9) at (1.75, 1) {};
		\node [style=none] (10) at (2.75, 1.75) {$A \ox (B \ox C)$};
		\node [style=none] (11) at (2.25, -3) {$(A \ox B) \ox C$};
		\node [style=none] (12) at (0.5, 1) {$A \ox B$};
		\node [style=none] (13) at (2.25, 0.75) {$C$};
		\node [style=none] (14) at (2.25, -2.25) {$B \ox C$};
		\node [style=none] (15) at (0, -2) {$A$};
	\end{pgfonlayer}
	\begin{pgfonlayer}{edgelayer}
		\draw (5.center) to (0);
		\draw [bend right, looseness=1.00] (0) to (2);
		\draw [in=-120, out=160, looseness=1.00] (0) to (1);
		\draw (3) to (4.center);
		\draw [in=-150, out=75, looseness=0.75] (1) to (3);
		\draw [in=-60, out=120, looseness=1.25] (2) to (1);
		\draw [in=-30, out=75, looseness=1.00] (2) to (3);
	\end{pgfonlayer}
\end{tikzpicture} ~~~~~~~~ (b)~~~~ \begin{tikzpicture}
	\begin{pgfonlayer}{nodelayer}
		\node [style=oa] (0) at (1.5, 0.75) {};
		\node [style=oa] (1) at (2, -1.75) {};
		\node [style=oa] (2) at (0.75, -0.25) {};
		\node [style=oa] (3) at (1, -3) {};
		\node [style=none] (4) at (1, -3.75) {};
		\node [style=none] (5) at (1.5, 1.5) {};
		\node [style=none] (6) at (1.8, -1.35) {};
		\node [style=none] (7) at (2.15, -1.35) {};
		\node [style=none] (8) at (1.35, -2.65) {};
		\node [style=none] (9) at (0.75, -2.65) {};
		\node [style=none] (10) at (2.25, -3.5) {$A \oa (B \oa C)$};
		\node [style=none] (11) at (2.5, 1.25) {$(A \oa B) \oa C$};
		\node [style=none] (12) at (0.75, 0.5) {$A$};
		\node [style=none] (13) at (2.25, -2.5) {$B \oa C$};
		\node [style=none] (14) at (1.25, -1.25) {$B$};
		\node [style=none] (15) at (2.5, -0) {$C$};
		\node [style=none] (16) at (0.25, -2.5) {$A$};
	\end{pgfonlayer}
	\begin{pgfonlayer}{edgelayer}
		\draw (5.center) to (0);
		\draw [bend right, looseness=1.00] (0) to (2);
		\draw [in=60, out=-20, looseness=1.00] (0) to (1);
		\draw (3) to (4.center);
		\draw [in=30, out=-105, looseness=0.75] (1) to (3);
		\draw [in=120, out=-60, looseness=1.25] (2) to (1);
		\draw [in=135, out=-130, looseness=1.00] (2) to (3);
	\end{pgfonlayer}
\end{tikzpicture} ~~~~~~~~ (c) ~~~~ \begin{tikzpicture}
	\begin{pgfonlayer}{nodelayer}
		\node [style=ox] (0) at (1.25, 0.75) {};
		\node [style=ox] (1) at (0.5, -1.25) {};
		\node [style=oa] (2) at (2, -0.25) {};
		\node [style=oa] (3) at (1.25, -2.5) {};
		\node [style=none] (4) at (1.25, -3.25) {};
		\node [style=none] (5) at (1.25, 1.5) {};
		\node [style=none] (6) at (1.5, -2.25) {};
		\node [style=none] (7) at (1, -2.25) {};
		\node [style=none] (8) at (1, 0.5) {};
		\node [style=none] (9) at (1.5, 0.5) {};
		\node [style=none] (10) at (2.5, 1.25) {$A \ox (B \oa C)$};
		\node [style=none] (11) at (0.25, -0) {$A$};
		\node [style=none] (12) at (2.55, 0.5) {$B \oa C$};
		\node [style=none] (13) at (1.25, -0.5) {$B$};
		\node [style=none] (14) at (2.5, -1) {$C$};
		\node [style=none] (15) at (0, -2) {$A \ox B$};
		\node [style=none] (16) at (2.55, -3.25) {$(A \ox B) \oa C$};
	\end{pgfonlayer}
	\begin{pgfonlayer}{edgelayer}
		\draw (5.center) to (0);
		\draw [bend left, looseness=1.00] (0) to (2);
		\draw [in=90, out=-150, looseness=1.00] (0) to (1);
		\draw (3) to (4.center);
		\draw [in=150, out=-75, looseness=1.00] (1) to (3);
		\draw [in=60, out=-120, looseness=1.25] (2) to (1);
		\draw [in=30, out=-75, looseness=1.00] (2) to (3);
	\end{pgfonlayer}
\end{tikzpicture} ~~~~~~~~ (d) ~~~~\begin{tikzpicture}
	\begin{pgfonlayer}{nodelayer}
		\node [style=ox] (0) at (1.3, 0.75) {};
		\node [style=ox] (1) at (2.05, -1.25) {};
		\node [style=oa] (2) at (0.55, -0.25) {};
		\node [style=oa] (3) at (1.3, -2.5) {};
		\node [style=none] (4) at (1.3, -3.25) {};
		\node [style=none] (5) at (1.3, 1.5) {};
		\node [style=none] (6) at (1.05, -2.25) {};
		\node [style=none] (7) at (1.55, -2.25) {};
		\node [style=none] (8) at (1.55, 0.5) {};
		\node [style=none] (9) at (1.05, 0.5) {};
		\node [style=none] (10) at (0, 1.25) {$(A \oa B) \oa C$};
		\node [style=none] (11) at (2.3, -0) {$C$};
		\node [style=none] (12) at (0, 0.5) {$A \oa B$};
		\node [style=none] (13) at (1.3, -0.5) {$B$};
		\node [style=none] (14) at (0.04999995, -1) {$A$};
		\node [style=none] (15) at (2.55, -2) {$B \ox C$};
		\node [style=none] (16) at (0, -3.25) {$A \oa (B \ox C)$};
	\end{pgfonlayer}
	\begin{pgfonlayer}{edgelayer}
		\draw (5.center) to (0);
		\draw [bend right, looseness=1.00] (0) to (2);
		\draw [in=90, out=-30, looseness=1.00] (0) to (1);
		\draw (3) to (4.center);
		\draw [in=30, out=-105, looseness=1.00] (1) to (3);
		\draw [in=120, out=-60, looseness=1.25] (2) to (1);
		\draw [in=150, out=-105, looseness=1.00] (2) to (3);
	\end{pgfonlayer}
\end{tikzpicture}  \]

$\begin{tikzpicture}
	\begin{pgfonlayer}{nodelayer}
		\node [style=oa] (0) at (1, -3) {};
		\node [style=none] (1) at (1, -3.75) {};
		\node [style=none] (2) at (1.25, -2.75) {};
		\node [style=none] (3) at (0.75, -2.75) {};
		\node [style=none] (7) at (0.5, -2) {};
		\node [style=none] (8) at (1.5, -2) {};
	\end{pgfonlayer}
	\begin{pgfonlayer}{edgelayer}
		\draw (0) to (1.center);
		\draw [in=-90, out=150, looseness=1.00] (0) to (7.center);
		\draw [in=-90, out=30, looseness=1.00] (0) to (8.center);
	\end{pgfonlayer}
\end{tikzpicture}$ is the $\oa$-introduction rule, $\begin{tikzpicture}
	\begin{pgfonlayer}{nodelayer}
		\node [style=ox] (0) at (1, -3) {};
		\node [style=none] (1) at (1, -3.75) {};
		\node [style=none] (2) at (1.25, -2.75) {};
		\node [style=none] (3) at (0.75, -2.75) {};
		\node [style=none] (7) at (0.5, -2) {};
		\node [style=none] (8) at (1.5, -2) {};
	\end{pgfonlayer}
	\begin{pgfonlayer}{edgelayer}
		\draw (0) to (1.center);
		\draw [in=-90, out=150, looseness=1.00] (0) to (7.center);
		\draw [in=-90, out=30, looseness=1.00] (0) to (8.center);
	\end{pgfonlayer}
\end{tikzpicture}$ is $\ox$-introduction rule, $\begin{tikzpicture}
	\begin{pgfonlayer}{nodelayer}
		\node [style=ox] (0) at (1, -2.75) {};
		\node [style=none] (1) at (1, -2) {};
		\node [style=none] (2) at (1.25, -3) {};
		\node [style=none] (3) at (0.75, -3) {};
		\node [style=none] (7) at (0.5, -3.75) {};
		\node [style=none] (8) at (1.5, -3.75) {};
	\end{pgfonlayer}
	\begin{pgfonlayer}{edgelayer}
		\draw (0) to (1.center);
		\draw [in=90, out=-150, looseness=1.00] (0) to (7.center);
		\draw [in=90, out=-30, looseness=1.00] (0) to (8.center);
	\end{pgfonlayer}
\end{tikzpicture}$ is the $\ox$-elimination rule, $\begin{tikzpicture}
	\begin{pgfonlayer}{nodelayer}
		\node [style=oa] (0) at (1, -2.75) {};
		\node [style=none] (1) at (1, -2) {};
		\node [style=none] (2) at (0.5, -3.75) {};
		\node [style=none] (3) at (1.5, -3.75) {};
	\end{pgfonlayer}
	\begin{pgfonlayer}{edgelayer}
		\draw (0) to (1.center);
		\draw [in=90, out=-150, looseness=1.00] (0) to (2.center);
		\draw [in=90, out=-30, looseness=1.00] (0) to (3.center);
	\end{pgfonlayer}
\end{tikzpicture}$ is the $\oa$-elimination rule.

The unitors are drawn as follows:
\[ (a) ~~~~ \begin{tikzpicture}
	\begin{pgfonlayer}{nodelayer}
		\node [style=circle] (0) at (0, -0) {$\top$};
		\node [style=none] (1) at (0, -2) {};
		\node [style=none] (2) at (0.75, -0) {};
		\node [style=none] (3) at (0.75, -2) {};
		\node [style=none] (4) at (1, -1) {$A$};
		\node [style=none] (5) at (0, -2.5) {$(u_\ox^L)^{-1}: A \to \top \ox A$};
	\end{pgfonlayer}
	\begin{pgfonlayer}{edgelayer}
		\draw (0) to (1.center);
		\draw (2.center) to (3.center);
	\end{pgfonlayer}
\end{tikzpicture} ~~~~~~~~ (b) ~~~~ \begin{tikzpicture}
	\begin{pgfonlayer}{nodelayer}
		\node [style=circle] (0) at (0, -0) {$\top$};
		\node [style=none] (1) at (0.75, 0.75) {};
		\node [style=none] (2) at (0.75, -2) {};
		\node [style=none] (3) at (1, -0.25) {$A$};
		\node [style=none] (4) at (0.25, -2.5) {$u_\ox^L: \top \ox A \to A$};
		\node [style=circle, scale=0.6] (5) at (0.75, -1.25) {};
		\node [style=none] (6) at (0, 0.75) {};
	\end{pgfonlayer}
	\begin{pgfonlayer}{edgelayer}
		\draw (1.center) to (2.center);
		\draw [dotted, bend right, looseness=1.25] (0) to (5);
		\draw (6.center) to (0);
	\end{pgfonlayer}
\end{tikzpicture} ~~~~~~~~ (c) ~~~ \begin{tikzpicture}
	\begin{pgfonlayer}{nodelayer}
		\node [style=circle] (0) at (0, -2.5) {$\bot$};
		\node [style=none] (1) at (0.75, -3.5) {};
		\node [style=none] (2) at (0.75, -0.5) {};
		\node [style=none] (3) at (1, -0.75) {$A$};
		\node [style=none] (4) at (0.5, -3.75) {$(u_\oa^L)^{-1}:  A \to \bot \oa A$};
		\node [style=circle, scale=0.6] (5) at (0.75, -1.5) {};
		\node [style=none] (6) at (0, -3.5) {};
	\end{pgfonlayer}
	\begin{pgfonlayer}{edgelayer}
		\draw (1.center) to (2.center);
		\draw [dotted, bend left, looseness=1.25, dotted] (0) to (5);
		\draw (6.center) to (0);
	\end{pgfonlayer}
\end{tikzpicture} ~~~~~~~ (d) ~~~ \begin{tikzpicture}
	\begin{pgfonlayer}{nodelayer}
		\node [style=circle] (0) at (0, -2.75) {$\bot$};
		\node [style=none] (1) at (0.75, -0.75) {};
		\node [style=none] (2) at (0.75, -3) {};
		\node [style=none] (3) at (1, -2.75) {$A$};
		\node [style=none] (4) at (0.25, -3.5) {$u_\oa^L:  \bot \oa A \to A$};
		\node [style=none] (5) at (0, -0.75) {};
	\end{pgfonlayer}
	\begin{pgfonlayer}{edgelayer}
		\draw (1.center) to (2.center);
		\draw (5.center) to (0);
	\end{pgfonlayer}
\end{tikzpicture}  \]

Diagram $(a)$ is called the left $\top$-introduction, $(b)$ is called the left $\top$-elimination, $(c)$ is the left $\bot$-introduction, and $(d)$ is the left $\bot$-elimination. The unit $\top$ is introduced, and the counit $\bot$ is eliminated using the thinning links which are shown using dotted wires in the diagrams. 

The following are a set of circuit equalities (which when oriented become reduction rewrite rules):
\[ \begin{tikzpicture}
	\begin{pgfonlayer}{nodelayer}
		\node [style=circle] (0) at (0, -1) {$\top$};
		\node [style=none] (1) at (0.75, 0.25) {};
		\node [style=none] (2) at (0.75, -3) {};
		\node [style=none] (3) at (1, -1.25) {$A$};
		\node [style=circle, scale = 0.4] (4) at (0.75, -2.25) {};
		\node [style=circle] (5) at (0, -0) {$\top$};
	\end{pgfonlayer}
	\begin{pgfonlayer}{edgelayer}
		\draw (1.center) to (2.center);
		\draw [dotted, bend right, looseness=1.25] (0) to (4);
		\draw (5) to (0);
	\end{pgfonlayer}
\end{tikzpicture} =  \begin{tikzpicture}
	\begin{pgfonlayer}{nodelayer}
		\node [style=none] (0) at (2.25, -1.5) {$A$};
		\node [style=none] (1) at (2, 0.25) {};
		\node [style=none] (2) at (2, -3) {};
	\end{pgfonlayer}
	\begin{pgfonlayer}{edgelayer}
		\draw (1.center) to (2.center);
	\end{pgfonlayer}
\end{tikzpicture}
~~~~~~~~
\begin{tikzpicture}
	\begin{pgfonlayer}{nodelayer}
		\node [style=circle] (0) at (0, -1.75) {$\bot$};
		\node [style=none] (1) at (0.75, -3) {};
		\node [style=none] (2) at (0.75, 0.25) {};
		\node [style=none] (3) at (1, -1.5) {$A$};
		\node [style=circle, scale=0.4] (4) at (0.75, -0.5) {};
		\node [style=circle] (5) at (0, -2.75) {$\bot$};
	\end{pgfonlayer}
	\begin{pgfonlayer}{edgelayer}
		\draw (1.center) to (2.center);
		\draw [dotted, bend left, looseness=1.25] (0) to (4);
		\draw (5) to (0);
	\end{pgfonlayer}
\end{tikzpicture} = \begin{tikzpicture}
	\begin{pgfonlayer}{nodelayer}
		\node [style=none] (0) at (2.25, -1.5) {$A$};
		\node [style=none] (1) at (2, 0.25) {};
		\node [style=none] (2) at (2, -3) {};
	\end{pgfonlayer}
	\begin{pgfonlayer}{edgelayer}
		\draw (1.center) to (2.center);
	\end{pgfonlayer}
\end{tikzpicture}
~~~~~~~~
\begin{tikzpicture}
	\begin{pgfonlayer}{nodelayer}
		\node [style=ox] (0) at (2, -0) {};
		\node [style=ox] (1) at (2, -1) {};
		\node [style=none] (2) at (1.5, -2) {};
		\node [style=none] (3) at (2.5, -2) {};
		\node [style=none] (4) at (1.5, 1) {};
		\node [style=none] (5) at (2.5, 1) {};
		\node [style=none] (6) at (1.25, 0.75) {$A$};
		\node [style=none] (7) at (2.75, 0.75) {$B$};
		\node [style=none] (8) at (1.25, -1.75) {$A$};
		\node [style=none] (9) at (2.75, -1.75) {$B$};
	\end{pgfonlayer}
	\begin{pgfonlayer}{edgelayer}
		\draw (0) to (1);
		\draw [in=90, out=-135, looseness=1.00] (1) to (2.center);
		\draw [in=90, out=-45, looseness=1.00] (1) to (3.center);
		\draw [in=-90, out=135, looseness=1.00] (0) to (4.center);
		\draw [in=-90, out=45, looseness=1.00] (0) to (5.center);
	\end{pgfonlayer}
\end{tikzpicture} = \begin{tikzpicture}
	\begin{pgfonlayer}{nodelayer}
		\node [style=none] (0) at (1.5, -2) {};
		\node [style=none] (1) at (2, -2) {};
		\node [style=none] (2) at (1.5, 1) {};
		\node [style=none] (3) at (2, 1) {};
		\node [style=none] (4) at (1.25, 0.75) {$A$};
		\node [style=none] (5) at (2.25, 0.75) {$B$};
		\node [style=none] (6) at (1.25, -1.75) {$A$};
		\node [style=none] (7) at (2.25, -1.75) {$B$};
	\end{pgfonlayer}
	\begin{pgfonlayer}{edgelayer}
		\draw (2.center) to (0.center);
		\draw (3.center) to (1.center);
	\end{pgfonlayer}
\end{tikzpicture}
~~~~~~~~
\begin{tikzpicture}
	\begin{pgfonlayer}{nodelayer}
		\node [style=oa] (0) at (2, -0) {};
		\node [style=oa] (1) at (2, -1) {};
		\node [style=none] (2) at (1.5, -2) {};
		\node [style=none] (3) at (2.5, -2) {};
		\node [style=none] (4) at (1.5, 1) {};
		\node [style=none] (5) at (2.5, 1) {};
		\node [style=none] (6) at (1.25, 0.75) {$A$};
		\node [style=none] (7) at (2.75, 0.75) {$B$};
		\node [style=none] (8) at (1.25, -1.75) {$A$};
		\node [style=none] (9) at (2.75, -1.75) {$B$};
	\end{pgfonlayer}
	\begin{pgfonlayer}{edgelayer}
		\draw (0) to (1);
		\draw [in=90, out=-135, looseness=1.00] (1) to (2.center);
		\draw [in=90, out=-45, looseness=1.00] (1) to (3.center);
		\draw [in=-90, out=135, looseness=1.00] (0) to (4.center);
		\draw [in=-90, out=45, looseness=1.00] (0) to (5.center);
	\end{pgfonlayer}
\end{tikzpicture} = \begin{tikzpicture}
	\begin{pgfonlayer}{nodelayer}
		\node [style=none] (0) at (1.5, -2) {};
		\node [style=none] (1) at (2, -2) {};
		\node [style=none] (2) at (1.5, 1) {};
		\node [style=none] (3) at (2, 1) {};
		\node [style=none] (4) at (1.25, 0.75) {$A$};
		\node [style=none] (5) at (2.25, 0.75) {$B$};
		\node [style=none] (6) at (1.25, -1.75) {$A$};
		\node [style=none] (7) at (2.25, -1.75) {$B$};
	\end{pgfonlayer}
	\begin{pgfonlayer}{edgelayer}
		\draw (2.center) to (0.center);
		\draw (3.center) to (1.center);
	\end{pgfonlayer}
\end{tikzpicture} \]
The following are also circuit equalities (and when oriented become expansion rules):
\[ \begin{tikzpicture}
	\begin{pgfonlayer}{nodelayer}
		\node [style=ox] (0) at (0, -0.75) {};
		\node [style=ox] (1) at (0, -2) {};
		\node [style=none] (2) at (0, -3) {};
		\node [style=none] (3) at (0, -0) {};
	\end{pgfonlayer}
	\begin{pgfonlayer}{edgelayer}
		\draw [bend right=60, looseness=1.50] (0) to (1);
		\draw [bend right=60, looseness=1.50] (1) to (0);
		\draw (3.center) to (0);
		\draw (1) to (2.center);
	\end{pgfonlayer}
\end{tikzpicture} = \begin{tikzpicture}
	\begin{pgfonlayer}{nodelayer}
		\node [style=none] (0) at (1, -0) {};
		\node [style=none] (1) at (1, -3) {};
		\node [style=none] (2) at (1.65, -1.75) {$A \ox B$};
	\end{pgfonlayer}
	\begin{pgfonlayer}{edgelayer}
		\draw (0.center) to (1.center);
	\end{pgfonlayer}
\end{tikzpicture}
~~~~~~~~
\begin{tikzpicture}
	\begin{pgfonlayer}{nodelayer}
		\node [style=oa] (0) at (0, -0.75) {};
		\node [style=oa] (1) at (0, -2) {};
		\node [style=none] (2) at (0, -3) {};
		\node [style=none] (3) at (0, -0) {};
	\end{pgfonlayer}
	\begin{pgfonlayer}{edgelayer}
		\draw [bend right=60, looseness=1.50] (0) to (1);
		\draw [bend right=60, looseness=1.50] (1) to (0);
		\draw (3.center) to (0);
		\draw (1) to (2.center);
	\end{pgfonlayer}
\end{tikzpicture} = \begin{tikzpicture}
	\begin{pgfonlayer}{nodelayer}
		\node [style=none] (0) at (1, -0) {};
		\node [style=none] (1) at (1, -3) {};
		\node [style=none] (2) at (1.65, -1.75) {$A \oa B$};
	\end{pgfonlayer}
	\begin{pgfonlayer}{edgelayer}
		\draw (0.center) to (1.center);
	\end{pgfonlayer}
\end{tikzpicture}
~~~~~~~~~
\begin{tikzpicture}
	\begin{pgfonlayer}{nodelayer}
		\node [style=circle] (0) at (0, -0.75) {$\top$};
		\node [style=circle] (1) at (-1, -1.5) {$\top$};
		\node [style=none] (2) at (-1, -3.25) {};
		\node [style=none] (3) at (0, 0) {};
		\node [style=circle, scale=0.4] (4) at (-1, -2.5) {};
	\end{pgfonlayer}
	\begin{pgfonlayer}{edgelayer}
		\draw (3.center) to (0);
		\draw (1) to (2.center);
		\draw [dotted, in=-90, out=30, looseness=1.25] (4) to (0);
	\end{pgfonlayer}
\end{tikzpicture} = \begin{tikzpicture}
	\begin{pgfonlayer}{nodelayer}
		\node [style=none] (0) at (1, -0) {};
		\node [style=none] (1) at (1, -3.25) {};
		\node [style=none] (2) at (1.25, -1.75) {$\top$};
	\end{pgfonlayer}
	\begin{pgfonlayer}{edgelayer}
		\draw (0.center) to (1.center);
	\end{pgfonlayer}
\end{tikzpicture}
~~~~~~~~~
\begin{tikzpicture}
	\begin{pgfonlayer}{nodelayer}
		\node [style=circle] (0) at (0, -2.5) {$\bot$};
		\node [style=circle] (1) at (-1, -1.75) {$\bot$};
		\node [style=none] (2) at (-1, 0) {};
		\node [style=none] (3) at (0, -3.25) {};
		\node [style=circle, scale=0.4] (4) at (-1, -0.75) {};
	\end{pgfonlayer}
	\begin{pgfonlayer}{edgelayer}
		\draw (3.center) to (0);
		\draw (1) to (2.center);
		\draw [dotted, in=90, out=-30, looseness=1.25] (4) to (0);
	\end{pgfonlayer}
\end{tikzpicture} = \begin{tikzpicture}
	\begin{pgfonlayer}{nodelayer}
		\node [style=none] (0) at (1, -0) {};
		\node [style=none] (1) at (1, -3.25) {};
		\node [style=none] (2) at (1.25, -1.75) {$\bot$};
	\end{pgfonlayer}
	\begin{pgfonlayer}{edgelayer}
		\draw (0.center) to (1.center);
	\end{pgfonlayer}
\end{tikzpicture} \]

As in linear logic, not all circuit diagrams constructed from these basic components represent a valid LDC circuit.  In his seminal paper on linear logic, \cite{Gir87}, Girard introduced a criterion for the correctness of his representation of proofs using proof nets based on switching links.  A valid proof structure must be connected and acyclic for all the switching link choices.  Using this correctness criterion has the disadvantage of requiring exponential time in the number of switching links.  Danos and Regnier \cite{DaR89} improved this situation significantly by providing an algorithm for correctness which takes linear time (see \cite{Gue99})
%S. Guerrini
%Correctness of multiplicative proof nets is linear
%14th Annual IEEE Symposium on Logic in Computer Science, LICS 1999, IEEE Computer Society, Trento, Italy (1999), pp. 454-463 
on the size of the circuit.  To verify the validity of the circuit diagrams of LDCs, Blute et.al. \cite{BCST96}, provided a boxing algorithm which was based on Danos and Regnier's more efficient algorithm which we now describe.

In order to verify that an LDC circuit is valid, circuit components are ``boxed''  using the below rules. The primitive generating maps are 
automatically boxed. 
\[ (a_1)~~~ \]

Double lines refer to multiple number of wires. $\ox$-introduction and $\oa$-elimination are boxed in $(a_1)$ and $(a_2)$ respectively. In $(b_1)$, it is shown how a box `eats' the $\ox$-elimination: in $(b_2)$ the dual rule shows a $\oa$-introduction being eaten. $(c)$ shows how boxes can be amalgamated when they are connected by a single wire. In $(e_1)$-$(e_4)$, it is shown how the thinning links can be boxed. By progressively enclosing the components of the circuit in boxes using these rules, if we end up with a single box (or a wire), precisely when the circuit is valid. 

As an example, we verify the validity of the left linear distributor:
\[ \begin{tikzpicture}
	\begin{pgfonlayer}{nodelayer}
		\node [style=ox] (0) at (1.25, 0.75) {};
		\node [style=ox] (1) at (0.5, -1.25) {};
		\node [style=oa] (2) at (2, -0.25) {};
		\node [style=oa] (3) at (1.25, -2.5) {};
		\node [style=none] (4) at (1.25, -3.25) {};
		\node [style=none] (5) at (1.25, 1.5) {};
	\end{pgfonlayer}
	\begin{pgfonlayer}{edgelayer}
		\draw (5.center) to (0);
		\draw [bend left, looseness=1.00] (0) to (2);
		\draw [in=90, out=-150, looseness=1.00] (0) to (1);
		\draw (3) to (4.center);
		\draw [in=150, out=-75, looseness=1.00] (1) to (3);
		\draw [in=60, out=-120, looseness=1.25] (2) to (1);
		\draw [in=30, out=-75, looseness=1.00] (2) to (3);
	\end{pgfonlayer}
\end{tikzpicture} \stackrel{a_1,a_2}{\Rightarrow} \begin{tikzpicture}
	\begin{pgfonlayer}{nodelayer}
		\node [style=ox] (0) at (1.25, 0.75) {};
		\node [style=ox] (1) at (0.5, -1.25) {};
		\node [style=oa] (2) at (2, -0.25) {};
		\node [style=oa] (3) at (1.25, -2.5) {};
		\node [style=none] (4) at (1.25, -3.25) {};
		\node [style=none] (5) at (1.25, 1.5) {};
		\node [style=none] (6) at (1.5, 0.25) {};
		\node [style=none] (7) at (2.5, 0.25) {};
		\node [style=none] (8) at (2.5, -1) {};
		\node [style=none] (9) at (1.5, -1) {};
		\node [style=none] (10) at (1.25, -0.5) {};
		\node [style=none] (11) at (0, -0.5) {};
		\node [style=none] (12) at (0, -1.75) {};
		\node [style=none] (13) at (1.25, -1.75) {};
	\end{pgfonlayer}
	\begin{pgfonlayer}{edgelayer}
		\draw (5.center) to (0);
		\draw [bend left, looseness=1.00] (0) to (2);
		\draw [in=90, out=-150, looseness=1.00] (0) to (1);
		\draw (3) to (4.center);
		\draw [in=150, out=-75, looseness=1.00] (1) to (3);
		\draw [in=60, out=-120, looseness=1.25] (2) to (1);
		\draw [in=30, out=-75, looseness=1.00] (2) to (3);
		\draw (9.center) to (8.center);
		\draw (8.center) to (7.center);
		\draw (7.center) to (6.center);
		\draw (6.center) to (9.center);
		\draw (11.center) to (12.center);
		\draw (12.center) to (13.center);
		\draw (13.center) to (10.center);
		\draw (10.center) to (11.center);
	\end{pgfonlayer}
\end{tikzpicture} \stackrel{c}{\Rightarrow} \begin{tikzpicture}
	\begin{pgfonlayer}{nodelayer}
		\node [style=ox] (0) at (1.25, 0.75) {};
		\node [style=ox] (1) at (0.5, -1.25) {};
		\node [style=oa] (2) at (2, -0.25) {};
		\node [style=oa] (3) at (1.25, -2.5) {};
		\node [style=none] (4) at (1.25, -3.25) {};
		\node [style=none] (5) at (1.25, 1.5) {};
		\node [style=none] (6) at (2.5, 0.25) {};
		\node [style=none] (7) at (2.5, -1.75) {};
		\node [style=none] (8) at (0, 0.25) {};
		\node [style=none] (9) at (0, -1.75) {};
	\end{pgfonlayer}
	\begin{pgfonlayer}{edgelayer}
		\draw (5.center) to (0);
		\draw [bend left, looseness=1.00] (0) to (2);
		\draw [in=90, out=-150, looseness=1.00] (0) to (1);
		\draw (3) to (4.center);
		\draw [in=150, out=-75, looseness=1.00] (1) to (3);
		\draw [in=60, out=-120, looseness=1.25] (2) to (1);
		\draw [in=30, out=-75, looseness=1.00] (2) to (3);
		\draw (7.center) to (6.center);
		\draw (8.center) to (9.center);
		\draw (8.center) to (6.center);
		\draw (9.center) to (7.center);
	\end{pgfonlayer}
\end{tikzpicture} \stackrel{b_1,b_2}{\Rightarrow} \begin{tikzpicture}
	\begin{pgfonlayer}{nodelayer}
		\node [style=ox] (0) at (1.25, 0.75) {};
		\node [style=ox] (1) at (0.5, -1.25) {};
		\node [style=oa] (2) at (2, -0.25) {};
		\node [style=oa] (3) at (1.25, -2.5) {};
		\node [style=none] (4) at (1.25, -3.25) {};
		\node [style=none] (5) at (1.25, 1.5) {};
		\node [style=none] (6) at (2.5, 1.25) {};
		\node [style=none] (7) at (2.5, -3) {};
		\node [style=none] (8) at (0, 1.25) {};
		\node [style=none] (9) at (0, -3) {};
	\end{pgfonlayer}
	\begin{pgfonlayer}{edgelayer}
		\draw (5.center) to (0);
		\draw [bend left, looseness=1.00] (0) to (2);
		\draw [in=90, out=-150, looseness=1.00] (0) to (1);
		\draw (3) to (4.center);
		\draw [in=150, out=-75, looseness=1.00] (1) to (3);
		\draw [in=60, out=-120, looseness=1.25] (2) to (1);
		\draw [in=30, out=-75, looseness=1.00] (2) to (3);
		\draw (7.center) to (6.center);
		\draw (8.center) to (9.center);
		\draw (8.center) to (6.center);
		\draw (9.center) to (7.center);
	\end{pgfonlayer}
\end{tikzpicture} \]
% Examples here
In the firts step the $\ox$-introduction and $\oa$-elimination are boxed. In the second step the boxes are amalgamated along the single wire joining them. In the third step,  the box absorbs the $\ox$-elimination and $\oa$-introduction.

In contrast, we now show that the reverse of the linear distributor is invalid as the boxing process gets stuck:
\[ \begin{tikzpicture}
	\begin{pgfonlayer}{nodelayer}
		\node [style=ox] (0) at (1.25, -2.5) {};
		\node [style=ox] (1) at (2, -0.25) {};
		\node [style=oa] (2) at (0.5, -1.25) {};
		\node [style=oa] (3) at (1.25, 0.75) {};
		\node [style=none] (4) at (1.25, 1.5) {};
		\node [style=none] (5) at (1.25, -3.25) {};
	\end{pgfonlayer}
	\begin{pgfonlayer}{edgelayer}
		\draw (5.center) to (0);
		\draw [in=-83, out=157, looseness=1.00] (0) to (2);
		\draw [in=-90, out=30, looseness=1.00] (0) to (1);
		\draw (3) to (4.center);
		\draw [in=-15, out=90, looseness=1.00] (1) to (3);
		\draw [in=-120, out=60, looseness=1.25] (2) to (1);
		\draw [in=-150, out=90, looseness=1.00] (2) to (3);
	\end{pgfonlayer}
\end{tikzpicture} \stackrel{a_1,a_2}{\Rightarrow} \begin{tikzpicture}
	\begin{pgfonlayer}{nodelayer}
		\node [style=ox] (0) at (1.25, -2.5) {};
		\node [style=ox] (1) at (2, -0.5) {};
		\node [style=oa] (2) at (0.5, -1.5) {};
		\node [style=oa] (3) at (1.25, 0.75) {};
		\node [style=none] (4) at (1.25, 1.5) {};
		\node [style=none] (5) at (1.25, -3.25) {};
		\node [style=none] (6) at (0.5, 1.25) {};
		\node [style=none] (7) at (2, 1.25) {};
		\node [style=none] (8) at (2, 0.25) {};
		\node [style=none] (9) at (0.5, 0.25) {};
		\node [style=none] (10) at (0.5, -2) {};
		\node [style=none] (11) at (2, -2) {};
		\node [style=none] (12) at (2, -3) {};
		\node [style=none] (13) at (0.5, -3) {};
	\end{pgfonlayer}
	\begin{pgfonlayer}{edgelayer}
		\draw (5.center) to (0);
		\draw [in=-83, out=157, looseness=1.00] (0) to (2);
		\draw [in=-90, out=30, looseness=1.00] (0) to (1);
		\draw (3) to (4.center);
		\draw [in=-15, out=90, looseness=1.00] (1) to (3);
		\draw [in=-120, out=60, looseness=1.25] (2) to (1);
		\draw [in=-150, out=90, looseness=1.00] (2) to (3);
		\draw (9.center) to (6.center);
		\draw (6.center) to (7.center);
		\draw (7.center) to (8.center);
		\draw (9.center) to (8.center);
		\draw (12.center) to (13.center);
		\draw (13.center) to (10.center);
		\draw (10.center) to (11.center);
		\draw (11.center) to (12.center);
	\end{pgfonlayer}
\end{tikzpicture}  \]

Before presenting the definition of linear functors, we briefly recall the definition of monoidal functors.
A functor $F: \X \to \Y$ between monoidal categories is a monoidal functor if it is equipped with natural transformations $m_\ox: F(A) \ox F(B) \to F(A \ox B)$ and $m_I: I \to F(I)$ such that the following diagrams commute:

\[
\xymatrix{
(F(A) \ox F(B)) \ox F(C) \ar[d]_{a_\ox} \ar[r]^{m_\ox \ox 1} & F(A \ox B) \ox F(C) \ar[r]^{m_\ox} & F((A \ox B) \ox C) \ar[d]^{F(a_\ox)} \\
F(A) \ox (F(B) \ox F(C)) \ar[r]_{1 \ox m_\ox} & F(A) \ox F(B \ox C) \ar[r]_{m_\ox} & F(A \ox (B \ox C)) 
} 
\]
\[
\xymatrix{
F(A) \ox I \ar[drr]^{u_\ox^R} \ar[d]_{1 \ox m_I} \\
F(A) \ox F(I) \ar[r]_{m_\ox} & F(A \ox I) \ar[r]_{F(u_\ox^L)} & F(A)}  ~~~~~~~~~~~~~~
\xymatrix{
I \ox F(A) \ar[drr]^{u_\ox^L} \ar[d]_{m_I \ox 1} \\
F(I) \ox F(A) \ar[r]_{m_\ox} & F(I \ox A) \ar[r]_{F(u_\ox^L)} & F(A)
}
\]

The first diagram for the monoidal functor is the associative law, and the other two diagrams are the right and the left unit laws respectively.

\begin{defi} \cite[Definition 1]{CS99}
Given linearly distributive categories $\X$ and $\Y$, a linear functor $F: \X \to \Y$ consists of 
\begin{enumerate}[(i)]
\item a pair of functors $F = (F_\ox, F_\oa)$: $(F_\ox, m_\ox, m_\top)$ which is monoidal with respect to $\ox$ and $(F_\oa, n_\oa, n_\bot)$ which is comonoidal with respect to $\oa$. We refer to $m_\ox$ and $n_\oa$ as {\bf tensor laxors}, and $m_\top$ and $n_\bot$ as {\bf unit laxors}. 

\item natural transformations:
\begin{align*}
\nu_\ox^R &: F_\ox(A \oa B) \to F_\oa(A) \oa F_\ox(B) \\
\nu_\ox^L &: F_\ox(A \oa B) \to F_\ox(A) \oa F_\oa(B) \\
\nu_\oa^R &: F_\ox(A) \ox F_\oa(B) \to F_\oa(A \ox B) \\
\nu_\oa^L &: F_\oa(A) \ox F_\ox(B) \to F_\oa( A \ox B)
\end{align*}
\end{enumerate}

such that the following coherence conditions hold:
\begin{enumerate}[{\bf \small [LF.1]}]
\item 
\begin{enumerate}[(a)]
\item $F_\ox(u^L_\oa) = \nu^R_\ox (n_\bot \oa 1) u^L_\oa$
\[ \xymatrix{
F_\ox( \bot \oa A) \ar[r]^{F_\ox(u_\oa^L)} \ar[d]_{\nu_\ox^R} & F_\ox(A) \\
F_\oa(\bot) \oa F_\ox(A) \ar[r]_{n_\bot \oa 1} & \bot \oa F_\ox(A) \ar[u]_{u_\oa^L}
} \]
\item $\nu_\ox^L ( 1 \oa n_\bot ) u_\oa^R = F_\ox( u_\oa^R) $
\item $(u_\ox^L)^{-1} (m_\top \ox 1) \nu_\oa^R = F_\oa((u_\ox^L)^{-1})$
\item $(u_\ox^R)^{-1} (m_\top \ox 1) \nu_\oa^L = F_\oa((u_\ox^R)^{-1}) $
\end{enumerate} 
\item 
\begin{enumerate}[(a)]
\item $F_\ox(a_\oa) \nu^R_\ox (1 \oa \nu^R_\ox) = \nu^R_\ox (n_\oa \oa 1) a_\oa$
\[ \xymatrix{
F_\ox((A \oa B) \oa C) \ar[r]^{F_\ox(a_\oa)} \ar[r]^{F_\ox(a_\oa)} \ar[d]_{\nu_\ox^R} & 
F_\ox(A \oa (B \oa C)) \ar[d]^{\nu_\ox^R} \\
F_\oa(A \oa B) \oa F_\ox(C) \ar[d]_{n_\oa \oa 1} &
F_\oa(A) \oa F_\ox (B \oa C) \ar[d]^{1 \oa \nu_\ox^R} \\
(F_\oa (A) \oa F_\oa(B)) \oa F_\ox(C) \ar[r]_{a_\oa} &
F_\oa(A) \oa (F_\oa(B) \oa F_\oa(C))
} \]
\item $F_\ox(a_\oa) \nu_\ox^L (1 \oa n_\oa ) = \nu^L_\oa (\nu^L \oa 1 ) a_\oa$
\item $(m_\ox \ox 1) \nu_\oa^R F_\oa(a_\ox) = a_\ox (1 \ox \nu_\oa^R) \nu_\oa^R$
\item $(\nu^R_\oa \ox 1) \nu_\oa^L F_\oa(a_\ox) = a_\ox (1 \ox m_\ox) \nu_\oa^L$
\end{enumerate}
\item 
\begin{enumerate}[(a)]
\item $F_\ox(a_\oa)\nu^R_\ox(1 \oa \nu^L_\ox) = \nu_\ox^L (\nu^R_\ox \oa 1) a_\oa$
\[ \xymatrix{
F_\ox((A \oa B) \oa C) \ar[r]^{F_\ox(a_\oa)} \ar[d]_{\nu_\ox^L} & 
F_\ox(A \oa (B \oa C))  \ar[d]^{\nu_\ox^R} \\
F_\ox(A \oa B) \oa F_\oa(C) \ar[d]_{\nu_\ox^R \oa 1} &
F_\oa(A) \oa F_\ox (B \oa C) \ar[d]^{1 \oa \nu_\ox^L} \\
(F_\oa (A) \oa F_\ox(B)) \oa F_\oa(C) \ar[r]_{a_\oa}&
F_\oa(A) \oa (F_\ox(B) \oa F_\oa(C))}
\]
\item $(\nu^R_\oa \ox 1) \nu^L_\oa F_\oa(a_\ox) = a_\ox (1 \ox \nu_\oa^L) \nu_\oa^R$
\end{enumerate}
\item 
\begin{enumerate}[(a)]
\item $(1 \ox \nu^R_\ox) \partial^L (\nu^R_\oa \oa 1) = m_\ox F_\ox(\partial^L) \nu^R_\ox$
\[ \xymatrix{
F_\ox(A) \ox F_\ox(B \oa C) \ar[r]^{1 \ox \nu_\ox^R} \ar[d]_{m_\ox} &
F_\ox(A) \ox (F_\oa(B) \oa F_\ox(C)) \ar[d]^{\partial^L} \\
F_\ox(A \ox (B \oa C)) \ar[d]_{F_\ox(\partial^L)} &
(F_\ox(A) \ox F_\oa(B)) \oa F_\ox(C) \ar[d]^{\nu_\oa^R \oa 1} \\
F_\ox((A \ox B) \oa C) \ar[r]_{\nu_\ox^R} &
F_\oa(A \oa B) \oa F_\ox(C)
} \]
\item $(\nu_\ox^L \ox 1) \partial^R (1 \oa \nu_\oa^L) = m_\ox F_\ox(\partial^R) \nu_\ox^L$
\item $(1 \ox \nu_\ox^L) \partial^L (\nu_\oa^L \oa 1) = \nu_\oa^L F_\oa(\partial^L) n_\oa$
\item $(\nu_\ox^R \ox 1) \partial^R (1 \oa \nu_\oa^R) = \nu_\oa^R F_\oa(\partial^R) n_\oa $
\end{enumerate}
\item 
\begin{enumerate}[(a)]
\item $(1 \ox \nu^L_\ox) \partial^L(m_\ox \oa 1) = m_\ox F_\ox(\partial^L) \nu^L_\ox$ 
\[ \xymatrix{
F_\ox(A) \ox F_\ox(B \oa C) \ar[r]^{1 \ox \nu_\ox^L} \ar[d]_{m_\ox} &
F_\ox(A) \ox (F_\ox(B) \oa F_\oa(C)) \ar[d]^{\partial^L} \\
F_\ox(A \ox (B \oa C)) \ar[d]_{F_\ox(\partial^L)} &
(F_\ox(A) \ox F_\ox(B)) \oa F_\oa(C) \ar[d]^{m_\ox \oa 1} \\
F_\ox((A \ox B) \oa C) \ar[r]_{\nu_\ox^L} &
F_\ox(A \ox B) \oa F_\oa(C)
}\]
\item $(\nu_\ox^R \ox 1) \partial^R (1 \oa m_\ox) = m_\ox F_\ox(\partial^R) \nu_\ox^R$
\item $(1 \ox n_\oa) \partial^L (\nu_\oa^R \oa 1) = \nu_\oa^R F_\oa(\partial^L) n_\oa $
\item $(n_\oa \ox 1) \partial^R (1 \oa \nu_\oa^L) = \nu_\oa^L F_\oa(\partial^R) n_ \oa $
\end{enumerate}
\end{enumerate}
\end{defi}

In the graphical calculus, functors are represented by linear functor boxes \cite{CS99}. A linear functor box can be monoidal or comonoidal. When the functor box is monoidal $(F_\ox)$, it has one principal output wire (represented by a port where the wire exits the box) and the other wires are auxiliary. When the box is comonoidal ($F_\oa$), it has one principal input wire with a port and the other wires are auxiliary. The functor boxes are subject to a very natural ``box eats box'' calculus described in \cite{CS99}.  A box can eat another box only when a ported wire meets an auxiliary wire. 

The linear strengths are drawn in the graphical calculus as follows:
\[
\nu_\oa^L = \begin{tikzpicture}
	\begin{pgfonlayer}{nodelayer}
		\node [style=none] (0) at (-2, 2) {};
		\node [style=none] (1) at (-1.5, 2) {};
		\node [style=none] (2) at (-2.25, 2) {};
		\node [style=none] (3) at (-2.25, 1) {};
		\node [style=none] (4) at (-0.75, 1) {};
		\node [style=none] (5) at (-0.75, 2) {};
		\node [style=none] (6) at (-2, 2.75) {};
		\node [style=none] (7) at (-1, 2.75) {};
		\node [style=none] (61) at (-2.75, 2.75) {$F_\oa(A)$};
		\node [style=none] (71) at (-0.25, 2.75) {$F_\ox(B)$};
		\node [style=none] (8) at (-1.5, 0.25) {};
		\node [style=none] (81) at (-2.25, 0) {$F_\oa(A \ox B) $};
		\node [style=ox] (9) at (-1.5, 1.5) {};
		\node [style=none] (10) at (-2, 1.25) {$F$};
	\end{pgfonlayer}
	\begin{pgfonlayer}{edgelayer}
		\draw [in=-90, out=-90, looseness=1.25] (0.center) to (1.center);
		\draw [bend right=15, looseness=1.00] (6.center) to (9);
		\draw [bend right=15, looseness=0.75] (9) to (7.center);
		\draw (2.center) to (3.center);
		\draw (3.center) to (4.center);
		\draw (4.center) to (5.center);
		\draw (5.center) to (2.center);
		\draw (9) to (8.center);
	\end{pgfonlayer}
\end{tikzpicture} ~~~~~~~~ \nu_\oa^R = \begin{tikzpicture}
	\begin{pgfonlayer}{nodelayer}
		\node [style=none] (0) at (-1.5, 2) {};
		\node [style=none] (1) at (-1, 2) {};
		\node [style=none] (2) at (-2.25, 2) {};
		\node [style=none] (3) at (-2.25, 1) {};
		\node [style=none] (4) at (-0.75, 1) {};
		\node [style=none] (5) at (-0.75, 2) {};
		\node [style=none] (6) at (-2, 2.75) {};
		\node [style=none] (7) at (-1, 2.75) {};
		\node [style=none] (61) at (-2.75, 2.75) {$F_\ox(A)$};
		\node [style=none] (71) at (-0.25, 2.75) {$F_\oa(B)$};
		\node [style=none] (8) at (-1.5, 0.25) {};
		\node [style=none] (81) at (-2.25, 0) {$F_\oa(A \ox B) $};
		\node [style=ox] (9) at (-1.5, 1.5) {};
		\node [style=none] (10) at (-2, 1.25) {$F$};
	\end{pgfonlayer}
	\begin{pgfonlayer}{edgelayer}
		\draw [in=-90, out=-90, looseness=1.25] (0.center) to (1.center);
		\draw [bend right=15, looseness=1.00] (6.center) to (9);
		\draw [bend right=15, looseness=0.75] (9) to (7.center);
		\draw (2.center) to (3.center);
		\draw (3.center) to (4.center);
		\draw (4.center) to (5.center);
		\draw (5.center) to (2.center);
		\draw (9) to (8.center);
	\end{pgfonlayer}
\end{tikzpicture}  ~~~~~~~~ \nu_\ox^L = \begin{tikzpicture}
	\begin{pgfonlayer}{nodelayer}
		\node [style=none] (0) at (-2.25, 1) {};
		\node [style=none] (1) at (-2.25, 2) {};
		\node [style=none] (2) at (-0.75, 2) {};
		\node [style=none] (3) at (-0.75, 1) {};
		\node [style=none] (4) at (-2, 0.25) {};
		\node [style=none] (5) at (-1, 0.25) {};
		\node [style=none] (6) at (-1.5, 2.75) {};
		\node [style=none] (41) at (-2.75, 0.25) {$F_\ox(A)$};
		\node [style=none] (51) at (-0.25, 0.25) {$F_\oa(B)$};
		\node [style=none] (61) at (-2, 3) {$F_\ox(A \oa B)$};
		\node [style=oa] (7) at (-1.5, 1.5) {};
		\node [style=none] (8) at (-2, 1.75) {$F$};
		\node [style=none] (9) at (-2, 1) {};
		\node [style=none] (10) at (-1.5, 1) {};
	\end{pgfonlayer}
	\begin{pgfonlayer}{edgelayer}
		\draw [bend left=15, looseness=1.00] (4.center) to (7);
		\draw [bend left=15, looseness=0.75] (7) to (5.center);
		\draw (0.center) to (1.center);
		\draw (1.center) to (2.center);
		\draw (2.center) to (3.center);
		\draw (3.center) to (0.center);
		\draw (7) to (6.center);
		\draw [in=90, out=90, looseness=1.25] (9.center) to (10.center);
	\end{pgfonlayer}
\end{tikzpicture} ~~~~~~~~~ \nu_\ox^R = \begin{tikzpicture}
	\begin{pgfonlayer}{nodelayer}
		\node [style=none] (0) at (-2.25, 1) {};
		\node [style=none] (1) at (-2.25, 2) {};
		\node [style=none] (2) at (-0.75, 2) {};
		\node [style=none] (3) at (-0.75, 1) {};
		\node [style=none] (4) at (-2, 0.25) {};
		\node [style=none] (5) at (-1, 0.25) {};
		\node [style=none] (41) at (-2.75, 0.25) {$F_\oa(A)$};
		\node [style=none] (51) at (-0.25, 0.25) {$F_\ox(B)$};
		\node [style=none] (61) at (-2, 3) {$F_\ox(A \oa B)$};
		\node [style=none] (6) at (-1.5, 2.75) {};
		\node [style=oa] (7) at (-1.5, 1.5) {};
		\node [style=none] (8) at (-2, 1.75) {$F$};
		\node [style=none] (9) at (-1.5, 1) {};
		\node [style=none] (10) at (-1, 1) {};
	\end{pgfonlayer}
	\begin{pgfonlayer}{edgelayer}
		\draw [bend left=15, looseness=1.00] (4.center) to (7);
		\draw [bend left=15, looseness=0.75] (7) to (5.center);
		\draw (0.center) to (1.center);
		\draw (1.center) to (2.center);
		\draw (2.center) to (3.center);
		\draw (3.center) to (0.center);
		\draw (7) to (6.center);
		\draw [in=90, out=90, looseness=1.25] (9.center) to (10.center);
	\end{pgfonlayer}
\end{tikzpicture} 
\]
\[
m_\ox = \begin{tikzpicture}
	\begin{pgfonlayer}{nodelayer}
		\node [style=none] (0) at (-2.25, 2) {};
		\node [style=none] (1) at (-2.25, 1) {};
		\node [style=none] (2) at (-0.75, 1) {};
		\node [style=none] (3) at (-0.75, 2) {};
		\node [style=none] (4) at (-2, 2.75) {};
		\node [style=none] (5) at (-1, 2.75) {};
		\node [style=none] (41) at (-2.75, 2.75) {$F_\ox(A)$};
		\node [style=none] (51) at (-0.25, 2.75) {$F_\ox(B)$};
		\node [style=none] (6) at (-1.5, 0.25) {};
		\node [style=none] (61) at (-2, 0) {$F_\ox(A \ox B)$};
		\node [style=ox] (7) at (-1.5, 1.5) {};
		\node [style=none] (8) at (-2, 1.25) {$F$};
		\node [style=none] (9) at (-1.75, 1) {};
		\node [style=none] (10) at (-1.25, 1) {};
	\end{pgfonlayer}
	\begin{pgfonlayer}{edgelayer}
		\draw [bend right=15, looseness=1.00] (4.center) to (7);
		\draw [bend right=15, looseness=0.75] (7) to (5.center);
		\draw (0.center) to (1.center);
		\draw (1.center) to (2.center);
		\draw (2.center) to (3.center);
		\draw (3.center) to (0.center);
		\draw (7) to (6.center);
		\draw [in=90, out=90, looseness=1.25] (9.center) to (10.center);
	\end{pgfonlayer}
\end{tikzpicture} ~~~~~ m_\top =  \begin{tikzpicture}
	\begin{pgfonlayer}{nodelayer}
	      \node [style=none] (8) at (0, 1) {};
		\node [style=circle] (0) at (0, -0) {$\top$};
		\node [style=none] (1) at (0.75, -0.5) {};
		\node [style=none] (2) at (2.75, -0.5) {};
		\node [style=none] (3) at (0.75, -1.5) {};
		\node [style=none] (4) at (2.75, -1.5) {};
		\node [style=none] (5) at (1.75, -2.5) {};
		\node [style=circle] (6) at (1.75, -1) {$\top$};
		\node [style=circle, scale=0.5] (7) at (1.75, -2) {};
	\end{pgfonlayer}
	\begin{pgfonlayer}{edgelayer}
		\draw (1.center) to (2.center);
		\draw (2.center) to (4.center);
		\draw (4.center) to (3.center);
		\draw (3.center) to (1.center);
		\draw (6) to (5.center);
		\draw [dotted, in=-165, out=-90, looseness=1.25] (0) to (7);
		\draw (0) to (8.center);
	\end{pgfonlayer}
\end{tikzpicture} = \begin{tikzpicture}
	\begin{pgfonlayer}{nodelayer}
	      \node [style=none] (8) at (2.75, 1) {};
		\node [style=circle] (0) at (2.75, -0) {$\top$};
		\node [style=none] (1) at (2, -0.5) {};
		\node [style=none] (2) at (0, -0.5) {};
		\node [style=none] (3) at (2, -1.5) {};
		\node [style=none] (4) at (0, -1.5) {};
		\node [style=none] (5) at (1, -2.5) {};
		\node [style=circle] (6) at (1, -1) {$\top$};
		\node [style=circle, scale=0.5] (7) at (1, -2) {};
	\end{pgfonlayer}
	\begin{pgfonlayer}{edgelayer}
		\draw (1.center) to (2.center);
		\draw (2.center) to (4.center);
		\draw (4.center) to (3.center);
		\draw (3.center) to (1.center);
		\draw (6) to (5.center);
		\draw (0) to (8.center);
		\draw [dotted, in=-15, out=-90, looseness=1.25] (0) to (7);
	\end{pgfonlayer}
\end{tikzpicture} \] \[ n_\oa =  \begin{tikzpicture}
	\begin{pgfonlayer}{nodelayer}
		\node [style=none] (0) at (-2.25, 1) {};
		\node [style=none] (1) at (-2.25, 2) {};
		\node [style=none] (2) at (-0.75, 2) {};
		\node [style=none] (3) at (-0.75, 1) {};
		\node [style=none] (4) at (-2, 0.25) {};
		\node [style=none] (5) at (-1, 0.25) {};
		\node [style=none] (41) at (-2.75, 0.25) {$F_\oa(A)$};
		\node [style=none] (51) at (-0.25, 0.25) {$F_\oa(B)$};
		\node [style=none] (61) at (-2, 3) {$F_\oa(A \oa B)$};
		\node [style=none] (6) at (-1.5, 2.75) {};
		\node [style=oa] (7) at (-1.5, 1.5) {};
		\node [style=none] (8) at (-2, 1.75) {$F$};
		\node [style=none] (9) at (-1.25, 2) {};
		\node [style=none] (10) at (-1.75, 2) {};
	\end{pgfonlayer}
	\begin{pgfonlayer}{edgelayer}
		\draw [bend left=15, looseness=1.00] (4.center) to (7);
		\draw [bend left=15, looseness=0.75] (7) to (5.center);
		\draw (0.center) to (1.center);
		\draw (1.center) to (2.center);
		\draw (2.center) to (3.center);
		\draw (3.center) to (0.center);
		\draw (7) to (6.center);
		\draw [in=-90, out=-90, looseness=1.25] (9.center) to (10.center);
	\end{pgfonlayer}
\end{tikzpicture}
~~~~~ n_\bot = 
\begin{tikzpicture}
	\begin{pgfonlayer}{nodelayer}
		\node [style=none] (0) at (-3, 2) {};
		\node [style=none] (1) at (-3, 1) {};
		\node [style=none] (2) at (-1, 2) {};
		\node [style=none] (3) at (-1, 1) {};
		\node [style=circle] (4) at (-2, 1.5) {$\bot$};
		\node [style=circle] (5) at (0, 1) {$\bot$};
		\node [style=none] (6) at (-2, 3) {};
		\node [style=circle, scale=0.5] (7) at (-2, 2.5) {};
		\node [style=none] (8) at (0, 0) {};
	\end{pgfonlayer}
	\begin{pgfonlayer}{edgelayer}
		\draw (0.center) to (1.center);
		\draw (1.center) to (3.center);
		\draw (3.center) to (2.center);
		\draw (2.center) to (0.center);
		\draw (6.center) to (4);
		\draw (8.center) to (5);
		\draw [dotted, bend left=45, looseness=1.25] (7) to (5);
	\end{pgfonlayer}
\end{tikzpicture} = \begin{tikzpicture}
	\begin{pgfonlayer}{nodelayer}
		\node [style=none] (0) at (0, 2) {};
		\node [style=none] (1) at (0, 1) {};
		\node [style=none] (2) at (-2, 2) {};
		\node [style=none] (3) at (-2, 1) {};
		\node [style=circle] (4) at (-1, 1.5) {$\bot$};
		\node [style=circle] (5) at (-3, 1) {$\bot$};
		\node [style=none] (6) at (-1, 3) {};
		\node [style=circle, scale=0.5] (7) at (-1, 2.5) {};
		\node [style=none] (8) at (-3, 0) {};
	\end{pgfonlayer}
	\begin{pgfonlayer}{edgelayer}
		\draw (0.center) to (1.center);
		\draw (1.center) to (3.center);
		\draw (3.center) to (2.center);
		\draw (2.center) to (0.center);
		\draw (6.center) to (4);
		\draw (8.center) to (5);
		\draw [dotted, bend right=45, looseness=1.25] (7) to (5);
	\end{pgfonlayer}
\end{tikzpicture} \]

When working in the categorical doctrine of {\em symmetric} LDCs we will expect the linear functors to preserve the symmetry.   Thus, a {\bf symmetric linear functor} is a linear functor $F= (F_\ox,F_\oa)$ 
which satisfies in addition:
\[ \xymatrix{F_\ox(A) \ox F_\ox(B) \ar[d]_{c_\ox} \ar[rr]^{m_\ox} & & F_\ox(A \ox B) \ar[d]^{F_\ox(c_\ox)} \\
                   F_\ox(B) \ox F_\ox(A) \ar[rr]_{m_\ox} & & F_\ox(B \ox A) }
    ~~~~~~
   \xymatrix{F_\oa(A \oa B)  \ar[d]_{F_\oa(c_\oa)} \ar[rr]^{n_\ox} & & F_\oa(A) \oa F_\oa(B)\ar[d]^{c_\oa} \\
                   F_\oa(B \ox A) \ar[rr]_{n_\oa} & &  F_\oa(B) \oa F_\oa(A) } \]
                   
Natural transformations between linear functors also break into two components linking respectively the tensor functors and, in the opposite  direction, the par functors:

\begin{defi} \cite[Definition 3]{CS99}
A {\bf linear (natural) transformation},  $\alpha: F \to G$,  between parallel linear functors $F,G: \X \to \Y$  consists of a pair of natural transformations $\alpha = (\alpha_\ox,\alpha_\oa)$ such that $\alpha_\ox: F_\ox \to G_\ox$ is a monoidal transformation and $\alpha_\oa: G_\oa \to F_\oa$ is a comonoidal transformation satisfying the following coherence conditions:
\begin{enumerate}[{\bf \small [LT.1]}]
\item $a_\ox \nu^R_\ox (a_\oa \oa 1) = \nu^R_\ox (1 \oa a_\ox)$
\[ 
\xymatrix{
F_\ox (A \oa B) \ar[rr]^{\alpha_\ox} \ar[d]_{\nu_\ox^R} & & G_\ox(A \oa B) \ar[d]^{\nu_\ox^R} \\
F_\oa(A) \oa F_\ox(B) \ar[dr]_{1 \oa \alpha_\ox} & & G_\oa(A) \oa G_\ox(B) \ar[ld]^{\alpha_\oa \oa 1} \\
& F_\oa(A) \oa G_\ox(B)&
}
\]
\item $\alpha_\ox \nu_\ox^L (1 \oa \alpha_\oa) = \nu_\ox^L (\alpha_\ox \oa 1)$
\item $(1 \ox \alpha_\ox) \nu_\oa^L (\alpha_\oa) = (\alpha_\oa \ox 1) \nu_\oa^L$
\item $(\alpha_\ox \ox 1) \nu_\oa^R \alpha_\oa = (1 \ox \alpha_\oa) \nu_\oa^R$
\end{enumerate}
\end{defi}

Conditions {\bf [LT.1]} - {\bf [LT.4]} are represented graphically as follows:
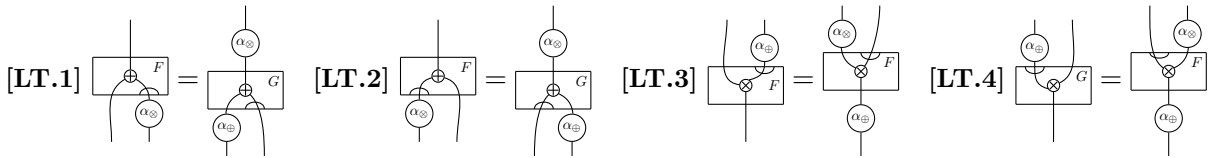
\[
\mbox{\small\bf [LT.1]}~ \begin{tikzpicture} %nat-a
	\begin{pgfonlayer}{nodelayer}
		\node [style=none] (0) at (-0.5, -0.75) {};
		\node [style=none] (1) at (-0.5, 0.25) {};
		\node [style=none] (2) at (-2.5, 0.25) {};
		\node [style=none] (3) at (-2.5, -0.75) {};
		\node [style=oa] (4) at (-1.5, -0.25) {};
		\node [style=none] (5) at (-1, -2) {};
		\node [style=none] (6) at (-2, -2) {};
		\node [style=none] (7) at (-1.5, 1.25) {};
		\node [style=circle, scale=2] (8) at (-1, -1.25) {};
		\node [style=none] (9) at (-1, -1.25) {$\alpha_\ox$};
		\node [style=none] (10) at (-0.75, -0) {$F$};
		\node [style=none] (11) at (-1.25, -0.75) {};
		\node [style=none] (12) at (-0.75, -0.75) {};
	\end{pgfonlayer}
	\begin{pgfonlayer}{edgelayer}
		\draw (0.center) to (1.center);
		\draw (1.center) to (2.center);
		\draw (2.center) to (3.center);
		\draw (0.center) to (3.center);
		\draw [in=90, out=-150, looseness=1.00] (4) to (6.center);
		\draw (4) to (7.center);
		\draw (5.center) to (8);
		\draw [in=-30, out=90, looseness=1.00] (8) to (4);
		\draw [bend left=90, looseness=1.25] (11.center) to (12.center);
	\end{pgfonlayer}
\end{tikzpicture} 
 = \begin{tikzpicture} %nat-b
	\begin{pgfonlayer}{nodelayer}
		\node [style=none] (0) at (-0.5, -0.75) {};
		\node [style=none] (1) at (-0.5, 0.25) {};
		\node [style=none] (2) at (-2.5, 0.25) {};
		\node [style=none] (3) at (-2.5, -0.75) {};
		\node [style=oa] (4) at (-1.5, -0.25) {};
		\node [style=none] (5) at (-1, -2) {};
		\node [style=none] (6) at (-2, -2) {};
		\node [style=none] (7) at (-2, -1.25) {$\alpha_\oa$};
		\node [style=circle, scale=2] (8) at (-1.5, 1) {};
		\node [style=none] (9) at (-1.5, 2) {};
		\node [style=none] (10) at (-1.5, 1) {$\alpha_\ox$};
		\node [style=none] (11) at (-0.75, -0) {$G$};
		\node [style=circle, scale=2] (12) at (-2, -1.25) {};
		\node [style=none] (13) at (-1.5, -0.75) {};
		\node [style=none] (14) at (-1, -0.75) {};
	\end{pgfonlayer}
	\begin{pgfonlayer}{edgelayer}
		\draw (0.center) to (1.center);
		\draw (1.center) to (2.center);
		\draw (2.center) to (3.center);
		\draw (0.center) to (3.center);
		\draw (9.center) to (8);
		\draw (8) to (4);
		\draw [in=90, out=-44, looseness=0.75] (4) to (5.center);
		\draw [bend left=90, looseness=1.25] (13.center) to (14.center);
		\draw (6.center) to (12);
		\draw [bend left, looseness=1.00] (12) to (4);
	\end{pgfonlayer}
\end{tikzpicture}  ~~~\mbox{\small\bf [LT.2]}~ \begin{tikzpicture}
	\begin{pgfonlayer}{nodelayer}
		\node [style=none] (0) at (-2.5, -0.75) {};
		\node [style=none] (1) at (-2.5, 0.25) {};
		\node [style=none] (2) at (-0.5, 0.25) {};
		\node [style=none] (3) at (-0.5, -0.75) {};
		\node [style=oa] (4) at (-1.5, -0.25) {};
		\node [style=none] (5) at (-2, -2) {};
		\node [style=none] (6) at (-1, -2) {};
		\node [style=none] (7) at (-1.5, 1.25) {};
		\node [style=circle, scale=2] (8) at (-2, -1.25) {};
		\node [style=none] (9) at (-2, -1.25) {$\alpha_\ox$};
		\node [style=none] (10) at (-0.75, 0) {$F$};
		\node [style=none] (11) at (-1.75, -0.75) {};
		\node [style=none] (12) at (-2.25, -0.75) {};
	\end{pgfonlayer}
	\begin{pgfonlayer}{edgelayer}
		\draw (0.center) to (1.center);
		\draw (1.center) to (2.center);
		\draw (2.center) to (3.center);
		\draw (0.center) to (3.center);
		\draw [in=90, out=-15, looseness=1.00] (4) to (6.center);
		\draw (4) to (7.center);
		\draw (5.center) to (8);
		\draw [in=-165, out=90, looseness=1.25] (8) to (4);
		\draw [bend right=90, looseness=1.25] (11.center) to (12.center);
	\end{pgfonlayer}
\end{tikzpicture} = \begin{tikzpicture}
	\begin{pgfonlayer}{nodelayer}
		\node [style=none] (0) at (-2.5, -0.75) {};
		\node [style=none] (1) at (-2.5, 0.25) {};
		\node [style=none] (2) at (-0.5, 0.25) {};
		\node [style=none] (3) at (-0.5, -0.75) {};
		\node [style=oa] (4) at (-1.5, -0.25) {};
		\node [style=none] (5) at (-2, -2) {};
		\node [style=none] (6) at (-1, -2) {};
		\node [style=none] (7) at (-1, -1.25) {$\alpha_\oa$};
		\node [style=circle, scale=2] (8) at (-1.5, 1) {};
		\node [style=none] (9) at (-1.5, 2) {};
		\node [style=none] (10) at (-1.5, 1) {$\alpha_\ox$};
		\node [style=none] (11) at (-0.75, -0) {$G$};
		\node [style=circle, scale=2] (12) at (-1, -1.25) {};
		\node [style=none] (13) at (-1.5, -0.75) {};
		\node [style=none] (14) at (-2, -0.75) {};
	\end{pgfonlayer}
	\begin{pgfonlayer}{edgelayer}
		\draw (0.center) to (1.center);
		\draw (1.center) to (2.center);
		\draw (2.center) to (3.center);
		\draw (0.center) to (3.center);
		\draw (9.center) to (8);
		\draw (8) to (4);
		\draw [in=90, out=-135, looseness=0.75] (4) to (5.center);
		\draw [bend right=90, looseness=1.25] (13.center) to (14.center);
		\draw (6.center) to (12);
		\draw [bend right, looseness=1.00] (12) to (4);
	\end{pgfonlayer}
\end{tikzpicture} ~~~ \mbox{\small\bf [LT.3]}~ \begin{tikzpicture}
	\begin{pgfonlayer}{nodelayer}
		\node [style=none] (0) at (-0.5, 0) {};
		\node [style=none] (1) at (-0.5, -1) {};
		\node [style=none] (2) at (-2.5, -1) {};
		\node [style=none] (3) at (-2.5, -0) {};
		\node [style=ox] (4) at (-1.5, -0.5) {};
		\node [style=none] (5) at (-1, 1.25) {};
		\node [style=none] (6) at (-2, 1.25) {};
		\node [style=none] (7) at (-1.5, -2) {};
		\node [style=circle, scale=2] (8) at (-1, 0.5) {};
		\node [style=none] (9) at (-1, 0.5) {$\alpha_\oa$};
		\node [style=none] (10) at (-0.75, -0.5) {$F$};
		\node [style=none] (11) at (-1.25, 0) {};
		\node [style=none] (12) at (-0.75, 0) {};
	\end{pgfonlayer}
	\begin{pgfonlayer}{edgelayer}
		\draw (0.center) to (1.center);
		\draw (1.center) to (2.center);
		\draw (2.center) to (3.center);
		\draw (0.center) to (3.center);
		\draw [in=-90, out=165, looseness=1.00] (4) to (6.center);
		\draw (4) to (7.center);
		\draw (5.center) to (8);
		\draw [in=15, out=-90, looseness=1.25] (8) to (4);
		\draw [bend right=90, looseness=1.25] (11.center) to (12.center);
	\end{pgfonlayer}
\end{tikzpicture} = \begin{tikzpicture}
	\begin{pgfonlayer}{nodelayer}
		\node [style=none] (0) at (-0.5, 0.75) {};
		\node [style=none] (1) at (-0.5, -0.25) {};
		\node [style=none] (2) at (-2.5, -0.25) {};
		\node [style=none] (3) at (-2.5, 0.75) {};
		\node [style=ox] (4) at (-1.5, 0.25) {};
		\node [style=none] (5) at (-1, 2) {};
		\node [style=none] (6) at (-2, 2) {};
		\node [style=none] (7) at (-2, 1.25) {$\alpha_\ox$};
		\node [style=circle, scale=2] (8) at (-1.5, -1) {};
		\node [style=none] (9) at (-1.5, -2) {};
		\node [style=none] (10) at (-1.5, -1) {$\alpha_\oa$};
		\node [style=none] (11) at (-0.75, 0.5) {$F$};
		\node [style=circle, scale=2] (12) at (-2, 1.25) {};
		\node [style=none] (13) at (-1.5, 0.75) {};
		\node [style=none] (14) at (-1, 0.75) {};
	\end{pgfonlayer}
	\begin{pgfonlayer}{edgelayer}
		\draw (0.center) to (1.center);
		\draw (1.center) to (2.center);
		\draw (2.center) to (3.center);
		\draw (0.center) to (3.center);
		\draw (9.center) to (8);
		\draw (8) to (4);
		\draw [in=-90, out=44, looseness=0.75] (4) to (5.center);
		\draw [bend right=90, looseness=1.25] (13.center) to (14.center);
		\draw (6.center) to (12);
		\draw [bend right, looseness=1.00] (12) to (4);
	\end{pgfonlayer}
\end{tikzpicture}~~~ \mbox{\small\bf [LT.4]}~ \begin{tikzpicture}
	\begin{pgfonlayer}{nodelayer}
		\node [style=none] (0) at (-2.5, 0) {};
		\node [style=none] (1) at (-2.5, -1) {};
		\node [style=none] (2) at (-0.5, -1) {};
		\node [style=none] (3) at (-0.5, 0) {};
		\node [style=ox] (4) at (-1.5, -0.5) {};
		\node [style=none] (5) at (-2, 1.25) {};
		\node [style=none] (6) at (-1, 1.25) {};
		\node [style=none] (7) at (-1.5, -2) {};
		\node [style=circle, scale=2] (8) at (-2, 0.5) {};
		\node [style=none] (9) at (-2, 0.5) {$\alpha_\oa$};
		\node [style=none] (10) at (-0.75, -0.25) {$G$};
		\node [style=none] (11) at (-1.75, 0) {};
		\node [style=none] (12) at (-2.25, 0) {};
	\end{pgfonlayer}
	\begin{pgfonlayer}{edgelayer}
		\draw (0.center) to (1.center);
		\draw (1.center) to (2.center);
		\draw (2.center) to (3.center);
		\draw (0.center) to (3.center);
		\draw [in=-90, out=15, looseness=1.00] (4) to (6.center);
		\draw (4) to (7.center);
		\draw (5.center) to (8);
		\draw [in=165, out=-90, looseness=1.25] (8) to (4);
		\draw [bend left=90, looseness=1.25] (11.center) to (12.center);
	\end{pgfonlayer}
\end{tikzpicture} = \begin{tikzpicture}
	\begin{pgfonlayer}{nodelayer}
		\node [style=none] (0) at (-2.5, 0.75) {};
		\node [style=none] (1) at (-2.5, -0.25) {};
		\node [style=none] (2) at (-0.5, -0.25) {};
		\node [style=none] (3) at (-0.5, 0.75) {};
		\node [style=ox] (4) at (-1.5, 0.25) {};
		\node [style=none] (5) at (-2, 2) {};
		\node [style=none] (6) at (-1, 2) {};
		\node [style=none] (7) at (-1, 1.25) {$\alpha_\ox$};
		\node [style=circle, scale=2] (8) at (-1.5, -1) {};
		\node [style=none] (9) at (-1.5, -2) {};
		\node [style=none] (10) at (-1.5, -1) {$\alpha_\oa$};
		\node [style=none] (11) at (-0.75, 0.5) {$F$};
		\node [style=circle, scale=2] (12) at (-1, 1.25) {};
		\node [style=none] (13) at (-1.5, 0.75) {};
		\node [style=none] (14) at (-2, 0.75) {};
	\end{pgfonlayer}
	\begin{pgfonlayer}{edgelayer}
		\draw (0.center) to (1.center);
		\draw (1.center) to (2.center);
		\draw (2.center) to (3.center);
		\draw (0.center) to (3.center);
		\draw (9.center) to (8);
		\draw (8) to (4);
		\draw [in=-90, out=136, looseness=0.75] (4) to (5.center);
		\draw [bend left=90, looseness=1.00] (13.center) to (14.center);
		\draw (6.center) to (12);
		\draw [bend left, looseness=1.00] (12) to (4);
	\end{pgfonlayer}
\end{tikzpicture} 
\]

An adjunction of linear functors,  $(\eta, \epsilon): F \dashv G$ is an adjunction in the usual sense 
(i.e. satisfying the triangle equalities) in the 2-category of LDCs with linear functors and 
linear natural transformations.   In particular, such an adjunction yields a pair of 
adjunctions: $(\eta_\ox, \epsilon_\ox): F_\ox \dashv G_\ox$ which is a monoidal adjunction, 
and $(\epsilon_\oa, \eta_\oa): G_\oa \dashv F_\oa$ which is a comonoidal adjunction.  
By Kelly's results \cite{Kel97}, a functor with a right adjoint is comonoidal if and only if its 
right adjoint is monoidal. This leads to the observation that:

\begin{lem}
\label{Lemma: strong adjunction}
If $(\eta, \epsilon): F \dashvv G$ is an adjunction of linear functors, then $F_\ox$ is iso-monoidal (or strong) with respect to $\ox$ and $F_\oa$ is iso-comonoidal making the 
linear functor $F$ strong.
\end{lem}
\begin{proof}
Since $(\eta_\ox, \epsilon_\ox): F_\ox \dashv G_\ox$ is a monoidal adjunction, the left adjoint $(F_\ox, m_\ox, m_\top)$ is a strong monoidal functor. Similarly, since $(\epsilon_\oa, \eta_\oa): G_\oa \dashv F_\oa$ is a comonoidal adjunction, the right adjoint $(F_\oa, n_\oa, n_\bot)$ is a strong comonoidal functor.
\end{proof}

A {\bf linear equivalence} is a linear adjunction in which the unit and counit are linear natural isomorphisms.

\subsection{Mix categories}

In this paper we shall be predominately concerned with LDCs which have a mix map:

\begin{defi} \cite{CS97a}
An LDC is a {\bf mix category} in case there is a {\bf mix map} ${\sf m}:\bot\to\top$ in $\X$ such that:
\[
\xymatrixcolsep{4pc}
\xymatrix{
A \ox B \ar[r]^{1 \ox u_\oa^{L^{-1}}} \ar[d]_{(u_\oa^R)^{-1} \ox 1} \ar@{.>}[ddrr]^{\mx_{A,B}} & A \ox (\bot \oa B) \ar[r]^{1 \ox (\m \oa 1)} & A \ox ( \top \oa B) \ar[d]^{\partial^L} \\
(A \oa \bot) \ox B \ar[d]_{\partial^R} & & ( A \ox \top ) \oa B  \ar[d]^{u_\ox^R \oa 1} \\
A \oa (\bot \ox B) \ar[r]_{1 \oa (\m \ox 1)} & A \oa (\top \ox B) \ar[r]_{1 \oa u_\ox^L} &  A \oa B
}
\]
\end{defi}

The map $\mx_{A,B}$ is natural and is called the {\bf mixor}. The coherence condition for the mix map has the following form in string diagrams (where the mix map is represented by an empty box):
\begin{align*}
\mx_{A,B}:=
\begin{tikzpicture}
	\begin{pgfonlayer}{nodelayer}
		\node [style=ox] (0) at (0, 0.2500001) {};
		\node [style=circ] (1) at (0.5000001, -0.2500001) {};
		\node [style=circ] (2) at (0, -1) {$\bot$};
		\node [style=map] (3) at (0, -1.75) {~};
		\node [style=circ] (4) at (0, -2.5) {$\top$};
		\node [style=circ] (5) at (-0.5000001, -3.25) {};
		\node [style=oa] (6) at (0, -3.75) {};
		\node [style=none] (7) at (0, 0.7499999) {};
		\node [style=none] (8) at (0, -4.25) {};
	\end{pgfonlayer}
	\begin{pgfonlayer}{edgelayer}
		\draw (7) to (0);
		\draw (0) to (1);
		\draw [in=45, out=-60, looseness=1.00] (1) to (6);
		\draw [in=120, out=-135, looseness=1.00] (0) to (5);
		\draw (5) to (6);
		\draw (6) to (8);
		\draw [densely dotted, in=-90, out=45, looseness=1.00] (5) to (4);
		\draw (4) to (3);
		\draw (3) to (2);
		\draw [densely dotted, in=-135, out=90, looseness=1.00] (2) to (1);
	\end{pgfonlayer}
\end{tikzpicture}
=
\begin{tikzpicture}
	\begin{pgfonlayer}{nodelayer}
		\node [style=circ] (0) at (-0.5000001, -0.2500001) {};
		\node [style=circ] (1) at (0, -1) {$\bot$};
		\node [style=map] (2) at (0, -1.75) {~};
		\node [style=circ] (3) at (0, -2.5) {$\top$};
		\node [style=circ] (4) at (0.5000001, -3.25) {};
		\node [style=none] (5) at (0, 0.7499999) {};
		\node [style=none] (6) at (0, -4.25) {};
		\node [style=oa] (7) at (0, -3.75) {};
		\node [style=ox] (8) at (0, 0.2500001) {};
	\end{pgfonlayer}
	\begin{pgfonlayer}{edgelayer}
		\draw [densely dotted, in=-90, out=150, looseness=1.00] (4) to (3);
		\draw (3) to (2);
		\draw (2) to (1);
		\draw [densely dotted, in=-45, out=90, looseness=1.00] (1) to (0);
		\draw (8) to (5);
		\draw (8) to (0);
		\draw [in=135, out=-120, looseness=1.00] (0) to (7);
		\draw (7) to (6);
		\draw (7) to (4);
		\draw [in=-45, out=60, looseness=1.00] (4) to (8);
	\end{pgfonlayer}
\end{tikzpicture}
\end{align*}

In a mix category, the associator, the distributor and the mix maps interact as follows. See Lemma 2, and proposition 3 in \cite{BCS00} for a proof.
\begin{equation*}
\mbox{\bf [mix.]}~~~~~~~ \xymatrix{
(A \oa B) \ox C \ar[r]^{\delta^R} \ar[d]_{\mx} \ar@{}[dr]|{(a)} &  A \oa (B \ox C) \ar[d]^{1 \oa \mx} \\
(A \oa B) \oa C \ar[r]_{a_\oa} & A \oa (B \oa C)
} ~~~~~~~~~~~~~
\xymatrix{
(A \ox B) \ox C \ar[r]^{\mx} \ar[d]_{a_\ox} \ar@{}[dr]|{(b)} & A \oa ( B \ox C ) \ar@{<-}[d]^{\delta^L} \\
A \ox (B \ox C) \ar[r]_{1 \ox \mx} & A \ox (B \oa C)
}
\end{equation*}
\[  ~~~~~~~~~~~~ \xymatrix{
C \ox (A \oa B) \ar[r]^{\delta^L} \ar[d]_{\mx} \ar@{}[dr]|{(c)} &  (C \ox A) \oa B \ar[d]^{ \mx \oa 1} \\
C \oa (A \oa B) \ar[r]_{a_\oa^{-1}} & (C \oa A) \oa B
} ~~~~~~~~~~~~~
\xymatrix{
A \ox (B \ox C) \ar[r]^{\mx} \ar[d]_{a_\ox^{-1}} \ar@{}[dr]|{(d)} & A \oa ( B \ox C ) \ar@{<-}[d]^{\delta^R} \\
(A \ox B) \ox C \ar[r]_{\mx \ox 1} & (A \oa B) \ox C
} \]

There are many examples of mix categories including Coherence spaces \cite{Gir87}, 
and Finiteness spaces \cite{Ehr05}.

When the mix map ${\sf m}$ is an isomorphism, then $\X$ is said to be an {\bf isomix category}. 
Recall that, when ${\sf m}$ is an isomorphism, the coherence requirement for the mixor is 
automatically satisfied (see \cite[Lemma 6.6]{CS97a}). Finiteness spaces \cite{Ehr05} and Chu Spaces \cite{Barr06}
 provide examples of isomix categories.

An isomix category, $(\X,\ox, \oa)$ always has two linear functors ${\sf Mx}_\downarrow: 
(\X,\ox, \ox) \to (\X,\ox,\oa)$  and ${\sf Mx}_\uparrow: (\X,\oa, \oa) \to (\X,\ox,\oa)$  
given by the  identity functor, that is $({\sf Mx}_\uparrow)_\ox = ({\sf Mx}_\uparrow)_\oa = 
{\sf Id} = ({\sf Mx}_\downarrow)_\ox = ({\sf Mx}_\downarrow)_\oa$.   
The linear strengths and monoidal maps are given by the inverse of the mix map and  the mixor.  
These mix functors take the degenerate linear structure on the tensor (respectively the par) 
and spread it out over both the tensor structures. 

\begin{lem} \label{mix-functor}
For any isomix category $\X$ the functors ${\sf Mx}_\downarrow: (\X,\ox,\ox) \to (\X,\ox, \oa)$ 
and ${\sf Mx}_\uparrow: (\X,\oa,\oa) \to (\X,\ox, \oa)$ are linear functors.
\end{lem}

\begin{proof}
We show that ${\sf Mx}_\downarrow: (\X,\ox, \ox) \to (\X,\ox,\oa)$ is a linear functor: 
the monoidal and comonoidal components of the functor are given by 
$(1, 1, 1)$ and  $( 1, \mx, \m^{-1})$ respectively. The linear strenghts are $\nu_\ox^L = \nu_\ox^R: A \ox B \to 
A \oa B := \mx$ and $\nu_\oa^L = \nu_\oa^R: A \oa B \to A \oa B := 1$.

First we show $( 1, \mx, \m^{-1}): (\X,\ox,\ox) \to (\X,\ox, \oa)$ is a monoidal functor:

\begin{itemize}
\item The associative law for monoidal functors, $(\mx \ox 1)~\mx~a_\oa = a_\ox~(1 \ox \mx)~\mx$, 
is satisfied:
\[
\begin{tikzpicture} %asso1
\begin{pgfonlayer}{nodelayer}
\node [style=circle, scale=0.5] (0) at (-0.5, 2) {};
\node [style=circle, scale=0.5] (1) at (0.5, 1.25) {};
\node [style=map] (2) at (0, 1.75) {};
\node [style=ox] (3) at (0, 2.5) {};
\node [style=oa] (4) at (0, 0.75) {};
\node [style=ox] (5) at (1, 3) {};
\node [style=oa] (6) at (0, -1) {};
\node [style=oa] (7) at (1, -1.75) {};
\node [style=none] (8) at (1, -3) {};
\node [style=none] (9) at (1, 3.75) {};
\node [style=none] (10) at (-0.75, -3) {};
\node [style=circle, scale=0.5] (11) at (1.5, -1.25) {};
\node [style=map] (12) at (0.75, -0.5) {};
\node [style=circle, scale=0.5] (13) at (0, -0) {};
\end{pgfonlayer}
\begin{pgfonlayer}{edgelayer}
\draw [dotted, in=90, out=-15, looseness=1.25] (0) to (2);
\draw [dotted, in=165, out=-90, looseness=1.25] (2) to (1);
\draw (9.center) to (5);
\draw [bend left=45, looseness=1.00] (5) to (7);
\draw [bend right=15, looseness=1.00] (5) to (3);
\draw [bend right=15, looseness=1.00] (6) to (7);
\draw (7) to (8.center);
\draw [in=90, out=-126, looseness=1.00] (6) to (10.center);
\draw [bend right=60, looseness=1.25] (3) to (4);
\draw (4) to (6);
\draw [bend left=60, looseness=1.25] (3) to (4);
\draw [dotted, in=90, out=0, looseness=1.50] (13) to (12);
\draw [dotted, in=165, out=-90, looseness=1.25] (12) to (11);
\end{pgfonlayer}
\end{tikzpicture} = \begin{tikzpicture} %asso2
\begin{pgfonlayer}{nodelayer}
\node [style=circle, scale=0.5] (0) at (0.25, 1) {};
\node [style=circle, scale=0.5] (1) at (1.75, -0.5) {};
\node [style=map] (2) at (1, 0.5) {};
\node [style=ox] (3) at (-0.25, 1.75) {};
\node [style=ox] (4) at (1, 3) {};
\node [style=oa] (5) at (1, -1.75) {};
\node [style=none] (6) at (1, -3) {};
\node [style=none] (7) at (1, 3.75) {};
\node [style=none] (8) at (-1, -3) {};
\node [style=circle, scale=0.5] (9) at (0.5, -0.75) {};
\node [style=map] (10) at (-0.25, -0) {};
\node [style=circle, scale=0.5] (11) at (-0.75, 1) {};
\end{pgfonlayer}
\begin{pgfonlayer}{edgelayer}
\draw [dotted, in=90, out=-15, looseness=1.25] (0) to (2);
\draw [dotted, in=165, out=-90, looseness=1.25] (2) to (1);
\draw (7.center) to (4);
\draw [bend left=45, looseness=1.00] (4) to (5);
\draw [bend right=15, looseness=1.00] (4) to (3);
\draw (5) to (6.center);
\draw [dotted, in=90, out=0, looseness=1.50] (11) to (10);
\draw [dotted, in=165, out=-90, looseness=1.25] (10) to (9);
\draw [in=90, out=-165, looseness=0.50] (3) to (8.center);
\draw [in=135, out=-45, looseness=1.00] (3) to (5);
\end{pgfonlayer}
\end{tikzpicture} 
  = \begin{tikzpicture} %asso3
\begin{pgfonlayer}{nodelayer}
\node [style=circle, scale=0.5] (0) at (0.25, 2.5) {};
\node [style=circle, scale=0.5] (1) at (1.5, -0.5) {};
\node [style=map] (2) at (1, 0.5) {};
\node [style=ox] (3) at (-0.25, 1.75) {};
\node [style=ox] (4) at (1, 3) {};
\node [style=oa] (5) at (1, -1.25) {};
\node [style=none] (6) at (1, -3) {};
\node [style=none] (7) at (1, 3.75) {};
\node [style=none] (8) at (-1, -3) {};
\node [style=circle, scale=0.5] (9) at (1, -2.25) {};
\node [style=map] (10) at (-0.25, -0) {};
\node [style=circle, scale=0.5] (11) at (-0.75, 1) {};
\end{pgfonlayer}
\begin{pgfonlayer}{edgelayer}
\draw [dotted, in=90, out=-15, looseness=1.25] (0) to (2);
\draw [dotted, in=165, out=-90, looseness=1.25] (2) to (1);
\draw (7.center) to (4);
\draw [bend left=45, looseness=1.00] (4) to (5);
\draw [bend right=15, looseness=1.00] (4) to (3);
\draw (5) to (6.center);
\draw [dotted, in=90, out=0, looseness=1.50] (11) to (10);
\draw [dotted, in=165, out=-90, looseness=1.25] (10) to (9);
\draw [in=90, out=-165, looseness=0.50] (3) to (8.center);
\draw [in=135, out=-45, looseness=1.00] (3) to (5);
\end{pgfonlayer}
\end{tikzpicture}
\]

\item The unit laws for monoidal functors hold.   Here is the pictorial proof of $(1 \ox \m^{-1}) \mx = u_\ox^L(u_\oa^L)^{-1}$, where
the filled rectangles represent $\m^{-1}$:
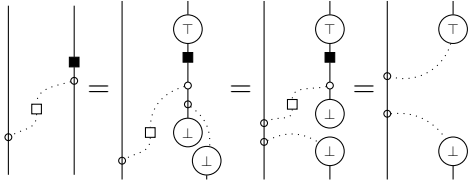
\[
\begin{tikzpicture}
\begin{pgfonlayer}{nodelayer}
\node [style=map, fill=black] (0) at (0, 2) {};
\node [style=circle, scale=0.5] (1) at (0, 1.5) {};
\node [style=circle, scale=0.5] (2) at (-1.75, -0) {};
\node [style=map] (3) at (-1, 0.75) {};
\node [style=none] (4) at (0, 3.5) {};
\node [style=none] (5) at (0, -1) {};
\node [style=none] (6) at (-1.75, -1) {};
\node [style=none] (7) at (-1.75, 3.5) {};
\end{pgfonlayer}
\begin{pgfonlayer}{edgelayer}
\draw [dotted, in=90, out=-165, looseness=1.00] (1) to (3);
\draw [dotted, in=30, out=-90, looseness=1.25] (3) to (2);
\draw (4.center) to (0);
\draw (7.center) to (6.center);
\draw (0) to (5.center);
\end{pgfonlayer}
\end{tikzpicture} = \begin{tikzpicture}
\begin{pgfonlayer}{nodelayer}
\node [style=map, fill=black] (0) at (0, 2.25) {};
\node [style=circle, scale=0.5] (1) at (0, 1.5) {};
\node [style=circle, scale=0.5] (2) at (-1.75, -0.5) {};
\node [style=map] (3) at (-1, 0.25) {};
\node [style=none] (4) at (0, 3.75) {};
\node [style=none] (5) at (0.5, -1) {};
\node [style=none] (6) at (-1.75, -1) {};
\node [style=none] (7) at (-1.75, 3.75) {};
\node [style=circle] (8) at (0, 3) {$\top$};
\node [style=circle] (9) at (0, 0.25) {$\bot$};
\node [style=circle] (10) at (0.5, -0.5) {$\bot$};
\node [style=circle, scale=0.5] (11) at (0, 1) {};
\end{pgfonlayer}
\begin{pgfonlayer}{edgelayer}
\draw [dotted, dotted, in=90, out=-165, looseness=1.00] (1) to (3);
\draw [dotted, dotted, in=30, out=-90, looseness=1.25] (3) to (2);
\draw (7.center) to (6.center);
\draw (8) to (0);
\draw (0) to (1);
\draw (1) to (9);
\draw (10) to (5.center);
\draw [bend left, looseness=1.00, dotted] (11) to (10);
\draw (4.center) to (8);
\end{pgfonlayer}
\end{tikzpicture} 
  = \begin{tikzpicture}
\begin{pgfonlayer}{nodelayer}
\node [style=map, fill=black] (0) at (0, 2.25) {};
\node [style=circle, scale=0.5] (1) at (0, 1.5) {};
\node [style=circle, scale=0.5] (2) at (-1.75, 0.5) {};
\node [style=map] (3) at (-1, 1) {};
\node [style=none] (4) at (0, 3.75) {};
\node [style=none] (5) at (0, -1) {};
\node [style=none] (6) at (-1.75, -1) {};
\node [style=none] (7) at (-1.75, 3.75) {};
\node [style=circle] (8) at (0, 3) {$\top$};
\node [style=circle] (9) at (0, 0.75) {$\bot$};
\node [style=circle] (10) at (0, -0.25) {$\bot$};
\node [style=circle, scale=0.5] (11) at (-1.75, -0) {};
\end{pgfonlayer}
\begin{pgfonlayer}{edgelayer}
\draw [dotted, dotted, in=90, out=-165, looseness=1.00] (1) to (3);
\draw [dotted, dotted, in=30, out=-90, looseness=1.25] (3) to (2);
\draw (7.center) to (6.center);
\draw (8) to (0);
\draw (0) to (1);
\draw (1) to (9);
\draw (10) to (5.center);
\draw [bend left, looseness=1.00, dotted] (11) to (10);
\draw (4.center) to (8);
\end{pgfonlayer}
\end{tikzpicture} = \begin{tikzpicture}
\begin{pgfonlayer}{nodelayer}
\node [style=circle, scale=0.5] (0) at (-1.75, 1.75) {};
\node [style=none] (1) at (0, 3.75) {};
\node [style=none] (2) at (0, -1) {};
\node [style=none] (3) at (-1.75, -1) {};
\node [style=none] (4) at (-1.75, 3.75) {};
\node [style=circle] (5) at (0, 3) {$\top$};
\node [style=circle] (6) at (0, -0.25) {$\bot$};
\node [style=circle, scale=0.5] (7) at (-1.75, 0.75) {};
\end{pgfonlayer}
\begin{pgfonlayer}{edgelayer}
\draw (4.center) to (3.center);
\draw (6) to (2.center);
\draw [dotted, bend left, looseness=1.00] (7) to (6);
\draw (1.center) to (5);
\draw [dotted, bend left=45, looseness=1.00] (5) to (0);
\end{pgfonlayer}
\end{tikzpicture}
\]
The other unit law holds similarly.
\end{itemize}

${\sf Mx}_\downarrow: (\X,\ox, \ox) \to (\X,\ox,\oa)$ satisfies all the coherence requirements of a linear functor:
{\bf [LF.1]}, {\bf [LF.2]}, and {\bf [LF.3]} hold because $({\sf Mx}_\downarrow)_\ox$ and $({\sf Mx}_\downarrow)_\oa$ are monoidal and comonoidal respectively, {\bf [LF.4]}(a) becomes 
$\mx a_\oa^{-1} = \partial^L (\mx \oa 1)$ and holds because:
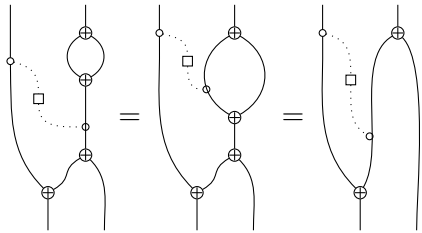
\[
\begin{tikzpicture} %dist1
\begin{pgfonlayer}{nodelayer}
\node [style=oa] (0) at (0, 2.25) {};
\node [style=oa] (1) at (0, 1) {};
\node [style=oa] (2) at (0, -1) {};
\node [style=oa] (3) at (-1, -2) {};
\node [style=none] (4) at (-1, -3) {};
\node [style=none] (5) at (0.5, -3) {};
\node [style=none] (6) at (-2, 3) {};
\node [style=none] (7) at (0, 3) {};
\node [style=map] (8) at (-1.25, 0.5) {};
\node [style=circle, scale=0.5] (9) at (0, -0.25) {};
\node [style=circle, scale=0.5] (10) at (-2, 1.5) {};
\end{pgfonlayer}
\begin{pgfonlayer}{edgelayer}
\draw (6.center) to (10);
\draw (7.center) to (0);
\draw [bend left=60, looseness=1.25] (0) to (1);
\draw [bend right=60, looseness=1.25] (0) to (1);
\draw (1) to (2);
\draw [in=90, out=-45, looseness=1.00] (2) to (5.center);
\draw [in=30, out=-150, looseness=1.50] (2) to (3);
\draw [in=-89, out=135, looseness=1.00] (3) to (10);
\draw (3) to (4.center);
\draw [in=90, out=-15, looseness=1.25,dotted] (10) to (8);
\draw [in=180, out=-90, looseness=1.25, dotted] (8) to (9);
\end{pgfonlayer}
\end{tikzpicture} 
  = \begin{tikzpicture}
\begin{pgfonlayer}{nodelayer}
\node [style=oa] (0) at (0, 2.25) {};
\node [style=oa] (1) at (0, -0) {};
\node [style=oa] (2) at (0, -1) {};
\node [style=oa] (3) at (-1, -2) {};
\node [style=none] (4) at (-1, -3) {};
\node [style=none] (5) at (0.5, -3) {};
\node [style=none] (6) at (-2, 3) {};
\node [style=none] (7) at (0, 3) {};
\node [style=map] (8) at (-1.25, 1.5) {};
\node [style=circle,scale=0.5] (9) at (-0.75, 0.75) {};
\node [style=circle, scale=0.5] (10) at (-2, 2.25) {};
\end{pgfonlayer}
\begin{pgfonlayer}{edgelayer}
\draw (6.center) to (10);
\draw (7.center) to (0);
\draw [bend left=60, looseness=1.25] (0) to (1);
\draw [bend right=60, looseness=1.25] (0) to (1);
\draw (1) to (2);
\draw [in=90, out=-45, looseness=1.00] (2) to (5.center);
\draw [in=30, out=-150, looseness=1.50] (2) to (3);
\draw [in=-89, out=135, looseness=1.00] (3) to (10);
\draw (3) to (4.center);
\draw [dotted, in=90, out=-15, looseness=1.25] (10) to (8);
\draw [dotted, in=180, out=-90, looseness=1.25] (8) to (9);
\end{pgfonlayer}
\end{tikzpicture} 
  = \begin{tikzpicture}
\begin{pgfonlayer}{nodelayer}
\node [style=oa] (0) at (0, 2.25) {};
\node [style=oa] (1) at (-1, -2) {};
\node [style=none] (2) at (-1, -3) {};
\node [style=none] (3) at (0.5, -3) {};
\node [style=none] (4) at (-2, 3) {};
\node [style=none] (5) at (0, 3) {};
\node [style=map] (6) at (-1.25, 1) {};
\node [style=circle, scale=0.5] (7) at (-0.75, -0.5) {};
\node [style=circle, scale=0.5] (8) at (-2, 2.25) {};
\end{pgfonlayer}
\begin{pgfonlayer}{edgelayer}
\draw (4.center) to (8);
\draw (5.center) to (0);
\draw [in=-89, out=135, looseness=1.00] (1) to (8);
\draw (1) to (2.center);
\draw [dotted, in=90, out=-15, looseness=1.25] (8) to (6);
\draw [dotted, in=150, out=-90, looseness=1.00] (6) to (7);
\draw [in=90, out=-45, looseness=0.50] (0) to (3.center);
\draw [in=60, out=-150, looseness=0.75] (0) to (1);
\end{pgfonlayer}
\end{tikzpicture}
\]
{\bf [LF.4]} (b) - (d) and {\bf [LF.5]} (a) - (d) are satisfied similarly. 

Thus, ${\sf Mx}_\downarrow$ is a linear functor. 

The proof that ${\sf Mx}_\uparrow$ is a linear functor is (linearly) dual.
\end{proof}

In fact, these linear functors and are examples of {\em  isomix Frobenius functors\/}, 
which we shall introduce formally in the Section~\ref{Sec: frobenius functors}.

\subsection{Compact LDCs}

A {\bf compact LDC} is an isomix category in which each mixor $\mx_{A,B}$ is an isomorphism.  
An important way in which compact LDCs arise is from the ``core'' of an isomix category:

%In any isomix category, $\X$, the core, $\Core(\X)$, forms a compact LDC as 
%here both $\mx_{A,B}: A \ox B \to A \oa B$ and ${\sf m}: \top \to \bot$ are isomorphisms.  

\begin{defi} \cite{BCS00}
The {\bf core} of a mix category, $\Core(\X) \subseteq \X$, is the full subcategory with objects $U$ such that the mixors
\[ U \ox (\_) \to^{\mx_{U,(\_)}} U \oa (\_) ~~~~\mbox{and}~~~~ (\_) \ox  U \to^{\mx_{(\_),U}} (\_) \oa U \]
are isomorphisms.
\end{defi}

\begin{prop} \cite[Proposition 3]{BCS00} 
If $\X$ is a mix-LDC and $A,B \in \Core(\X)$ then $A \oa B$ and $A \ox B \in \Core(\X)$ (and $A \oa B \simeq A \ox B$).  If $\X$ is an isomix-LDC, then $\top, \bot \in \Core(\X)$.  
\end{prop}

\begin{cor} \label{compact-mix-functor}
When $\X$ is a compact LDC,  the mix functors, ${\sf Mx}_\downarrow$ and ${\sf Mx}_\uparrow$, are linear isomorphisms. Consequently, compact LDCs are linearly equivalent to monoidal categories.
\end{cor}

We shall denote the inverse of ${\sf Mx}_\downarrow$ by 
${\sf Mx}^{*}_\downarrow: (\X,\ox,\oa) \to (\X,\oa,\oa)$: this is the identity functor as a mere functor, 
strict on the par structure, and on the tensor structure having as the unit laxor ${\sf m}$ and as the tensor laxor ${\sf mx}^{-1}$.   Similarly, we shall denote the inverse of ${\sf Mx}_\uparrow$ by ${\sf Mx}^{*}_\uparrow$.

%%%%%%%%%%%%%%%%%%%%%%%%%%%%%%%%%%%%%%%%%%%%%%%%

\subsection{Linear duals}

A key notion in the theory of LDCs is the notion of a linear adjoint \cite{CKS00}.  Here we shall refer to linear adjoints as ``linear duals'' in order to avoid any confusion with an adjunction of linear functors.   

\begin{defi} Suppose $\mathbb{X}$ is a LDC and $A,B \in\X$, then $B$ is {\bf left linear dual}  (or left linear adjoint) to $A$ -- or $A$ is {\bf right linear dual} (right linear adjoint) to $B$ -- written $(\eta, \epsilon): B \dashv \!\!\!\!\! \dashv  A$, if there exists $\eta: \top \rightarrow B \oa A$ and $\epsilon: A \ox B \rightarrow \bot$ such that the following diagrams commute:

\[
\xymatrix{
B \ar[r]^{(u_\ox^L)^{-1}} \ar@{=}[d] 
& \top \ox B \ar[r]^{\eta \ox 1} 
& (B \oa A) \ox B \ar[d]^{\partial_R} \\
B 
& B \oa \bot \ar[l]^{u_\oa^R} 
& B \oa (A \ox B) \ar[l]^{1 \oa \epsilon}
}
~~~~~
\xymatrix{
A \ar[r]^{(u_\ox^R)^{-1}} \ar@{=}[d] 
& A \ox \top  \ar[r]^{1 \ox \eta} 
& A  \ox  (B \oa A)\ar[d]^{\partial_L} \\
A
& \bot \oa A \ar[l]^{u_\oa^L} 
& (A \ox B) \oa A   \ar[l]^{ \epsilon \oa 1} }
\]
\end{defi}

The commuting diagrams are called often referred to as ``snake diagrams'' because of their shape when drawn in string calculus:
\[
\begin{tikzpicture}
	\begin{pgfonlayer}{nodelayer}
		\node [style=circle] (0) at (-3, 1.25) {$\eta$};
		\node [style=circle] (1) at (-2, 0) {$\epsilon$};
		\node [style=none] (2) at (-3.5, -0.5) {};
		\node [style=none] (3) at (-1.5, 2) {};
		\node [style=none] (4) at (-1.25, 1.75) {$B$};
		\node [style=none] (5) at (-3.75, -0.25) {$B$};
	\end{pgfonlayer}
	\begin{pgfonlayer}{edgelayer}
		\draw [style=none, in=-90, out=45, looseness=1.00] (1) to (3.center);
		\draw [style=none, in=150, out=-30, looseness=1.00] (0) to (1);
		\draw [style=none, in=90, out=-135, looseness=1.00] (0) to (2.center);
	\end{pgfonlayer}
\end{tikzpicture}  = \begin{tikzpicture}
  \draw (0,2.5) -- (0,0);
\end{tikzpicture} ~~~~~~~~~~
\begin{tikzpicture}%snake 1
	\begin{pgfonlayer}{nodelayer}
		\node [style=circle] (0) at (-3, 0.25) {$\epsilon$};
		\node [style=circle] (1) at (-2, 1.5) {$\eta$};
		\node [style=none] (2) at (-3.5, 2) {};
		\node [style=none] (3) at (-1.5, -0.5) {};
		\node [style=none] (4) at (-3.25, 2) {$A$};
		\node [style=none] (5) at (-1.75, -0.5) {$A$};
	\end{pgfonlayer}
	\begin{pgfonlayer}{edgelayer}
		\draw [style=none, in=90, out=-30, looseness=1.00] (1) to (3.center);
		\draw [style=none, in=-150, out=15, looseness=1.25] (0) to (1);
		\draw [style=none, in=-90, out=165, looseness=1.00] (0) to (2.center);
	\end{pgfonlayer}
\end{tikzpicture} =
\begin{tikzpicture}
  \draw (0,2.5) -- (0,0);
\end{tikzpicture} 
\]

\begin{lem} \cite{BCS00}
\begin{enumerate}[(i)]
\item In an LDC if $(\eta,\epsilon): B \dashvv A$ and $(\eta',\epsilon'): C \dashvv A$, then $B$ and $C$ are isomorphic;
\item In a symmetric LDC $(\eta, \epsilon): B \dashvv A$ if and only if $(\eta c_\oa, c_\ox \epsilon): A \dashvv B$;
\item In a mix-LDC if $B \in \Core(\X)$ and $B \dashvv A$, then $A \in \Core(\X)$.
\end{enumerate}
\end{lem}

\begin{lem} \cite{CKS00}
\label{Lemma: linear adjoints}
Linear functors preserve linear duals: when $F: \X \to \Y$ is a linear functor and $(\eta, \epsilon): A \dashvv B \in \X$, then $F_\ox(A) \dashvv F_\oa(B)$ and $F_\oa(A) \dashvv F_\ox(B)$.
\end{lem}
\begin{proof}
The unit and counit of the adjunction $(\eta', \epsilon'): F_\ox(A) \dashvv F_\oa(B)$ is given as follows:
\[ \eta' := \top \xrightarrow{m_\top} F_\ox(\top) \xrightarrow{F_\ox(\eta)} F_\ox( A \oa B) \xrightarrow{\nu_\ox^L} F_\ox(A) \oa F_\oa(B) = \begin{tikzpicture}
	\begin{pgfonlayer}{nodelayer}
		\node [style=circle] (0) at (0, -0.75) {$\eta$};
		\node [style=none] (1) at (-1, -2.75) {};
		\node [style=none] (2) at (1, -2.75) {};
		\node [style=circle] (3) at (-1.75, -0) {$\top$};
		\node [style=none] (4) at (-1.75, 1) {};
		\node [style=none] (5) at (-1, -2) {};
		\node [style=circle, scale=0.4] (6) at (-1, -2) {};
		\node [style=none] (7) at (-1, -0) {};
		\node [style=none] (8) at (1, -0) {};
		\node [style=none] (9) at (1, -1.25) {};
		\node [style=none] (10) at (-1, -1.25) {};
		\node [style=none] (11) at (-0.75, -1.25) {};
		\node [style=none] (12) at (-0.25, -1.25) {};
		\node [style=none] (13) at (0.75, -0.25) {$F$};
	\end{pgfonlayer}
	\begin{pgfonlayer}{edgelayer}
		\draw [bend right, looseness=1.00] (0) to (1.center);
		\draw [bend left, looseness=1.00] (0) to (2.center);
		\draw (4.center) to (3);
		\draw [dotted, bend right=15, looseness=1.00] (3) to (5.center);
		\draw (7.center) to (8.center);
		\draw (8.center) to (9.center);
		\draw (9.center) to (10.center);
		\draw (10.center) to (7.center);
		\draw [bend left=90, looseness=1.25] (11.center) to (12.center);
	\end{pgfonlayer}
\end{tikzpicture}\]
\[ \epsilon' := F_\oa(B) \ox F_\ox(A) \xrightarrow{\nu_\oa^L} F_\oa(B \ox A) \xrightarrow{F_\oa(\epsilon)} F_\oa(\bot) \xrightarrow{n_\bot} \bot = \begin{tikzpicture}
	\begin{pgfonlayer}{nodelayer}
		\node [style=circle] (0) at (0, 0.25) {$\epsilon$};
		\node [style=none] (1) at (-1, 2.25) {};
		\node [style=none] (2) at (1, 2.25) {};
		\node [style=circle] (3) at (-1.75, -1) {$\bot$};
		\node [style=none] (4) at (-1.75, -2) {};
		\node [style=none] (5) at (-1, 1.5) {};
		\node [style=circle, scale=0.4] (6) at (-1, 1.5) {};
		\node [style=none] (7) at (-1, -0.75) {};
		\node [style=none] (8) at (1, -0.75) {};
		\node [style=none] (9) at (1, 0.75) {};
		\node [style=none] (10) at (-1, 0.75) {};
		\node [style=none] (11) at (-0.75, 0.75) {};
		\node [style=none] (12) at (-0.25, 0.75) {};
		\node [style=none] (13) at (0.75, -0.5) {$F$};
	\end{pgfonlayer}
	\begin{pgfonlayer}{edgelayer}
		\draw [bend left, looseness=1.00] (0) to (1.center);
		\draw [bend right, looseness=1.00] (0) to (2.center);
		\draw (4.center) to (3);
		\draw [dotted, in=-120, out=92, looseness=1.00] (3) to (5.center);
		\draw (7.center) to (8.center);
		\draw (8.center) to (9.center);
		\draw (9.center) to (10.center);
		\draw (10.center) to (7.center);
		\draw [bend right=90, looseness=1.25] (11.center) to (12.center);
	\end{pgfonlayer}
\end{tikzpicture}\]

The unit and counit of the other adjunction is given similarly, however using the right linear strengths ($\nu_\ox^R$ and $\nu_\oa^R$).
\end{proof}

An LDC in which every object has a chosen left and right linear dual, respectively $(\eta{*},\epsilon{*}): A^{*}  \dashvv A$ and $({*}\eta,{*}\epsilon): A  \dashvv \!~^{*}A$, is a {\bf $*$-autonomous category}.   In the symmetric case a left linear dual gives a right linear dual using the symmetry: thus, it is standard to assume the existence of just the left dual with the right being the same object with the unit and counit given by symmetry (as above). 

Just as compact LDCs are linearly equivalent to monoidal categories so compact $*$-autonomous categories are linearly equivalent to compact closed categories. The equivalence is given by ${\sf Mx}_\uparrow$ which spreads the par onto  two tensor structures (or, indeed, by ${\sf Mx}_\downarrow$ which shows how to spread out a compact closed 
structure on the tensor).

In a symmetric $*$-autonomous category the left dual of an object is always canonically isomorphic to the right dual.  Moreover, even in non-symmetric $*$-autonomous categories, it is often the case that the two duals are coherently isomorphic:

\begin{defi}\cite{EggMcCurd12}
A {\bf cyclor} in a $*$-autonomous category $(\X, \ox, \top, \oa, \bot, ~^{*}(\_), (\_)^*)$ is a natural isomorphism $A^* \to^{\psi} \!~^{*}A$ satisfying the following coherence conditions:
\[
\xymatrix{
\bot^* \ar[rr]^{\psi} \ar[dr]_{m_\top^{-1}} & \ar@{}[d]|{\mbox{\tiny \bf [C.1]}} & ~^{*}\bot \ar[ld]^{m_\top^{-1}} \\
& \top &
} ~~~~~~~~ \xymatrix{
 & A \ar[ld]_{\eta} \ar[dr]^{\eta} \ar@{}[d]|{\mbox{\tiny \bf [C.2]}}  & \\
(^{*}A)^* \ar[r]_{\psi^*} & (A^*)^* \ar[r]_{\psi_{A^*}} & ~^{*}(A^*)
} ~~~~~~~~ \xymatrix{
\top^* \ar[rr]^{\psi} \ar[dr]_{n_\bot} & \ar@{}[d]|{\mbox{\tiny \bf [C.3]}} & ~^{*}\top \ar[ld]^{n_\bot} \\
& \bot &
}
\]

\[
\xymatrixcolsep{4pc}
\xymatrix{
(A \ox B)^* \ar[r]^{\psi} \ar[d]_{n_\oa} \ar@{}[dr]|{\mbox{\tiny \bf [C.4]}} & ~^{*}(A \ox B) \ar[d]^{n_\oa} \\
(B^* \oa A^*) \ar[r]_{\psi \oa \psi} & ~^{*}B \oa ~^{*}A
} ~~~~~~~~~ \xymatrix{
(A \oa B)^* \ar[r]^{\psi} \ar[d]_{m_\ox^{-1}} \ar@{}[dr]|{\mbox{\tiny \bf [C.5]}}  & ~^{*}(A \oa B) \ar[d]^{m_\ox^{-1}} \\
(B^* \ox A^*) \ar[r]_{\psi \ox \psi} & ~^{*}B \ox ~^{*}A
} \]
A $*$-autonomous category with a cyclor is said to be {\bf cyclic}.
\end{defi}

In {\mbox{ \bf [C.2]}}, $\eta: A \to (^*A)^*$ is drawn as follows: $ \begin{tikzpicture}
	\begin{pgfonlayer}{nodelayer}
		\node [style=none] (0) at (-0.5, 1) {};
		\node [style=none] (1) at (0, 3.5) {};
		\node [style=none] (2) at (0, 2.25) {};
		\node [style=none] (3) at (-1.25, 2.25) {};
		\node [style=none] (4) at (-1.25, 2.75) {};
		\node [style=none] (5) at (-0.5, 2.75) {};
		\node [style=none] (6) at (0.25, 3.25) {$A$};
		\node [style=none] (7) at (0, 1.25) {$(^*A)^*$};
		\node [style=none] (8) at (-1.75, 2.5) {$(^*A)^*$};
	\end{pgfonlayer}
	\begin{pgfonlayer}{edgelayer}
		\draw (1.center) to (2.center);
		\draw [bend left=90, looseness=1.75] (2.center) to (3.center);
		\draw (3.center) to (4.center);
		\draw [bend left=90, looseness=1.50] (4.center) to (5.center);
		\draw (5.center) to (0.center);
	\end{pgfonlayer}
\end{tikzpicture} $. The map $\eta: A \to ~^*(A^*)$ is given similarly.

The coherence conditions are not independent of each other: being cyclic is equivalent to any one of the following four pairs of conditions: ({\bf [C.1]}, {\bf [C.5]}), ({\bf [C.2]}, {\bf [C.5]}),  ({\bf [C.4]}, {\bf [C.2]}) and ({\bf [C.4]}, {\bf [C.3]}) (see \cite{EggMcCurd12}).

Condition {\bf [C.2]} which is used extensively in Section~\ref{daggers-duals-conjugation}:
\[
\begin{tikzpicture}
	\begin{pgfonlayer}{nodelayer}
		\node [style=circle] (0) at (-2, 1) {$\psi_{A^*}$};
		\node [style=circle] (1) at (-0.75, 1) {$\psi_A$};
		\node [style=none] (2) at (-2, 1.75) {};
		\node [style=none] (3) at (-0.75, 1.75) {};
		\node [style=none] (4) at (-2, -1) {};
		\node [style=none] (5) at (-0.75, -0) {};
		\node [style=none] (6) at (0.5, -0) {};
		\node [style=none] (7) at (0.5, 3) {};
		\node [style=none] (8) at (-1.5, 2.75) {$\eta*$};
		\node [style=none] (9) at (0, -0.75) {$*\epsilon$};
		\node [style=none] (10) at (-2.5, 1.75) {$A^{**}$};
		\node [style=none] (11) at (-2.7, -0.7) {$~^*(A^*)$};
		\node [style=none] (12) at (-0.5, 2) {$A^*$};
		\node [style=none] (13) at (-0.45, 0.25) {$~^*A$};
		\node [style=none] (14) at (0.75, 2.5) {$A$};
	\end{pgfonlayer}
	\begin{pgfonlayer}{edgelayer}
		\draw (0) to (4.center);
		\draw (2.center) to (0);
		\draw [bend left=90, looseness=2.00] (2.center) to (3.center);
		\draw (3.center) to (1);
		\draw (1) to (5.center);
		\draw [bend right=90, looseness=1.50] (5.center) to (6.center);
		\draw (6.center) to (7.center);
	\end{pgfonlayer}
\end{tikzpicture}  = \begin{tikzpicture}
	\begin{pgfonlayer}{nodelayer}
		\node [style=none] (0) at (-2, 0.25) {};
		\node [style=none] (1) at (-0.75, 0.25) {};
		\node [style=none] (2) at (-2, 3) {};
		\node [style=none] (3) at (-0.75, 2) {};
		\node [style=none] (4) at (0.5, 2) {};
		\node [style=none] (5) at (0.5, -1) {};
		\node [style=none] (6) at (-0.25, 2.75) {$*\eta$};
		\node [style=none] (7) at (-1.5, -0.75) {$\epsilon*$};
		\node [style=none] (8) at (1, -0.7) {$~^*(A^*)$};
		\node [style=none] (9) at (-0.5, 1.15) {$A^*$};
		\node [style=none] (10) at (-2.25, 2.75) {$A$};
	\end{pgfonlayer}
	\begin{pgfonlayer}{edgelayer}
		\draw [bend right=90, looseness=2.00] (0.center) to (1.center);
		\draw [bend left=90, looseness=1.50] (3.center) to (4.center);
		\draw (4.center) to (5.center);
		\draw (3.center) to (1.center);
		\draw (2.center) to (0.center);
	\end{pgfonlayer}
\end{tikzpicture}
\]

%  You must get the equallities in the right order so you can see the reasoning!!!!!
The requirement {\bf [C.2]} implies:
\[
 \text{\bf [C.2]$^{-1}$} 
  \begin{tikzpicture}
	\begin{pgfonlayer}{nodelayer}
		\node [style=circle] (0) at (0, 1) {$\psi^{-1}_A$};
		\node [style=circle] (1) at (2, 1) {$\psi^{-1}_{A^*}$};
		\node [style=none] (2) at (0, -0) {};
		\node [style=none] (3) at (2, -0) {};
		\node [style=none] (4) at (2, 3) {};
		\node [style=none] (5) at (0, 2) {};
		\node [style=none] (6) at (-1, 2) {};
		\node [style=none] (7) at (-1, -1.25) {};
		\node [style=none] (8) at (1, -1.1) {$\epsilon*$};
		\node [style=none] (9) at (-0.5, 2.75) {$*\eta$};
		\node [style=none] (10) at (2.5, 0.25) {$A^{**}$};
		\node [style=none] (11) at (-0.5, 0.25) {$A^*$};
		\node [style=none] (12) at (2.5, 2.75) {$~^*(A^*)$};
		\node [style=none] (13) at (-0.5, 1.75) {$~^*A$};
		\node [style=none] (14) at (-1.25, -0.75) {$A$};
	\end{pgfonlayer}
	\begin{pgfonlayer}{edgelayer}
		\draw (0) to (2.center);
		\draw [bend right=90, looseness=1.50] (2.center) to (3.center);
		\draw (3.center) to (1);
		\draw (4.center) to (1);
		\draw (5.center) to (0);
		\draw [bend right=90, looseness=2.00] (5.center) to (6.center);
		\draw (6.center) to (7.center);
	\end{pgfonlayer}
\end{tikzpicture} =  \left(
\begin{tikzpicture}
	\begin{pgfonlayer}{nodelayer}
		\node [style=circle] (0) at (-2, 1) {$\psi_{A^*}$};
		\node [style=circle] (1) at (-0.75, 1) {$\psi_A$};
		\node [style=none] (2) at (-2, 1.75) {};
		\node [style=none] (3) at (-0.75, 1.75) {};
		\node [style=none] (4) at (-2, -1) {};
		\node [style=none] (5) at (-0.75, -0) {};
		\node [style=none] (6) at (0.5, -0) {};
		\node [style=none] (7) at (0.5, 3) {};
		\node [style=none] (8) at (-1.5, 2.75) {$\eta*$};
		\node [style=none] (9) at (0, -0.75) {$*\epsilon$};
		\node [style=none] (10) at (-2.5, 1.75) {$A^{**}$};
		\node [style=none] (11) at (-2.7, -0.7) {$~^*(A^*)$};
		\node [style=none] (12) at (-0.5, 2) {$A^*$};
		\node [style=none] (13) at (-0.45, 0.25) {$~^*A$};
		\node [style=none] (14) at (0.75, 2.5) {$A$};
	\end{pgfonlayer}
	\begin{pgfonlayer}{edgelayer}
		\draw (0) to (4.center);
		\draw (2.center) to (0);
		\draw [bend left=90, looseness=2.00] (2.center) to (3.center);
		\draw (3.center) to (1);
		\draw (1) to (5.center);
		\draw [bend right=90, looseness=1.50] (5.center) to (6.center);
		\draw (6.center) to (7.center);
	\end{pgfonlayer}
\end{tikzpicture} \right)^{-1} = \left( \begin{tikzpicture}
	\begin{pgfonlayer}{nodelayer}
		\node [style=none] (0) at (-2, 0.25) {};
		\node [style=none] (1) at (-0.75, 0.25) {};
		\node [style=none] (2) at (-2, 3) {};
		\node [style=none] (3) at (-0.75, 2) {};
		\node [style=none] (4) at (0.5, 2) {};
		\node [style=none] (5) at (0.5, -1) {};
		\node [style=none] (6) at (-0.25, 2.75) {$*\eta$};
		\node [style=none] (7) at (-1.5, -0.75) {$\epsilon*$};
		\node [style=none] (8) at (1, -0.7) {$~^*(A^*)$};
		\node [style=none] (9) at (-0.5, 1.15) {$A^*$};
		\node [style=none] (10) at (-2.25, 2.75) {$A$};
	\end{pgfonlayer}
	\begin{pgfonlayer}{edgelayer}
		\draw [bend right=90, looseness=2.00] (0.center) to (1.center);
		\draw [bend left=90, looseness=1.50] (3.center) to (4.center);
		\draw (4.center) to (5.center);
		\draw (3.center) to (1.center);
		\draw (2.center) to (0.center);
	\end{pgfonlayer}
\end{tikzpicture} \right)^{-1}
= \begin{tikzpicture}
	\begin{pgfonlayer}{nodelayer}
		\node [style=none] (0) at (-2, 0.25) {};
		\node [style=none] (1) at (-0.75, 0.25) {};
		\node [style=none] (2) at (-2, 3) {};
		\node [style=none] (3) at (-0.75, 2) {};
		\node [style=none] (4) at (0.5, 2) {};
		\node [style=none] (5) at (0.5, -1) {};
		\node [style=none] (6) at (-0.25, 2.75) {$\eta*$};
		\node [style=none] (7) at (-1.5, -0.75) {$*\epsilon$};
		\node [style=none] (8) at (-2.75, 2.75) {$~^*(A^*)$};
		\node [style=none] (9) at (-0.5, 1) {$A^*$};
		\node [style=none] (10) at (0.75, -0.75) {$A$};
	\end{pgfonlayer}
	\begin{pgfonlayer}{edgelayer}
		\draw [bend right=90, looseness=2.00] (0.center) to (1.center);
		\draw [bend left=90, looseness=1.50] (3.center) to (4.center);
		\draw (4.center) to (5.center);
		\draw (3.center) to (1.center);
		\draw (2.center) to (0.center);
	\end{pgfonlayer}
\end{tikzpicture}
\]

Symmetric $*$-autonomous categories always have a canonical cyclor:

\[
\begin{tikzpicture}
	\begin{pgfonlayer}{nodelayer}
		\node [style=none] (0) at (-0.25, 2.25) {};
		\node [style=none] (1) at (-0.25, 0.5) {};
		\node [style=none] (12) at (-0.85, 1.7) {$*\eta$};
		\node [style=none] (2) at (-1.5, 0.5) {};
		\node [style=none] (3) at (-1.5, 1) {};
		\node [style=none] (4) at (-0.5, 1) {};
		\node [style=none] (34) at (-0.85, -0.4) {$\epsilon*$};
		\node [style=none] (5) at (-0.5, -1) {};
		\node [style=none] (6) at (0, 2) {$A^{*}$};
		\node [style=none] (7) at (-0.15, -0.75) {${~^*A}$};
	\end{pgfonlayer}
	\begin{pgfonlayer}{edgelayer}
		\draw (0.center) to (1.center);
		\draw [bend left=90, looseness=1.75] (1.center) to (2.center);
		\draw (2.center) to (3.center);
		\draw [bend left=90, looseness=1.75] (3.center) to (4.center);
		\draw (4.center) to (5.center);
	\end{pgfonlayer}
\end{tikzpicture} 
\]

We shall use the cyclor in Section~\ref{daggers-duals-conjugation} to show how conjugation and dagger are related in the presence of dualization.

%%%%%%%%%%%%%%%%%%%%%%%%%%%%%%%%%%%%%%%%%%%%%%%%

\section{Frobenius functors and daggers}
%%%%%%%%%%%%%%%%%%%%%%%%%%%%%%%%%%%%%%%%%%%%%%%%
We shall be interested here in linear functors between LDCs called Frobenius functors 
which come in various flavours, including mix functors and isomix functors, 
as illustrated in Figure \ref{linear-functor-family}.  These functors are directly related 
to the Frobenius monoidal functors of \cite{DP08} and they are referred to as degenerate 
linear functors in \cite{BPS12}.  Furthermore, we have already seen two rather basic examples, 
namely, ${\sf Mx}_\uparrow$ and ${\sf Mx}_\downarrow$.

\begin{figure}[ht]
\begin{center}
\begin{tikzpicture} [scale=1.5]
	\begin{pgfonlayer}{nodelayer}
		\node [style=none] (0) at (-3, 2) {};
		\node [style=none] (1) at (-3, -0) {};
		\node [style=none] (2) at (1, -0) {};
		\node [style=none] (3) at (1, 2) {};
		\node [style=none] (4) at (-4, 3) {};
		\node [style=none] (5) at (-4, -1) {};
		\node [style=none] (6) at (2, -1) {};
		\node [style=none] (7) at (2, 3) {};
		\node [style=none] (8) at (-5, -2) {};
		\node [style=none] (9) at (3, -2) {};
		\node [style=none] (10) at (3, 4) {};
		\node [style=none] (11) at (-5, 4) {};
		\node [style=none] (12) at (-6, 5) {};
		\node [style=none] (13) at (4, 5) {};
		\node [style=none] (14) at (4, -3) {};
		\node [style=none] (15) at (-6, -3) {};
		\node [style=none] (16) at (-1, 4.5) {Linear functors};
		\node [style=none] (17) at (-1, 3.5) {Frobenius functors};
		\node [style=none] (18) at (-1, 2.5) {Mix functors};
		\node [style=none] (19) at (-1, 1.5) {Isomix functors};
		\node [style=none] (34) at (-1, 1) {(Normal functors)};
		\node [style=none] (20) at (-1, 0.5) {$m_\top F(\m^{-1})n_\bot = \m^{-1}$};
		\node [style=none] (21) at (-1, -0.5) {$n_\bot\m m_\top = F(\m)$};
		\node [style=none] (22) at (-1, -1.5) {$F_\ox = F_\oa; m_\ox = \nu_\oa^L = \nu_\oa^R; n_\oa = \nu_\ox^L = \nu_\ox^R$};
		\node [style=none] (23) at (-1, -2.75) {};
		\node [style=none] (24) at (0, 6) {};
		\node [style=none] (25) at (0, -4) {};
		\node [style=none] (26) at (5.75, 6) {};
		\node [style=none] (27) at (5.75, -4) {};
		\node [style=none] (28) at (-0.9999999, 7) {};
		\node [style=none] (29) at (-0.9999999, -5) {};
		\node [style=none] (30) at (6.75, -5) {};
		\node [style=none] (31) at (6.75, 7) {};
		\node [style=none] (32) at (3, 6.75) {cyclic functors};
		\node [style=none] (33) at (3, 5.5) {symmetric functors};
	\end{pgfonlayer}
	\begin{pgfonlayer}{edgelayer}
		\draw (0.center) to (1.center);
		\draw (1.center) to (2.center);
		\draw (2.center) to (3.center);
		\draw (3.center) to (0.center);
		\draw (4.center) to (5.center);
		\draw (5.center) to (6.center);
		\draw (6.center) to (7.center);
		\draw (7.center) to (4.center);
		\draw (11.center) to (10.center);
		\draw (10.center) to (9.center);
		\draw (9.center) to (8.center);
		\draw (8.center) to (11.center);
		\draw (12.center) to (13.center);
		\draw (13.center) to (14.center);
		\draw (12.center) to (15.center);
		\draw (15.center) to (14.center);
		\draw[dotted]  (28.center) to (29.center);
		\draw[dotted] (29.center) to (30.center);
		\draw[dotted]  (30.center) to (31.center);
		\draw[dotted]  (31.center) to (28.center);
		\draw[dotted]  (24.center) to (25.center);
		\draw[dotted]  (25.center) to (27.center);
		\draw[dotted]  (27.center) to (26.center);
		\draw[dotted] (26.center) to (24.center);
	\end{pgfonlayer}
\end{tikzpicture}
\end{center}
\caption{Linear functor family}
\label{linear-functor-family}
\end{figure}

Frobenius functors preserve linear duals and with an additional coherence condition they preserve the mix map.  The coherence requirements for a dagger on an LDC are implied by requiring that the dagger functor be a Frobenius involutive equivalence.  Once the dagger  is understood we can consider $\dagger$-mix categories and their functors which we shall take to be mix Frobenius functors with a further requirement concerning the preservation of the dagger.

\subsection{Frobenius functors}
\label{Sec: frobenius functors}

\begin{defi}
A {\bf Frobenius functor} is a linear functor $F$ such that:
\begin{enumerate}[{\bf \small [FLF.1]}]
\item $F_\ox = F_\oa $
\item $m_\ox = \nu_\oa^R = \nu_\oa^L $ 
\item $n_\oa = \nu_\ox^L = \nu_\ox^R$
\end{enumerate}
\end{defi}

The left and right linear strengths of $\ox$ and $\oa$ coinciding with the $m_\ox$ and $n_\oa$ respectively means that in the diagrammatic calculus, ports can be moved around freely:
\[
\begin{tikzpicture}
	\begin{pgfonlayer}{nodelayer}
		\node [style=none] (0) at (-2, 2) {};
		\node [style=none] (1) at (-1.5, 2) {};
		\node [style=none] (2) at (-2.25, 2) {};
		\node [style=none] (3) at (-2.25, 1) {};
		\node [style=none] (4) at (-0.75, 1) {};
		\node [style=none] (5) at (-0.75, 2) {};
		\node [style=none] (6) at (-2, 2.75) {};
		\node [style=none] (7) at (-1, 2.75) {};
		\node [style=none] (61) at (-2.75, 2.75) {$F_\oa(A)$};
		\node [style=none] (71) at (-0.25, 2.75) {$F_\ox(B)$};
		\node [style=none] (8) at (-1.5, 0.25) {};
		\node [style=none] (81) at (-2.25, 0) {$F_\oa(A \ox B) $};
		\node [style=ox] (9) at (-1.5, 1.5) {};
		\node [style=none] (10) at (-2, 1.25) {$F$};
	\end{pgfonlayer}
	\begin{pgfonlayer}{edgelayer}
		\draw [in=-90, out=-90, looseness=1.25] (0.center) to (1.center);
		\draw [bend right=15, looseness=1.00] (6.center) to (9);
		\draw [bend right=15, looseness=0.75] (9) to (7.center);
		\draw (2.center) to (3.center);
		\draw (3.center) to (4.center);
		\draw (4.center) to (5.center);
		\draw (5.center) to (2.center);
		\draw (9) to (8.center);
	\end{pgfonlayer}
\end{tikzpicture} = \begin{tikzpicture}
	\begin{pgfonlayer}{nodelayer}
		\node [style=none] (0) at (-1.5, 2) {};
		\node [style=none] (1) at (-1, 2) {};
		\node [style=none] (2) at (-2.25, 2) {};
		\node [style=none] (3) at (-2.25, 1) {};
		\node [style=none] (4) at (-0.75, 1) {};
		\node [style=none] (5) at (-0.75, 2) {};
		\node [style=none] (6) at (-2, 2.75) {};
		\node [style=none] (7) at (-1, 2.75) {};
		\node [style=none] (61) at (-2.75, 2.75) {$F_\ox(A)$};
		\node [style=none] (71) at (-0.25, 2.75) {$F_\oa(B)$};
		\node [style=none] (8) at (-1.5, 0.25) {};
		\node [style=none] (81) at (-2.25, 0) {$F_\oa(A \ox B) $};
		\node [style=ox] (9) at (-1.5, 1.5) {};
		\node [style=none] (10) at (-2, 1.25) {$F$};
	\end{pgfonlayer}
	\begin{pgfonlayer}{edgelayer}
		\draw [in=-90, out=-90, looseness=1.25] (0.center) to (1.center);
		\draw [bend right=15, looseness=1.00] (6.center) to (9);
		\draw [bend right=15, looseness=0.75] (9) to (7.center);
		\draw (2.center) to (3.center);
		\draw (3.center) to (4.center);
		\draw (4.center) to (5.center);
		\draw (5.center) to (2.center);
		\draw (9) to (8.center);
	\end{pgfonlayer}
\end{tikzpicture} = \begin{tikzpicture}
	\begin{pgfonlayer}{nodelayer}
		\node [style=none] (0) at (-2.25, 2) {};
		\node [style=none] (1) at (-2.25, 1) {};
		\node [style=none] (2) at (-0.75, 1) {};
		\node [style=none] (3) at (-0.75, 2) {};
		\node [style=none] (4) at (-2, 2.75) {};
		\node [style=none] (5) at (-1, 2.75) {};
		\node [style=none] (41) at (-2.75, 2.75) {$F_\ox(A)$};
		\node [style=none] (51) at (-0.25, 2.75) {$F_\ox(B)$};
		\node [style=none] (6) at (-1.5, 0.25) {};
		\node [style=none] (61) at (-2, 0) {$F_\ox(A \ox B)$};
		\node [style=ox] (7) at (-1.5, 1.5) {};
		\node [style=none] (8) at (-2, 1.25) {$F$};
		\node [style=none] (9) at (-1.75, 1) {};
		\node [style=none] (10) at (-1.25, 1) {};
	\end{pgfonlayer}
	\begin{pgfonlayer}{edgelayer}
		\draw [bend right=15, looseness=1.00] (4.center) to (7);
		\draw [bend right=15, looseness=0.75] (7) to (5.center);
		\draw (0.center) to (1.center);
		\draw (1.center) to (2.center);
		\draw (2.center) to (3.center);
		\draw (3.center) to (0.center);
		\draw (7) to (6.center);
		\draw [in=90, out=90, looseness=1.25] (9.center) to (10.center);
	\end{pgfonlayer}
\end{tikzpicture}  ~~~~~~~~ \begin{tikzpicture}
	\begin{pgfonlayer}{nodelayer}
		\node [style=none] (0) at (-2.25, 1) {};
		\node [style=none] (1) at (-2.25, 2) {};
		\node [style=none] (2) at (-0.75, 2) {};
		\node [style=none] (3) at (-0.75, 1) {};
		\node [style=none] (4) at (-2, 0.25) {};
		\node [style=none] (5) at (-1, 0.25) {};
		\node [style=none] (6) at (-1.5, 2.75) {};
		\node [style=none] (41) at (-2.75, 0.25) {$F_\ox(A)$};
		\node [style=none] (51) at (-0.25, 0.25) {$F_\oa(B)$};
		\node [style=none] (61) at (-2, 3) {$F_\ox(A \oa B)$};
		\node [style=oa] (7) at (-1.5, 1.5) {};
		\node [style=none] (8) at (-2, 1.75) {$F$};
		\node [style=none] (9) at (-2, 1) {};
		\node [style=none] (10) at (-1.5, 1) {};
	\end{pgfonlayer}
	\begin{pgfonlayer}{edgelayer}
		\draw [bend left=15, looseness=1.00] (4.center) to (7);
		\draw [bend left=15, looseness=0.75] (7) to (5.center);
		\draw (0.center) to (1.center);
		\draw (1.center) to (2.center);
		\draw (2.center) to (3.center);
		\draw (3.center) to (0.center);
		\draw (7) to (6.center);
		\draw [in=90, out=90, looseness=1.25] (9.center) to (10.center);
	\end{pgfonlayer}
\end{tikzpicture} = \begin{tikzpicture}
	\begin{pgfonlayer}{nodelayer}
		\node [style=none] (0) at (-2.25, 1) {};
		\node [style=none] (1) at (-2.25, 2) {};
		\node [style=none] (2) at (-0.75, 2) {};
		\node [style=none] (3) at (-0.75, 1) {};
		\node [style=none] (4) at (-2, 0.25) {};
		\node [style=none] (5) at (-1, 0.25) {};
		\node [style=none] (41) at (-2.75, 0.25) {$F_\oa(A)$};
		\node [style=none] (51) at (-0.25, 0.25) {$F_\ox(B)$};
		\node [style=none] (61) at (-2, 3) {$F_\ox(A \oa B)$};
		\node [style=none] (6) at (-1.5, 2.75) {};
		\node [style=oa] (7) at (-1.5, 1.5) {};
		\node [style=none] (8) at (-2, 1.75) {$F$};
		\node [style=none] (9) at (-1.5, 1) {};
		\node [style=none] (10) at (-1, 1) {};
	\end{pgfonlayer}
	\begin{pgfonlayer}{edgelayer}
		\draw [bend left=15, looseness=1.00] (4.center) to (7);
		\draw [bend left=15, looseness=0.75] (7) to (5.center);
		\draw (0.center) to (1.center);
		\draw (1.center) to (2.center);
		\draw (2.center) to (3.center);
		\draw (3.center) to (0.center);
		\draw (7) to (6.center);
		\draw [in=90, out=90, looseness=1.25] (9.center) to (10.center);
	\end{pgfonlayer}
\end{tikzpicture} = \begin{tikzpicture}
	\begin{pgfonlayer}{nodelayer}
		\node [style=none] (0) at (-2.25, 1) {};
		\node [style=none] (1) at (-2.25, 2) {};
		\node [style=none] (2) at (-0.75, 2) {};
		\node [style=none] (3) at (-0.75, 1) {};
		\node [style=none] (4) at (-2, 0.25) {};
		\node [style=none] (5) at (-1, 0.25) {};
		\node [style=none] (41) at (-2.75, 0.25) {$F_\oa(A)$};
		\node [style=none] (51) at (-0.25, 0.25) {$F_\oa(B)$};
		\node [style=none] (61) at (-2, 3) {$F_\oa(A \oa B)$};
		\node [style=none] (6) at (-1.5, 2.75) {};
		\node [style=oa] (7) at (-1.5, 1.5) {};
		\node [style=none] (8) at (-2, 1.75) {$F$};
		\node [style=none] (9) at (-1.25, 2) {};
		\node [style=none] (10) at (-1.75, 2) {};
	\end{pgfonlayer}
	\begin{pgfonlayer}{edgelayer}
		\draw [bend left=15, looseness=1.00] (4.center) to (7);
		\draw [bend left=15, looseness=0.75] (7) to (5.center);
		\draw (0.center) to (1.center);
		\draw (1.center) to (2.center);
		\draw (2.center) to (3.center);
		\draw (3.center) to (0.center);
		\draw (7) to (6.center);
		\draw [in=-90, out=-90, looseness=1.25] (9.center) to (10.center);
	\end{pgfonlayer}
\end{tikzpicture}
\]
\[ \nu_\oa^L = \nu_\oa^R = m_\ox  ~~~~~~~~~~~~~~~~~~~~~~~~~~~~~~~~~~~~  \nu_\ox^L = \nu_\ox^R = n_\oa \]
This implies that the ports can be omitted in the circuits.

A Frobenius functor is {\bf symmetric} if as a linear functor it preserves the symmetries of the tensor and par.  

\begin{lem}
\label{Lemma: Frobenius}
Suppose $\X$ and $\Y$ are LDCs. The following are equivalent:
\begin{enumerate}[(a)]
\item $F: \X \to \Y$ is a Frobenius linear functor.
\item $F$ is $\ox$-monoidal and $\oa$-comonoidal such that 
\[ 
\hspace{-1cm}
\xymatrixcolsep{2pc}
\xymatrix{
F(A) \ox F(B \oa C) \ar[r]^{1 \ox n_\oa} \ar[d]_{m_\ox} \ar@{}[dr]|{\tiny{\bf [F.1]}} & F(A) \ox (F(B) \oa F(C)) \ar[d]^{\partial^L} \\
F(A \ox (B \oa C)) \ar[d]_{F(\partial^L)} & (F(A) \ox F(B)) \oa F(C) \ar[d]^{m_\ox \oa 1} \\
F((A \oa B) \oa C) \ar[r]_{n_\oa} & F(A \oa B) \oa F(C)
}\] \[ \xymatrix{
F( A \oa B) \ox F(C) \ar[r]^{n_\oa \ox 1}  \ar[d]_{m_\ox} \ar@{}[dr]|{\tiny{\bf [F.2]}} & (F(A) \oa F(B)) \ox F(C) \ar[d]^{\partial^R} \\
F( (A \oa B) \ox C) \ar[d]_{F(\partial^R)} & F(A) \oa (F(B) \ox F(C)) \ar[d]^{1 \oa m_\ox} \\
F(A \oa (B \ox C)) \ar[r]_{n_\oa} & F(A) \oa F(B \ox C) 
}
\]
\end{enumerate}
\end{lem}
\begin{proof}
For (a)$\Rightarrow$(b), fix $F := F_\ox = F_\oa$,  then $F$ is $\ox$-monoidal and $\oa$-comonoidal. Conditions {\bf \small [F.1]} and {\bf \small [F.2]} are given by {\bf  \small [LF.5]-(a)} and {\bf \small [LF.5]-(b)}. For the other direction, define $F_\ox = F_\oa := F$. Then it is straightforward to check that all the axioms of Frobenius linear functors are satisfied by $(F_\ox, F_\oa)$.
\end{proof}

Conditions {\bf \small  [F.1]} and {\bf \small [F.2]} in Lemma \ref{Lemma: Frobenius} are diagrammatically represented as follows:
\[
{\bf [F.1]}~~~~~~~~
\begin{tikzpicture}
	\begin{pgfonlayer}{nodelayer}
		\node [style=oa] (0) at (-3, 2) {};
		\node [style=ox] (1) at (-4, 1) {};
		\node [style=none] (2) at (-4.5, 3) {};
		\node [style=none] (3) at (-3, 3) {};
		\node [style=none] (4) at (-2.5, -0) {};
		\node [style=none] (5) at (-4, -0) {};
		\node [style=none] (6) at (-4.75, 2.5) {};
		\node [style=none] (7) at (-4.75, 0.5) {};
		\node [style=none] (8) at (-2.25, 0.5) {};
		\node [style=none] (9) at (-2.25, 2.5) {};
		\node [style=none] (10) at (-2.5, 2.25) {$F$};
	\end{pgfonlayer}
	\begin{pgfonlayer}{edgelayer}
		\draw [in=15, out=-150, looseness=1.00] (0) to (1);
		\draw [in=-90, out=135, looseness=1.00] (1) to (2.center);
		\draw (0) to (3.center);
		\draw [in=90, out=-45, looseness=1.00] (0) to (4.center);
		\draw (1) to (5.center);
		\draw (6.center) to (9.center);
		\draw (9.center) to (8.center);
		\draw (8.center) to (7.center);
		\draw (7.center) to (6.center);
	\end{pgfonlayer}
\end{tikzpicture} = \begin{tikzpicture}
	\begin{pgfonlayer}{nodelayer}
		\node [style=oa] (0) at (-3, 2.5) {};
		\node [style=ox] (1) at (-4, 1) {};
		\node [style=none] (2) at (-4.5, 3.5) {};
		\node [style=none] (3) at (-3, 3.5) {};
		\node [style=none] (4) at (-2.5, -0) {};
		\node [style=none] (5) at (-4, -0) {};
		\node [style=none] (6) at (-4.75, 1.5) {};
		\node [style=none] (7) at (-4.75, 0.5) {};
		\node [style=none] (8) at (-3.25, 0.5) {};
		\node [style=none] (9) at (-3.25, 1.5) {};
		\node [style=none] (10) at (-3.75, 2) {};
		\node [style=none] (11) at (-2.25, 2) {};
		\node [style=none] (12) at (-2.25, 3) {};
		\node [style=none] (13) at (-3.75, 3) {};
		\node [style=none] (14) at (-2.5, 2.75) {$F$};
		\node [style=none] (15) at (-3.5, 1.25) {$F$};
	\end{pgfonlayer}
	\begin{pgfonlayer}{edgelayer}
		\draw [in=60, out=-135, looseness=1.00] (0) to (1);
		\draw [in=-90, out=135, looseness=1.00] (1) to (2.center);
		\draw (0) to (3.center);
		\draw [in=90, out=-45, looseness=1.00] (0) to (4.center);
		\draw (1) to (5.center);
		\draw (6.center) to (9.center);
		\draw (9.center) to (8.center);
		\draw (8.center) to (7.center);
		\draw (7.center) to (6.center);
		\draw (13.center) to (10.center);
		\draw (10.center) to (11.center);
		\draw (11.center) to (12.center);
		\draw (12.center) to (13.center);
	\end{pgfonlayer}
\end{tikzpicture} 
~~~~~~~~~~~~~~~~~~~~~~~~~~~~ {\bf [F.2]}~~~~~~~~ \begin{tikzpicture}
	\begin{pgfonlayer}{nodelayer}
		\node [style=oa] (0) at (-3.25, 1) {};
		\node [style=ox] (1) at (-4, 1.75) {};
		\node [style=none] (2) at (-4, 3) {};
		\node [style=none] (3) at (-2.75, 3) {};
		\node [style=none] (4) at (-3.25, -0) {};
		\node [style=none] (5) at (-4.5, -0) {};
		\node [style=none] (6) at (-4.75, 2.5) {};
		\node [style=none] (7) at (-4.75, 0.5) {};
		\node [style=none] (8) at (-2.25, 0.5) {};
		\node [style=none] (9) at (-2.25, 2.5) {};
		\node [style=none] (10) at (-2.5, 2.25) {$F$};
	\end{pgfonlayer}
	\begin{pgfonlayer}{edgelayer}
		\draw [in=-45, out=165, looseness=1.25] (0) to (1);
		\draw [in=-90, out=90, looseness=1.00] (1) to (2.center);
		\draw [in=-90, out=45, looseness=1.00] (0) to (3.center);
		\draw [in=90, out=-90, looseness=1.00] (0) to (4.center);
		\draw [in=90, out=-150, looseness=0.75] (1) to (5.center);
		\draw (6.center) to (9.center);
		\draw (9.center) to (8.center);
		\draw (8.center) to (7.center);
		\draw (7.center) to (6.center);
	\end{pgfonlayer}
\end{tikzpicture} = \begin{tikzpicture}
	\begin{pgfonlayer}{nodelayer}
		\node [style=oa] (0) at (-3.25, 0.5) {};
		\node [style=ox] (1) at (-4, 1.75) {};
		\node [style=none] (2) at (-4, 3) {};
		\node [style=none] (3) at (-2.75, 3) {};
		\node [style=none] (4) at (-3.25, -0.5) {};
		\node [style=none] (5) at (-4.5, -0.5) {};
		\node [style=none] (6) at (-2.5, 1) {};
		\node [style=none] (7) at (-2.5, -0) {};
		\node [style=none] (8) at (-3.75, -0) {};
		\node [style=none] (9) at (-3.75, 1) {};
		\node [style=none] (10) at (-2.5, 2.25) {};
		\node [style=none] (11) at (-3.25, 2.25) {};
		\node [style=none] (12) at (-3.25, 1.25) {};
		\node [style=none] (13) at (-4.5, 2.25) {};
		\node [style=none] (14) at (-4.5, 1.25) {};
		\node [style=none] (15) at (-3.4, 1.5) {$F$};
		\node [style=none] (16) at (-2.75, 0.25) {$F$};
	\end{pgfonlayer}
	\begin{pgfonlayer}{edgelayer}
		\draw [in=-45, out=135, looseness=1.25] (0) to (1);
		\draw [in=-90, out=90, looseness=1.00] (1) to (2.center);
		\draw [in=-90, out=45, looseness=1.00] (0) to (3.center);
		\draw [in=90, out=-90, looseness=1.00] (0) to (4.center);
		\draw [in=90, out=-150, looseness=0.75] (1) to (5.center);
		\draw (6.center) to (9.center);
		\draw (9.center) to (8.center);
		\draw (8.center) to (7.center);
		\draw (7.center) to (6.center);
		\draw (11.center) to (13.center);
		\draw (13.center) to (14.center);
		\draw (14.center) to (12.center);
		\draw (12.center) to (11.center);
	\end{pgfonlayer}
\end{tikzpicture}
\]

Frobenius functors compose: the composition is defined as the usual composition of linear functors \cite{CS97}. 

It is immediate from Lemma \ref{Lemma: linear adjoints} that Frobenius functors preserve linear duals.  In fact if $F: \X \to \Y$ is a Frobenius functor 
and $A \dashvv B$ is a linear dual, as the duals $F_\ox(A) \dashvv F_\oa(B)$ and $F_\oa(A) \dashvv F_\ox(B)$ now coincide,  we just obtain the 
one dual $F(A) \dashvv F(B)$.   In the case when the Frobenius functor is between cyclic $*$-autonomous categories we expect the 
functor to be  {\bf cyclor-preserving} in the following sense:
\[ \mbox{\bf [CFF]} ~~~~~\xymatrix{ F(X^{*}) \ar[d]_{\cong} \ar[rr]^{F(\psi)} & & F(\!\!~^{*}X) \ar[d]^{\cong} \\
                    F(X)^{*} \ar[rr]_{\psi} & & \!\!~^{*}F(X) } \]
where the left and right vertical arrows are respectively the maps:
\[ (u^R_\ox)^{-1} (\eta* \ox 1) \partial^R (1 \oa (m^F_\ox F(\epsilon*) n^F_\bot)) u^R_\oa 
       ~~~\mbox{and}~~~(u^R)^{-1} (1 \ox *\eta) \partial^L (m^F_\ox \oa 1)((F(*\epsilon)n_\bot^F) \oa 1) u^L_\oa \]
The cyclor preserving condition may be pictorially represented as follows:
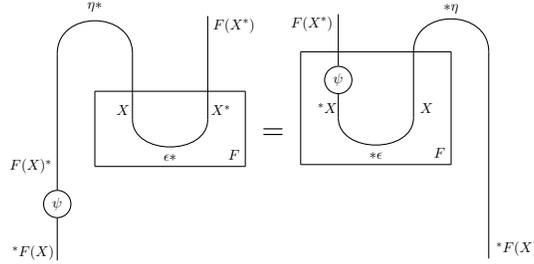
\[ %newcoh1

\]

\begin{lem}
\label{Lemma: cyclic Frob}
Suppose $F$ is a cyclor preserving Frobenius linear functor, then
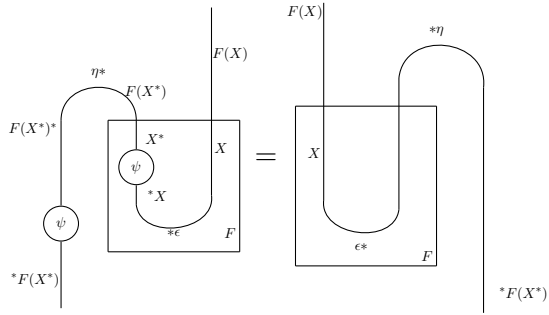
\[
%cyclor1
%
$
\end{proof}

\begin{defi}
Suppose $\X$ and $\Y$ are mix categories. $F: \X \to \Y$ is a {\bf mix functor} if it is a Frobenius functor such that 
\[
\mbox{\bf{[mix-FF]}}~~~~~
\xymatrix{
F(\bot) \ar@/_2pc/[rrr]_{F(\m)} \ar[r]^{n_\bot} & \bot \ar[r]^{\m} & \top \ar[r]^{m_\top} & F(\top) \\
}
\]
\end{defi}

This is diagrammatically represented as follows:
\[
%
 \]

\begin{lem}
\label{Lemma: Mix Frobenius linear functor}
Mix functors preserve the mix map:
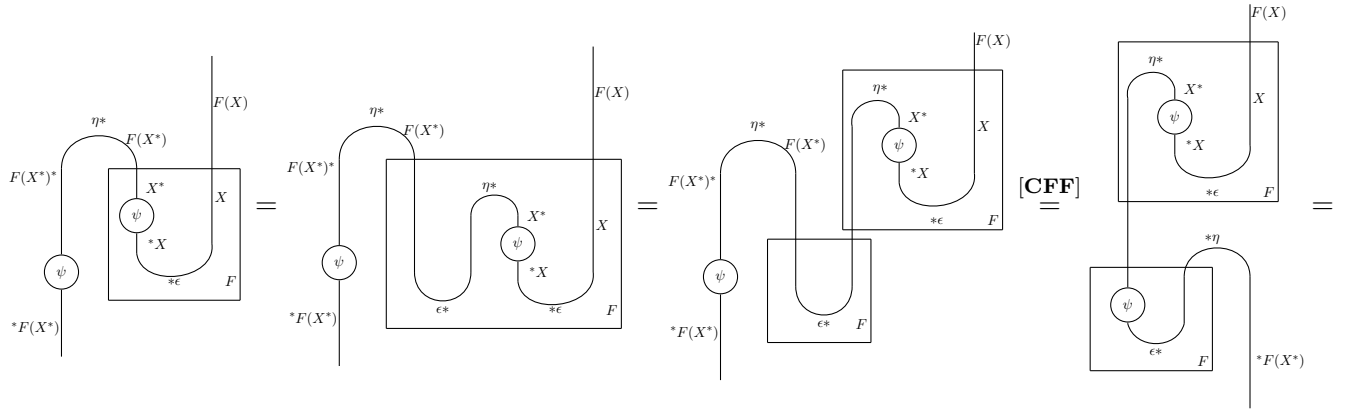
\[
\xymatrix{
F(A) \ox F(B) \ar[r]^{\mx} \ar[d]_{m_\ox} \ar[r]^{\mx} & F(A) \oa F(B) \\
F(A \ox B) \ar[r]_{F(\mx)} \ar[r]_{F(\mx)} & F(A \oa B) \ar[u]_{n_\oa}
}
\]
\end{lem}
\begin{proof}
\[
%
 
\]
\end{proof}

Linear natural isomorphisms between Frobenius functors $(\alpha_\ox, \alpha_\oa): F \to G$ often take a special form with $\alpha_\ox = \alpha_\oa^{-1}$: this allows the coherence 
requirements to be simplified.  The next results describe some basic circumstances in which this happens:
 
\begin{lem}
\label{Lemma: Frobenius linear transformation}
Suppose $F: \X \to \Y$, and $G: \X \to \Y$ are Frobenius linear functors and $\alpha := (\alpha_\ox, \alpha_\oa): F \Rightarrow G$ is a linear natural transformation.  Then, the following are equivalent:
\begin{enumerate}[(i)]
\item One of {\bf [nat.1](a)} or {\bf [nat.1](b)} holds, and one of  $\alpha_\ox$ or $\alpha_\oa$ is an isomorphism. \[
\mbox{\bf \small [nat.1]} ~~~
\xymatrixcolsep{5pc}
\xymatrix{
\top \ar[r]^{m_\top} \ar[dr]_{m_\top} & G(\top) \ar[d]^{\alpha_\oa} \ar@{}[dl]|(.35){\tiny{(a)}} \\
& F(\top)
} ~~~~~ \text{ or }~~~~~  \xymatrix{
F(\bot) \ar[r]^{\alpha_\ox} \ar[dr]_{n_\bot} & G(\bot) \ar[d]^{n_\bot} \ar@{}[dl]|(.35){\tiny{(b)}} \\
& \bot
}
\]

\item One of {\bf [nat.1](a)} or {\bf [nat.1](b)} holds and one of
\[ \mbox{\bf \small [nat.2]}~~~~ 
\xymatrix{
G(A) \ox F(B) \ar[r]^{1 \ox \alpha_\ox} \ar[d]_{\alpha_\oa \ox 1} \ar@{}[ddr]|{\tiny{(a)}} & G(A) \ox G(B) \ar[dd]^{m_\ox^G} \\
F(A) \ox F(B) \ar[d]_{m_\ox^F} & \\
F(A \ox B) \ar[r]_{\alpha_\ox} & G(A \ox B)
} ~~~ \text{ or }~~~ \xymatrix{
F(A) \ox G(B) \ar[r]^{\alpha_\ox \ox 1} \ar[d]_{1 \ox \alpha_\oa} \ar@{}[ddr]|{\tiny{(b)}}  & G(A) \ox G(B) \ar[dd]^{m_\ox^G} \\
F(A) \ox F(B) \ar[d]_{m_\ox^F} & \\
F(A \ox B) \ar[r]_{\alpha_\ox} & G(A \ox B)
}
\]
\[
\text{or}~~~~~ 
\xymatrix{
G(A \oa B) \ar[r]^{n_\oa^G} \ar[d]_{\alpha_\oa} \ar@{}[ddr]|{\tiny{(c)}}  & G(A) \oa G(B) \ar[dd]^{1 \oa \alpha_\oa} \\
F(A \oa B) \ar[d]_{n_\oa^F} & \\
F(A) \oa F(B) \ar[r]_{\alpha_\ox \oa 1} & G(A) \oa F(B)
} ~~~\text{or}~~~
\xymatrix{
G(A \oa B) \ar[r]^{n_\oa^G} \ar[d]_{\alpha_\oa} \ar@{}[ddr]|{\tiny{(d)}} & G(A) \oa G(B) \ar[dd]^{\alpha_\oa \oa 1} \\
F(A \oa B) \ar[d]_{n_\oa} & \\
F(A) \oa F(B) \ar[r]_{1 \ox \alpha_\ox} & F(A) \oa G(B)
}
\] holds.
\item $\alpha_\ox^{-1} = \alpha_\oa$
\item $\alpha' := (\alpha_\oa, \alpha_\ox): G \Rightarrow F$ is a linear transformation.
\end{enumerate}
\end{lem}

Conditions {\bf [nat.2]} are as follows in the graphical calculus:
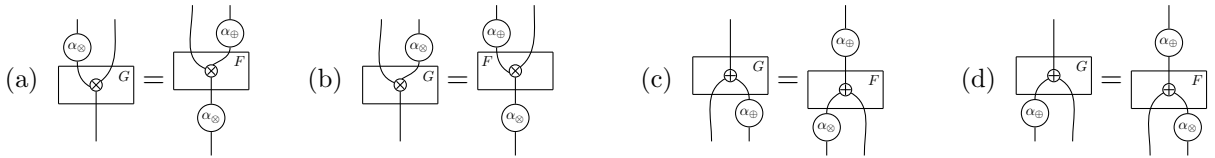
\[ \mbox{\small (a)}~~

\]

\begin{proof}
\begin{description}
\item[(i) $\Rightarrow$ (iii)] Here is the proof assuming {\bf \small [nat.1](a)} that $\alpha_\ox \alpha_\oa = 1$:
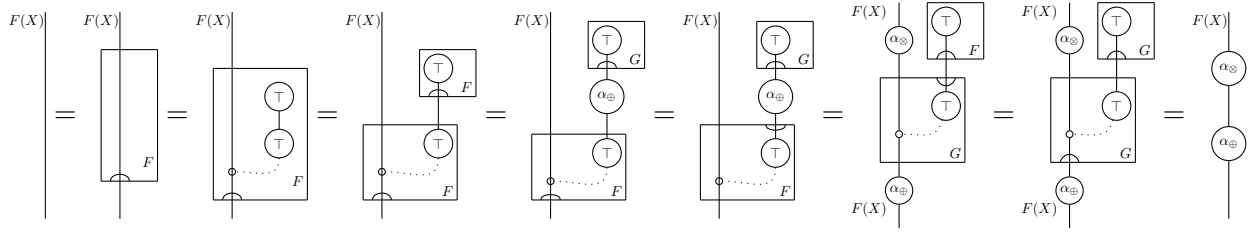
\[ \hspace{-0.75cm}
%
$

\smallskip

if either $\alpha_\ox$ or $\alpha_\oa$ are isomorphisms, then $\alpha_\oa \alpha_\ox =1$.
\item[(ii) $\Rightarrow$ (iii)]
The assumption of {\bf [nat.1](a)} or {\bf (b)} yields, as above, that $\alpha_\ox \alpha_\oa =1$. 
Using {\bf [nat.2](c)} for example gives $\alpha_\oa \alpha_\ox =1$:
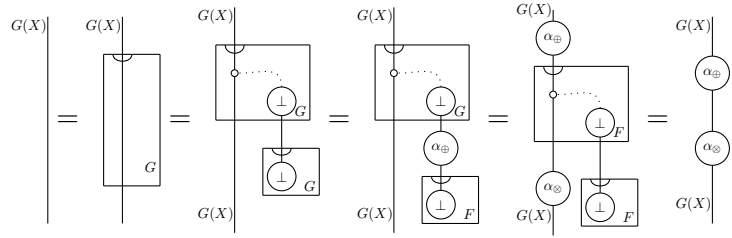
\[
%
\]
Since, $\alpha_\ox \alpha_\oa = 1$ and $\alpha_\ox \alpha_\oa = 1$ we have $\alpha_\ox = \alpha_\oa^{-1}$. The other combinations of rules are used in similar fashion.
\item[(iii) $\Rightarrow$ (iv)]
If $\alpha_\ox = \alpha_\oa^{-1}$, then 

$(\alpha_\oa \ox \alpha_\oa) m_\ox^F = m_\ox^G \alpha_\ox: G(A) \ox G(B) \to F(A \ox B)$
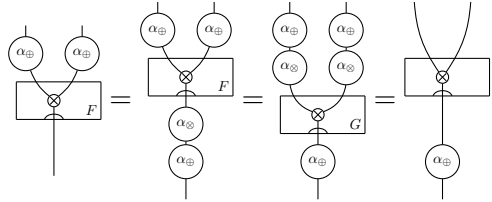
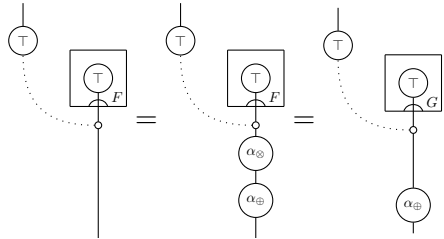
\[
%
\]
Thus, $\alpha_\ox$ is comonoidal. Similarly, it can be proven that $\alpha_\oa$ is monoidal. The axioms {\bf \small [LT.4] (a)}-{\bf \small (d)} for a linear transformation are satisfied for $(\alpha_\oa, \alpha_\ox)$ because $\alpha_\oa = \alpha_\ox^{-1}$. 

\item[(iv) $\Rightarrow$ ((i) and (ii))] The axioms {\bf \small [nat.1]} and
  {\bf \small [nat.2]} are given by the fact that $(\alpha_\oa, \alpha_\ox)$ is
  a linear transformation. \qedhere
\end{description}
\end{proof}

Frobenius functors between isomix categories are especially important in the development of dagger linearly distributive categories and they often satisfy an additional property:

\begin{defi}
A Frobenius functor between isomix categories is an {\bf isomix functor} in case it is a mix functor which satisfies, in addition, the following diagram:
\[ \mbox{\bf{[isomix-FF]}}~~~~~\xymatrix{ \top \ar@/^/[rrr]^{{\sf m}^{-1}} \ar[dr]_{m_\top} & & & \bot \\ & F(\top) \ar[r]_{F({\sf m}^{-1})} & F(\bot) \ar[ur]_{n_\bot} } \]
\end{defi}

Recall that a linear functor is {\bf normal} in case both $m_\top$ and $n_\bot$ are isomorphisms.  We observe:

\begin{lem} 
\label{Lemma: isomix functor}
For a mix Frobenius functor, $F: \X \to \Y$, between isomix categories the following are equivalent:
\begin{enumerate}[(i)]
\item $n_\bot: F(\bot) \to \bot$  or $m_\top: \top \to F(\top)$ is an isomorphism;
\item $F$ is a normal functor;
\item $F$ is an isomix functor.
\end{enumerate}
\end{lem}

\begin{proof}~
\begin{description}
\item[(i) $\Rightarrow$ (ii)]  Note that, as $F$ is a mix functor $F({\sf m}) = n_\bot {\sf m}~m_\top$.  As the mix map ${\sf m}$ is an isomorphism so is $F({\sf m})$ which implies that if  $n_\bot$ is an isomorphism then $m_\top$ must be an isomorphism and vice versa.  Thus, $F$ will be a normal functor. 
\item[(ii) $\Rightarrow$ (iii)] If $F$ is normal then $n_\bot$ and $m_\top$ are isomorphisms and so 
\[ \infer={m_\top F({\sf m}^{-1}) n_\bot = {\sf m}^{-1}}{\infer={F({\sf m}^{-1}) = m_\top^{-1} {\sf m}^{-1} n_\bot^{-1}}{F({\sf m}) = n_\bot {\sf m}~m_\top}} \]
\item[(iii) $\Rightarrow$ (i)] The mix-preservation for $F$ makes  $n_\bot$ a section (and $m_\top$ a retraction) while the isomix-preservation makes $m_\bot$ a retraction (and $m_\top$ a section). 
This means $n_\bot$ is an isomorphism ($m_\top$ is an isomorphism). \qedhere
\end{description}
\end{proof}

\begin{cor} \label{Corollary: normal-nat-iso}
$\alpha := (\alpha_\ox, \alpha_\oa)$ is a linear natural isomorphism between isomix Frobenius linear functors if and only if $\alpha_\ox = \alpha_\oa^{-1}$.
\end{cor}

\begin{proof}
Note  that if we can establish {\bf [nat.1](a)} or {\bf (b)} then we can prove that $\alpha_\ox\alpha_\oa =1$ and, as $\alpha_\ox$ is an isomorphism it follows that $\alpha_\oa\alpha_\ox =1$.
Thus, it suffices to show that {\bf [nat.1](a)} holds:
\[ m_\top \alpha_\oa G({\sf m}^{-1}) n_\bot = m_\top F({\sf m}^{-1})\alpha_\oa  n_\bot = m_\top F({\sf m}^{-1})  n_\bot = {\sf m}^{-1} = m_\top G({\sf m}^{-1}) n_\bot  \]
However, as $G({\sf m}^{-1}) n_\bot$ is an isomorphism, it follows that $m_\top \alpha_\oa = m_\top$.
\end{proof}

Lemma \ref{Lemma: Frobenius linear transformation} and Corollary \ref{Corollary: normal-nat-iso} are generalizations of \cite[Proposition 7]{DP08}. \cite[Proposition 7]{DP08} states the following:

Let $\X$ and $\Y$ be monoidal categories and $(\eta, \epsilon): A \dashvv B \in \X$. If $F,G: \X \to \Y$ are Frobenius monoidal functors with a natural transformation $\alpha: F \Rightarrow G$ which is both monoidal and comonoidal, then $\alpha_A$ is invertible. 

In Lemma \ref{Lemma: Frobenius linear transformation}, when $A \dashvv B \in \X$, then $\alpha_\oa$ is defined as follows:
\[
\alpha_\oa: G(A) \to F(A) = \begin{tikzpicture}
	\begin{pgfonlayer}{nodelayer}
		\node [style=circle] (0) at (-0.25, -0) {$\alpha_\ox$};
		\node [style=none] (1) at (1.25, -2.75) {};
		\node [style=none] (2) at (-0.25, -1) {};
		\node [style=none] (3) at (-1.75, -1) {};
		\node [style=none] (4) at (-1.75, 2.75) {};
		\node [style=circle] (5) at (-1, -1.75) {$\epsilon$};
		\node [style=none] (6) at (-0.25, 1) {};
		\node [style=none] (7) at (1.25, 1) {};
		\node [style=circle] (8) at (0.5, 1.75) {$\eta$};
		\node [style=none] (9) at (1.5, 2.5) {};
		\node [style=none] (10) at (1.5, 1) {};
		\node [style=none] (11) at (-0.5, 1) {};
		\node [style=none] (12) at (-0.5, 2.5) {};
		\node [style=none] (13) at (-2, -0.75) {};
		\node [style=none] (14) at (0, -0.75) {};
		\node [style=none] (15) at (0, -2.5) {};
		\node [style=none] (16) at (-2, -2.5) {};
		\node [style=none] (17) at (-2.25, 2.5) {$G(A)$};
		\node [style=none] (18) at (1.75, -2.5) {$F(A)$};
		\node [style=none] (19) at (0.25, 0.65) {$F(B)$};
		\node [style=none] (20) at (0.5, -0.5) {$G(B)$};
		\node [style=none] (21) at (1.25, 2.25) {$F$};
		\node [style=none] (22) at (-0.25, -2.25) {$G$};
	\end{pgfonlayer}
	\begin{pgfonlayer}{edgelayer}
		\draw (0) to (2.center);
		\draw (3.center) to (4.center);
		\draw [in=0, out=-90, looseness=1.25] (2.center) to (5);
		\draw [in=180, out=-90, looseness=1.25] (3.center) to (5);
		\draw (11.center) to (10.center);
		\draw (12.center) to (11.center);
		\draw (10.center) to (9.center);
		\draw (9.center) to (12.center);
		\draw (7.center) to (1.center);
		\draw [in=90, out=0, looseness=1.25] (8) to (7.center);
		\draw (6.center) to (0);
		\draw [in=180, out=90, looseness=1.25] (6.center) to (8);
		\draw (13.center) to (16.center);
		\draw (16.center) to (15.center);
		\draw (15.center) to (14.center);
		\draw (14.center) to (13.center);
	\end{pgfonlayer}
\end{tikzpicture}
\]

For these special linear isomorphisms with $\alpha_\ox = \alpha_\oa^{-1}$ we can simplify the coherence requirements:

\begin{lem} \label{simplifying-coherences}
Suppose $F$ and  $G$ are Frobenius functors and $\alpha: F \to G$ is a natural isomorphism then:
\begin{enumerate}[(i)]
\item If $\alpha: F \to G$ is $\ox$-monoidal and $\oa$-comonoidal then $(\alpha,\alpha^{-1})$ is a linear transformation;
\item If $F$ and $G$ are strong Frobenius functors and $\alpha$ is $\ox$-monoidal and $\oa$-monoidal then $(\alpha,\alpha^{-1})$ is a linear transformation.
\end{enumerate}
\end{lem}

\begin{proof}~
\begin{enumerate}[{\em (i)}]
\item If $\alpha$ is $\ox$-monoidal and $\oa$-comonoidal then so is $\alpha^{-1}$ supporting the possibility that it is a component of a linear transformation.
Considering  {\bf [LT.1]} we show that $(\alpha,\alpha^{-1})$ satisfies this requirement as:
\[ \xymatrix{ F(A \oa B) \ar[rr]^{\alpha_\ox}  \ar[d]_{n_\oa = \nu^R_\ox} & & G(A \oa B) \ar[d]^{n_\oa = \nu^R_\ox} \\
                    F(A) \oa F(B) \ar[dr]_{1 \oa \alpha} \ar[rr]^{\alpha \oa \alpha} & & G(A) \oa G(B)  \ar[dl]^{\alpha^{-1} \oa 1} \\
                    & F(A) \oa G(B) } \]
The remaining requirements follow in a similar manner.   
\item  When the laxors for the functors are isomorphisms then being monoidal implies being comonoidal. \qedhere
\end{enumerate}
\end{proof}

%%%%%%%%%%%%%%%%%%%%%%%%%%%%%%%%%%%%%%%%%%%%%%%%

\subsection{Dagger mix categories}

\label{Section: dagger LDC}

Conventionally, in categorical quantum mechanics a dagger is defined as a contravariant functor which is an involution that is stationary on objects $(A^\dagger = A)$. Before proceeding to define the dagger functor for LDCs, the notion of the opposite LDC and whence the notion of a contravariant linear functors have to be developed.  For LDCs we cannot expect the dagger to be stationary on objects, however, it is still possible that it can act like an involution.

If $(\X, \ox, \top, \oa, \bot)$ is a linear distributive category, the {\bf opposite linear distributive category} is $(\X, \ox, \top, \oa, \bot)^{\op} := (\X^{\op}, \oa, \bot, \ox, \top)$ where $\X^{\op}$ is the usual opposite category with the monoidal structures are flipped as follows: 
\[\ox^{\op} := \oa ~~~~~~~ \top^{\op} := \bot ~~~~~~~ \oa^{\op} := \ox ~~~~~~~ \bot^{\op} := \top\]
$(\_)^\op$ is an endo functor for the category of LDCs and linear functors. It is an involution:  

$(\X, \ox, \top, \oa, \bot)^{{\op}~{\op}} = (\X, \ox, \top, \oa, \bot) $. 

Let $(F_\ox, F_\oa): (\X, \ox, \top, \oa, \bot)^{\op} \to (\X, \ox, \top, \oa, \bot)$ be a linear functor. The opposite linear functor  $(F_\ox, F_\oa)^{\op}: (\X, \ox, \top, \oa, \bot) \to  (\X, \ox, \top, \oa, \bot)^{\op}$ given by the pair of opposite functors $(F_\oa^{\op}, F_\ox^{\op})$. Observe that $F^{\op}$ is  a mix Frobenius linear functor if and only if $F$ is.

\begin{defi}
\label{Definition: daggerLDC succinct}
A {\bf dagger linearly distributive category} ($\dagger$-LDC), is an LDC, $\X$, with a contravariant Frobenius linear functor $(\_)^\dagger: \X^{\op} \to \X$ which is a linear involutive equivalence   $(\_)^\dagger ~\dashvv ~ (\_)^{\dagger^{\op}}: \X^{\op} \to \X$.
\end{defi}

First note that saying this is an {\bf involutive} equivalence asserts that the unit and counit of the equivalence are the same (although one is in the opposite category).  Thus, the adjunction expands to take the form $(\imath,\imath): (\_)^\dagger ~\dashvv ~ (\_)^{\dagger^{\op}}: \X^{\op} \to \X$.  However, the unit and counit are linear natural transformations so $\imath$ expands to $\imath = (\imath_\ox,\imath_\oa)$.  As the dagger functor is 
a left adjoint, it is strong and, thus, is normal.  Furthermore, as the unit of an equivalence, $\imath$ is a linear natural isomorphism.  This means  $\imath = (\imath_\ox,\imath_\oa)$ satisfies the requirements of Lemma \ref{Lemma: Frobenius linear transformation}, implying that $\imath_\ox^{-1}  = \imath_\oa$.  Simplifying notation we shall set $\iota:= \imath_\oa$  so the unit linear transformation is $\imath := (\iota^{-1},\iota)$. We then can simplify the requirements of $\imath$ to the map  
$\iota: A \to (A^\dagger)^\dagger$ which we refer to as the {\bf involutor}.

A {\bf symmetric  $\dagger$-LDC} is a $\dagger$-LDC which is a symmetric LDC for which the dagger is a symmetric linear functor. A {\bf cyclic $\dagger$-$*$-autonomous category} is a 
$\dagger$-LDC with chosen left are right duals and a cyclor which is preserved by the dagger.  A $\dagger$-{\bf mix} category is a $\dagger$-LDC for which $(\_)^\dagger: \X^{\op} \to \X$ 
is a mix functor.  As the dagger functor is strong (and so normal) if the category is an isomix category then being {\bf $\dagger$-mix} already implies that 
the dagger is an isomix functor.  Thus, a {\bf $\dagger$-\bf isomix} category is a $\dagger$-mix category which happens to be an isomix category.

In the remainder of the section, we unfold the definition of a $\dagger$-isomix category and give the coherence requirements explicitly.

\begin{prop}
\label{Definition: daggerLDC elaborate}
A dagger linearly distributive category is an LDC with a functor $(\_)^\dag:\X^\op\to \X$ and 
natural isomorphisms 
\begin{align*}
\text{ \bf laxors: }  A^\dag \ox B^\dag &\xrightarrow{ \lambda_\ox} (A\oa B)^\dag ~~~~~ A^\dag \oa B^\dag \xrightarrow{ \lambda_\oa} (A\ox B)^\dag \\
\top &\xrightarrow{\lambda_\top} \bot^\dag ~~~~~~~~~~~~~~~~~~~~~ \bot \xrightarrow{\lambda_\bot} \top^\dag \\
\text{ \bf involutor: }  A &\xrightarrow{\iota} (A^\dag)^\dag 
\end{align*}
such that the following coherences hold:
\begin{enumerate}[{\bf [$\dagger$-ldc.1]}]
\item Interaction of $\lambda_\ox, \lambda_\oa$  with associators:
\[
\begin{tabular}{cc}
\xymatrix{
A^\dag \ox (B^\dag \ox C^\dag)                \ar@{->}[r]^{{a_\ox}^{-1}}    \ar@{->}[d]_{1 \ox  \lambda_\ox}   
  & (A^\dag \ox B^\dag) \ox C^\dag              \ar@{->}[d]^{\lambda_\ox \ox 1}  \\
A^\dag \ox ( B \oa C)^\dag                     \ar@{->}[d]_{\lambda_\ox}  
  & (A \oa B)^\dag \ox C^\dag                   \ar@{->}[d]^{\lambda_\ox}    \\
(A \oa (B \oa C))^\dag                          \ar@{->}[r]_{a_\oa^\dag}
  & ( (A\oa B) \oa C)^\dag
} &\xymatrix{
A^\dag \oa (B^\dag \oa C^\dag)                \ar@{->}[r]^{a_\oa^{-1}}    \ar@{->}[d]_{1 \oa \lambda_\oa }  
  & (A^\dag \oa B^\dag) \oa C^\dag              \ar@{->}[d]^{\lambda_\oa \oa 1}  \\
   A^\dag \oa (B \ox C)^\dag                 \ar@{->}[d]_{\lambda_\oa}  
  & (A \ox B)^\dag \oa C^\dag                    \ar@{->}[d]^{\lambda_\oa}    \\
(A \ox (B\ox C))^\dag                          \ar@{->}[r]_{a_\ox^\dag}
  & ((A\ox B) \ox C)^\dag
}
\end{tabular}
\]

\item Interaction of $\lambda_\top, \lambda_\bot$ with unitors: 
\begin{center}
\begin{tabular}{cc}
\xymatrix{
\top \ox A^\dag                            \ar@{->}[rr]^{\lambda_\top\ox 1} \ar@{->}[d]_{u_\ox^R}  
&  & \bot^\dag\ox A^\dag           \ar@{->}[d]^{\lambda_\ox}\\
A^\dag                                     
&  & (\bot \oa A)^\dag                      \ar@{<-}[ll]^{(u_\oa^R)^\dag}\\
} &
\xymatrix{
\bot \oa A^\dag                            \ar@{->}[rr]^{\lambda_\bot\oa 1} \ar@{->}[d]_{u_\oa^R} 
& & \top^\dag\oa A^\dag                    \ar@{->}[d]^{\lambda_\oa}\\
A^\dag                                     
&  & (\top \ox A)^\dag                      \ar@{<-}[ll]^{(u_\ox^R)^\dag}\\
} 
\end{tabular}
\end{center}
and two symmetric diagrams for $u_\ox^L$ and $u_\oa^L$ must also be satisfied.

\item Interaction of $\lambda_\ox, \lambda_\oa$ with linear distributors:
\begin{center}
\begin{tabular}{cc}
\xymatrix{
A^\dag \ox(B^\dag\oa C^\dag)              \ar@{->}[r]^{\partial^L} \ar@{->}[d]_{1\ox\lambda_\oa}
  & (A^\dag \ox B^\dag)\oa C^\dag         \ar@{->}[d]_{\lambda_\ox\oa 1}\\
A^\dag \ox (B\ox C)^\dag                  \ar@{->}[d]_{\lambda_\ox}
  & (A\oa B)^\dag \oa C^\dag              \ar@{->}[d]^{\lambda_\oa}\\
(A\oa (B\ox C))^\dag                      \ar@{->}[r]_{(\partial^R)^\dag}
  & ((A\oa B)\ox C)^\dag
} &
\xymatrix{
(A^\dag \oa B^\dag) \ox C^\dag     \ar@{->}[r]^{\partial^R} \ar@{->}[d]_{\lambda_\oa\ox 1} 
  & A^\dag \oa (B^\dag \ox C^\dag) \ar@{->}[d]_{1\oa \lambda_\ox}\\
(A\ox B)^\dag \ox C^\dag           \ar@{->}[d]_{\lambda_\ox}
  & A^\dag \oa (B\oa C)^\dag       \ar@{->}[d]^{\lambda_\oa}\\
((A\ox B)\oa C)^\dag               \ar@{->}[r]_{(\partial^L)^\dag}
  & (A\ox (B\oa C))^\dag
} 
\end{tabular}
\end{center}

\item Interaction of $\iota: A \rightarrow A^{\dagger\dagger}$ with $\lambda_\ox$, $\lambda_\oa$:
\begin{center}
\begin{tabular}{cc}
\xymatrix{
A\oa B                          \ar@{->}[r]^{\iota} \ar@{->}[d]_{\iota \oa \iota} 
  & ((A\oa B)^\dag)^\dag        \ar@{->}[d]^{\lambda_\ox^\dag}\\
(A^\dag)^\dag   \oa (B^\dag)^\dag                \ar@{->}[r]_{\lambda_\oa}
  & (A^\dag\ox B^\dag)^\dag
} &  \xymatrix{
A\ox B                          \ar@{->}[r]^{\iota} \ar@{->}[d]_{\iota \ox \iota} 
  & ((A\ox B)^\dag)^\dag        \ar@{->}[d]^{\lambda_\oa^\dag}\\ 
(A^\dag)^\dag  \ox (B^\dag)^\dag                \ar@{->}[r]_{\lambda_\ox}
  & (A^\dag\oa B^\dag)^\dag
}
\end{tabular}
\end{center}
\item Interaction of $\iota: A \rightarrow A^{\dagger\dagger}$ with $\lambda_\top$, $\lambda_\bot$:
\begin{center}
\begin{tabular}{cc}
$\begin{matrix} \xymatrix{
&\bot                   \ar@{->}[r]^{\iota} \ar@{->}[dr]_{\lambda_\bot}  
  & (\bot^\dag)^\dag   \ar@{->}[d]^{\lambda_\top^\dag}\\
&{}
  & \top^\dag 
} \end{matrix}$ & 
 $\begin{matrix} \xymatrix{
&\top                   \ar@{->}[r]^{\iota} \ar@{->}[dr]_{\lambda_\top} 
  & (\top^\dag)^\dag   \ar@{->}[d]^{\lambda_\bot^\dag} \\
&{}
  & \bot^\dag 
} \end{matrix}$
\end{tabular}
\end{center}
\item $\iota_{A^\dagger} = (\iota_A^{-1})^\dagger: A^\dagger \to A^{\dagger\dagger\dagger}$  
\end{enumerate}
\end{prop}

The structure is presented using strong monoidal laxors:  to form a linear functor the laxor $\lambda_\oa$ needs to be reversed by taking its inverse. Then, we have $\nu_\ox^l = \nu_\ox^r := \lambda_\oa^{-1}$ and $\nu_\oa^l = \nu_\oa^r := \lambda_\ox$.  Once this adjustment is made all the required coherences for $\dagger$ to be a linear functor are present.

Note that {\bf [$\dagger$-ldc.6]} equivalently expresses the triangle identities of the adjunction $(\iota, \iota):: \dagger^{\op} \dashv \dagger : \X^{\op} \to \X$.   

The coherences for the involutor asserts that it is a monoidal transformation for both the tensor and par: by Lemma \ref{simplifying-coherences} (ii) this suffices to show that it is a linear transformation.

A {\bf symmetric $\dagger$-LDC} is a $\dagger$-LDC which is a symmetric LDC and for which the following additional diagrams commute:
\begin{enumerate}[{\bf [$\dagger$-ldc.7]}]
\item Interaction of $\lambda_\ox , \lambda_\oa$ with symmetry maps: 
\[
\begin{tabular}{cc}
\xymatrix{
A^\dag \ox B^\dag                          \ar@{->}[r]^{\lambda_\ox} \ar@{->}[d]_{c_\ox}         
  & (A\oa B)^\dag                          \ar@{->}[d]^{c_\oa^\dag}\\
B^\dag \ox A^\dag                          \ar@{->}[r]_{\lambda_\ox}
  & (B\oa A)^\dag\\
} & \xymatrix{
A^\dag \oa  B^\dag                          \ar@{->}[r]^{\lambda_\oa} \ar@{->}[d]_{c_\oa}  
  & (A\ox B)^\dag                          \ar@{->}[d]^{c_\ox^\dag}\\
B^\dag \oa A^\dag                          \ar@{->}[r]_{\lambda_\oa}
  & (B\ox A)^\dag\\
}
\end{tabular}
\]
\end{enumerate}

A {\bf $\dagger$-mix category} is a $\dagger$-LDC which has a mix map and satisfies the following additional coherence:

\[ \mbox{\bf [$\dagger$-\text{mix}]}  ~~~~\begin{array}[c]{c} 
\xymatrix{
\bot                 \ar@{->}[r]^{{\sf m}} \ar@{->}[d]_{\lambda_\bot}  
  & \top             \ar@{->}[d]^{\lambda_\top}\\
\top^\dag            \ar@{->}[r]_{{\sf m}^\dag}
  & \bot^\dag
} 
\end{array} \]
If ${\sf m}$ is an isomorphism, then $\X$ is a {\bf $\dagger$-isomix category} and, since $(\_)^\dagger$ is normal, $(\_)^\dagger$ is an isomix Frobenius functor.

\begin{lem}
\label{lemma: mixdagger}
Suppose $\X$ is a $\dagger$-mix category then the following diagram commutes:
\[
\xymatrix{ A^\dag \ox B^\dag  \ar[r]^{\mx} \ar[d]_{\lambda_\ox}&  A^\dag \oa B^\dag \ar[d]^{\lambda_\oa} \\
                (A \oa B)^\dag \ar[r]_{\mx^\dag}  & (A \ox B)^\dag }
\]
\end{lem}
\begin{proof}
The proof follows directly from Lemma \ref{Lemma: Mix Frobenius linear functor}.
\end{proof}

 With respect to its applications to quantum theory, this article primarily focuses on $\dagger$-isomix categories. As we will see in Section~\ref{Sec: unitary}, the notion of unitary objects and unitary isomorphisms is supported only within a $\dagger$-isomix category.

It is useful to observe that objects in the core are closed under taking the dagger and duals.

\begin{lem}
\label{Lemma: mixdagger}
Suppose $\X$ is a $\dagger$-mix category and $A \in \Core(\X)$ then $A^\dagger \in$ $\Core(\X)$.
\end{lem}
\begin{proof}
The natural transformation $A^\dagger \ox X \xrightarrow{\mx} A^\dagger \oa X$ is an isomorphism as follows:
\[
\xymatrix{
A^\dagger \ox X \ar[r]^{1 \ox \iota} \ar[d]_{\mx} \ar@{}[dr]|{\scalebox{0.95}{\tiny\bf (nat. {\sf mx})}}
& A^\dagger \ox X^{\dagger\dagger} \ar[r]^{\lambda_\ox} \ar[d]_{\mx} \ar@{}[dr]|{\scalebox{.95}{\tiny\bf {lem. \ref{lemma: mixdagger}}}}
& (A \oa X^\dagger)^\dagger \ar[d]^{\mx^\dagger} \\
A^\dagger \oa X \ar[r]_{1 \oa \iota} 
& A^\dagger \oa X^{\dagger \dagger} \ar[r]_{\lambda_\oa}
&(A \ox X^\dagger)^\dagger
}  \]
\end{proof}

\begin{lem}
Let $\X$ be $\dagger$-LDC. If $A \dashvv B$ then $B^\dagger \dashvv A^\dagger$.
\end{lem}
\begin{proof}
The statement  follows from Lemma \ref{Lemma: linear adjoints}: Frobenius functors preserve linear adjoints. Explicitly, if $(\eta,\epsilon): A \dashvv B$ then $(\lambda_\top\epsilon^\dag\lambda_\oa^{-1},\lambda_\ox\eta^\dagger \lambda_\bot^{-1}): B^\dagger \dashvv A^\dagger$. 
\end{proof}

Suppose $\X$ is a $\dagger$-$*$-autonomous category and $(\eta*, \epsilon*): A^* \dashvv A$, then $((\epsilon*)^\dagger, (\eta*)^\dagger): A^\dagger \dashvv (A^*)^\dagger$, where $((\epsilon*)^\dagger, (\eta*)^\dagger) :=   (\lambda_\top\epsilon*^\dag\lambda_\oa^{-1},\lambda_\ox\eta*^\dagger \lambda_\bot^{-1})$. We draw $(\epsilon*)^\dagger$ and $(*\epsilon)^\dagger$ as dagger cups, and $(\eta*)^\dagger$ and $(*\eta)^\dagger$ as dagger caps which are pictorially represented as follows:
 \[
 
 \]

A $\dagger$-$*$-autonomous category is a {\bf cyclic} $\dagger$-$*$-autonomous category when the dagger preserves the cyclor in the following sense. 
\[ %last1
%
 \]

\begin{lem} 
\label{Lemma: cyclic dagger}
In a cyclic, $\dagger$-$*$-autonomous category,
\[ %last3
%
\]
\end{lem}
\begin{proof}
Proved by direct application of Lemma \ref{Lemma: cyclic Frob}.
\end{proof}

%%%%%%%%%%%%%%%%%%%%%%%%%%%%%%%%%%%%%%%%%%%
\subsection{Dagger functor box}

Suppose $\X$ is a $\dagger$-LDC and  $f: A \rightarrow B \in \X$. Then, the map $f^\dagger: B^\dagger \rightarrow A^\dagger$ is graphically depicted as follows:
\[ %
 \]
The rectangle is a functor box for the $\dag$-functor.  Notice how we use vertical mirroring to express the contravariance of the $\dag$-functor. By the functoriality of $(\_)^\dagger$, we have:
%
.

These contravariant functor boxes compose... contravariantly. Given maps $f: A \to B$ and $g: B \to C$:
\[
%
\]

The following are the representations of the basic natural isomorphisms of a $\dagger$-LDC:
\[
\begin{tabular}{L L}
\lambda_\top: \top \rightarrow \bot^\dagger = 
\begin{tikzpicture}
	\begin{pgfonlayer}{nodelayer}
		\node [style=none] (0) at (-2, 0) {};
		\node [style=none] (1) at (-2, -1) {};
		\node [style=none] (2) at (-1, -1) {};
		\node [style=none] (3) at (-1, 0) {};
		\node [style=circle] (4) at (-1.5, -0.5) {$\bot$};
		\node [style=none] (5) at (-1.25, 0) {};
		\node [style=none] (6) at (-1.75, -1) {};
		\node [style=none] (7) at (-1.5, -2) {};
		\node [style=none] (8) at (-1.5, -1) {};
		\node [style=circle, scale=0.4] (9) at (-1.5, -1.5) {};
		\node [style=circle] (10) at (-2.5, -0.5) {$\top$};
		\node [style=none] (11) at (-2.5, 0.5) {};
	\end{pgfonlayer}
	\begin{pgfonlayer}{edgelayer}
		\draw [style=none] (0.center) to (1.center);
		\draw [style=none] (1.center) to (2.center);
		\draw [style=none] (2.center) to (3.center);
		\draw [style=none] (3.center) to (0.center);
		\draw [style=none] (4) to (5.center);
		\draw [style=none] (7.center) to (8.center);
		\draw [style=none] (11.center) to (10);
		\draw [densely dotted, bend right=45, looseness=1.25] (10) to (9);
	\end{pgfonlayer}
\end{tikzpicture}  & 
\lambda_\top^{-1}: \bot^\dagger \rightarrow \top =
\begin{tikzpicture}
	\begin{pgfonlayer}{nodelayer}
		\node [style=none] (0) at (-2, 0) {};
		\node [style=none] (1) at (-2, -1) {};
		\node [style=none] (2) at (-1, -1) {};
		\node [style=none] (3) at (-1, 0) {};
		\node [style=circle] (4) at (-1.5, -0.5) {$\bot$};
		\node [style=none] (5) at (-1.25, 0) {};
		\node [style=none] (6) at (-1.75, -1) {};
		\node [style=none] (7) at (-1.5, -1) {};
		\node [style=circle] (8) at (-0.5, -0.25) {$\top$};
		\node [style=none] (9) at (-1.5, 0) {};
		\node [style=none] (10) at (-1.5, 0.5) {};
		\node [style=none] (11) at (-0.5, -2) {};
		\node [style=circle, scale=0.4] (12) at (-0.5, -1.25) {};
	\end{pgfonlayer}
	\begin{pgfonlayer}{edgelayer}
		\draw [style=none] (0.center) to (1.center);
		\draw [style=none] (1.center) to (2.center);
		\draw [style=none] (2.center) to (3.center);
		\draw [style=none] (3.center) to (0.center);
		\draw [style=none] (4) to (6.center);
		\draw [style=none] (10.center) to (9.center);
		\draw [style=none] (8) to (11.center);
		\draw [densely dotted, bend right=45, looseness=1.00] (7.center) to (12);
		\draw [densely dotted] (4) to (5.center);
	\end{pgfonlayer}
\end{tikzpicture}  \\
\end{tabular}
\]

Dagger boxes interact with involutor $A \xrightarrow{\iota} A^{\dagger\dagger}$ as follows:
$$
%
$$

%the issue is more nuanced than this, i think
It is important to note that one may not have a legal proof net inside a $\dagger$-box. This complicates the correctness criterion. However, the required correctness criterion is discussed in \cite{MP05}.

%%%%%%%%%%%%%%%%%%%%%%%%%%%%%%%%%%%%%%%%%%%
\subsection{Functors for $\dagger$-linearly distributive categories}

Clearly the functors and transformations between $\dagger$-LDCs must ``preserve'' the dagger in some sense.  Precisely we have:

\begin{defi}
	$F: \X \to \Y$ is a {\bf $\dagger$-linear functor} between $\dagger$-LDCs when $F$ is a linear functor equipped with a linear natural isomorphism 
	$\rho^F= (\rho_\ox^F: F_\ox(A^\dagger) \to F_\oa(A)^\dagger ,\rho_\oa^F:  F_\ox(A)^\dagger \to F_\oa(A^\dagger))$ called the {\bf preservator}, 
	such that  the following diagrams commute:
	\[ 
	\xymatrix{
		F_\ox(X) \ar[r]^{\iota} \ar[d]_{F_\ox(\iota)} \ar@{}[dr]|{\mbox{\tiny {\bf [$\dagger$-LF.1]}}} & 
		F_\ox(X)^{\dagger \dagger} \ar@{<-}[d]^{(\rho^F_\oa)^\dagger} \\
		F_\ox(X^{\dagger \dagger}) \ar[r]_{\rho^F_\ox} & F_\oa(X^\dagger)^\dagger
	} ~~~~~~~~~ \xymatrix{
		F_\oa(X) \ar[r]^{\iota} \ar[d]_{F_\oa(\iota)}  \ar@{}[dr]|{\mbox{\tiny {\bf [$\dagger$-LF.2]}}} & F_\oa(X)^{\dagger \dagger} \ar[d]^{(\rho^F_\ox)^\dagger} \\
		F_\oa(X^{\dagger \dagger}) \ar@{<-}[r]_{\rho^F_\oa} & F_\ox(X^\dagger)^\dagger
	}
	\]
\end{defi}

In case that $F$ is a normal mix functor between $\dagger$-isomix categories, then by Lemma \ref{Lemma: isomix functor}, $F$ is an isomix functor and, therefore by Corollary \ref{Corollary: normal-nat-iso},  the preservators become inverses, $\rho^F_\ox = (\rho^F_\oa)^{-1}$.  
This means the squares {\bf [$\dagger$-LF.1]} and {\bf [$\dagger$-LF.2]} coincide to give a single condition for the tensor preservator:
\[ \xymatrix{
	F(X) \ar[r]^{\iota} \ar[d]_{F(\iota)} \ar@{}[dr]|{\mbox{\tiny {\bf [$\dagger$-isomix]}}} & 
	F(X)^{\dagger \dagger} \ar@{->}[d]^{(\rho^F_\ox)^\dagger} \\
	F(X^{\dagger \dagger}) \ar[r]_{\rho^F_\ox} & F(X^\dagger)^\dagger
} \]

%%%%%%%%%%%%%%%%

In case when $F$ is an isomix functor, by Lemma \ref{Lemma: Frobenius linear transformation}, $\rho := \rho_\ox$ is monoidal on $\ox$ and comonoidal on $\oa$:

\[
\bf{[P.1]} ~~~~~~~
\xymatrix{
	F(A^\dagger) \ox F(B^\dagger) \ar[r]^{\rho \ox \rho} \ar[d]_{m_\ox} \ar@{}[ddr]|{(a)} & F(A)^\dagger \ox F(B)^\dagger \ar[d]^{\lambda_\ox} \\
	F(A^\dagger \ox B^\dagger) \ar[d]_{F(\lambda_\ox)} & (F(A) \ox F(B))^\dagger \ar[d]^{n_\oa^\dagger} \\
	F((A \oa B)^\dagger) \ar[r]_{\rho} & (F(A\oa B))^\dagger
}~~~~~~~~~ \xymatrix{
	\top \ar[d]_{m_\top} \ar[dr]^{\lambda_\top} \ar@{}[ddr]|{(b)} & \\
	F(\top) \ar[d]_{F(\lambda_\top)} & \bot^\dagger \ar[dr]^{n_\bot^\dagger} \\
	F(\bot^\dagger) \ar[rr]_{\rho} & & (F(\bot))^\dagger
}
\]

\[
\bf{[P.2]} ~~~~~~~
\xymatrix{
	F((A \ox B)^\dagger) \ar[r]^{\rho} \ar[d]_{F(\lambda_\oa^{-1})} \ar@{}[ddr]|{(a)} & F(A \ox B)^\dagger \ar[d]^{m_\ox^\dagger} \\
	F(A^\dagger \oa B^\dagger) \ar[d]_{n_\oa^F} & (F(A) \ox F(B))^\dagger \ar[d]^{\lambda_\oa^{-1}} \\
	F(A^\dagger) \oa F(B^\dagger) \ar[r]_{\rho \oa \rho} & F(A)^\dagger \oa F(B)^\dagger
} ~~~~~~~~~ \xymatrix{
	F(\top^\dagger) \ar[dd]_{\rho} \ar[dr]^{F(\lambda_\bot^{-1})} \ar@{}[ddr]|{(b)} & \\
	& F(\bot) \ar[dr]^{n_\bot} & \\
	F(\top)^\dagger \ar[r]_{m_\top^\dagger} & \top^\dagger \ar[r]_{\lambda_\bot^{-1}} & \bot
}
\]

Pictorial representation of {\bf[P.2]-(a)} is as follows:
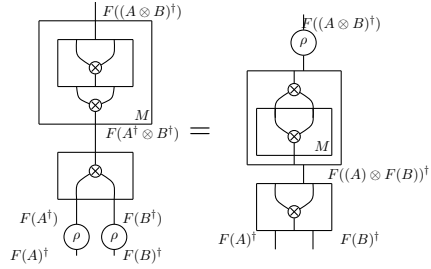
\[
\begin{tikzpicture} %rho-comon-b
\begin{pgfonlayer}{nodelayer}
\node [style=none] (0) at (1.75, 3.25) {};
\node [style=none] (1) at (2.25, -0) {};
\node [style=none] (2) at (1.25, -0) {};
\node [style=none] (3) at (0.25, -0) {};
\node [style=none] (4) at (0.75, 2.25) {};
\node [style=ox] (5) at (1.75, 1.5) {};
\node [style=none] (6) at (2.75, 2.25) {};
\node [style=none] (7) at (2.75, 1) {};
\node [style=none] (8) at (3.25, 2.75) {};
\node [style=none] (9) at (1.25, 2.25) {};
\node [style=circle] (10) at (2.25, -3) {$\rho$};
\node [style=none] (11) at (1.75, 2.25) {};
\node [style=none] (12) at (1.75, 1) {};
\node [style=none] (13) at (2.25, -3.5) {};
\node [style=none] (14) at (0.75, 1) {};
\node [style=none] (15) at (3, 0.25) {$M$};
\node [style=none] (16) at (3.25, -0) {};
\node [style=circle] (17) at (1.25, -3) {$\rho$};
\node [style=none] (18) at (1.25, -3.5) {};
\node [style=none] (19) at (2.25, 2.25) {};
\node [style=none] (20) at (0.25, 2.75) {};
\node [style=none] (21) at (1.25, -0.75) {};
\node [style=ox] (22) at (1.75, -1.25) {};
\node [style=none] (23) at (2.25, -0.75) {};
\node [style=none] (24) at (2.75, -2) {};
\node [style=none] (25) at (0.75, -0.75) {};
\node [style=none] (26) at (0.75, -2) {};
\node [style=none] (27) at (1.25, -2) {};
\node [style=none] (28) at (1.75, -0.75) {};
\node [style=none] (29) at (2.75, -0.75) {};
\node [style=none] (30) at (2.25, -2) {};
\node [style=none] (31) at (2.25, 1) {};
\node [style=none] (32) at (1.25, 1) {};
\node [style=none] (33) at (1.75, -0) {};
\node [style=ox] (34) at (1.75, 0.5) {};
\node [style=none] (35) at (3, 3) {$F((A \ox B)^\dagger)$};
\node [style=none] (36) at (3, -0.25) {$F(A^\dagger \ox B^\dagger)$};
\node [style=none] (37) at (0.25, -2.5) {$F(A^\dagger)$};
\node [style=none] (38) at (3, -2.5) {$F(B^\dagger)$};
\node [style=none] (39) at (0, -3.5) {$F(A)^\dagger$};
\node [style=none] (40) at (3, -3.5) {$F(B)^\dagger$};
\end{pgfonlayer}
\begin{pgfonlayer}{edgelayer}
\draw [bend right, looseness=1.50] (9.center) to (5);
\draw [bend right, looseness=1.50] (5) to (19.center);
\draw (5) to (12.center);
\draw (4.center) to (6.center);
\draw (6.center) to (7.center);
\draw (14.center) to (7.center);
\draw (14.center) to (4.center);
\draw (20.center) to (3.center);
\draw (3.center) to (16.center);
\draw (16.center) to (8.center);
\draw (8.center) to (20.center);
\draw (17) to (18.center);
\draw (10) to (13.center);
\draw (11.center) to (0.center);
\draw [bend left, looseness=1.50] (27.center) to (22);
\draw [bend left, looseness=1.50] (22) to (30.center);
\draw (22) to (28.center);
\draw (26.center) to (24.center);
\draw (24.center) to (29.center);
\draw (25.center) to (29.center);
\draw (25.center) to (26.center);
\draw [bend right, looseness=1.50] (32.center) to (34);
\draw [bend right, looseness=1.50] (34) to (31.center);
\draw (34) to (33.center);
\draw (33.center) to (28.center);
\draw (27.center) to (17);
\draw (30.center) to (10);
\end{pgfonlayer}
\end{tikzpicture} = \begin{tikzpicture} %rho-comon-a
\begin{pgfonlayer}{nodelayer}
\node [style=none] (0) at (1.5, 0.5) {};
\node [style=none] (1) at (-0.25, -0.75) {};
\node [style=none] (2) at (-0.25, 1.75) {};
\node [style=none] (3) at (1.25, -0.75) {};
\node [style=ox] (4) at (1, -0) {};
\node [style=none] (5) at (0.5, 0.5) {};
\node [style=none] (6) at (2.25, 1.75) {};
\node [style=none] (7) at (1.5, 0.75) {};
\node [style=none] (8) at (1.25, -1.25) {};
\node [style=none] (9) at (1.25, 1.75) {};
\node [style=ox] (10) at (1, 1.25) {};
\node [style=none] (11) at (2, 0.75) {};
\node [style=none] (12) at (0.5, 1.75) {};
\node [style=none] (13) at (1, -0.75) {};
\node [style=none] (14) at (2, -0.5) {};
\node [style=none] (15) at (0, -0.5) {};
\node [style=circle] (16) at (1.25, 2.5) {$\rho$};
\node [style=none] (17) at (0, 0.75) {};
\node [style=none] (18) at (2.25, -0.75) {};
\node [style=none] (19) at (1, 1.75) {};
\node [style=none] (20) at (1.75, -0.25) {$M$};
\node [style=none] (21) at (0.5, 0.75) {};
\node [style=none] (22) at (1.25, 3.25) {};
\node [style=none] (23) at (0, -1.25) {};
\node [style=none] (24) at (0.5, -1.25) {};
\node [style=ox] (25) at (1, -2) {};
\node [style=none] (26) at (1, -2.5) {};
\node [style=none] (27) at (2, -1.25) {};
\node [style=none] (28) at (0, -2.5) {};
\node [style=none] (29) at (1.5, -1.25) {};
\node [style=none] (30) at (2, -2.5) {};
\node [style=none] (31) at (0.5, -2.5) {};
\node [style=none] (32) at (0.5, -3) {};
\node [style=none] (33) at (1.5, -3) {};
\node [style=none] (34) at (1.5, -2.5) {};
\node [style=none] (35) at (2.25, 3) {$F((A \ox B)^\dagger)$};
\node [style=none] (36) at (3.25, -1) {$F((A) \ox F(B))^\dagger$};
\node [style=none] (37) at (-0.5, -2.75) {$F(A)^\dagger$};
\node [style=none] (38) at (2.75, -2.75) {$F(B)^\dagger$};
\node [style=none] (39) at (2.25, 2.25) {};
\end{pgfonlayer}
\begin{pgfonlayer}{edgelayer}
\draw [bend left, looseness=1.50] (5.center) to (10);
\draw [bend left, looseness=1.50] (10) to (0.center);
\draw (1.center) to (18.center);
\draw (18.center) to (6.center);
\draw (2.center) to (6.center);
\draw (2.center) to (1.center);
\draw [bend right, looseness=1.50] (21.center) to (4);
\draw [bend right, looseness=1.50] (4) to (7.center);
\draw (4) to (13.center);
\draw (17.center) to (11.center);
\draw (11.center) to (14.center);
\draw (15.center) to (14.center);
\draw (15.center) to (17.center);
\draw (16) to (9.center);
\draw (19.center) to (10);
\draw (3.center) to (8.center);
\draw (22.center) to (16);
\draw [bend right, looseness=1.50] (24.center) to (25);
\draw [bend right, looseness=1.50] (25) to (29.center);
\draw (25) to (26.center);
\draw (23.center) to (27.center);
\draw (27.center) to (30.center);
\draw (28.center) to (30.center);
\draw (28.center) to (23.center);
\draw (31.center) to (32.center);
\draw (34.center) to (33.center);
\end{pgfonlayer}
\end{tikzpicture}
\]

%%%%%%%%%%%%%%%%

%%%%%%%%%%%%%%%%
For linear natural transformations $(\beta_\ox, \beta_\oa): F \to G$ between $\dagger$-linear functors, we demand that $\beta_\ox$ and $\beta_\oa$ are related by:

\[ \xymatrix{F_\ox(A^\dagger) \ar[d]_{\rho^F_\ox} \ar[rr]^{\beta_\ox} && G_\ox(A^\dagger)  \ar[d]^{\rho^G_\ox} \\
	(F_\oa(X))^\dagger \ar[rr]_{\beta_\oa^\dagger} & &  (G_\oa(X))^\dagger}
~~~~~~
\xymatrix{(G_\ox(X))^\dagger \ar[d]_{\rho^G_\oa} \ar[rr]^{\beta_\ox^\dagger} &&  (F_\ox(X))^\dagger \ar[d]^{\rho^F_\oa} \\
	G_\oa(A^\dagger) \ar[rr]_{\beta_\oa} & &  F_\oa(A^\dagger)} \]
Notice that this means that $\beta_\ox$ is completely determined by $\beta_\oa$ in the following sense:
\[ \xymatrix{ F_\ox(A) \ar[d]_{F_\ox(\iota)}  \ar[rr]^{\beta_\ox} && G_\ox(A) \ar[d]^{G_\ox(\iota)} \\
	F_\ox(A^{\dagger\dagger})  \ar[d]_{\rho^F_\ox} \ar[rr]^{\beta_\ox} & & G_\ox(A^{\dagger\dagger}) \ar[d]^{\rho^G_\ox}  \\
	F_\oa(A^\dagger)^\dagger \ar[rr]_{\beta_\oa^\dagger} && G_\oa(A^\dagger)^\dagger } \]
Because the vertical maps are isomorphisms, this diagram can be used to express $\beta_\ox$ in terms of $\beta_\oa$.  Similarly, $\beta_\oa$ can be expressed in terms 
of $\beta_\ox$.  Thus, it is possible to express the coherences in terms of just one of these transformations.

%%%%%%%%%%%%%%%%%%%%%%%%%%%%%%%%%%%%%%%%%%%%%%%%

\subsection{Examples of $\dagger$-LDCs}
\label{subsection: Dagger LDC examples}

In this section, we discuss some basic examples of $\dagger$-isomix categories. The first example is a compact LDC. The example of Finiteness spaces is non-compact and all these examples are $*$-autonomous category. These categories all have a non-stationary dagger functor.  More examples of $\dagger$-isomix categories can be found in Section~\ref{Sec: MUC examples}.

%%%%%%%%%%%%%%%%%%%%%%%%%%%%%%%%%%%%%%%%%%%%%%%%%%%%%%%
\subsubsection{Finite dimensional framed vector spaces}
\label{subsection:fdfv}

In this section we describe the category of ``framed'' finite dimensional vector spaces, where a frame in this context is just a 
choice of basis.  Thus, the objects in this category are vector spaces with a chosen basis while the maps, ignoring the basis, are simply homomorphisms of 
the vector spaces.

The category of finite dimensional framed vectors spaces, ${\sf FFVec}_K$, is a monoidal category defined as follows:
\begin{description}
\item[Objects]  The objects are pairs $(V,{\cal V})$ where $V$ is a finite dimensional $K$-vector space and ${\cal V} = \{ v_1,...,v_n \}$ is a basis. 
\item[Maps]  A map $(V,{\cal V}) \xrightarrow{f} (W,{\cal W})$ is a linear map  $V \xrightarrow{f} W$ in ${\sf FdVec}_K$.
\item[Tensor] $(V,{\cal V}) \ox (W,{\cal W}) = (V \ox W,\{ v \ox w | v \in {\cal V}, w \in {\cal W} \})$ where $V \ox W$ is the usual tensor product.  The unit is 
$(K,\{ e \})$ where $e$ is the unit of the field $K$.
\end{description}

To define the ``dagger'' we  must first choose a conjugation $\overline{(\_)}: K \to K$ (see more details in Section~\ref{Sec: conjugation}), that is a field homomorphism with $k = \overline{(\overline{k})}$. The 
canonical example being conjugation of the complex numbers, however, the conjugation can be arbitrarily chosen -- so could also, for example, be the identity.  
This conjugation then can be extended to a (covariant) functor:
\[ \overline{(\_)}: {\sf FFVec}_K \to {\sf FFVec}_K; \begin{array}[c]{c} \xymatrix{(V,{\cal V}) \ar[d]^{f} \\ (W,{\cal W})} \end{array} 
                 \mapsto \begin{array}[c]{c} \xymatrix{\overline{(V,{\cal V}) } \ar[d]^{\overline{f}} \\ \overline{(W,{\cal W}) }} \end{array} \]
where $\overline{(V,{\cal V})}$ is the vector space with the same basis but with the conjugate action $c~\overline{\cdot}~v = \overline{c} \cdot  v$.  The conjugate 
homomorphism, $\overline{f}$, is then the same underlying map which is homomorphism between the conjugate spaces.

${\sf FFVec}_K$ is also a compact closed category with $(V,{\cal B})^{*} = (V^{*}, \{ \widetilde{b_i} | b_i \in {\cal B} \})$ where 
\[ V^{*} = V \multimap K~~~~\mbox{and}~~~~\widetilde{b_i}: V \to K; \left(\sum_j \beta_j \cdot b_j \right) \mapsto  \beta_i \]
This makes $(\_)^{*}: {\sf FFVec}_K^{\rm op} \to {\sf FFVec}_K$ a contravariant functor whose action is determined by precomposition.   Finally, we 
define the ``dagger'' to be the composite $(V,{\cal B})^\dagger = \overline{(V,{\cal B})^{*}}$.

This is a compact LDC with tensor and par being identified (so the linear distribution is the associator) and is isomix.  We must show that it is a $\dagger$-LDC.  
Towards this aim we define the required natural transformations on the basis:

\[ \lambda_\ox = \lambda_\oa: (V,{\cal V})^\dag \ox (W,{\cal W})^\dag \to ((V,{\cal V}) \ox (W,{\cal W}))^\dag;  \widetilde{v_i} \ox \widetilde{w_j} \mapsto \widetilde{v_i \ox w_j} \]

\[ \lambda_\top = \lambda_\bot: (K,\{ e\}) \to (K,\{ e\})^\dag; k \mapsto \overline{k} \]
\[ \iota: (V,{\cal V}) \to ((V,{\cal V})^\dag)^\dag; v \mapsto \lambda f. f(v) \]
Note that the last transformation is given in a basis independent manner. Importantly, it may also be given in a basis dependent manner as
$\iota(v_i) = \widetilde{\widetilde{v_i}}$ as the behaviour of these two maps is the  same when applied to the basis of $(V,{\cal V})^\dag$ namely 
the elements $\widetilde{v_j}$:
\[ \iota(v_i)(\widetilde{v_j}) = (\lambda f. f(v_i)) \widetilde{v_j} = \widetilde{v_j} v_i = \partial_{i,j} = \widetilde{\widetilde{v_i}}(\widetilde{v_j})  \]
Also note that $\widetilde{v_i \ox w_j} = (\widetilde{v_i} \ox \widetilde{w_j}) u_\ox$, where $u_\ox: K \ox K \to K$ is the multiplication of the field.
With these definitions in hand it is straightforward to check that this gives a $\dagger$-LDC by checking the required coherences on basis elements.
To demonstrate the technique consider the coherence {\bf [$\dagger$-ldc.4]}:
\[ \xymatrix{A \oa B \ar[d]_{\iota \oa \iota} \ar[rr]^\iota & & ((A \oa B)^\dag)^\dag \ar[d]^{\lambda_\ox^\dag} \\
                   (A^\dag)^\dag \oa (B^\dag)^\dag  \ar[rr]_{\lambda_\oa} & & (A^\dag \ox B^\dag)^\dag} \]
 We must show (identifying tensor and par) that $\lambda_\ox^\dag (\iota(a_i \ox b_j)) = \lambda_\ox(\iota \ox \iota(a_i \ox b_j))$.  Now the result  is a 
 higher-order term so it suffices to show the evaluations on basis elements are the same.  This means we need to show:
 $\lambda_\ox^\dag (\iota(a_i \ox b_j))(\widetilde{a_p} \ox \widetilde{b_q}) = \lambda_\ox(\iota \ox \iota(a_i \ox b_j))(\widetilde{a_p} \ox \widetilde{b_q})$
 \begin{eqnarray*}
(\lambda_\ox(\iota \ox \iota(a_i \ox b_j)))(\widetilde{a_p} \ox \widetilde{b_q}) & = & (\lambda_\ox(\widetilde{\widetilde{a_i}} \ox \widetilde{\widetilde{b_j}}))(\widetilde{a_p} \ox \widetilde{b_q}) \\
&  = & (\widetilde{\widetilde{a_i} \ox \widetilde{b_j}}) (\widetilde{a_p} \ox \widetilde{b_q}) \\
& = & (\widetilde{a_p} \ox \widetilde{b_q})(\widetilde{\widetilde{a_i}} \ox \widetilde{\widetilde{b_j}}) u_\ox ~~~\mbox{(diagrammatic order)}\\
& = & \partial_{p,i} \partial_{q,j} \\
(\lambda_\ox^\dag (\iota(a_i \ox b_j)))(\widetilde{a_p} \ox \widetilde{b_q}) 
& = & (\lambda_\ox^\dag (\widetilde{\widetilde{a_i \ox b_j}}))(\widetilde{a_p} \ox \widetilde{b_q})  \\ 
& = & (\widetilde{a_p} \ox \widetilde{b_q}) \lambda_\ox \widetilde{\widetilde{a_i \ox b_j}} ~~~\mbox{(diagrammatic order)}\\
& = & \widetilde{a_p \ox b_q} \widetilde{\widetilde{a_i \ox b_j}} \\
& = & \partial_{p,i} \partial_{q,j} 
\end{eqnarray*}

Thus, $ {\sf FFVec}_K$ is a compact $\dagger$-isomix category where the $\dagger$ functor shifts objects i.e., $A \neq A^\dagger$.

\subsubsection{Finiteness matrices}
\label{Sec: Finiteness matrices}

Finiteness spaces were introduced by Ehrhard, \cite{Ehr05}, as a model of linear logic.  The type system can be used to produce a typed system for infinite 
dimensional matrix multiplication in which no sums become infinite.   This system of infinite dimensional matrices forms an isomix $*$-autonomous category.  
If these matrices have entries in the complex numbers then there is a natural notion of conjugation and this gives a $\dagger$-isomix category. We shall see that 
by taking the core of this category, which is the category of finite dimensional matrices, one obtains a basic example of a 
mixed unitary category as described in Section~\ref{Sec: unitary} of this paper.  This example is also explored further in the sequels to this paper \cite{CS19,CoS20}.

\begin{defi}
A {\bf finiteness space} is a pair $X:=(|X|,\mathcal{F})$ with $|X|$ a set, called the {\bf web}, and $\mathcal{F}$ be a subset of $\mathcal{P} (|X|)$ such that $\mathcal{F} = \mathcal{F}^{\perp\perp}$ where
\[\mathcal{F}^\perp:=\{ x' \subseteq X |\forall x \in \mathcal{F}.  x\cap x'\text{ is finite} \}. \]
The elements of $\mathcal{F}$ are called the {\bf finitary sets} of the finiteness space of $X$, and the elements of $\mathcal{F}^\perp$ are called {\bf cofinitary sets}.
\end{defi}

Observe that if the web of $X$ is finite, then $\mathcal{F}$ is forced to be the whole powerset of $X$.  Finiteness spaces organize themselves into a (symmetric) $*$-autonomous category, ${\sf FRel}$:

\begin{description}
\item[Objects] Finiteness spaces $X:=(|X|, \mathcal{F})$.
\item[Maps] A map $R:(|X|, \mathcal{F}) \to (|Y|, \mathcal{G})$ is a {\bf finiteness relation} that is a relation $R: |X| \to |Y|$ so that:
\begin{itemize}
\item For all $x \in \mathcal{F}$, $xR \in \mathcal{G}$ where $xR = \{ y | y \in Y, \text{ there exists } b \in x \text{ such that } bRy \}$
\item For all $y \in \mathcal{G}^\perp$, $Ry \in \mathcal{F}^\perp$ where $Ry = \{ x | x \in X, \text{ there exists } t \in y \text{ such that } xRt \}$
\end{itemize}
\item[Composition and identities] same as in sets and relations.
\item[Tensor] 
Given finiteness spaces $X = (|X|, \mathcal{F})$, and $Y=(|Y|, \mathcal{G})$, 
$X \ox Y := (|X| \x |Y|,  \mathcal{F} \x \mathcal{G}) \text{ where} $
\begin{align*}
\mathcal{F} \x \mathcal{G} &:=~ \downarrow\{ u \times v | u \in \mathcal{F}, v \in \mathcal{G} \} \\
& ~= ~ \{ w | \exists  u \in \mathcal{F},  v \in \mathcal{G}. w \subseteq u \times v \}
\end{align*}
This is well-defined as $\mathcal{F} \x \mathcal{G} = (\mathcal{F} \x \mathcal{G})^{\perp\perp}$ is proven in \cite{Ehr05} Lemma 2.

Given maps $R: X_1 \to Y_1$ and $S: X_2 \to Y_2$:
\[ R \ox S := \{((x_1, x_2), (y_1,y_2))| (x_1,y_1)\in R, (x_2,y_2) \in S \} \]
\item[Monoidal tensor unit] $\top := (\{*\}, \mathcal{P}(\{*\}))$
\item[Dualizing functor] $(X, \mathcal{F})^* := (X, \mathcal{F}^\perp)$
\end{description}

Finiteness spaces with finiteness relations, ${\sf FRel}$, form an isomix $*$-autonomous category, where the objects in the core are precisely the objects whose webs are finite sets.
There is a faithful underlying structure preserving functor $U: {\sf FRel} \to {\sf Rel}$.

The category of finiteness matrices over the complex numbers, ${\sf FMat}(\C)$, can be built from the category of finiteness spaces:
\begin{description}
\item[Objects] Finiteness spaces $X = (|X|, \mathcal{F})$;
\item[Maps]  $M: X \to Y$, matrices $M: |X| \x |Y| \to \C; (x,y) \mapsto M_{x,y}$ such that ${\sf supp}(M) = \{ (x,y) | M_{x,y} \not= 0 \}$ is a finiteness relation;
\item[Composition] If $M: X \to Y$ and $N: Y \to Z$ then $(MN)_{x,z} = \sum_{y \in Y} M_{x,y}N_{y,z}$ where for fixed $x$ and $z$ both $M_{x,y}$ and 
$N_{y,z}$ are non-zero for only finitely many $y$ as this is the intersection of a finitary set of $Y$ with a cofinitary set of $Y$.  Identity maps are given as usual 
by diagonal matrices.
\item[Tensor]  If $M: X \to Y$ and $N: X' \to Y'$ then $M \ox N: X \ox X' \to Y \ox Y'; ((x,y),(x',y')) \mapsto M_{x,y}N_{x',y'}$;
\item[Dual and dagger] If $M: X \to Y$ then $M^{*}: Y^{*} \to X^{*}; (y,x) \mapsto M_{x,y}$ while $M^\dagger: Y^{*} \to X^{*}; (y,x) \mapsto \overline{M_{x,y}}$.
\end{description}

The coherence maps for the isomix $*$-autonomous of ${\sf FMat}(\C)$ are simply the characteristic functions of the relations which give the coherence maps in ${\sf FRel}$.
The support gives a colax 2-functor from ${\sf supp}: {\sf FMat}(\C) \to {\sf FRel}$ where ${\sf supp}(M N) \subseteq  {\sf supp}(M) {\sf supp}(N)$.  

%%%%%%%%%%%%%%%%%%%%%%%%%%%%%%%%%%%%%%%%%%%%%%%%
\subsubsection{Category of abstract state spaces}
\label{Sec: Asp}

This model is inspired by the category of convex operational models \cite{BaW11}. 
The following is a way to construct new $\dagger$-isomix categories from an exisiting one.

\begin{defi}
	Let $\X$ be a $\dagger$-isomix category. Define $\Asp(\X)$ as follows:
	\begin{description}
		\item[Objects] $(A, e_A:A \to \bot, u_A: \top \to A)$
		\item[Arrows]  $f: A \to B \in \X$ such that the following diagram commutes:
		
		$
		\xymatrix{
			& \top \ar[dl]_{u_A}  \ar[dr]^{u_B} & \\
			A \ar[rr]^{f} \ar[dr]_{e_A}  & & B \ar[ld]^{e_B} \\
			& \bot &
		}
		$
	\end{description}
	Identity arrow and composition are inherited directly from $\X$. $\Asp(\X)$ is a LDC:
	\begin{description}
		\item[$\ox$ on objects] $(A, e_A, u_A) \ox (B, e_B, u_B) := (A \ox B, e', u')$ where, $e' := \mx (e_A \oa e_B) u_\oa$ and $u' := u_\ox^{-1} (u_A \ox u_B)$. The unit of $\ox$ is given by $(\top, \m^{-1}: \top \to \bot, 1_\top)$.
		\item[$\oa$ on objects] $(A, e_A, u_A) \oa (B, e_B, u_B) := (A \oa B, e', u')$ where, $e' := (e_A \oa e_B) u_\oa$ and $u' := u_\ox^{-1} (u_A \ox u_B) \mx$. The unit of $\oa$ is $(\bot, 1_\bot, \m^{-1}: \top \to \bot)$
	\end{description}
\end{defi}

$\Asp(\X)$ is also $\dagger$-isomix category with \[ (A, e, u)^\dagger := (A^\dagger, u^\dagger \lambda_\bot^{-1}, \lambda_\top e^\dagger)\] All the basic natural isomorphisms are inherited from $\X$. Hence, $\Asp(\X)$ is a $\dagger$-isomix category.

%%%%%%%%%%%%%%%%%%%%%%%%%%%%%%%%%%%%%%%%%%%%%%%%
\section{Daggers, duals, and conjugation}
\label{daggers-duals-conjugation}

%%%%%%%%%%%%%%%%%%%%%%%%%%%%%%%%%%%%%%%%%%%%%%%%

The goal of this section is to review the interaction of the dualizing, conjugation and dagger functors. 
In dagger compact closed categories, the dagger functor $(\_)^\dagger$, and the dualizing functor $(\_)^*$ 
commute with each other and their composite gives the conjugate functor $(\_)_*$. Similary, $(\_)_*$ and $(\_)^*$ 
when composed gives the dagger functor.  Our aim is to  generalize these interactions to $\dagger$-LDCs and to achieve 
this at a reasonable level of abstraction.   To achieve this we shall need the notion which we here refer to as ``conjugation'' 
but was investigated by Egger in \cite{Egg11} under the moniker of ``involution'' (which clashes with our usage).  

%%%%%%%%%%%%%%%%%%%%%%%%%%%%%%%%%%%%%%%%%%%%%%%%%%%%%%%%%%%%%%%%%%%%%%
\subsection{Duals}

The reverse of an LDC, $\X$, written $\X^\rev := (\X, \ox, \top, \oa, \bot)^{\sf rev}$ 

$= (\X, \ox^{\sf rev}, \top, \oa^{\sf rev}, \bot)$ where,
\[ A \ox^{\sf rev} B := B \ox A ~~~~~~~~~~~ A \oa^{\sf rev} B := B \oa A \]
and the associators and distributors are adjusted accordingly.  Similar to the opposite of an LDC, we have $(\X^{\sf rev})^{\sf rev} = \X$.  

In a $*$-autonomous category, taking the left (or right) linear dual of an object extends to a Frobenius linear functor as follows:
\begin{align*}
(\_)^*: (\X^{\sf op})^{\sf rev} &\to \X  ; ~~ A \mapsto A^* ; ~~
\begin{tikzpicture}
	\begin{pgfonlayer}{nodelayer}
		\node [style=none] (0) at (0, 2) {};
		\node [style=none] (1) at (-0.5, 1) {};
		\node [style=none] (2) at (0.5, 1) {};
		\node [style=none] (3) at (-0.5, -0) {};
		\node [style=none] (4) at (0.5, -0) {};
		\node [style=none] (5) at (0, 1) {};
		\node [style=none] (6) at (0, -1) {};
		\node [style=none] (7) at (0, -0) {};
		\node [style=none] (8) at (0, 0.5) {$f$};
		\node [style=none] (9) at (0.25, 1.75) {$A$};
		\node [style=none] (10) at (0.25, -0.75) {$B$};
	\end{pgfonlayer}
	\begin{pgfonlayer}{edgelayer}
		\draw (1.center) to (3.center);
		\draw (3.center) to (4.center);
		\draw (4.center) to (2.center);
		\draw (2.center) to (1.center);
		\draw (0.center) to (5.center);
		\draw (7.center) to (6.center);
	\end{pgfonlayer}
\end{tikzpicture} \mapsto  \begin{tikzpicture}
	\begin{pgfonlayer}{nodelayer}
		\node [style=none] (0) at (0, 1.5) {};
		\node [style=none] (1) at (-0.5, 1) {};
		\node [style=none] (2) at (0.5, 1) {};
		\node [style=none] (3) at (-0.5, -0) {};
		\node [style=none] (4) at (0.5, -0) {};
		\node [style=none] (5) at (0, 1) {};
		\node [style=none] (6) at (0, -0.5) {};
		\node [style=none] (7) at (0, -0) {};
		\node [style=none] (8) at (-1, 1.5) {};
		\node [style=none] (9) at (-1, -1) {};
		\node [style=none] (10) at (1, -0.5) {};
		\node [style=none] (11) at (1, 2) {};
		\node [style=none] (12) at (0, 0.5) {$f$};
		\node [style=none] (13) at (1.25, 1.75) {$B^*$};
		\node [style=none] (14) at (-1.25, -0.75) {$A^*$};
	\end{pgfonlayer}
	\begin{pgfonlayer}{edgelayer}
		\draw (1.center) to (3.center);
		\draw (3.center) to (4.center);
		\draw (4.center) to (2.center);
		\draw (2.center) to (1.center);
		\draw (0.center) to (5.center);
		\draw (7.center) to (6.center);
		\draw [bend right=90, looseness=1.25] (0.center) to (8.center);
		\draw (8.center) to (9.center);
		\draw [bend right=90, looseness=1.25] (6.center) to (10.center);
		\draw (10.center) to (11.center);
	\end{pgfonlayer}
\end{tikzpicture}
\end{align*}

The $(\_)^*$ functor is both contravariant and, op-monoidal and op-comonoidal:

\[ m_{\ox}:  A^* \ox B^* \to (B \oa A)^* :=  \begin{tikzpicture}
\begin{pgfonlayer}{nodelayer}
\node [style=none] (0) at (-0.25, 2) {};
\node [style=ox] (1) at (-0.25, 1) {};
\node [style=none] (2) at (-0.5, 0.25) {};
\node [style=none] (3) at (0.25, -0) {};
\node [style=none] (4) at (-1.25, 0.25) {};
\node [style=none] (5) at (-2, -0) {};
\node [style=oa] (6) at (-1.5, 1) {};
\node [style=none] (7) at (-1.5, 1.5) {};
\node [style=none] (8) at (-2.25, 1.5) {};
\node [style=none] (9) at (-2.25, -1) {};
\node [style=none] (10) at (0.5, 1.75) {$B^* \ox A^*$};
\node [style=none] (11) at (-0.8, 0.5) {$B^*$};
\node [style=none] (12) at (0.5, 0.5) {$A^*$};
\node [style=none] (13) at (-3.2, -0.75) {$(A \oa B)^* $};
\end{pgfonlayer}
\begin{pgfonlayer}{edgelayer}
\draw [in=90, out=-45, looseness=1.00] (1) to (3.center);
\draw [in=90, out=-135, looseness=1.00] (1) to (2.center);
\draw [bend left=90, looseness=1.25] (2.center) to (4.center);
\draw [in=-60, out=90, looseness=1.00] (4.center) to (6);
\draw (6) to (7.center);
\draw [bend right=90, looseness=1.25] (7.center) to (8.center);
\draw (8.center) to (9.center);
\draw [in=90, out=-135, looseness=0.75] (6) to (5.center);
\draw [in=-90, out=-90, looseness=1.00] (5.center) to (3.center);
\draw (1) to (0.center);
\end{pgfonlayer}
\end{tikzpicture} 
~~~~~~~~~~~~~~~ m_\top: \top \to \bot^* :=  \begin{tikzpicture} %pic2
\begin{pgfonlayer}{nodelayer}
\node [style=circle] (0) at (0, 1.75) {$\top$};
\node [style=none] (1) at (0, 3) {};
\node [style=circle, scale=0.5] (2) at (-1, 1) {};
\node [style=circle] (3) at (-1, -0.25) {$\bot$};
\node [style=none] (4) at (-2, -1) {};
\node [style=none] (5) at (-2, 1.5) {};
\node [style=none] (6) at (-1, 1.5) {};
\node [style=none] (7) at (0.25, 2.75) {$\top$};
\node [style=none] (8) at (-2.5, -0.75) {$\bot^*$};
\end{pgfonlayer}
\begin{pgfonlayer}{edgelayer}
\draw (1.center) to (0);
\draw (3) to (6.center);
\draw [bend left=90, looseness=2.00] (5.center) to (6.center);
\draw (5.center) to (4.center);
\draw [in=0, out=-90, looseness=1.50, dotted] (0) to (2);
\end{pgfonlayer}
\end{tikzpicture} \]
\[ n_{\oa}:  (A \ox B)^* \to B^* \oa A^* :=  \begin{tikzpicture} %pic3
\begin{pgfonlayer}{nodelayer}
\node [style=none] (0) at (2, 2) {};
\node [style=ox] (1) at (1, 0.25) {};
\node [style=none] (2) at (2, -0.5) {};
\node [style=none] (3) at (1, -0.5) {};
\node [style=none] (4) at (0.5, 0.75) {};
\node [style=none] (5) at (1.5, 0.75) {};
\node [style=oa] (6) at (-0.75, 0.25) {};
\node [style=none] (7) at (-0.75, -1.75) {};
\node [style=none] (8) at (-0.25, 0.75) {};
\node [style=none] (9) at (-1.25, 0.75) {};
\node [style=none] (10) at (2, 2.25) {$(A \ox B)^*$};
\node [style=none] (11) at (-0.75, -2) {$B^* \oa A^*$};
\end{pgfonlayer}
\begin{pgfonlayer}{edgelayer}
\draw (0.center) to (2.center);
\draw [bend left=90, looseness=1.50] (2.center) to (3.center);
\draw (3.center) to (1);
\draw (6) to (7.center);
\draw [in=150, out=-90, looseness=1.00] (9.center) to (6);
\draw [in=-90, out=11, looseness=1.00] (6) to (8.center);
\draw [in=165, out=-90, looseness=1.25] (4.center) to (1);
\draw [in=-90, out=15, looseness=1.25] (1) to (5.center);
\draw [bend left=90, looseness=2.00] (8.center) to (4.center);
\draw [bend left=90, looseness=1.25] (9.center) to (5.center);
\end{pgfonlayer}
\end{tikzpicture} ~~~~~~~~~~~~~~~ n_\bot: \top^* \to \bot :=  \begin{tikzpicture}
\begin{pgfonlayer}{nodelayer}
\node [style=circle] (0) at (-2, 1) {$\top$};
\node [style=none] (1) at (-1, 2) {};
\node [style=none] (2) at (-2, -0) {};
\node [style=none] (3) at (-1, -0) {};
\node [style=circle] (4) at (0, -1) {$\bot$};
\node [style=circle, scale=0.5] (5) at (-1, 1) {};
\node [style=none] (6) at (0, -2) {};
\node [style=none] (7) at (-0.6, 1.75) {$\top^*$};
\end{pgfonlayer}
\begin{pgfonlayer}{edgelayer}
\draw (0) to (2.center);
\draw [bend right=90, looseness=1.75] (2.center) to (3.center);
\draw (3.center) to (1.center);
\draw [in=90, out=-15, looseness=1.25, dotted] (5) to (4);
\draw (6.center) to (4);
\end{pgfonlayer}
\end{tikzpicture} \]

These maps are op-monoidal and op-comonoidal laxors, hence are isomorphisms, which satisfy the obvious coherences. Thus, $(\_)^*$ is a strong Frobenius linear functor. 

In the rest of the section, we will write $(\X^\op)^\rev$ as $\X^{\op\rev}$.

\begin{lem} If $\X$ is an isomix category, then $(\_)^*: \X^{\op\rev} \to \X$ is an isomix functor. \end{lem}
\begin{proof}
Because,  $(\_)^*$ is a strong Frobenius functor, by Lemma \ref{Lemma: isomix functor}, it suffices  to prove that $(\_)^*$ preserves mix, i.e.,  $(\_)^*$  is a mix functor i.e., we need  to show that $n_\bot \!~\m~ m_\top : \top^* \to \bot^* =  \m^* $. The proof is as follows:
\[ n_\bot \m m_\top = \begin{tikzpicture}
	\begin{pgfonlayer}{nodelayer}
		\node [style=circle] (0) at (-3, 2.75) {$\top$};
		\node [style=none] (1) at (-2, 2) {};
		\node [style=none] (2) at (-2, 3) {};
		\node [style=circle] (3) at (-1, 1.25) {$\bot$};
		\node [style=map] (4) at (-1, 0.5) {};
		\node [style=circle] (5) at (-1, -0.25) {$\top$};
		\node [style=circle] (6) at (-2, -2) {$\bot$};
		\node [style=none] (7) at (-3, -1.25) {};
		\node [style=none] (8) at (-3, -2.25) {};
		\node [style=circle, scale=0.5] (9) at (-2, 1.75) {};
		\node [style=circle, scale=0.5] (10) at (-2, -0.75) {};
		\node [style=none] (11) at (-3, 2) {};
		\node [style=none] (12) at (-2, -1.25) {};
	\end{pgfonlayer}
	\begin{pgfonlayer}{edgelayer}
		\draw (1.center) to (2.center);
		\draw (7.center) to (8.center);
		\draw [bend left, looseness=1.00, dotted] (9) to (3);
		\draw [bend left=45, looseness=1.00, dotted] (5) to (10);
		\draw (3) to (4);
		\draw (4) to (5);
		\draw [bend right=90, looseness=3.25] (11.center) to (1.center);
		\draw (0) to (11.center);
		\draw [bend left=90, looseness=3.75] (7.center) to (12.center);
		\draw (12.center) to (6);
	\end{pgfonlayer}
\end{tikzpicture} =\begin{tikzpicture}
	\begin{pgfonlayer}{nodelayer}
		\node [style=circle] (0) at (1.25, 3) {$\top$};
		\node [style=none] (1) at (0.25, 3) {};
		\node [style=circle] (2) at (-1, 1.25) {$\bot$};
		\node [style=map] (3) at (-1, 0.5) {};
		\node [style=circle] (4) at (-1, -0.25) {$\top$};
		\node [style=circle] (5) at (-2, -2.25) {$\bot$};
		\node [style=none] (6) at (-3, -2.25) {};
		\node [style=circle, scale=0.5] (7) at (0.25, 2) {};
		\node [style=circle, scale=0.5] (8) at (-2, -0.75) {};
		\node [style=none] (9) at (0.25, 1.5) {};
		\node [style=none] (10) at (1.25, 1.5) {};
		\node [style=none] (11) at (-2, -0.25) {};
		\node [style=none] (12) at (-3, -0.25) {};
	\end{pgfonlayer}
	\begin{pgfonlayer}{edgelayer}
		\draw [dotted, bend right, looseness=1.25] (7) to (2);
		\draw [dotted, bend left=45, looseness=1.00] (4) to (8);
		\draw (2) to (3);
		\draw (3) to (4);
		\draw (1.center) to (9.center);
		\draw [bend right=90, looseness=2.00] (9.center) to (10.center);
		\draw (10.center) to (0);
		\draw (5) to (11.center);
		\draw (12.center) to (6.center);
		\draw [bend left=90, looseness=1.75] (12.center) to (11.center);
	\end{pgfonlayer}
\end{tikzpicture} = \begin{tikzpicture} %m3
	\begin{pgfonlayer}{nodelayer}
		\node [style=circle] (0) at (1.25, 3) {$\top$};
		\node [style=none] (1) at (0.25, 2.75) {};
		\node [style=circle] (2) at (-1, 1.25) {$\bot$};
		\node [style=map] (3) at (-1, 0.5) {};
		\node [style=circle] (4) at (-1, -0.25) {$\top$};
		\node [style=circle] (5) at (-2, -2.25) {$\bot$};
		\node [style=none] (6) at (-3, -2) {};
		\node [style=circle, scale=0.5] (7) at (0.25, 2) {};
		\node [style=circle, scale=0.5] (8) at (-2, -0.75) {};
		\node [style=none] (9) at (0.25, -1.75) {};
		\node [style=none] (10) at (1.25, -1.75) {};
		\node [style=none] (11) at (-2, 2.5) {};
		\node [style=none] (12) at (-3, 2.5) {};
	\end{pgfonlayer}
	\begin{pgfonlayer}{edgelayer}
		\draw [dotted, bend right, looseness=1.25] (7) to (2);
		\draw [dotted, bend left=45, looseness=1.00] (4) to (8);
		\draw (2) to (3);
		\draw (3) to (4);
		\draw (1.center) to (9.center);
		\draw [bend right=90, looseness=2.00] (9.center) to (10.center);
		\draw (10.center) to (0);
		\draw (5) to (11.center);
		\draw (12.center) to (6.center);
		\draw [bend left=90, looseness=1.75] (12.center) to (11.center);
	\end{pgfonlayer}
\end{tikzpicture} = \begin{tikzpicture} %m3
	\begin{pgfonlayer}{nodelayer}
		\node [style=circle] (0) at (1.25, 3) {$\top$};
		\node [style=none] (1) at (0.25, 2.75) {};
		\node [style=circle] (2) at (-0.75, 1.25) {$\bot$};
		\node [style=map] (3) at (-0.75, 0.5) {};
		\node [style=circle] (4) at (-0.75, -0.25) {$\top$};
		\node [style=circle] (5) at (-2, -2.25) {$\bot$};
		\node [style=none] (6) at (-3, -2) {};
		\node [style=circle, scale=0.5] (7) at (-2, 2) {};
		\node [style=circle, scale=0.5] (8) at (0.25, -0.75) {};
		\node [style=none] (9) at (0.25, -1.75) {};
		\node [style=none] (10) at (1.25, -1.75) {};
		\node [style=none] (11) at (-2, 2.5) {};
		\node [style=none] (12) at (-3, 2.5) {};
	\end{pgfonlayer}
	\begin{pgfonlayer}{edgelayer}
		\draw [dotted, bend left, looseness=1.25] (7) to (2);
		\draw [dotted, bend right=45, looseness=1.00] (4) to (8);
		\draw (2) to (3);
		\draw (3) to (4);
		\draw (1.center) to (9.center);
		\draw [bend right=90, looseness=2.00] (9.center) to (10.center);
		\draw (10.center) to (0);
		\draw (5) to (11.center);
		\draw (12.center) to (6.center);
		\draw [bend left=90, looseness=1.75] (12.center) to (11.center);
	\end{pgfonlayer}
\end{tikzpicture} = \begin{tikzpicture}
	\begin{pgfonlayer}{nodelayer}
		\node [style=map] (0) at (-0.75, 1) {};
		\node [style=none] (1) at (-1.75, -1.75) {};
		\node [style=none] (2) at (-0.75, -0.25) {};
		\node [style=none] (3) at (0.25, -0.25) {};
		\node [style=none] (4) at (-0.75, 2.25) {};
		\node [style=none] (5) at (-1.75, 2.25) {};
		\node [style=none] (6) at (0.25, 3.25) {};
	\end{pgfonlayer}
	\begin{pgfonlayer}{edgelayer}
		\draw [bend right=90, looseness=2.00] (2.center) to (3.center);
		\draw (5.center) to (1.center);
		\draw [bend left=90, looseness=1.75] (5.center) to (4.center);
		\draw (4.center) to (0);
		\draw (0) to (2.center);
		\draw (6.center) to (3.center);
	\end{pgfonlayer}
\end{tikzpicture} = \m^*
\]
\end{proof}

\begin{lem}
$(\eta, \epsilon) :: (\_)^* \dashvv ~ {^*(\_)^{\op\rev}} : \X^{\op \rev} \to \X$  
\[ 
\eta_\ox: X \to ~^*(X^*) := \begin{tikzpicture} %snakea
	\begin{pgfonlayer}{nodelayer}
		\node [style=none] (0) at (1, -2) {};
		\node [style=none] (1) at (1, 0.75) {};
		\node [style=none] (2) at (0, 0.75) {};
		\node [style=none] (12) at (0.5, 1.75) {$*\eta$};		
		\node [style=none] (3) at (0, -0.25) {};
		\node [style=none] (4) at (-1, -0.25) {};
		\node [style=none] (12) at (-0.5, -1.5) {$\epsilon*$};		
		\node [style=none] (5) at (-1, 2.5) {};
		\node [style=none] (6) at (-1.25, 2) {$X$};
		\node [style=none] (7) at (-0.35, 0.5) {$X^*$};
		\node [style=none] (8) at (1.6, -1.75) {$^*(X^*)$};
	\end{pgfonlayer}
	\begin{pgfonlayer}{edgelayer}
		\draw (4.center) to (5.center);
		\draw [bend right=90, looseness=3.25] (4.center) to (3.center);
		\draw (3.center) to (2.center);
		\draw [bend left=90, looseness=2.50] (2.center) to (1.center);
		\draw (1.center) to (0.center);
	\end{pgfonlayer}
\end{tikzpicture} \in \X
~~~~~~~~~~~~ \eta_\oa := \eta_\ox^{-1} \]
\[
\epsilon_\oa: X \to  (^*X)^*:= \begin{tikzpicture}
	\begin{pgfonlayer}{nodelayer}
		\node [style=none] (0) at (-0.75, -2) {};
		\node [style=none] (1) at (-0.75, 0.75) {};
		\node [style=none] (2) at (0.25, 0.75) {};
		\node [style=none] (12) at (-0.5, 1.75) {$\eta*$};
		\node [style=none] (3) at (0.25, -0.25) {};
		\node [style=none] (4) at (1.25, -0.25) {};
		\node [style=none] (34) at (0.65, -1.5) {$*\epsilon$};
		\node [style=none] (5) at (1.25, 2.5) {};
		\node [style=none] (6) at (1.65, 2) {$X$};
		\node [style=none] (7) at (0.6, 0.5) {$^*X$};
		\node [style=none] (8) at (-1.25, -1.75) {$(^*X)^*$};
	\end{pgfonlayer}
	\begin{pgfonlayer}{edgelayer}
		\draw (4.center) to (5.center);
		\draw [bend left=90, looseness=3.25] (4.center) to (3.center);
		\draw (3.center) to (2.center);
		\draw [bend right=90, looseness=2.50] (2.center) to (1.center);
		\draw (1.center) to (0.center);
	\end{pgfonlayer}
\end{tikzpicture} \in \X ~~~~~~~~~~~ \epsilon_\ox := \epsilon_\oa^{-1}
\]
is a linear equivalence of Frobenius linear functors.
\end{lem}
\begin{proof}
The proof is straightforward in the graphical calculus.
\end{proof}

For a cyclic $*$-autonomous category, we can straighten out this equivalence to be a dualizing involutive equivalence (i.e. so that the unit and counit are equal):

\begin{lem}
$(\eta', \epsilon') :: (\_)^* \dashvv  ((\_)^{*})^{\op\rev}: \X^{\op\rev} \to \X$  where $\eta'_\ox = {\eta'_\oa}^{-1} := \eta_\ox \psi^{-1}$, $\epsilon'_\ox = \epsilon'_\oa:= \epsilon \psi^*$ and $\eta' = \epsilon'$.
\end{lem}
\begin{proof}

The unit and counit are drawn as follows:

\[
\eta_\ox' =  \begin{tikzpicture}
	\begin{pgfonlayer}{nodelayer}
		\node [style=none] (0) at (-3, 3) {};
		\node [style=none] (1) at (-3, 1) {};
		\node [style=none] (2) at (-2, 1) {};
		\node [style=none] (3) at (-2, 2) {};
		\node [style=none] (4) at (-1, 2) {};
		\node [style=circle] (5) at (-1, -0) {$\psi^{-1}$};
		\node [style=none] (6) at (-1, -1) {};
		\node [style=none] (7) at (-2.5, 0.25) {$\epsilon*$};
		\node [style=none] (8) at (-1.5, 2.75) {$~^*\eta$};
		\node [style=none] (9) at (-3.25, 2.75) {$X$};
		\node [style=none] (10) at (-0.5, 1.75) {$~^*(X^*)$};
		\node [style=none] (11) at (-0.5, -0.7) {$X^{**}$};
		\node [style=none] (12) at (-2.35, 1.5) {$X^*$};
	\end{pgfonlayer}
	\begin{pgfonlayer}{edgelayer}
		\draw (0.center) to (1.center);
		\draw [bend right=90, looseness=1.75] (1.center) to (2.center);
		\draw (2.center) to (3.center);
		\draw [bend left=90, looseness=2.00] (3.center) to (4.center);
		\draw (4.center) to (5);
		\draw (5) to (6.center);
	\end{pgfonlayer}
\end{tikzpicture} \in \X ~~~~~~~~~~ \epsilon_\ox' = \begin{tikzpicture}
	\begin{pgfonlayer}{nodelayer}
		\node [style=none] (0) at (-0.75, -0.5) {};
		\node [style=none] (1) at (-0.75, 2.25) {};
		\node [style=none] (2) at (0.25, 2.25) {};
		\node [style=none] (3) at (0.25, 1) {};
		\node [style=none] (4) at (1.25, 1) {};
		\node [style=none] (5) at (1.25, 3) {};
		\node [style=none] (6) at (0.75, 0.25) {$*\epsilon$};
		\node [style=none] (7) at (-0.25, 3) {$\eta*$};
		\node [style=none] (8) at (1.5, 2.5) {$X^*$};
		\node [style=none] (9) at (-1.75, -0.5) {};
		\node [style=circle, scale=2] (10) at (-1.75, 0.5) {};
		\node [style=none] (11) at (-1.75, 1.25) {};
		\node [style=none] (12) at (-2.75, 1.25) {};
		\node [style=none] (13) at (-2.75, -1.25) {};
		\node [style=none] (14) at (-1.25, -1.25) {$\epsilon*$};
		\node [style=none] (15) at (-2.25, 2) {$\eta*$};
		\node [style=none] (16) at (-1.75, 0.5) {$\psi$};
		\node [style=none] (17) at (-3.25, -0.75) {$X^{**}$};
		\node [style=none] (18) at (-0.25, -0) {$~^*(X^*)$};
	\end{pgfonlayer}
	\begin{pgfonlayer}{edgelayer}
		\draw [bend left=90, looseness=2.00] (1.center) to (2.center);
		\draw (2.center) to (3.center);
		\draw [bend right=90, looseness=2.00] (3.center) to (4.center);
		\draw (4.center) to (5.center);
		\draw (1.center) to (0.center);
		\draw [bend right=90, looseness=1.75] (9.center) to (0.center);
		\draw (11.center) to (10);
		\draw (10) to (9.center);
		\draw [bend right=90, looseness=1.75] (11.center) to (12.center);
		\draw (12.center) to (13.center);
	\end{pgfonlayer}
\end{tikzpicture} = %dualizor1
 \begin{tikzpicture}
	\begin{pgfonlayer}{nodelayer}
		\node [style=none] (0) at (-0.75, 3) {};
		\node [style=none] (1) at (-0.75, -0) {};
		\node [style=none] (2) at (-1.75, -0) {};
		\node [style=none] (3) at (-1.75, 2) {};
		\node [style=none] (4) at (-2.75, 2) {};
		\node [style=none] (5) at (-2.75, -1) {};
		\node [style=none] (6) at (-1.25, -0.75) {$*\epsilon$};
		\node [style=none] (7) at (-2.25, 2.75) {$\eta^*$};
		\node [style=none] (8) at (-0.5, 2.75) {$X$};
		\node [style=none] (9) at (-3.25, -0.5) {$X^{**}$};
		\node [style=circle] (10) at (-1.75, 1) {$\psi$};
		\node [style=none] (11) at (-2.25, 0.25) {$~^*X$};
		\node [style=none] (12) at (-2.1, 1.7) {$X^*$};
	\end{pgfonlayer}
	\begin{pgfonlayer}{edgelayer}
		\draw (0.center) to (1.center);
		\draw [bend left=90, looseness=1.75] (1.center) to (2.center);
		\draw [bend right=90, looseness=2.00] (3.center) to (4.center);
		\draw (4.center) to (5.center);
		\draw (3.center) to (10);
		\draw (10) to (2.center);
	\end{pgfonlayer}
\end{tikzpicture} \stackrel{{\bf \tiny [C.2]}}{=} %dualizor1
 \begin{tikzpicture}
	\begin{pgfonlayer}{nodelayer}
		\node [style=none] (0) at (-3, 3) {};
		\node [style=none] (1) at (-3, 1) {};
		\node [style=none] (2) at (-2, 1) {};
		\node [style=none] (3) at (-2, 2) {};
		\node [style=none] (4) at (-1, 2) {};
		\node [style=circle, scale=2.3] (5) at (-1, -0) {};
		\node [style=none] (13) at (-1, -0) {$\psi^{-1}$};
		\node [style=none] (6) at (-1, -1) {};
		\node [style=none] (7) at (-2.5, 0.25) {$\epsilon*$};
		\node [style=none] (8) at (-1.5, 2.75) {$~^*\eta$};
		\node [style=none] (9) at (-3.25, 2.75) {$X$};
		\node [style=none] (10) at (-0.5, 1.75) {$~^*(X^*)$};
		\node [style=none] (11) at (-0.5, -0.7) {$X^{**}$};
		\node [style=none] (12) at (-2.35, 1.5) {$X^*$};
	\end{pgfonlayer}
	\begin{pgfonlayer}{edgelayer}
		\draw (0.center) to (1.center);
		\draw [bend right=90, looseness=1.75] (1.center) to (2.center);
		\draw (2.center) to (3.center);
		\draw [bend left=90, looseness=2.00] (3.center) to (4.center);
		\draw (4.center) to (5);
		\draw (5) to (6.center);
	\end{pgfonlayer}
\end{tikzpicture} \in \X
\]

The cyclor is a linear transformation which is an isomorphism as it is monoidal with respect to both tensor and par and adjoints are determined only upto isomorphism. It remains to check that the triangle identities hold:
\[
\eta_{\ox'_{X^*}}(\epsilon_\ox')^* = 
\begin{tikzpicture} %triangle1
	\begin{pgfonlayer}{nodelayer}
		\node [style=circle] (0) at (0.75, 0.25) {$\psi^{-1}$};
		\node [style=none] (1) at (0.75, -2.5) {};
		\node [style=none] (2) at (0.75, 2.25) {};
		\node [style=none] (3) at (-0.25, 2.25) {};
		\node [style=none] (4) at (-0.25, 1) {};
		\node [style=none] (5) at (-1.25, 1) {};
		\node [style=none] (6) at (-1.25, 3) {};
		\node [style=none] (7) at (1.5, -0.5) {$X^{***}$};
		\node [style=none] (8) at (-0.75, 0.25) {$\epsilon*$};
		\node [style=none] (9) at (0.25, 3) {$*\eta$};
		\node [style=none] (10) at (-1.5, 2.5) {$X^*$};
		\node [style=none] (11) at (-1, -2.5) {};
		\node [style=circle] (12) at (-1, -1.5) {$\psi^{-1}$};
		\node [style=none] (13) at (-1, -0.75) {};
		\node [style=none] (14) at (-2, -0.75) {};
		\node [style=none] (15) at (-2, -3.25) {};
		\node [style=none] (16) at (-1.5, -0.1) {$*\eta$};
		\node [style=none] (17) at (0, -3.35) {$\epsilon*$};
		\node [style=none] (18) at (-2.5, -3) {};
	\end{pgfonlayer}
	\begin{pgfonlayer}{edgelayer}
		\draw (2.center) to (0);
		\draw (0) to (1.center);
		\draw [bend right=90, looseness=2.00] (2.center) to (3.center);
		\draw (3.center) to (4.center);
		\draw [bend left=90, looseness=2.00] (4.center) to (5.center);
		\draw (5.center) to (6.center);
		\draw [bend right=75, looseness=1.25] (11.center) to (1.center);
		\draw (11.center) to (12);
		\draw (12) to (13.center);
		\draw [bend right=90, looseness=1.50] (13.center) to (14.center);
		\draw (14.center) to (15.center);
	\end{pgfonlayer}
\end{tikzpicture} = \begin{tikzpicture} %triangle2
	\begin{pgfonlayer}{nodelayer}
		\node [style=none] (0) at (0, -2.5) {};
		\node [style=none] (1) at (0, 2.5) {};
		\node [style=none] (2) at (1, 2.5) {};
		\node [style=none] (3) at (1, 1.25) {};
		\node [style=none] (4) at (2, 1.25) {};
		\node [style=none] (5) at (2, 3.25) {};
		\node [style=none] (6) at (-0.5, 1.25) {$X^{***}$};
		\node [style=none] (7) at (1.5, 0.5) {$*\epsilon$};
		\node [style=none] (8) at (0.5, 3.25) {$\eta*$};
		\node [style=none] (9) at (2.25, 2.75) {$X^*$};
		\node [style=none] (10) at (-1, -2.5) {};
		\node [style=circle] (11) at (-1, -1.5) {$\psi^{-1}$};
		\node [style=none] (12) at (-1, -0.75) {};
		\node [style=none] (13) at (-2, -0.75) {};
		\node [style=none] (14) at (-2, -3.25) {};
		\node [style=none] (15) at (-1.5, -0) {$*\eta$};
		\node [style=none] (16) at (-0.5, -3) {$\epsilon*$};
		\node [style=circle] (17) at (1, 2) {$\psi$};
		\node [style=none] (18) at (-2.25, -2.75) {$X^*$};
	\end{pgfonlayer}
	\begin{pgfonlayer}{edgelayer}
		\draw [bend left=90, looseness=2.00] (1.center) to (2.center);
		\draw [bend right=90, looseness=2.00] (3.center) to (4.center);
		\draw (4.center) to (5.center);
		\draw [bend right=75, looseness=1.25] (10.center) to (0.center);
		\draw (10.center) to (11);
		\draw (11) to (12.center);
		\draw [bend right=90, looseness=1.50] (12.center) to (13.center);
		\draw (13.center) to (14.center);
		\draw (1.center) to (0.center);
		\draw (2.center) to (17);
		\draw (17) to (3.center);
	\end{pgfonlayer}
\end{tikzpicture} = \begin{tikzpicture}
	\begin{pgfonlayer}{nodelayer}
		\node [style=none] (0) at (-1, -1) {};
		\node [style=none] (1) at (0, -1) {};
		\node [style=none] (2) at (0, 3.25) {};
		\node [style=none] (3) at (-0.5, -1.75) {$*\epsilon$};
		\node [style=none] (4) at (0.25, 2.75) {$X^*$};
		\node [style=circle] (5) at (-1, 1.5) {$\psi^{-1}$};
		\node [style=none] (6) at (-1, 2.25) {};
		\node [style=none] (7) at (-2, 2.25) {};
		\node [style=none] (8) at (-2, -3.25) {};
		\node [style=none] (9) at (-1.5, 3) {$*\eta$};
		\node [style=circle] (10) at (-1, -0.25) {$\psi$};
		\node [style=none] (11) at (-2.25, -2.75) {$X^*$};
	\end{pgfonlayer}
	\begin{pgfonlayer}{edgelayer}
		\draw [bend right=90, looseness=2.00] (0.center) to (1.center);
		\draw (1.center) to (2.center);
		\draw (5) to (6.center);
		\draw [bend right=90, looseness=1.50] (6.center) to (7.center);
		\draw (7.center) to (8.center);
		\draw (10) to (0.center);
		\draw (5) to (10);
	\end{pgfonlayer}
\end{tikzpicture} = ~ \begin{tikzpicture} \draw (2, 3.25) -- (2,-3.5); \end{tikzpicture} ~ = ~ 1
\] The other triangle identity holds similarly.
\end{proof}

The equality of $\eta'$ and $\epsilon'$ is immediate from {\bf [C.2]} for cyclors with the map $\eta'=\epsilon'$ being the {\bf dualizor}.   In the symmetric case, the dualizor of this 
equivalence may be drawn as:
 \[
 %dualizor1
\begin{tikzpicture}
	\begin{pgfonlayer}{nodelayer}
		\node [style=none] (0) at (-3, 3) {};
		\node [style=none] (1) at (-3, 1) {};
		\node [style=none] (2) at (-2, 1) {};
		\node [style=none] (3) at (-2, 2) {};
		\node [style=none] (4) at (-1, 2) {};
		\node [style=circle] (5) at (-1, -0) {$\psi^{-1}$};
		\node [style=none] (6) at (-1, -1) {};
		\node [style=none] (7) at (-2.5, 0.25) {$\epsilon*$};
		\node [style=none] (8) at (-1.5, 2.75) {$*\eta$};
		\node [style=none] (9) at (-3.25, 2.75) {$A$};
		\node [style=none] (10) at (-0.5, 1.75) {$~^*(A^*)$};
		\node [style=none] (11) at (-0.5, -0.7) {$A^{**}$};
		\node [style=none] (12) at (-2.25, 1.5) {$A^*$};
	\end{pgfonlayer}
	\begin{pgfonlayer}{edgelayer}
		\draw (0.center) to (1.center);
		\draw [bend right=90, looseness=1.75] (1.center) to (2.center);
		\draw (2.center) to (3.center);
		\draw [bend left=90, looseness=2.00] (3.center) to (4.center);
		\draw (4.center) to (5);
		\draw (5) to (6.center);
	\end{pgfonlayer}
\end{tikzpicture} =
\begin{tikzpicture}
\begin{pgfonlayer}{nodelayer}
\node [style=none] (0) at (-1.75, 2) {};
\node [style=none] (1) at (-1.75, -0) {};
\node [style=none] (2) at (0, -0) {};
\node [style=none] (12) at (-0.75, -1) {$\epsilon*$};
\node [style=none] (3) at (0, 1) {};
\node [style=none] (4) at (-1, 1) {};
\node [style=none] (34) at (-0.5, 1.75) {$\eta*$};
\node [style=none] (5) at (-1, -2.75) {};
\node [style=none] (6) at (-2, 1.75) {$A$};
\node [style=none] (8) at (-1.5, -2.25) {$A^{**}$};
\end{pgfonlayer}
\begin{pgfonlayer}{edgelayer}
\draw (0.center) to (1.center);
\draw [bend right=90, looseness=1.50] (1.center) to (2.center);
\draw (2.center) to (3.center);
\draw [bend right=90, looseness=1.75] (3.center) to (4.center);
\draw (4.center) to (5.center);
\end{pgfonlayer}
\end{tikzpicture} = %dualizor2
 \begin{tikzpicture}
	\begin{pgfonlayer}{nodelayer}
		\node [style=none] (0) at (-0.75, 3) {};
		\node [style=none] (1) at (-0.75, -0) {};
		\node [style=none] (2) at (-1.75, -0) {};
		\node [style=none] (3) at (-1.75, 2) {};
		\node [style=none] (4) at (-2.75, 2) {};
		\node [style=none] (5) at (-2.75, -1) {};
		\node [style=none] (6) at (-1.25, -0.75) {$*\epsilon$};
		\node [style=none] (7) at (-2.25, 2.75) {$\eta*$};
		\node [style=none] (8) at (-0.5, 2.75) {$A$};
		\node [style=none] (9) at (-3.25, -0.5) {$A^{**}$};
		\node [style=circle] (10) at (-1.75, 1) {$\psi$};
		\node [style=none] (11) at (-2.25, 0.25) {$~^*A$};
		\node [style=none] (12) at (-2, 1.7) {$A^*$};
	\end{pgfonlayer}
	\begin{pgfonlayer}{edgelayer}
		\draw (0.center) to (1.center);
		\draw [bend left=90, looseness=1.75] (1.center) to (2.center);
		\draw [bend right=90, looseness=2.00] (3.center) to (4.center);
		\draw (4.center) to (5.center);
		\draw (3.center) to (10);
		\draw (10) to (2.center);
	\end{pgfonlayer}
\end{tikzpicture}
\]
%%%%%%%%%%%%%%%%%%%%%%%%%%%%%%%%%%%%%%%%%%%%%%%%%%%%%%%%%%%%%%%%%%%%%%

\subsection{Conjugation}
\label{Sec: conjugation}

Recall the following structure from Egger \cite{Egg11}:

\begin{defi}
A {\bf conjugation} for a monoidal category $(X, \otimes, I)$ consists of a functor $\overline{(\_)}: \X^\rev \to \X$ with natural isomorphisms:
\[ \bar{A} \otimes \bar{B} \to^{\chi} \bar{B \otimes A} ~~~~~~~~~~~~~ \bar{\bar{A}} \to^{\varepsilon} A  \]
called respectively the (tensor reversing) {\bf conjugating laxor} and the {\bf conjugator}
such that 
\[\bar{\bar{\bar{A}}} \to^{\bar{\varepsilon_A} = \varepsilon_{\bar{A}} } \bar{A} \] 
and  
\[
\xymatrix{
(\bar{A} \ox \bar{B}) \ox \bar{C} \ar[rr]^{a_\ox} \ar[d]_{\chi \ox 1}  \ar@{}[ddrr]|{\bf [CF.1]_\ox}  & & \bar{A} \ox (\bar{B} \ox \bar{C}) \ar[d]^{1 \ox \chi} \\
\bar{(B \ox A)} \ox \bar{C} \ar[d]_{\chi} & & \bar{A} \ox \bar{(C \ox B)} \ar[d]^{\chi} \\
\bar{C \ox (B \ox A)} \ar[rr]_{\bar{a_\otimes^{-1}}} & & \bar{(C \ox B) \ox A}
} ~~~~~~~~ \xymatrix{
\bar{\bar{A}} \ox \bar{\bar{B}} \ar[rr]^{\chi} \ar[dd]_{\varepsilon \ox \varepsilon}  \ar@{}[ddrr]|{\bf [CF.2]_\ox}  & & \bar{\bar{B} \ox \bar{A}} \ar[dd]^{\bar{\chi}} \\ 
& & \\
A \ox B & & \bar{\bar{A \ox B}} \ar[ll]^{\varepsilon}
}
\]
\end{defi}

A monoidal category is {\bf conjugative} when it has a conjugation functor.

A symmetric monoidal category, which is conjugative, is {\bf symmetric conjugative} in case it satisfies the additional coherence:
\[ 
\xymatrix{
\bar{A} \otimes \bar{B} \ar[d]_{c_\otimes} \ar[rr]^{\chi}  \ar@{}[drr]|{\bf [CF.3]_\ox}  &  & \bar{B \otimes A} \ar[d]^{\bar{c_\otimes}} \\
\bar{B} \otimes \bar{A} \ar[rr]_{\chi} &  & \bar{A \otimes B} 
}
\]

Notice that we have not specified any coherence for the unit $I$. This is because the expected coherences are automatic:

\begin{lem}
\label{Lemma: unit conjugate} \cite[Lemma 2.3]{Egg11}
For every conjugative monoidal category, there exists a unique isomorphism $I \xrightarrow{\chi^{\!\!\!\circ}} \bar{I}$ such that 
\[
\xymatrix{
I \otimes \bar{A} \ar[rr]^{\chi^{\!\!\!\circ} \otimes 1} \ar[d]_{u_\otimes}  \ar@{}[drr]|{\bf [CF.4]_\top}  & & \bar{I} \otimes \bar{A} \ar[d]^{\chi} \\
\bar{A} \ar[rr]_{\bar{u_\otimes^{-1}}} & & \bar{A \otimes I } } ~~~~~ 
\xymatrix{ \bar{A} \ox I \ar[rr]^{1 \otimes \chi^{\!\!\!\circ}} \ar[d]_{u_\otimes}  \ar@{}[drr]|{\bf [CF.5]_\top}  & & \bar{A} \ox \bar{I} \ar[d]^{\chi} \\
\bar{A} \ar[rr]_{\bar{u_\otimes^{-1}}} & & \bar{I \otimes A} } ~~~~~
\xymatrix{
I \ar[rr]^{\chi^{\!\!\!\circ}} \ar@{=}[d]  \ar@{}[drr]|{\bf [CF.6]_\top}  & & \bar{I} \ar[d]^{\bar{\chi^{\!\!\!\circ}}} \\
I \ar[rr]_{\varepsilon^{-1}} & & \bar{\bar{I}}
}
\]
\end{lem}

\begin{defi} \cite{Egg11}
A {\bf conjugative LDC} is a linearly distributive category $(\X, \ox, \top, \oa, \bot)$ together with a conjugating functor $\bar{(\_)}: \X \to \X$ and natural isomorphisms:
\[ \bar{A} \ox \bar{B} \xrightarrow{\chi_\ox} \bar{B \ox A} ~~~~~~~~~~~~ \bar{A \oa B} \xrightarrow{\chi_\oa} \bar{B} \oa \bar{A} ~~~~~~~~~~~~ \bar{\bar{A}} \xrightarrow{\varepsilon} A \] 
\end{defi}

such that $(\X, \ox, \top, \chi_\ox, \varepsilon)$ and $(\X, \oa, \bot, \chi_\oa^{-1}, \varepsilon)$ are conjugative (symmetric) monoidal categories with respect to the conjugating functor and  the following diagrams commute:

\[ 
\xymatrix{
\bar{B \oa C} \ox \bar{A} \ar[rr]^{\chi_\oa \ox 1} \ar[d]_{\chi_\ox}  \ar@{}[ddrr]|{\bf [CF.7]}  & & (\bar{C} \oa \bar{B}) \ox \bar{A} \ar[d]^{\partial} \\
\bar{(A \ox (B \oa C))} \ar[d]_{\bar{\partial}} & & \bar{C} \oa ( \bar{B} \ox \bar{A}) \ar[d]^{1 \oa \chi_\ox} \\
\bar{((A \ox B) \oa C)} \ar[rr]_{\chi_\oa} & & \bar{C} \oa \bar{A \ox B}
} ~~~~~~~~~~ 
\xymatrix{
\bar{A} \ox \bar{C \oa B} \ar[rr]^{\chi_\ox} \ar[d]_{1 \ox \chi_\oa} \ar@{}[ddrr]|{\bf [CF.8]}   & & \bar{(C \oa B) \ox A} \ar[d]^{\bar{\partial}} \\
\bar{A} \ox (\bar{B} \oa \bar{C}) \ar[d]_{\partial} & & \bar{C \oa (B \ox A)} \ar[d]^{\chi_\oa} \\
(\bar{A} \ox \bar{B}) \oa \bar{C} \ar[rr]_{\chi_\ox \oa 1} & & \bar{(B \ox A)} \oa \bar{C}
} 
\]

Note, by Lemma \ref{Lemma: unit conjugate}, there exists canonical isomorphisms $\top \xrightarrow{\chi^{\!\!\!\circ}_\top} \bar{\top}$ and $\bot \xrightarrow{\chi^{\!\!\!\circ}_\bot} \bar{\bot}$, hence conjugation is a normal functor.  However, the conjugation is not necessarily a mix functor when $\X$ is a mix category.   For
conjugation to be a mix functor, the following extra condition must be satisfied:
\[ {\bf \small [CF.9]}~~~~~ \xymatrix{\overline{\bot} \ar[rd]_{(\chi^{\!\!\!\circ}_\bot)^{-1}} \ar@/^/[rrr]^{\overline{{\sf m}}} & & & \overline{\top} \\ & \bot \ar[r]_{{\sf m}} & \top \ar[ur]_{\chi^{\!\!\!\circ}_\top} } \]

\begin{prop}
A conjugative LDC is precisely a LDC, $\X$, with a Frobenius adjoint $(\epsilon^{-1}, \epsilon): \overline{(\_)} \dashv \overline{(\_)}^\rev: \X^\rev \to \X$ where $\epsilon := (\varepsilon,\varepsilon^{-1})$.  Furthermore, if $\X$ is an isomix category and conjugation is a mix functor then conjugation is an isomix equivalence.
\end{prop}

\begin{proof}
It is clear that $\overline{(\_)}$ is a strong Frobenius functor so being mix implies isomix.   Also, $\varepsilon$ is clearly monoidal for tensor and par.  The triangle equalities give $ \overline{\varepsilon^{-1}} \varepsilon = 1: \overline{A} \to \overline{A}$ thus $\varepsilon = \overline{\varepsilon}$.
\end{proof}

Clearly conjugation should flip left duals into right duals:

\begin{lem}
\label{Lemma: involutive linear adjoint}
If $B \dashvv A$ is a linear dual then $\bar{A} \dashvv \bar{B}$ is a linear dual.
\end{lem}

\begin{proof}
Suppose $(\eta, \varepsilon): B \dashvv A$. Then, $(\chi^{\!\!\!\circ}_\top \bar{\eta} \chi_\oa,\chi_\ox \bar{\varepsilon} \chi^{\!\!\!\circ}_\bot): \bar{A} \dashvv \bar{B}$.
\end{proof}

When a $*$-autonomous category is cyclic one expects that conjugation should interact with the cyclor in a coherent fashion:

\begin{defi}\cite{EggMcCurd12}
\label{Defn: conjugative cyclic}
A {\bf conjugative cyclic $*$-autonomous category} is a conjugative $*$-autonomous category together with a cyclor $A^* \to^{\psi} \!\!~^{*}\!A$ such that 
\[
\xymatrix{
(\bar{A})^* \ar[rr]^{\psi} \ar[d]_{\simeq}  & & ^{*}(\bar{A}) \ar[d]^{\simeq} \\
\bar{(^{*}A)} \ar[rr]_{\bar{\psi^{-1}}} & & \bar{(A^*)}
}
\]
which gives a map $\sigma: (\overline{A})^{*} \to \overline{(A^{*})}$.
\end{defi}

The above condition is drawn as follows:
\[ \sigma =
\begin{tikzpicture} %cca1
	\begin{pgfonlayer}{nodelayer}
		\node [style=none] (0) at (-1, 0.75) {$\overline{\psi^{-1}}$};
		\node [style=circle, scale=2.5] (0) at (-1, 0.75) {};
		\node [style=none] (1) at (-1, -0.5) {};
		\node [style=none] (2) at (-1, 2.25) {};
		\node [style=none] (3) at (0, 2.25) {};
		\node [style=none] (4) at (0, 1) {};
		\node [style=none] (5) at (1, 1) {};
		\node [style=none] (6) at (1, 3) {};
		\node [style=none] (7) at (1.5, 2.75) {$(\overline{A})^*$};
		\node [style=none] (8) at (0.5, 0.2) {$\epsilon^*$};
		\node [style=none] (9) at (-0.5, 3.1) {$\overline{~^*\eta}$};
		\node [style=none] (10) at (-1.5, -0.25) {$\overline{A^*}$};
		\node [style=none] (12) at (-0.25, 1.5) {$\overline{A}$};
		\node [style=none] (11) at (-1.5, 1.75) {$\overline{~^*A}$};
	\end{pgfonlayer}
	\begin{pgfonlayer}{edgelayer}
		\draw (2.center) to (0);
		\draw (0) to (1.center);
		\draw [bend left=90, looseness=2.00] (2.center) to (3.center);
		\draw (3.center) to (4.center);
		\draw [bend right=90, looseness=2.00] (4.center) to (5.center);
		\draw (5.center) to (6.center);
	\end{pgfonlayer}
\end{tikzpicture} = \begin{tikzpicture} %cca1
	\begin{pgfonlayer}{nodelayer}
		\node [style=circle] (0) at (-1, 2.5) {$\psi$};
		\node [style=none] (1) at (-1, 3.25) {};
		\node [style=none] (2) at (-1, 0.5) {};
		\node [style=none] (3) at (0, 0.5) {};
		\node [style=none] (4) at (0, 1.75) {};
		\node [style=none] (5) at (1, 1.75) {};
		\node [style=none] (6) at (1, -0.25) {};
		\node [style=none] (7) at (-1.75, 3.25) {$(\overline{A})^*$};
		\node [style=none] (8) at (0.5, 2.6) {$\overline{\eta^*}$};
		\node [style=none] (9) at (-0.5, -0.35) {$~^*\epsilon$};
		\node [style=none] (10) at (1.5, 0.25) {$\overline{A^*}$};
		\node [style=none] (11) at (-1.75, 1.5) {$~^*(\overline{A})$};
		\node [style=none] (12) at (-0.25, 1.25) {$\overline{A}$};
	\end{pgfonlayer}
	\begin{pgfonlayer}{edgelayer}
		\draw (2.center) to (0);
		\draw (0) to (1.center);
		\draw [bend right=90, looseness=2.00] (2.center) to (3.center);
		\draw (3.center) to (4.center);
		\draw [bend left=90, looseness=2.00] (4.center) to (5.center);
		\draw (5.center) to (6.center);
	\end{pgfonlayer}
\end{tikzpicture}
\]

When the $*$-autonomous category is symmetric, conjugation automatically preserves the canonical cyclor.

\begin{lem}
\label{Lemma: varepsi monoidal}
In a conjugative $*$-autonomous category, 
\[ %epsi1 
\begin{tikzpicture} 
	\begin{pgfonlayer}{nodelayer}
		\node [style=circle] (0) at (0, 1) {$\varepsilon$};
		\node [style=none] (1) at (0, -0) {};
		\node [style=none] (2) at (-1.75, -0) {};
		\node [style=none] (3) at (-1.75, 2) {};
		\node [style=none] (4) at (0, 2) {};
		\node [style=none] (5) at (-1, 3) {$\overline{\overline{\eta*}}$};
		\node [style=none] (6) at (-2.25, 0.3) {$\overline{\overline{X^*}}$};
		\node [style=none] (7) at (0.5, 1.7) {$\overline{\overline{X}}$};
		\node [style=none] (8) at (0.5, 0.25) {$X$};
	\end{pgfonlayer}
	\begin{pgfonlayer}{edgelayer}
		\draw (4.center) to (0);
		\draw (0) to (1.center);
		\draw [bend left=90, looseness=1.50] (3.center) to (4.center);
		\draw (3.center) to (2.center);
	\end{pgfonlayer}
\end{tikzpicture} = %epsi2
\begin{tikzpicture}
	\begin{pgfonlayer}{nodelayer}
		\node [style=circle] (0) at (-1.75, 1) {$\varepsilon^{-1}$};
		\node [style=none] (1) at (-1.75, -0) {};
		\node [style=none] (2) at (0, -0) {};
		\node [style=none] (3) at (0, 2) {};
		\node [style=none] (4) at (-1.75, 2) {};
		\node [style=none] (5) at (-1, 3) {$\eta*$};
		\node [style=none] (6) at (-2.5, 0.3) {$\overline{\overline{X^*}}$};
		\node [style=none] (7) at (0.5, 0.5) {$X$};
		\node [style=none] (8) at (-2.25, 1.7) {$X^*$};
	\end{pgfonlayer}
	\begin{pgfonlayer}{edgelayer}
		\draw (4.center) to (0);
		\draw (0) to (1.center);
		\draw [bend right=90, looseness=1.50] (3.center) to (4.center);
		\draw (3.center) to (2.center);
	\end{pgfonlayer}
\end{tikzpicture}
\]
\[
\chi_\top^{\!\!\!\circ}~ \overline{\chi_\top^{\!\!\!\circ}} ~ \overline{\overline{\eta*}} ~ \overline{\chi_\oa}^{-1} \chi_\oa^{-1} (1 \oa \varepsilon) = \eta* (\varepsilon^{-1} \oa 1)  : \top \to \overline{\overline{X^*}} \oa X 
\]
\end{lem}
\begin{proof}
\begin{align*}
\chi_\top^{\!\!\!\circ} ~ \overline{\chi_\top^{\!\!\!\circ}} ~ \overline{\overline{\eta*}} ~ \overline{\chi}^{-1} \chi^{-1} (1 \oa \varepsilon) &=  \chi_\top^{\!\!\!\circ} ~ \overline{\chi_\top^{\!\!\!\circ}} ~ \overline{\overline{\eta*}} ~ \overline{\chi}^{-1} \chi^{-1} (\varepsilon \varepsilon^{-1} \oa \varepsilon) \\
&=\chi_\top^{\!\!\!\circ} ~ \overline{\chi_\top^{\!\!\!\circ}}  ~ \overline{\overline{\eta*}} ~ \overline{\chi}^{-1} \chi^{-1} (\varepsilon  \oa \varepsilon) (\varepsilon^{-1} \oa 1) \\
&\stackrel{{\bf \small [CF.2]_\oa}}{=} \chi^{\!\!\!\circ} ~ \overline{\chi_\top^{\!\!\!\circ}} ~ \overline{\overline{\eta*}} \varepsilon (\varepsilon^{-1} \oa 1) \\
&\stackrel{{\bf \small nat.}}{=} \chi_\top^{\!\!\!\circ} ~ \overline{\chi_\top^{\!\!\!\circ}} ~  \varepsilon \eta\!* (\varepsilon^{-1} \oa 1) \\
&\stackrel{{\bf \small [CF.6]_\top}}{=} \eta\!* (\varepsilon^{-1} \oa 1) \qedhere
\end{align*}
\end{proof}

%%%%%%%%%%%%%%%%%%%%%%%%%%%%%%%%%%%%%%%%%%%%%%%%%%%%%%%%%%%%%%%%%%%%%%

\subsection{Dagger and conjugation}

The interaction of the dagger and conjugation for cyclic $*$-autonomous categories in the presence of the dualizing functor is illustrated by the following diagram:
\[
\xymatrixcolsep{5pc}
\xymatrixrowsep{5pc}
\xymatrix{
\X^{\op} \ar@/^1pc/[rr]^{(\_)^\dagger} \ar@/^1pc/[dr]|{((\_)^*)^\rev} \ar@{}[rr]|{\bot} &~ \ar@{}[d]|{\cong} & 
\X  \ar@/^1pc/[ll]|{((\_)^{\dagger})^{\op}} \ar@/_1pc/[dl]|{\overline{(\_)}^\rev} \ar@{}[dl]|{\dashv} \\
 & \X^\rev \ar@/^1pc/[ul]^{(\_)^{*^{\op}}} \ar@/_1pc/[ur]_{\overline{(\_)}} \ar@{}[ul]|{\dashv} &
 }
\]
Specifically we have: 

\begin{thm}
Every cyclic, $\dagger$-$*$-autonomous category is a conjugative $*$-autonomous category.
\end{thm}
\begin{proof}
Let $\X$ be a cyclic, $\dagger$-$*$-autonomous category. Then composing adjoints gives the equivalence  $(\_)^{\dagger^*} \dashv (\_)^{*^\dagger}$. To build a conjugation, however,
we need an equivalence between the same functors: to obtain such an equivalence we use the natural equivalence $\omega: (\_)^{\dagger*} \to (\_)^{*\dagger} $ from the cyclor preserving condition for Frobenius linear functors.  A conjugative equivalence, in addition, requires that the unit and counit of the equivalence be inverses of each other.  
The unit and counit of the equivalence are given by {\em (a)} and {\em (b)} respectively;
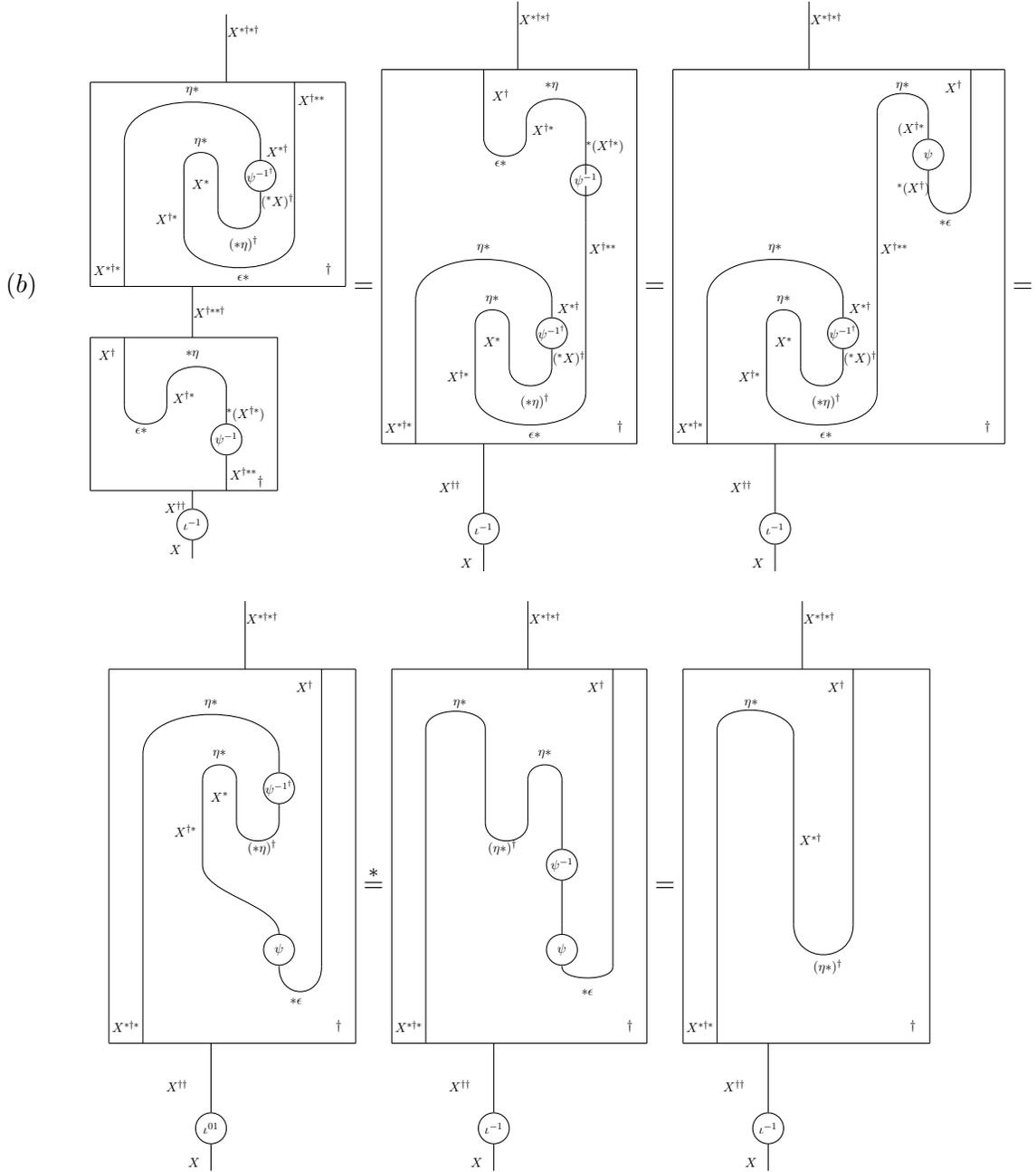
\[ \mbox{\em (a)}~~~~~~~~~  \xymatrix{
\X^{\sf rev} \ar@{=}[rrrr] \ar[dr]_{(\_)^{*^\op}} & 
&  
\ar@{}[d]|{\Downarrow ~ \eta_\ox'}&  
& 
\X^{\sf rev} \ar@{<-}[ld]^{(\_)^{*^\op}}  
& \\
&
 \X^{\op} \ar@{=}[rr] \ar[dr]_{\dagger} &  
 \ar@{}[d]|{\Downarrow ~ \iota^{-1}} & 
 \X^{\sf op} \ar@{<-}[ld]^{\dagger^\op} \ar@{}[rr]_{\omega}  \ar@{}[rr]^{\Longrightarrow}  & 
 & 
 \X^{\sf oprev} \ar@{<-}@/^1pc/[llld]^{(\_)^{*^{\sf oprev}}} \ar@/_1pc/[ul]_{\dagger^\rev}  \\
& & \X & &
}
\]

\[
\mbox{\em (b)}~~~~~~~~~~~~~~ \xymatrix{
& & &
\X^{\rev}  \ar@{<-}@/_1pc/[llld]_{\dagger^\rev} \ar@{<-}[dl]_{(\_)^{*^{\rev}}} \ar[dr]^{(\_)^{*^{\sf op}}} \ar@{}[d]|{\Downarrow ~ \epsilon'_\ox} & &\\
\X^{\sf op rev} \ar@{<-}[rd]_{(\_)^{*^{\sf oprev}}} \ar@{}[rr]_{\omega^{-1}}  \ar@{}[rr]^{\Longrightarrow}  &  & 
\X^{\sf op} \ar@{=}[rr] \ar@{<-}[dl]^{\dagger^\op} & \ar@{}[d]|{\Downarrow ~ \iota^{-1}} &
\X^{\sf op} \ar[dr]^{\dagger} \\
& \X \ar@{=}[rrrr]&  & & & \X }
\] where the isomorphism $\omega: (\_)^{\dagger*} \to (\_)^{*\dagger}$ is from the cyclor preserving condition, {\bf [CFF]}, for Frobenius linear functors:
\[
\omega :=
\]

It remains to show that the unit and the counit maps are inverses of each other in $\X$:

\[ (a) ~~~~~ %dbig2
 %
 \stackrel{\ref{Lemma:  cyclic dagger}}{=} %big8
%
\]
 $(*)$ holds because $\dagger$ preserves the cyclor. Thus, $(a)$ and $(b)$ are inverses of each other.
\end{proof}

Next, we show that a conjugation functor together with a dualizing functors gives a $\dagger$:

\begin{thm}
\label{Theorem: conjugation+dualizing}
Every cyclic, conjugative $*$-autonomous category  is also a $\dagger$-$*$-autonomous category.
\end{thm}
\begin{proof}
Let $\X$ be a cyclic, conjugative $*$-autonomous category then  $\bar{(\_)^*} \dashv \bar{(\_)}^*$ is an equivalence. To build a dagger we need an equivalence on the same 
functor: we obtain this by using the natural equivalence $\sigma: \bar{(\_)^*} \to \bar{(\_)}^*$ from Definition \ref{Defn: conjugative cyclic}.  An involutive equivalence, in addition, 
requires the unit and counit of the (contravariant) equivalence to be the same map (which we called the involutor, $\iota$). We show that this is the case: 

The unit and counit of the equivalence is given by {\em (a)} and {\em (b)} respectively;
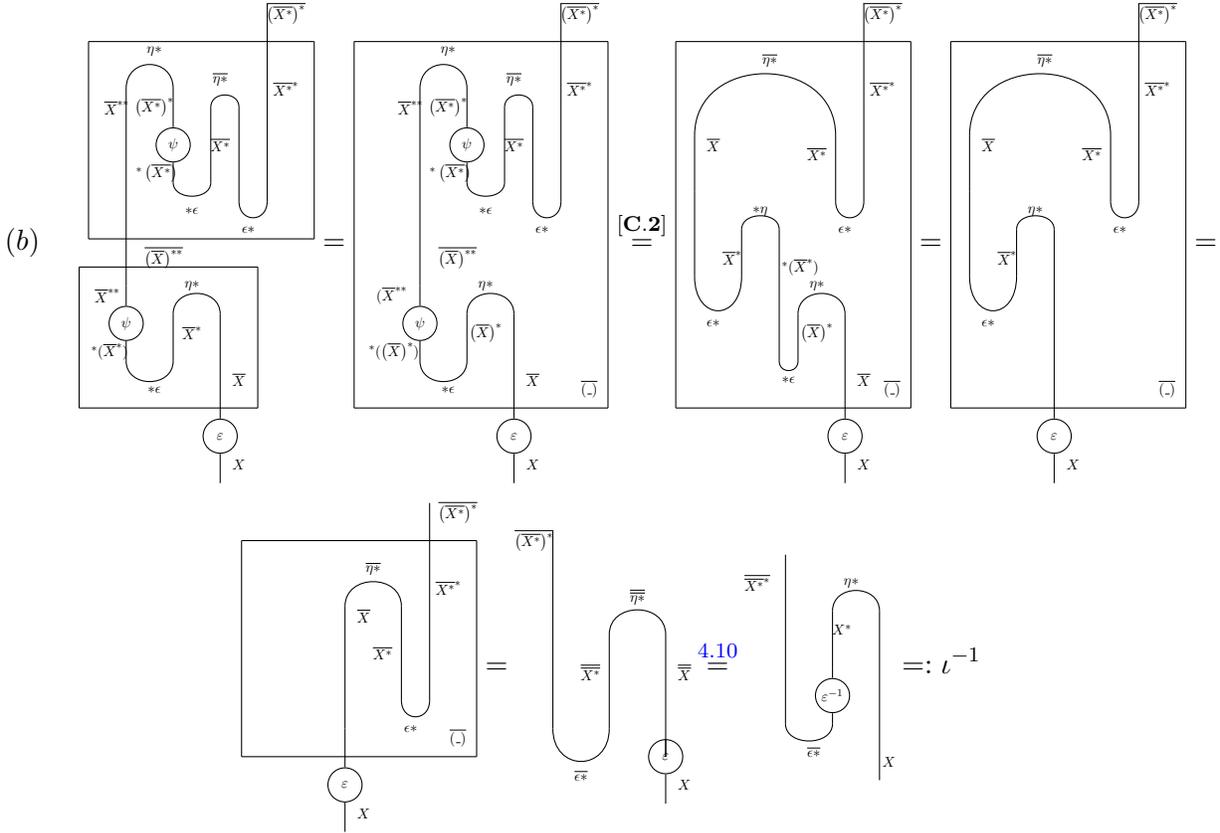
\[ \mbox{\em (a)}~~~~~~~~~  \xymatrix{
\X^{\op} \ar@{=}[rrrr] \ar[dr]_{(\_)^{*^{\sf rev}}} & 
&  
\ar@{}[d]|{\Downarrow ~ \eta_\ox'}&  
& 
\X^{\op} \ar@{<-}[ld]^{(\_)^{*^{\sf op}}}  
& \\
&
 \X^{\sf rev} \ar@{=}[rr] \ar[dr]_{\bar{(\_)}} &  
 \ar@{}[d]|{\Downarrow ~ \varepsilon^{-1}} & 
 \X^{\sf rev} \ar@{<-}[ld]^{\bar{(\_)}^{\sf rev}} \ar@{}[rr]_{\sigma}  \ar@{}[rr]^{\Longrightarrow}  & 
 & 
 \X^{\sf oprev} \ar@{<-}@/^1pc/[llld]^{(\_)^{*^{\sf oprev}}} \ar@/_1pc/[ul]_{\bar{(\_)}^\op}  \\
& & \X & &
}
\]

\[
\mbox{\em (b)}~~~~~~~~~~~~~~ \xymatrix{
& & &
\X^{\op}  \ar@{<-}@/_1pc/[llld]_{\bar{(\_)}^\op} \ar@{<-}[dl]_{(\_)^{*^{\op}}} \ar[dr]^{(\_)^{*^{\sf rev}}} \ar@{}[d]|{\Downarrow ~ \epsilon'_\ox} & &\\
\X^{\sf op rev} \ar@{<-}[rd]_{(\_)^{*^{\sf oprev}}} \ar@{}[rr]_{\sigma^{-1}}  \ar@{}[rr]^{\Longrightarrow}  &  & 
\X^{\sf rev} \ar@{=}[rr] \ar@{<-}[dl]^{\bar{(\_)}^{\sf rev}} & \ar@{}[d]|{\Downarrow ~ \varepsilon} &
\X^{\sf rev} \ar[dr]^{\bar{(\_)}} \\
& \X \ar@{=}[rrrr]&  & & & \X }
\] where $\sigma: \bar{A}^* \to \bar{A^*}$ is given in Definition \ref{Defn: conjugative cyclic}.  Below we show that the unit and counit coincide in $\X$.

\[
(a)~~~~~ %big1 
 \stackrel{\ref{Lemma: varepsi monoidal}}{=}
%big8
%
 =: \iota^{-1}
\] 
\end{proof}

Observe that for composition of the dualizing functor and the conjugation functor to yield a dagger, and vice versa, a $*$-autonomous category is required to be cyclic with the cyclor being preserved by the conjugation (see Definition \ref{Defn: conjugative cyclic}) and the dagger (see just before Lemma \ref{Lemma: cyclic dagger}). 

\subsection{Examples}

In this section, we cover examples of $\dagger$-isomix categories where the $\dagger$ is given by conjugation and the dualizing functor. 
\subsubsection{Category of a group with conjugation}

\begin{defi}
A {\bf group with conjugation} is a group $(G, ., e)$  together with a function $\overline{(\_)}: G \to G$ such that, for all $g \in G$, $\overline{\overline{g}} = g$, and for all $g, h \in G$, $\overline{g . h} = \overline{h} \overline{g}$, and $\overline e = e$.
\end{defi}

Let $(G,., e)$ be a group with conjugation. The discrete category $\D{ (G,., e})$ whose objects are the elements of the group is a monoidal category with the tensor product given by $g \ox h := g.h$, and the monoidal unit $e$. Moreover, $\D{(G,.,e)}$ is a compact closed category where $g^* := g^{-1}$ and it has a trivial conjugative cyclor (See Definition \ref{Defn: conjugative cyclic}). Thus, $\D{(G,.,e})$ is a compact $\dagger$-isomix-$*$-autonomous category with $g^\dagger := \overline{g^*}$ gives an example of how the conjugation gives rise to a dagger.

Here are some examples of groups with conjugation and the discrete categories given by them:
\begin{itemize}
\item Suppose we fix the group to be $(\C, +, 0)$ where the objects are complex numbers and the tensor product is addition. The dual and conjugation of complex numbers are given as follows: $(a+ib)^* = -a - ib$ and $\overline{a + ib} := a - ib$. Hence, \[(a + ib)^\dagger := \overline{(a+ib)^*} = \overline{(-a-ib)} = -a + ib\] 
\item Consider the multiplicative group $(\C^*, ., 1)$ where the objects are non-zero complex numbers and the tensor product is given by multiplication. The dualizing and the conjugation functors are given as follows: 
\[ (a+ib)^* = c+id, \text{ where } ac-bd=1 \text { and }  ad+bc=0 \]  \[ \overline{a + ib} := a-ib \] $(a+ib)^\dagger$ is given by $\overline{(a+ib)^*}$.
\item Suppose the group is fixed to be $\D{(P(x), +, 0})$ where $P(x)$ is a polynomial ring. 

$\D{(P(x), +, 0})$ is a conjugative compact closed category: $P(x)^* = -P(x)$ and $\overline{P(x)} = P(-x)$. Then, $P(x)^\dagger = -P(-x)$.
\item Consider the general linear group of degree 2, $(\mathbb{M}_2, . ,I_2)$ over complex numbers. Then, the discrete category $\D(\mathbb{M}_2, ., I_2)$ has a dualizing functor given by matrix inverse and conjugation is given by conjugate transpose: $\overline{\left(
\begin{matrix}
a+ib & m+in \\
c+id & p+iq
\end{matrix}
\right)} :=  \left(
\begin{matrix}
a-ib & c-id \\
m-in & p-iq
\end{matrix}
\right)
$. Then,  $\D(\mathbb{M}_2, .,I_2)$ is a $\dagger$-isomix $*$-autonomous category with:
 \[ \left(
  \begin{matrix}
 a+ib & m+in \\
 c+id & p+iq 
 \end{matrix}
 \right)^\dagger := \overline{\left(
 \begin{matrix}
a+ib & m+in \\
c+id & p+iq 
 \end{matrix}
 \right)^*} =  \left(
 \begin{matrix}
a-ib & c-id \\
m-in & p-iq 
\end{matrix}
 \right)^{-1} \] 
\end{itemize}

%The example in subsection \ref{subsection: discrete example} is an example for composition of conjugation and dualizing functor producing a $\dagger$ functor for LDCs. In subsection \ref{subsection: discrete example}, the discrete category built from an Abelian group is a conjugative monoidal category with the conjugation functor given by $\overline{(\_)}: \D{(G,., e})^{\rev} \to \D{(G,., e})$ given by the group involution. Then,  $\D{(G,.,e})$ is a compact $\dagger$ isomix $*$-autonomous category with $g^\dagger := \overline{g^*}$. 

%%%%%%%%%%%%%%%%%%%%%%%%%%%%%%%%%%%%%%%%%%%%%%%%%%%%%%%
\subsubsection{Chu Spaces}
\label{Section: Chu}

Applications of Chu Spaces to represent quantum systems have been studied in \cite{Abr12}, \cite{Abr13}. 
In this section we show that the Chu construction over a closed conjugative monoidal category,
 which  has pullbacks, produces a $\dagger$-isomix LDC, ${\sf Chu}_\X(I)$.  To get 
 the $*$-autonomous category and $\dagger$-structure on ${\sf Chu}_\X(I)$ we shall start 
 by explaining how one can produce  conjugative structure on the ${\sf Chu}$ category.  
 To achieve this we iteratively develop the structure of this category, starting with a 
 conjugative closed monoidal category, $\X$, which is not necessarily  symmetric.  
 Note that the fact that it is conjugative means that it is both left and right closed which 
 allows us to consider the non-commutative ${\sf Chu}$ construction: in this regard  
 we shall follow J\"urgen Koslowski's construction \cite{Jur06} using simplified ``Chu-cells'' on the 
 same dualizing object to obtain not a ${*}$-linear bicategory but a cyclic  $*$-autonomous category.  
 Furthermore, we shall choose a dualizing object which is conjugative in order  to obtain a conjugative 
 cyclic $*$-autonomous  category.

A conjugative object is an object $D$ of $\X$ with an isomorphism $d: \overline{D} \to D$ such that $\overline{d}d = \varepsilon: D \to \overline{\overline{D}}$.  
We can then define ${\sf Chu}_\X(D)$ as follows:

\begin{description}
\item[Objects] $(A, B, \psi_0, \psi_1)$ where $\psi_0: A \ox B \to D$ and $\psi_1: B \ox A \to D$  in $\X$ (these are the simplified Chu cells). %What are the simplified chu cells?
\item[Arrows] $(f,g): (A, B, \psi_0,\psi_1) \to (A', B', \psi_0',\psi_1')$ where $f: A \to A'$ and $g: B' \to B$ and the following diagrams commutes:
\[ \xymatrix{
& A \ox B' \ar[dl]_{1 \ox g} \ar[dr]^{f \ox 1} & \\
A \ox B \ar[dr]_{\psi_0} & & A' \ox B' \ar[ld]^{\psi'_0} \\
& D &}
~~~~~\xymatrix{
& B' \ox A \ar[dl]_{g \ox 1} \ar[dr]^{1 \ox f} & \\
B \ox A \ar[dr]_{\psi_1} & & B' \ox A' \ar[ld]^{\psi'_1} \\
& D &}
\]
\item[Compositon] $(f,g)(f',g') := (ff', g'g)$. Composition is well-defined as:
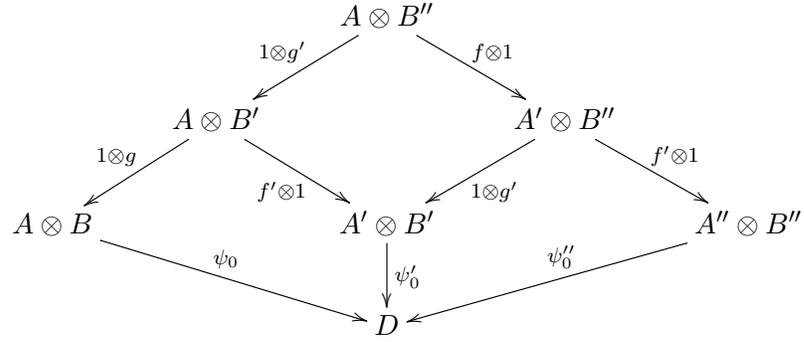
\[
\xymatrix{
&& A \ox B'' \ar[dl]_{1 \ox g'} \ar[dr]^{f \ox 1}   && \\
& A \ox B' \ar[dr]_{f \ox 1}  \ar[dl]_{1 \ox g}  & & A' \ox B'' \ar[dl]^{1 \ox g'}  \ar[dr]^{f' \ox 1} \\
A \ox B \ar[drr]^{\psi_0} & & A' \ox B' \ar[d]^{\psi'_0} && A'' \ox B'' \ar[lld]_{\psi_0''} \\
& & D & & 
}
\]
and similarly for the reverse Chu-maps $\psi_1$, $\psi'_1$ and $\psi_1''$.
The {\bf identity maps} are $(1_A,1_B): (A, B, \psi_0,\psi_1) \to (A, B, \psi_0,\psi_1)$ as expected.

\item[Tensor product $\ox$] $(A, B, \psi_0, \psi_1) \ox (A', B', \psi_0', \psi_1') := (A \ox A', E, \gamma_0, \gamma_1)$, where $E$ is the pullback in the following diagram:
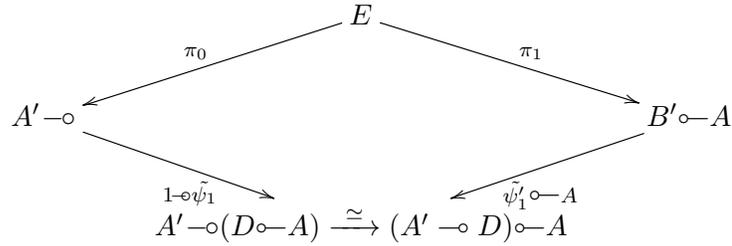
\[
\xymatrix{
& E \ar[rd]^{\pi_1} \ar[ld]_{\pi_0} & \\ %Ask Robin how to add a pullback corner
A' \multimap B \ar[rd]_{1 \multimap \tilde{\psi_1}} & & B' \poppilol A \ar[ld]^{\tilde{\psi_1'} \poppilol A}  \\
& A' \multimap (D \poppilol A) \to^{\simeq} (A' \multimap D) \poppilol A &
}
\]
with
\[
\infer{B \ox A \to^{\psi_1} B}{B \to^{\tilde{\psi_1}} D \poppilol A} ~~~~~~~~~~ \infer{A' \ox B' \to D}{B' \to^{\tilde{\psi_1'}} A' \multimap D}
\]
and, 
\[
\gamma_0 := (A \ox A') \ox E \to^{1 \ox \pi_0} (A \ox A') \ox (A' \multimap B) \to^{a_\ox} A \ox ( A' \ox A' \multimap B) \to^{1 \ox eval_{\multimap}} A \ox B \to^{\psi_0'} D
\]
\[
\gamma_1 := E \ox (A \ox A')  \to^{\pi_1 \ox 1} (B' \poppilol A) \ox (A \ox A')  \to^{a_\ox^{-1}} (B' \poppilol A \ox A) \ox A' \to^{eval_{\poppilol} \ox 1} B' \ox A' \to^{\psi_1'} D
\]
The tensor unit is $(I, D, u_\ox^l, u_\ox^r)$.
\end{description}

It is standard that ${\sf Chu}_\X(D)$  is a (non-commutative) $*$-autonomous category. Furthermore, it is cyclic because 
$$~^{*}(A,B,\psi_0,\psi_1) = (A,B,\psi_0,\psi_1)^{*} = (B,A,\psi_1,\psi_0).$$ 
In addition, ${\sf Chu}_\X(D)$  is conjugative with 
$$\overline{(A,B,\psi_0,\psi_1)} := (\overline{A},\overline{B},\chi \overline{\psi_1}d,\chi \overline{\psi_0}d)$$
and $\overline{(f,g)} = (\overline{f},\overline{g})$.
Finally, being conjugative cyclic $*$-autonomous implies that one has a dagger!

In the case that $\X$ is a symmetric monoidal closed category we may recapture the usual Chu construction \cite{Barr06}, which we denote ${\sf Chus}_\X(D)$.  Consider the full subcategory 
of Chu-objects with special Chu-cells of the form $(A,B, \psi,c_\otimes \psi)$ in which the symmetry map is used to obtain the second cell, this gives an inclusion 
${\sf Chus}_\X(D) \to {\sf Chu}_\X(D)$.  

We observe that $\X$ is symmetric conjugative when this subcategory is closed to the conjugation:

\begin{lem}
If $\X$  is a conjugative symmetric monoidal closed category  and $d: \overline{D} \to D$ is an involutive object, then ${\sf Chus}_\X(D)$ is a conjugative symmetric $*$-autonomous  category.
\end{lem}
\begin{proof}
It suffices to observe that the Chu-cells of $\overline{(A,B,\psi,c_\otimes\psi)}$ have the right form.  Using the coherence of the 
involution with symmetry, the first Chu-cell of this object has  $\chi \overline{c_\otimes \psi}d = c_\ox \chi \overline{\psi}d$
which is exactly the symmetry map applied to the second Chu-cell of the object as desired.
\end{proof}

To obtain an isomix category one can choose D = I. ${\sf Chus}_\X(I)$ is an isomix category because the unit for tensor and par are the same (namely $\top = \bot = (I,I, u_\ox^l = u_\ox^r)$).  The tensor unit is always a conjugative object since $(\chi^{\!\!\!\circ})^{-1}: \overline{I} \to I$; therefore, this is immediately a conjugative symmetric $*$-autonomous category.  Composing the conjugation with the dualizing functor gives us a dagger.

%%%%%%%%%%%%%%%%%%%%%%%%%%%%%%%%%%%%%%%%%%%%%%%%%%%%%%
\subsubsection{Category of Hopf Modules of a $*$-autonomous category}
\label{Sec: HModx}
In the previous example, starting from a conjugative closed monoidal category with pullbacks, we built a $\dagger$-$*$-autonomous category using the Chu construction. In this example\footnote{We thank J-S. P. Lemay for bringing our attention to this example.}, we start with any symmetric $*$-autonomous category, $\X$, and build the category of modules over a Hopf Algebra which is in turn a $\dagger$-$*$-autonomous category.  

First of all, it has been already proven in \cite{PaS09}, that the category of Hopf modules over a $\ox$-Hopf algebra in any symmetric $*$-autonomous category is also a $*$-autonomous category. Then we note that, whenever the Hopf Algebra is cocommutative, the resulting $*$-autonomous category has a conjugation functor. One can construct the dagger functor by composing the conjugation functor and dualizing functor as in Theorem \ref{Theorem: conjugation+dualizing}. We establish some basic definitions before describing the category of modules over a Hopf Algerba, ${\mbox{\bf H-Mod}}_\X$.

\begin{defi}
A {\bf bialgebra}  in a symmetric monoidal category is a 4-tuple 
$$(\nabla: B \ox B \to B, e:  I \to B, \Delta : B \to B \ox B, u: B \to I)$$
such that $(A, \nabla,e)$ is a monoid and $(A, \Delta, u)$  is a comonoid and $\nabla$ and $e$ are coalgebra homomorphisms with respect to the comultiplication and the counit.
\end{defi}

Note that instead of requiring that  $\nabla$ and $e$ are coalgebra homomorphisms, one could equivalently require $\Delta$ and $u$ are algebra homomorphims with respect to the multiplication and the unit.

 The components of a bialgebra are graphically depicted as follows:
 
 \[
 \begin{tikzpicture}
	\begin{pgfonlayer}{nodelayer}
		\node [style=none] (0) at (-3, 1) {};
		\node [style=none] (1) at (-2.5, 1) {};
		\node [style=none] (2) at (-2.75, 0.75) {};
		\node [style=none] (3) at (-3.25, 1.5) {};
		\node [style=none] (4) at (-2.25, 1.5) {};
		\node [style=none] (5) at (-2.75, 0.25) {};
	\end{pgfonlayer}
	\begin{pgfonlayer}{edgelayer}
		\draw (0.center) to (1.center);
		\draw (1.center) to (2.center);
		\draw (0.center) to (2.center);
		\draw [in=-90, out=30, looseness=1.00] (1.center) to (4.center);
		\draw [in=-90, out=150, looseness=1.00] (0.center) to (3.center);
		\draw (2.center) to (5.center);
	\end{pgfonlayer}
\end{tikzpicture} : A \ox A \to A ~~~~~~~ \begin{tikzpicture}
	\begin{pgfonlayer}{nodelayer}
		\node [style=none] (0) at (-3, 0.75) {};
		\node [style=none] (1) at (-2.5, 0.75) {};
		\node [style=none] (2) at (-2.75, 1) {};
		\node [style=none] (3) at (-3.25, 0.25) {};
		\node [style=none] (4) at (-2.25, 0.25) {};
		\node [style=none] (5) at (-2.75, 1.5) {};
	\end{pgfonlayer}
	\begin{pgfonlayer}{edgelayer}
		\draw (0.center) to (1.center);
		\draw (1.center) to (2.center);
		\draw (0.center) to (2.center);
		\draw [in=90, out=-30, looseness=1.00] (1.center) to (4.center);
		\draw [in=90, out=-150, looseness=1.00] (0.center) to (3.center);
		\draw (2.center) to (5.center);
	\end{pgfonlayer}
\end{tikzpicture} : A \to A \ox A ~~~~~~~~ \begin{tikzpicture}
	\begin{pgfonlayer}{nodelayer}
		\node [style=none] (0) at (-3, -0) {};
		\node [style=none] (1) at (-2.5, -0) {};
		\node [style=none] (2) at (-2.75, 0.25) {};
		\node [style=none] (3) at (-2.75, 1.5) {};
	\end{pgfonlayer}
	\begin{pgfonlayer}{edgelayer}
		\draw (0.center) to (1.center);
		\draw (1.center) to (2.center);
		\draw (0.center) to (2.center);
		\draw (2.center) to (3.center);
	\end{pgfonlayer}
\end{tikzpicture}: A \to I ~~~~~~~~~~ \begin{tikzpicture}
	\begin{pgfonlayer}{nodelayer}
		\node [style=none] (0) at (-3, 1.5) {};
		\node [style=none] (1) at (-2.5, 1.5) {};
		\node [style=none] (2) at (-2.75, 1.25) {};
		\node [style=none] (3) at (-2.75, 0) {};
	\end{pgfonlayer}
	\begin{pgfonlayer}{edgelayer}
		\draw (0.center) to (1.center);
		\draw (1.center) to (2.center);
		\draw (0.center) to (2.center);
		\draw (2.center) to (3.center);
	\end{pgfonlayer}
\end{tikzpicture} : I \to A
 \]
This gives a succinct graphical depiction of the coalgebra homomorphism laws; namely:
\[
\begin{tikzpicture}
	\begin{pgfonlayer}{nodelayer}
		\node [style=none] (0) at (-2, 2) {};
		\node [style=none] (1) at (-2.25, 1.75) {};
		\node [style=none] (2) at (-1.75, 1.75) {};
		\node [style=none] (3) at (-2, 2.75) {};
		\node [style=none] (4) at (-0.5, 2.75) {};
		\node [style=none] (5) at (-0.75, 1.75) {};
		\node [style=none] (6) at (-0.5, 2) {};
		\node [style=none] (7) at (-0.25, 1.75) {};
		\node [style=none] (8) at (-1.75, 0.25) {};
		\node [style=none] (9) at (-0.5, -0.75) {};
		\node [style=none] (10) at (-2.25, 0.25) {};
		\node [style=none] (11) at (-2, -0) {};
		\node [style=none] (12) at (-2, -0.75) {};
		\node [style=none] (13) at (-0.5, -0) {};
		\node [style=none] (14) at (-0.25, 0.25) {};
		\node [style=none] (15) at (-0.75, 0.25) {};
	\end{pgfonlayer}
	\begin{pgfonlayer}{edgelayer}
		\draw (0.center) to (1.center);
		\draw (1.center) to (2.center);
		\draw (2.center) to (0.center);
		\draw (3.center) to (0.center);
		\draw (6.center) to (5.center);
		\draw (5.center) to (7.center);
		\draw (7.center) to (6.center);
		\draw (4.center) to (6.center);
		\draw (11.center) to (10.center);
		\draw (10.center) to (8.center);
		\draw (8.center) to (11.center);
		\draw (12.center) to (11.center);
		\draw (13.center) to (15.center);
		\draw (15.center) to (14.center);
		\draw (14.center) to (13.center);
		\draw (9.center) to (13.center);
		\draw [in=15, out=-165, looseness=1.00] (5.center) to (8.center);
		\draw [bend right, looseness=1.25] (1.center) to (10.center);
		\draw [in=150, out=-30, looseness=1.00] (2.center) to (15.center);
		\draw [bend left, looseness=1.00] (7.center) to (14.center);
	\end{pgfonlayer}
\end{tikzpicture} = \begin{tikzpicture}
	\begin{pgfonlayer}{nodelayer}
		\node [style=none] (0) at (-2, -1.5) {};
		\node [style=none] (1) at (-2.25, -1.75) {};
		\node [style=none] (2) at (-1.75, -1.75) {};
		\node [style=none] (3) at (-2, -0.75) {};
		\node [style=none] (4) at (-1.75, 0.25) {};
		\node [style=none] (5) at (-2.25, 0.25) {};
		\node [style=none] (6) at (-2, -0) {};
		\node [style=none] (7) at (-2, -0.75) {};
		\node [style=none] (8) at (-2.75, 1) {};
		\node [style=none] (9) at (-1.25, 1) {};
		\node [style=none] (10) at (-2.75, -2.5) {};
		\node [style=none] (11) at (-1.25, -2.5) {};
	\end{pgfonlayer}
	\begin{pgfonlayer}{edgelayer}
		\draw (0.center) to (1.center);
		\draw (1.center) to (2.center);
		\draw (2.center) to (0.center);
		\draw (3.center) to (0.center);
		\draw (6.center) to (5.center);
		\draw (5.center) to (4.center);
		\draw (4.center) to (6.center);
		\draw (7.center) to (6.center);
		\draw [bend right=45, looseness=0.75] (4.center) to (9.center);
		\draw [bend right, looseness=1.25] (8.center) to (5.center);
		\draw [bend right, looseness=0.75] (1.center) to (10.center);
		\draw [bend left, looseness=1.00] (2.center) to (11.center);
	\end{pgfonlayer}
\end{tikzpicture} ~~~~~~~ \begin{tikzpicture} %bialg-2a
	\begin{pgfonlayer}{nodelayer}
		\node [style=none] (0) at (-1, 0.75) {};
		\node [style=none] (1) at (-1.25, 1) {};
		\node [style=none] (2) at (-0.75, 1) {};
		\node [style=none] (3) at (-1, -0) {};
		\node [style=none] (4) at (-1.25, -0.25) {};
		\node [style=none] (5) at (-0.75, -0.25) {};
		\node [style=none] (6) at (-1.75, -1) {};
		\node [style=none] (7) at (-0.25, -1) {};
	\end{pgfonlayer}
	\begin{pgfonlayer}{edgelayer}
		\draw (1.center) to (2.center);
		\draw (2.center) to (0.center);
		\draw (0.center) to (1.center);
		\draw (4.center) to (5.center);
		\draw (5.center) to (3.center);
		\draw (3.center) to (4.center);
		\draw (0.center) to (3.center);
		\draw [in=71, out=-135, looseness=1.00] (4.center) to (6.center);
		\draw [bend left, looseness=0.75] (5.center) to (7.center);
	\end{pgfonlayer}
\end{tikzpicture} =   \begin{tikzpicture}
	\begin{pgfonlayer}{nodelayer}
		\node [style=none] (0) at (-0.75, 0.75) {};
		\node [style=none] (1) at (-1, 1) {};
		\node [style=none] (2) at (-0.5, 1) {};
		\node [style=none] (3) at (-0.75, -1) {};
		\node [style=none] (4) at (0, 0.75) {};
		\node [style=none] (5) at (0.25, 1) {};
		\node [style=none] (6) at (0, -1) {};
		\node [style=none] (7) at (-0.25, 1) {};
	\end{pgfonlayer}
	\begin{pgfonlayer}{edgelayer}
		\draw (1.center) to (2.center);
		\draw (2.center) to (0.center);
		\draw (0.center) to (1.center);
		\draw (0.center) to (3.center);
		\draw (7.center) to (5.center);
		\draw (5.center) to (4.center);
		\draw (4.center) to (7.center);
		\draw (4.center) to (6.center);
	\end{pgfonlayer}
\end{tikzpicture} ~~~~~~~~ \begin{tikzpicture} %bialg2a
	\begin{pgfonlayer}{nodelayer}
		\node [style=none] (0) at (-1, -0.75) {};
		\node [style=none] (1) at (-1.25, -1) {};
		\node [style=none] (2) at (-0.75, -1) {};
		\node [style=none] (3) at (-1, 0) {};
		\node [style=none] (4) at (-1.25, 0.25) {};
		\node [style=none] (5) at (-0.75, 0.25) {};
		\node [style=none] (6) at (-1.75, 1) {};
		\node [style=none] (7) at (-0.25, 1) {};
	\end{pgfonlayer}
	\begin{pgfonlayer}{edgelayer}
		\draw (1.center) to (2.center);
		\draw (2.center) to (0.center);
		\draw (0.center) to (1.center);
		\draw (4.center) to (5.center);
		\draw (5.center) to (3.center);
		\draw (3.center) to (4.center);
		\draw (0.center) to (3.center);
		\draw [in=-71, out=135, looseness=1.00] (4.center) to (6.center);
		\draw [bend right, looseness=0.75] (5.center) to (7.center);
	\end{pgfonlayer}
\end{tikzpicture} = \begin{tikzpicture}
	\begin{pgfonlayer}{nodelayer}
		\node [style=none] (0) at (-1, -0.75) {};
		\node [style=none] (1) at (-1.25, -1) {};
		\node [style=none] (2) at (-0.75, -1) {};
		\node [style=none] (3) at (-1, 1) {};
		\node [style=none] (4) at (0, -0.75) {};
		\node [style=none] (5) at (0.25, -1) {};
		\node [style=none] (6) at (0, 1) {};
		\node [style=none] (7) at (-0.25, -1) {};
	\end{pgfonlayer}
	\begin{pgfonlayer}{edgelayer}
		\draw (1.center) to (2.center);
		\draw (2.center) to (0.center);
		\draw (0.center) to (1.center);
		\draw (0.center) to (3.center);
		\draw (7.center) to (5.center);
		\draw (5.center) to (4.center);
		\draw (4.center) to (7.center);
		\draw (4.center) to (6.center);
	\end{pgfonlayer}
\end{tikzpicture} ~~~~~~~~ \begin{tikzpicture}
	\begin{pgfonlayer}{nodelayer}
		\node [style=none] (0) at (-1, -0.75) {};
		\node [style=none] (1) at (-1.25, -1) {};
		\node [style=none] (2) at (-0.75, -1) {};
		\node [style=none] (3) at (-1, 0.75) {};
		\node [style=none] (4) at (-1.25, 1) {};
		\node [style=none] (5) at (-0.75, 1) {};
	\end{pgfonlayer}
	\begin{pgfonlayer}{edgelayer}
		\draw (1.center) to (2.center);
		\draw (2.center) to (0.center);
		\draw (0.center) to (1.center);
		\draw (0.center) to (3.center);
		\draw (4.center) to (3.center);
		\draw (3.center) to (5.center);
		\draw (5.center) to (4.center);
	\end{pgfonlayer}
\end{tikzpicture} = I
\]

\begin{defi} An {\bf antipode} for a bialgebra $(B, \nabla, \bialgunitmap{0.8}, \Delta, \bialgcounitmap{0.8})$ is an endomorphism $s: B \to B$ such that 
\[
\begin{tikzpicture} %hopf-1
	\begin{pgfonlayer}{nodelayer}
		\node [style=none] (0) at (-2, 2.25) {};
		\node [style=none] (1) at (-2.25, 2) {};
		\node [style=none] (2) at (-1.75, 2) {};
		\node [style=none] (3) at (-2, 2.75) {};
		\node [style=none] (4) at (-1.75, 0.5) {};
		\node [style=none] (5) at (-2.25, 0.5) {};
		\node [style=none] (6) at (-2, 0.25) {};
		\node [style=none] (7) at (-2, -0.5) {};
		\node [style=circle, scale=1.5] (8) at (-1.5, 1.25) {};
		\node [style=none] (9) at (-1.5, 1.25) {$s$};
	\end{pgfonlayer}
	\begin{pgfonlayer}{edgelayer}
		\draw (0.center) to (1.center);
		\draw (1.center) to (2.center);
		\draw (2.center) to (0.center);
		\draw (3.center) to (0.center);
		\draw (6.center) to (5.center);
		\draw (5.center) to (4.center);
		\draw (4.center) to (6.center);
		\draw (7.center) to (6.center);
		\draw [bend left, looseness=1.00] (2.center) to (8);
		\draw [bend left, looseness=0.75] (8) to (4.center);
		\draw [bend right=45, looseness=1.00] (1.center) to (5.center);
	\end{pgfonlayer}
\end{tikzpicture} = \begin{tikzpicture}
	\begin{pgfonlayer}{nodelayer}
		\node [style=none] (0) at (-1.75, 2.25) {};
		\node [style=none] (1) at (-1.5, 2) {};
		\node [style=none] (2) at (-2, 2) {};
		\node [style=none] (3) at (-1.75, 2.75) {};
		\node [style=none] (4) at (-2, 0.5) {};
		\node [style=none] (5) at (-1.5, 0.5) {};
		\node [style=none] (6) at (-1.75, 0.25) {};
		\node [style=none] (7) at (-1.75, -0.5) {};
		\node [style=circle, scale=1.5] (8) at (-2.25, 1.25) {};
		\node [style=none] (9) at (-2.25, 1.25) {$s$};
	\end{pgfonlayer}
	\begin{pgfonlayer}{edgelayer}
		\draw (0.center) to (1.center);
		\draw (1.center) to (2.center);
		\draw (2.center) to (0.center);
		\draw (3.center) to (0.center);
		\draw (6.center) to (5.center);
		\draw (5.center) to (4.center);
		\draw (4.center) to (6.center);
		\draw (7.center) to (6.center);
		\draw [bend right, looseness=1.00] (2.center) to (8);
		\draw [bend right, looseness=0.75] (8) to (4.center);
		\draw [bend left=45, looseness=1.00] (1.center) to (5.center);
	\end{pgfonlayer}
\end{tikzpicture} = \begin{tikzpicture}
	\begin{pgfonlayer}{nodelayer}
		\node [style=none] (0) at (1, 1.5) {};
		\node [style=none] (1) at (1, 0.5) {};
		\node [style=none] (2) at (1, -0.5) {};
		\node [style=none] (3) at (1, 2.75) {};
		\node [style=none] (4) at (0.75, 0.75) {};
		\node [style=none] (5) at (1.25, 0.75) {};
		\node [style=none] (6) at (1, 0.5) {};
		\node [style=none] (7) at (1.25, 1.25) {};
		\node [style=none] (8) at (0.75, 1.25) {};
		\node [style=none] (9) at (1, 1.5) {};
	\end{pgfonlayer}
	\begin{pgfonlayer}{edgelayer}
		\draw (3.center) to (0.center);
		\draw (1.center) to (2.center);
		\draw (4.center) to (5.center);
		\draw (5.center) to (6.center);
		\draw (6.center) to (4.center);
		\draw (8.center) to (7.center);
		\draw (7.center) to (9.center);
		\draw (9.center) to (8.center);
	\end{pgfonlayer}
\end{tikzpicture}
\]

A {\bf Hopf algebra} is a bialgebra with an antipode. An {\bf involutive Hopf algebra} is a hopf algebra where the antipode is self-inverse.
\end{defi}

A standard example of a Hopf algebra is a group algebra over a field:  for all group elements $g$, $\nabla : g \mapsto g \ox g$, $\bialgcounitmap{0.8}: g \mapsto 1$, $\Delta: g \ox h \mapsto gh$ and $s: g \mapsto g^{-1}$.

\begin{lem}
\label{Lemma: involutive s}
Suppose $\X$ is a symmetric monoidal category, then:

\begin{enumerate}[(i)]

\item \cite[Theroem 3.5]{Blu96} If $H$ is a commutative or a cocommutative Hopf Algebra in $\X$, then $s^2 = 1$ where $s$ is the antipode: 
so it is an involutive Hopf algebra.

\item \cite[Lemma 2.11]{Lem19} If $H$ is a commutative Hopf Algebra, then $s$ is a monoid homomorphism. If $H$ is a cocommutative 
Hopf Algebra, then $s$ is  a comonoid homomorphism.

\end{enumerate}

\end{lem}

\begin{defi}
A {\bf left module} for a bialgebra $(B, \nabla, u, \Delta, e)$ is a tuple $(M, a_M^l:B \ox M \to M)$ such that $a_M^l$ is a $B$-action i.e., the following diagram commutes:

\[
\xymatrix{
M \ar[r]^{u_\ox^l} \ar@{=}[dr] & \top \ox M \ar[d]^{\bialgunitmap{0.8}} \\ 
& M
}
\]
\end{defi}

We graphically depict $a_m^l$ as follows:
\[
\begin{tikzpicture}[xscale=-1]
	\begin{pgfonlayer}{nodelayer}
		\node [style=none] (0) at (-2.75, 1.25) {};
		\node [style=none] (1) at (-2.5, 1.25) {};
		\node [style=none] (2) at (-2.75, 1) {};
		\node [style=none] (3) at (-2.75, 0.25) {};
		\node [style=none] (4) at (-2, 2) {};
		\node [style=none] (5) at (-2.75, 2) {};
	\end{pgfonlayer}
	\begin{pgfonlayer}{edgelayer}
		\draw (0.center) to (1.center);
		\draw (1.center) to (2.center);
		\draw (0.center) to (2.center);
		\draw (2.center) to (3.center);
		\draw (5.center) to (0.center);
		\draw [in=-105, out=30, looseness=1.25] (1.center) to (4.center);
	\end{pgfonlayer}
\end{tikzpicture}: B \ox M \to M
\]

giving the graphical presentation of the module laws:

\[
\begin{tikzpicture}%leftaction-rule1
	\begin{pgfonlayer}{nodelayer}
		\node [style=none] (0) at (-1.5, 2.25) {};
		\node [style=none] (1) at (-0.5, 2.25) {};
		\node [style=none] (2) at (0.25, 2.25) {};
		\node [style=none] (3) at (-1, 0.5) {};
		\node [style=none] (4) at (0, -0.25) {};
		\node [style=none] (5) at (0.25, -0.25) {};
		\node [style=none] (6) at (0.25, -0.5) {};
		\node [style=none] (7) at (0.25, -1) {};
		\node [style=none] (8) at (-1.25, 0.75) {};
		\node [style=none] (9) at (-0.75, 0.75) {};
		\node [style=none] (10) at (-1, 0.5) {};
	\end{pgfonlayer}
	\begin{pgfonlayer}{edgelayer}
		\draw (2.center) to (5.center);
		\draw (6.center) to (5.center);
		\draw (5.center) to (4.center);
		\draw (4.center) to (6.center);
		\draw [in=-90, out=150, looseness=1.25] (4.center) to (3.center);
		\draw (6.center) to (7.center);
		\draw (8.center) to (9.center);
		\draw (9.center) to (10.center);
		\draw (10.center) to (8.center);
		\draw [bend right=15, looseness=1.00] (0.center) to (8.center);
		\draw [bend right=15, looseness=1.00] (9.center) to (1.center);
	\end{pgfonlayer}
\end{tikzpicture} = \begin{tikzpicture}
	\begin{pgfonlayer}{nodelayer}
		\node [style=none] (0) at (-1.75, 2) {};
		\node [style=none] (1) at (-1.5, 2) {};
		\node [style=none] (2) at (-1.5, 1.75) {};
		\node [style=none] (3) at (-1.5, 3) {};
		\node [style=none] (4) at (-2.25, 3) {};
		\node [style=none] (5) at (-0.75, -0) {};
		\node [style=none] (6) at (-0.75, 3) {};
		\node [style=none] (7) at (-0.75, 0.25) {};
		\node [style=none] (8) at (-1.5, 1.25) {};
		\node [style=none] (9) at (-1, 0.25) {};
		\node [style=none] (10) at (-0.75, -0.5) {};
	\end{pgfonlayer}
	\begin{pgfonlayer}{edgelayer}
		\draw (3.center) to (1.center);
		\draw (2.center) to (1.center);
		\draw (1.center) to (0.center);
		\draw (0.center) to (2.center);
		\draw [in=-90, out=150, looseness=1.25] (0.center) to (4.center);
		\draw (6.center) to (7.center);
		\draw (5.center) to (7.center);
		\draw (7.center) to (9.center);
		\draw (9.center) to (5.center);
		\draw [in=-90, out=150, looseness=1.25] (9.center) to (8.center);
		\draw (2.center) to (8.center);
		\draw (5.center) to (10.center);
	\end{pgfonlayer}
\end{tikzpicture} ~~~~\text{ and }~~ \begin{tikzpicture} %leftactionrule-1b
	\begin{pgfonlayer}{nodelayer}
		\node [style=none] (0) at (-0.75, 1) {};
		\node [style=none] (1) at (-0.75, 3) {};
		\node [style=none] (2) at (-0.75, 1.25) {};
		\node [style=none] (3) at (-1, 1.25) {};
		\node [style=none] (4) at (-0.75, -0) {};
		\node [style=none] (5) at (-1.75, 3) {};
		\node [style=none] (6) at (-1.25, 3) {};
		\node [style=none] (7) at (-1.5, 2.75) {};
	\end{pgfonlayer}
	\begin{pgfonlayer}{edgelayer}
		\draw (1.center) to (2.center);
		\draw (0.center) to (2.center);
		\draw (2.center) to (3.center);
		\draw (3.center) to (0.center);
		\draw (0.center) to (4.center);
		\draw (5.center) to (6.center);
		\draw (6.center) to (7.center);
		\draw (7.center) to (5.center);
		\draw [bend right, looseness=1.00] (7.center) to (3.center);
	\end{pgfonlayer}
\end{tikzpicture} = \begin{tikzpicture}
	\begin{pgfonlayer}{nodelayer}
		\node [style=none] (0) at (-0.75, 1.75) {};
		\node [style=none] (1) at (-0.75, 3.25) {};
		\node [style=none] (2) at (-0.75, 1.75) {};
		\node [style=none] (3) at (-0.75, 0.25) {};
	\end{pgfonlayer}
	\begin{pgfonlayer}{edgelayer}
		\draw (1.center) to (2.center);
		\draw (0.center) to (2.center);
		\draw (0.center) to (3.center);
	\end{pgfonlayer}
\end{tikzpicture}
\] 

\begin{defi}
Let $\X$ be a $*$-autonomous category and $H$ be a Hopf $\ox$-algebra in $\X$. The category of left H-modules in $\X$, ${\mbox{\bf H-Mod}}_\X$ has:
\begin{description}
\item[Objects]  Left $H$-modules $(A, a_A^l:H \ox A\to A)$:
\item[Arrows] A module homomorphism $(A, a_A^L:H \ox A\to A) \to^{f}(B, a_B^L:H \ox B\to B)$ is a map $A \to^{f} B$ such that the following diagram commutes:
\[
\xymatrix{
H \ox A \ar[r]^{ a_A^L} \ar[d]_{1 \ox f} & A \ar[d]^{f} \\
H \ox B \ar[r]_{ a_B^L} & B
}
\]

This is graphically depicted as follows:
\[
 \begin{tikzpicture} %act2
	\begin{pgfonlayer}{nodelayer}
		\node [style=none] (0) at (0, 1) {};
		\node [style=none] (1) at (0, 0.75) {};
		\node [style=none] (2) at (-0.25, 1) {};
		\node [style=none] (3) at (0, 2) {};
		\node [style=none] (4) at (0, -0.5) {};
		\node [style=none] (5) at (-0.75, 2) {};
		\node [style=none] (6) at (-0.75, 1.5) {};
		\node [style=circle, scale=1.5] (7) at (0, 0.25) {};
		\node [style=none] (8) at (0, 0.25) {$f$};
	\end{pgfonlayer}
	\begin{pgfonlayer}{edgelayer}
		\draw (2.center) to (0.center);
		\draw (0.center) to (1.center);
		\draw (1.center) to (2.center);
		\draw [bend left, looseness=1.00] (2.center) to (6.center);
		\draw (6.center) to (5.center);
		\draw (3.center) to (0.center);
		\draw (1.center) to (7);
		\draw (7) to (4.center);
	\end{pgfonlayer}
\end{tikzpicture} =
\begin{tikzpicture} %act1
	\begin{pgfonlayer}{nodelayer}
		\node [style=none] (0) at (0, 0.5) {};
		\node [style=none] (1) at (0, 0.25) {};
		\node [style=none] (2) at (-0.25, 0.5) {};
		\node [style=none] (3) at (0, 2) {};
		\node [style=none] (4) at (0, -0.25) {};
		\node [style=none] (5) at (-0.75, 2) {};
		\node [style=none] (6) at (-0.75, 1.5) {};
		\node [style=circle, scale=1.5] (7) at (0, 1.25) {};
		\node [style=none] (8) at (0, 1.25) {$f$};
	\end{pgfonlayer}
	\begin{pgfonlayer}{edgelayer}
		\draw (2.center) to (0.center);
		\draw (0.center) to (1.center);
		\draw (1.center) to (2.center);
		\draw (1.center) to (4.center);
		\draw (3.center) to (7);
		\draw (7) to (0.center);
		\draw [bend left, looseness=1.00] (2.center) to (6.center);
		\draw (6.center) to (5.center);
		\draw[fill=black] (0.center) -- (1.center) -- (2.center) -- (0.center);
	\end{pgfonlayer}
\end{tikzpicture} \]

\end{description}
\end{defi}

Observe that any left action is indeed a module homomorphism.

\begin{thmC} \cite{PaS09}
Let $\X$ be symmetric $*$-autonomous category and $H$ be a $\ox$-Hopf Algebra in $\X$ with bijective antipode ($s^2 = 1$). Then, ${\mbox{\bf H-Mod}}_\X$ is a $*$-autonomous category. If the Hopf Algebra, $H$, is cocommutative, then ${\mbox{\bf H-Mod}}_\X$ is a symmetric $*$-autonomous category.
\end{thmC}
\begin{proof} (Sketch)
The monoidal product $\ox$ for ${\mbox{\bf H-Mod}}_\X$ is defined as follows:
\[ (A, \leftaction{0.4}{white}) \ox (B, \leftaction{0.4}{black}) := (A \ox B, \leftaction{0.4}{gray} ) 
\text{ where, } \leftaction{0.8}{gray} := \begin{tikzpicture}
	\begin{pgfonlayer}{nodelayer}
		\node [style=none] (0) at (-0.25, -0.5) {};
		\node [style=none] (1) at (-0.25, -0.75) {};
		\node [style=none] (2) at (-0.5, -0.5) {};
		\node [style=none] (3) at (-0.25, 0.25) {};
		\node [style=none] (4) at (-1, 0.75) {};
		\node [style=none] (5) at (-1, -0) {};
		\node [style=none] (6) at (-0.25, -1.25) {};
		\node [style=none] (7) at (1, -0.75) {};
		\node [style=none] (8) at (1, 0.25) {};
		\node [style=none] (9) at (1, -1.25) {};
		\node [style=none] (10) at (0.25, -0) {};
		\node [style=none] (11) at (1, -0.5) {};
		\node [style=none] (12) at (-0.5, 0.75) {};
		\node [style=none] (13) at (0.75, -0.5) {};
		\node [style=none] (14) at (0.25, 1.75) {};
		\node [style=none] (15) at (1, 1.75) {};
		\node [style=none] (16) at (-0.75, 1.75) {};
		\node [style=none] (17) at (-0.75, 1) {};
		\node [style=none] (18) at (-1, 0.75) {};
		\node [style=none] (19) at (-0.5, 0.75) {};
	\end{pgfonlayer}
	\begin{pgfonlayer}{edgelayer}
		\draw (2.center) to (0.center);
		\draw (0.center) to (1.center);
		\draw (1.center) to (2.center);
		\draw [bend left, looseness=1.00] (2.center) to (5.center);
		\draw (5.center) to (4.center);
		\draw (3.center) to (0.center);
		\draw (1.center) to (6.center);
		\draw (13.center) to (11.center);
		\draw (11.center) to (7.center);
		\draw (7.center) to (13.center);
		\draw [bend left, looseness=1.00] (13.center) to (10.center);
		\draw [in=-45, out=90, looseness=1.00] (10.center) to (12.center);
		\draw (8.center) to (11.center);
		\draw (7.center) to (9.center);
		\draw [in=-90, out=90, looseness=1.00] (3.center) to (14.center);
		\draw (15.center) to (8.center);
		\draw (18.center) to (19.center);
		\draw (19.center) to (17.center);
		\draw (17.center) to (18.center);
		\draw (16.center) to (17.center);
		\draw[fill=black] (0.center) --  (1.center) --  (2.center) --  (0.center);
	\end{pgfonlayer}
\end{tikzpicture} \]

 The unit of $\ox$ is given by $(\top, H \ox \top \to^{u_\ox^R} H \to^{e} \top)$, the left action is drawn as 
 $
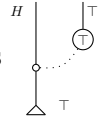
\begin{tikzpicture} %act6
	\begin{pgfonlayer}{nodelayer}
		\node [style=circle, scale=0.5] (0) at (-1, 0.25) {};
		\node [style=circle, scale=1.5] (1) at (0.25, 1) {};
		\node [style=none] (2) at (0.25, 1) {$\top$};
		\node [style=none] (3) at (-1, 2) {};
		\node [style=none] (4) at (0.25, 2) {};
		\node [style=none] (5) at (-0.25, -0.75) {$\top$};
		\node [style=none] (6) at (-1.5, 1.75) {$H$};
		\node [style=none] (7) at (0.5, 1.75) {$\top$};
		\node [style=none] (8) at (-1.25, -1) {};
		\node [style=none] (9) at (-0.75, -1) {};
		\node [style=none] (10) at (-1, -0.75) {};
		\node [style=none] (11) at (-1, -1) {};
	\end{pgfonlayer}
	\begin{pgfonlayer}{edgelayer}
		\draw (3.center) to (0);
		\draw [bend left, looseness=1.25, dotted] (1) to (0);
		\draw (4.center) to (1);
		\draw (8.center) to (9.center);
		\draw (9.center) to (10.center);
		\draw (8.center) to (10.center);
		\draw (0) to (10.center);
	\end{pgfonlayer}
\end{tikzpicture}$

 The par product is defined as: $(A, \leftaction{0.4}{white}) \oa (B, \leftaction{0.4}{black}) := (A \oa B, \leftaction{0.4}{gray} ) $
where,
\begin{align*}
 \leftaction{0.8}{gray} &:= H \ox (A \oa B) \to^{\Delta \ox 1} (H \ox H) \ox (A \oa B) \to^{a_\ox} H \ox (H \ox (A \oa B)) \\
  &\to^{c_\ox} (H \ox (A \oa B)) \ox H \to^{\partial^L \ox 1} ((H \ox A) \oa B) \ox H \to^{\partial^R} (H \ox A) \oa (B \ox H) \\
  &\to^{1 \oa c_\ox} (H \ox A) \oa (H \ox B) \to^{\leftaction{0.4}{white} \oa \leftaction{0.4}{black}} A \oa B
\end{align*}

and  the unit of $\oa$ is
\[
\bot := ( \bot, H \ox \bot \to^{u \ox \bot} \top \ox \bot \to^{u_\ox} \bot)
\]

All the basic natural isomorphisms are inherited directly from $\X$ and they are module homomorphisms. Thus, {\bf HMod}$_\X$ is a LDC.

The dualizing functor $(\_)^*$ is given as follows:$(A, \leftaction{0.4}{white}: H \ox A \to A)^* := (A^*, \leftaction{0.4}{white}*: H \ox A^* \to A^*) \ \text{ where, }$ \[ \leftaction{0.4}{white}* := 
\begin{tikzpicture} %act-tensor5
	\begin{pgfonlayer}{nodelayer}
		\node [style=none] (0) at (1.5, -0.25) {};
		\node [style=none] (1) at (1.5, -1) {};
		\node [style=none] (2) at (1.5, -0.5) {};
		\node [style=none] (3) at (2.25, -1) {};
		\node [style=none] (4) at (0.75, 0.5) {};
		\node [style=none] (5) at (0.75, 0.25) {};
		\node [style=none] (6) at (1.25, -0.25) {};
		\node [style=none] (7) at (1.5, 0.25) {};
		\node [style=none] (8) at (-0.25, 0.25) {};
		\node [style=none] (9) at (2.25, 3) {};
		\node [style=none] (10) at (-0.25, -1.25) {};
		\node [style=none] (11) at (0.75, 3) {};
		\node [style=none] (12) at (1, 2.75) {$H$};
		\node [style=none] (13) at (2.5, 2.75) {$A^*$};
		\node [style=none] (14) at (-0.25, -3) {};
		\node [style=none] (15) at (0, -2.75) {$A^*$};
		\node [style=circle, scale=1.5] (16) at (0.75, 2) {};
		\node [style=none] (17) at (0.75, 2) {$s$};
	\end{pgfonlayer}
	\begin{pgfonlayer}{edgelayer}
		\draw (6.center) to (0.center);
		\draw (0.center) to (2.center);
		\draw (2.center) to (6.center);
		\draw [bend left, looseness=1.00] (6.center) to (5.center);
		\draw (5.center) to (4.center);
		\draw (2.center) to (1.center);
		\draw [bend left=90, looseness=2.00] (8.center) to (7.center);
		\draw [bend right=90, looseness=2.00] (1.center) to (3.center);
		\draw (9.center) to (3.center);
		\draw (8.center) to (10.center);
		\draw (11.center) to (16);
		\draw (16) to (4.center);
		\draw (14.center) to (10.center);
		\draw (7.center) to (0.center);
	\end{pgfonlayer}
\end{tikzpicture}\]
 Equationally, 

\begin{align*}
\leftaction{0.4}{white}* &:= H \ox A^* \to^{s \ox 1} H \ox A^* \to^{u_\ox^{-1} \ox 1} (H \ox \top) \ox A^* \to^{1 \ox \eta \ox 1} (H \ox (A^* \oa A)) \ox A^* \\
&\to^{c_\ox \ox 1}  ((A^* \oa A) \ox H) \ox A^* \to^{\partial \ox 1} (A^* \oa (A \ox H)) \ox A^* \to^{1 \ox c_\ox \ox 1} (A^* \oa (H \ox A)) \ox A^* \\
& \to^{(1 \oa \leftaction{0.3}{white}) \ox 1} (A^* \oa A) \ox A^* \to^{\partial} A^* \oa (A \ox A^*) \to^{1 \oa \epsilon} A \oa \bot \to^{u_\oa^R} A^*
\end{align*}

The cups and caps are inherited directly from $\X$, hence the snake diagrams hold. The antipode in the definition of $\leftaction{0.4}{white}*: H \ox A^* \to A^*$ makes the cup and cap module morphisms.

Suppose $(A, \leftaction{0.4}{white}) \to^{f} (B, \leftaction{0.4}{black})$ is as a module morphism, then $f^* := B^* \to^{f^*} A^* \in \X$ which is also a module morphism. Thus, {\bf H-Mod}$_\X$ is a monoidal category with a dualizing functor, hence a $*$-autonomous category.

If $H$ is cocommutative, then $(A, \leftaction{0.4}{white}) \ox (B, \leftaction{0.4}{black}) \to^{c_\otimes} (B, \leftaction{0.4}{black}) \ox (A, \leftaction{0.4}{white})$ is a module homomorphism.

In that case, {\bf H-Mod}$_\X$ is a symmetric $*$-autonomous category.
\end{proof}

Futhermore, we can show that the category of Hopf modules is conjugative.

\begin{lem}
Let $\X$ be a symmetric $*$-autonomous category. ${\mbox{\bf H-Mod}}_\X$, the category of modules over a cocommutative Hopf Algebra H is a conjugative symmetric $*$-autonomous category.
\end{lem}
\begin{proof}
We already know that ${\mbox{\bf H-Mod}}_\X$ is a symmetric $*$-autonomous category. We define the conjugation functor $\bar{(\_)}: {\mbox{\bf H-Mod}}_\X \to {\mbox{\bf H-Mod}}_\X$ as follows:

\begin{itemize}
\item 
$\overline{(A, \leftaction{0.4}{white})} := (A, \overline{\leftaction{0.4}{white}})$ where,
$ \overline{\leftaction{0.5}{white}} := 
\begin{tikzpicture}[scale=1]
	\begin{pgfonlayer}{nodelayer}
		\node [style=none] (0) at (0, 1) {};
		\node [style=none] (1) at (0, 0.75) {};
		\node [style=none] (2) at (-0.25, 1) {};
		\node [style=none] (3) at (0, 2) {};
		\node [style=none] (4) at (0, 0.5) {};
		\node [style=circle] (5) at (-0.75, 1.5) {};
		\node [style=none] (6) at (-0.75, 2) {};
		\node [style=none] (7) at (-0.75, 1.5) {$s$};
	\end{pgfonlayer}
	\begin{pgfonlayer}{edgelayer}
		\draw (2.center) to (0.center);
		\draw (0.center) to (1.center);
		\draw (1.center) to (2.center);
		\draw (1.center) to (4.center);
		\draw (3.center) to (0.center);
		\draw [bend right, looseness=1.00] (5) to (2.center);
		\draw (5) to (6.center);
	\end{pgfonlayer}
\end{tikzpicture}
$
\item Suppose $f: (A, \leftaction{0.4}{white}) \to (B, \leftaction{0.4}{black})$, then $\overline{f} := f$
\end{itemize}

The  basic natural isomorphisms are given by: \[ \overline{(B, \leftaction{0.4}{black})} \ox \overline{(A, \leftaction{0.4}{white})}  \to^{\chi}  \overline{(A, \leftaction{0.4}{white}) \ox (B, \leftaction{0.4}{black})} := B \ox A \to^{(c_\ox)_{B,A}} A \ox B\]  \[(A, \overline{\overline{\leftaction{0.4}{white}}}) \to^{\varepsilon} (A, \leftaction{0.4}{white}) := 1\]

 The natural isormorphisms satisfy all the coherences of conjugative symmetric $*$-autonomous category.
\end{proof}

\begin{lem}
\label{Lemma conjugative}
Suppose $\X$ is a symmetric (iso)mix $*$-autonomous category, then {\bf H-Mod}$_\X$, the category of Hopf modules over a cocommutative Hopf Algebra H is a (iso)mix conjugative symmetric $*$-autonomous category.
\end{lem}
\begin{proof}~
The mix map $\m: \bot \to \top$ is inherited directly from $\X$.
\end{proof}

\begin{cor}
Suppose $\X$ is a symmetric (iso)mix $*$-autonomous, then {\bf H-Mod}$_\X$, the category of modules over a cocommutative Hopf Algebra H is a symmetric $\dagger$ (iso)mix $*$-autonomous category.
\end{cor}
\begin{proof}
From Lemma \ref{Lemma conjugative},  ${\mbox{\bf H-Mod}}_\X$ is an  (iso)mix conjugative symmetric $*$-autonomous category. Then, by Theorem \ref{Theorem: conjugation+dualizing} one can construct a dagger functor by composing the conjugation and the dualizing functor as follows:  $(\_)^\dagger := \overline{(\_)^*}:{\mbox{\bf H-Mod}}_\X^{\op} \to {\mbox{\bf H-Mod}}_\X$. Therefore,

$(A, \leftaction{0.4}{white})^\dagger := (A^*, \overline{ \leftaction{0.4}{white}^*})$ where,
\[
(A, \overline{\leftaction{0.4}{white}}^*) := 
\begin{tikzpicture} %dagger1
	\begin{pgfonlayer}{nodelayer}
		\node [style=none] (0) at (1.5, -1) {};
		\node [style=none] (1) at (2.25, -1) {};
		\node [style=none] (2) at (1.5, 0.5) {};
		\node [style=none] (3) at (1.5, 0.75) {};
		\node [style=none] (4) at (0.25, 0.75) {};
		\node [style=none] (5) at (2.25, 3.5) {};
		\node [style=none] (6) at (0.25, -1.25) {};
		\node [style=none] (7) at (0.75, 3.5) {};
		\node [style=none] (8) at (1, 3.25) {$H$};
		\node [style=none] (9) at (2.5, 3.25) {$A^*$};
		\node [style=none] (10) at (0.25, -1.75) {};
		\node [style=none] (11) at (0.5, -1.5) {$A^*$};
		\node [style=circle, scale=1.5] (12) at (0.75, 2) {};
		\node [style=none] (13) at (0.75, 2.25) {};
		\node [style=none] (14) at (0.75, 2) {$s$};
		\node [style=circle, scale=1.5] (15) at (0.75, 2.75) {};
		\node [style=none] (16) at (0.75, 2.75) {$s$};
		\node [style=none] (17) at (1.5, 0.5) {};
		\node [style=none] (18) at (1.5, 0.25) {};
		\node [style=none] (19) at (1.25, 0.5) {};
	\end{pgfonlayer}
	\begin{pgfonlayer}{edgelayer}
		\draw [in=-90, out=90, looseness=1.00] (2.center) to (3.center);
		\draw [bend left=90, looseness=1.75] (4.center) to (3.center);
		\draw [bend right=90, looseness=2.00] (0.center) to (1.center);
		\draw (5.center) to (1.center);
		\draw (4.center) to (6.center);
		\draw (13.center) to (12);
		\draw (7.center) to (15);
		\draw (15) to (13.center);
		\draw (6.center) to (10.center);
		\draw (19.center) to (17.center);
		\draw (17.center) to (18.center);
		\draw (18.center) to (19.center);
		\draw [bend right=15, looseness=1.25] (12) to (19.center);
		\draw (18.center) to (0.center);
	\end{pgfonlayer}
\end{tikzpicture} =
\begin{tikzpicture} %dagger2
	\begin{pgfonlayer}{nodelayer}
		\node [style=none] (0) at (1.5, -1) {};
		\node [style=none] (1) at (2.25, -1) {};
		\node [style=none] (2) at (1.5, 0.5) {};
		\node [style=none] (3) at (1.5, 1) {};
		\node [style=none] (4) at (-0.25, 1) {};
		\node [style=none] (5) at (2.25, 3) {};
		\node [style=none] (6) at (0.5, 2.75) {$H$};
		\node [style=none] (7) at (2.5, 2.75) {$A^*$};
		\node [style=none] (8) at (-0.25, -1.75) {};
		\node [style=none] (9) at (0, -1.5) {$A^*$};
		\node [style=none] (10) at (0.75, 3) {};
		\node [style=none] (11) at (1.25, 0.5) {};
		\node [style=none] (12) at (1.5, 0.5) {};
		\node [style=none] (13) at (1.5, 0.25) {};
		\node [style=none] (14) at (2.75, 3.75) {};
	\end{pgfonlayer}
	\begin{pgfonlayer}{edgelayer}
		\draw [in=-90, out=90, looseness=1.00] (2.center) to (3.center);
		\draw [bend left=90, looseness=2.00] (4.center) to (3.center);
		\draw [bend right=90, looseness=2.00] (0.center) to (1.center);
		\draw (5.center) to (1.center);
		\draw (11.center) to (12.center);
		\draw (11.center) to (13.center);
		\draw (13.center) to (12.center);
		\draw (0.center) to (13.center);
		\draw [in=-90, out=135, looseness=1.00] (11.center) to (10.center);
		\draw (4.center) to (8.center);
	\end{pgfonlayer}
\end{tikzpicture}
\]
\end{proof}

Thus, one can generate a $\dagger$-isomix category from a symmetric isomix $*$-autonomous category by choosing the Hopf modules over any cocommutative $\ox$- Hopf Algebra. 

%%%%%%%%%%%%%%%%%%%%%%%%%%%%%%%%%%%%%%%%%%%%%%%%

\section{Unitary structure and mixed unitary categories}
\label{Sec: unitary}

%%%%%%%%%%%%%%%%%%%%%%%%%%%%%%%%%%%%%%%%%%%%%%%%

The objective of this section is to introduce mixed unitary categories (MUCs) and their morphisms.   A mixed unitary category consists of 
a {\em unitary category\/}, $\U$, with a $\dagger$-isomix Frobenius functor $M: \U \to \C$ into a ``large'' $\dagger$ isomix category $\C$.   
We refer to $\U$ as the {\em unitary core\/} of the MUC.  The unitary core is to be regarded as providing the analogue of scalars for the larger 
category much as a field provides scalars for an algebra over that field.  

The section starts by describing the general notion of unitary structure in a $\dagger$-isomix category. This allows the definition of a unitary category as 
a compact $\dagger$-isomix category in which all objects have unitary structure satisfying certain coherence conditions.   We then show how to extract 
a unitary category from any compact $\dagger$-isomix category using {\em pre-unitary objects\/}.  This is a useful construction in practice.  However, it can 
just deliver a trivial unitary category -- trivial in the sense that all objects are isomorphic to the units.  This means that, in applying the construction, it is 
important to identify non-trivial pre-unitary objects to ensure that one is getting something worthwhile out.

Next we show, using the isomix functors ${\sf Mx}_\uparrow$ (or ${\sf Mx}_\downarrow$) that unitary categories are $\dagger$-linearly equivalent to 
$\dagger$-monoidal categories and, furthermore, that closed unitary categories are equivalent to $\dagger$-compact closed categories.  This provides an 
explicit connection from MUCs to the standard notions from categorical quantum mechanics.   One contribution of this more general perspective is that through the 
constructions in this section one can obtain examples not only of mixed unitary categories but also of $\dagger$-monoidal and $\dagger$-compact closed 
categories which might otherwise have been difficult to realize.

The final subsection introduces mixed unitary categories (MUCs).  These form the basis for our approach to infinite dimensional categorical quantum mechanics. 
A MUC has a unitary core which is a model of classical categorical quantum mechanics extended by a larger setting in which infinite dimensional objects can 
be modelled.

%%%%%%%%%%%%%%%%%%%%%%%%%%%%%%%%%%%%%%%%%%%%%%%%%

\subsection{Unitary structure}

The notion of unitary maps is central to both quantum information theory as well as quantum mechanics since the evolution of a closed quantum system is described by such maps. Categorically, within a $\dagger$-category, a unitary map is an isomorphism $f: A \to B$ such that $f^{-1} =  f^\dagger$. This definition of unitary isomorphism cannot be used directly within the framework of $\dagger$-LDCs since the types of $f^{-1}: B \to A$ and  $f^\dagger: B^\dagger \to A^\dagger$ are different. It is therefore
apparent that one can only ask to have unitary isomorphisms  between certain objects, which we call ``unitary objects'':

\begin{defi}
A  $\dagger$-isomix category, $\X$ has {\bf unitary structure} in case there is an essentially small class of objects $\mathcal{U}$, called the {\bf unitary objects} of $\X$ such that
\begin{enumerate}[{\bf [U.1]}]
\item for all $A \in \mathcal{U}$, $A \in  \Core(\X)$, and $A$ is equipped with an isomorphism, $\varphi_A: A \to A^\dag$, called the {\bf unitary structure map} of $A$
\item $\mathcal{U}$ is closed to $(\_)^\dag$ so that for all $A \in \mathcal{U}$, $\varphi_{A^\dag} = ((\varphi_A)^{-1})^\dag$ 
\item for all $A \in \mathcal{U}$, the following diagram commutes:
 \[   \xymatrix{  A   \ar[d]_{\varphi_A} \ar[drrr]^{\iota}  & \\ A^\dag \ar[rrr]_{\varphi_{A^\dag}}  & & & (A^\dag)^\dag  } \]
\item $\bot, \top \in \mathcal{U}$ such that:
\[ \xymatrixcolsep{2pc}
\xymatrix{
\bot \ar[r]^{\varphi_\bot} \ar[d]_{\lambda_\bot} \ar[dr]^{\m} & \bot^\dagger \ar[d]^{\lambda_\top^{-1}}  \\
\top^\dagger \ar[r]_{\varphi_\top^{-1}} & \top
}
\]
\item If $A , B \in \mathcal{U}$, then $A \ox B$ and $A \oa B \in \mathcal{U}$ such that :
\[ (a) ~~~~~ \xymatrixcolsep{3pc}
\xymatrix{
A \ox B \ar[r]^{\varphi_A \ox \varphi_B}_{\simeq} \ar@/_2pc/[rrr]_{\mx}&
 A^\dagger \ox B^\dagger \ar[r]^{\lambda_\oa}_{\simeq} & 
 (A \oa B) ^\dagger \ar[r]^{\varphi_{A \oa B}^{-1}} _{\simeq} &
A \oa B
}
\]
\[ (b) ~~~~~ \xymatrixcolsep{3pc}
\xymatrix{
A \ox B \ar[r]^{\varphi_{A \ox B}}_{\simeq} \ar@/_2pc/[rrr]_{\mx}&
 (A \ox B)^\dagger \ar[r]^{\lambda_\ox^{-1}}_{\simeq} & 
 A^\dagger \oa B^\dagger \ar[r]^{\varphi_A^{-1} \oa \varphi_B^{-1}} _{\simeq} &
A \oa B
}
\]
 \end{enumerate}
\end{defi}

%%%%%%%%%%%%%%%%%%%%%%%%%%%%

\begin{lem}
\label{Lemma: square root tensor unitary}
When $A$ and $B$ is a unitary object in a $\dagger$-isomix category then, $\varphi_{A^{\dagger\dagger}} = (\varphi_A)^{\dagger \dagger}: A^{\dagger\dagger} \to A^{\dagger \dagger \dagger}$.
\end{lem}
\begin{proof}~
\[ \varphi_{(A^\dagger)^{\dagger}} = ((\varphi_{A^\dagger})^{-1})^{\dagger} = ((((\varphi_A)^{-1})^\dagger)^{-1})^\dagger = ((((\varphi_A)^{-1})^{-1})^\dagger)^\dagger = ((\varphi_A)^\dagger)^\dagger  \qedhere \]
\end{proof}

%%%%%%%%%%%%%%%%%%

Often we shall want the unitary objects to have linear adjoints (or duals) but we shall need the analogue of $\dagger$-duals $(\eta^\dagger = c_\ox \epsilon$ and $\epsilon^\dagger = \eta c_\ox)$ from categorical quantum mechanics:

\begin{defi} \label{unitary-duals}
A {\bf unitary linear duality} $(\eta, \epsilon): A \dashvv_{~u} B$ between unitary objects  $A$ and $B$ is a linear duality satisfying in addition:
\[
\begin{matrix}
\xymatrix{ \\
{\bf [Udual.]} \\
}~~~
\xymatrix{
\top \ar@{}[ddrr]|{(a)} \ar[rr]^{\eta} \ar[d]_{\lambda_\top}  & & A \oa B \ar[d]^{\varphi_A \oa \varphi_B} \\
\bot^\dagger \ar[d]_{\epsilon^\dag} & & A^\dagger \oa B^\dagger \ar[d]^{c_\oa} \\ 
(B \ox A)^\dag \ar[rr]_{\lambda_\oa^{-1}} & & B^\dagger \oa A^\dagger} 
~~~~~ & \text{(or)} & ~~~~~~
\xymatrix{
A \ox B \ar@{}[ddrr]|{(b)} \ar[rr]^{\varphi_A \ox \varphi_B} \ar[d]_{c_\ox} & & A^\dag \ox B^\dag \ar[d]^{\lambda_\ox} \\
B \ox A \ar[d]_{\epsilon} & & (A \oa B)^\dagger \ar[d]^{\eta^\dagger} \\
\bot \ar[rr]_{\lambda_\bot} & & \top^\dagger } 
\end{matrix}
\]
\end{defi}

Observe that ${\bf [Udual.]} (a) \Leftrightarrow (b)$. In a compact $\dagger$-LDC, $\top \dashvv_{~u} \bot$. {\bf [Udual] (a)} is shown diagrammatically as follows:
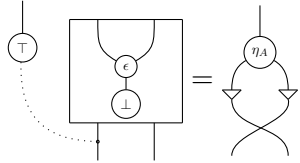
\[\begin{tikzpicture}
	\begin{pgfonlayer}{nodelayer}
		\node [style=none] (0) at (-2, 4) {};
		\node [style=none] (1) at (1, 4) {};
		\node [style=none] (2) at (-2, 2) {};
		\node [style=none] (3) at (1, 2) {};
		\node [style=circle] (4) at (-0.5, 2.5) {$\epsilon$};
		\node [style=none] (5) at (-1.25, 4) {};
		\node [style=none] (6) at (0.25, 4) {};
		\node [style=none] (7) at (-1.25, 2) {};
		\node [style=none] (8) at (0.25, 2) {};
		\node [style=none] (9) at (-1.25, 1) {};
		\node [style=none] (10) at (0.25, 1) {};
	\end{pgfonlayer}
	\begin{pgfonlayer}{edgelayer}
		\draw [in=150, out=-90, looseness=1.25] (5.center) to (4);
		\draw [in=30, out=-90, looseness=1.25] (6.center) to (4);
		\draw (0.center) to (1.center);
		\draw (1.center) to (3.center);
		\draw (3.center) to (2.center);
		\draw (2.center) to (0.center);
		\draw (7.center) to (9.center);
		\draw (8.center) to (10.center);
	\end{pgfonlayer}
\end{tikzpicture}  = \begin{tikzpicture}
	\begin{pgfonlayer}{nodelayer}
		\node [style=circle] (0) at (0.5, 5.75) {$\eta$};
		\node [style=none] (1) at (-0.5, 4.75) {};
		\node [style=none] (2) at (0, 4.75) {};
		\node [style=none] (3) at (-0.25, 4.5) {};
		\node [style=none] (4) at (1, 4.75) {};
		\node [style=none] (5) at (1.5, 4.75) {};
		\node [style=none] (6) at (1.25, 4.5) {};
		\node [style=none] (7) at (-0.25, 3) {};
		\node [style=none] (8) at (1.25, 3) {};
		\node [style=none] (9) at (-0.25, 4.75) {};
		\node [style=none] (10) at (1.25, 4.75) {};
	\end{pgfonlayer}
	\begin{pgfonlayer}{edgelayer}
		\draw (1.center) to (2.center);
		\draw (2.center) to (3.center);
		\draw (3.center) to (1.center);
		\draw (4.center) to (5.center);
		\draw (5.center) to (6.center);
		\draw (6.center) to (4.center);
		\draw [in=90, out=-150, looseness=1.00] (0) to (9.center);
		\draw [in=90, out=-30, looseness=1.00] (0) to (10.center);
		\draw [in=90, out=-75, looseness=0.75] (6.center) to (7.center);
		\draw [in=90, out=-90, looseness=1.00] (3.center) to (8.center);
	\end{pgfonlayer}
\end{tikzpicture} \]

\begin{lem}
Suppose $(\eta_1, \epsilon_1): V_1 \dashvv_{~u} U_1$ and $(\eta_2, \epsilon_2): V_2 \dashvv_{~u} U_2$. Then, $(V_1 \otimes V_2) \dashvv_{~u} (U_1 \oa U_2)$.
\end{lem}
\begin{proof}
Define $(\eta', \epsilon'): (V_1 \otimes V_2) \dashvv_{~u} (U_1 \oa U_2)$ so that  
$\eta' = \begin{tikzpicture} %opluseta
	\begin{pgfonlayer}{nodelayer}
		\node [style=circle] (0) at (-4, 3) {$\eta_1$};
		\node [style=circle] (1) at (-2, 3) {$\eta_2$};
		\node [style=ox] (2) at (-4, 1.75) {};
		\node [style=oa] (3) at (-2, 1.75) {};
		\node [style=none] (4) at (-4, 1) {};
		\node [style=none] (5) at (-2, 1) {};
	\end{pgfonlayer}
	\begin{pgfonlayer}{edgelayer}
		\draw [style=none, in=15, out=-165, looseness=1.00] (1) to (2);
		\draw [style=none, bend left, looseness=1.25] (1) to (3);
		\draw [style=none, in=180, out=-15, looseness=1.00] (0) to (3);
		\draw [style=none, bend left=45, looseness=1.25] (2) to (0);
		\draw [style=none] (2) to (4.center);
		\draw [style=none] (3) to (5.center);
	\end{pgfonlayer}
\end{tikzpicture} ~~~~~~~ 
\epsilon' = \begin{tikzpicture} %oplusepsi
	\begin{pgfonlayer}{nodelayer}
		\node [style=circle] (0) at (-4, 1) {$\epsilon_1$};
		\node [style=circle] (1) at (-2, 1) {$\epsilon_2$};
		\node [style=oa] (2) at (-4, 2.25) {};
		\node [style=ox] (3) at (-2, 2.25) {};
		\node [style=none] (4) at (-4, 3) {};
		\node [style=none] (5) at (-2, 3) {};
	\end{pgfonlayer}
	\begin{pgfonlayer}{edgelayer}
		\draw [style=none, in=-15, out=165, looseness=1.00] (1) to (2);
		\draw [style=none, bend right, looseness=1.25] (1) to (3);
		\draw [style=none, in=180, out=15, looseness=1.00] (0) to (3);
		\draw [style=none, bend right=45, looseness=1.25] (2) to (0);
		\draw [style=none] (2) to (4.center);
		\draw [style=none] (3) to (5.center);
	\end{pgfonlayer}
\end{tikzpicture}
$. This is easily checked to be a unitary linear adjoint.
\end{proof}

We can now define what it means for an isomorphism to be unitary:

\begin{defi}
Suppose $A$ and $B$ are unitary objects. An isomorphism $A\xrightarrow{f} B$ is said to be a {\bf unitary isomorphism} if the following diagram commutes:
\[  \xymatrix{A   \ar[r]^{\varphi_A}    \ar[d]_{f} \ar[r]^{\varphi_A} & A^\dag \\ B  \ar[r]_{\varphi_B} & B^\dag  \ar[u]_{f^\dag}  }  \]
\end{defi}

Observe that $\varphi$ is ``twisted'' natural for all unitary isomorphisms, thus, unitary isomorphisms compose and contain the identity maps. In a category in which the unitary structure maps are identity morphisms, one recovers the usual notion of unitary isomorphisms.

Our next objective is to show that all the coherence isomorphisms between unitary objects are unitary maps. First a warm up: %too poetic ... poetry is good!

\begin{lem}
\label{lemma:MUCProperties}
In a $\dagger$-isomix category with unitary structure:
\begin{enumerate}[(i)]
\item If $f$ is a unitary isomorphism, then so is $f^\dagger$;
\item If $f$ and $g$ are unitary, then so are $f \ox g$ and $f \oa g$;
\item Unitary isomorphisms are closed under composition.
\end{enumerate}
\end{lem}

\begin{proof}~
\begin{enumerate}[{\em (i)}]
\item Recall that $\varphi_{A^\dag} = (\varphi_A^{-1})^\dag$,  then $f^\dagger$ is unitary because 
\[ \xymatrix{B^\dag \ar[d]_{(\varphi_B^{-1})^\dag = \varphi_{B^\dag}} \ar[rr]^{f^\dag} & & A^\dag \ar[d]^{(\varphi_A^{-1})^\dag = \varphi_{A^\dag}} \\
   B^{\dag\dag}  & & A^{\dag\dag} \ar[ll]^{f^{\dag\dag}}} \]
is just the dagger functor applied to the unitary diagram of $f$.
\item Suppose $f$ and $g$ are unitary morphisms, then:
\[
\xymatrix{
A \ox B \ar@{->}[rrr]^{\varphi_{A \ox B}} \ar[ddd]_{f \ox g} \ar[dr]_{\mx}  \ar@{}[dddr]|{\mbox{\tiny \bf (nat. $\mx$)}~~~}
&   \ar@{}[dr]|{\mbox{\tiny {\bf [U.5(b)]}}} &  &  (A \ox B)^\dagger \ar@{}[lddd]|{~~~~~\mbox{\tiny \bf (nat. $\lambda_\oa)$}} \\
& A \oa B \ar[r]^{\varphi_A \oa \varphi_B} \ar[d]_{f \oa g} 
& A^\dagger \oa B^\dagger \ar[ur]_{\lambda_\oa}  & 
\\ & A' \oa B' \ar[r]_{\varphi_{A'} \ox \varphi_{B'}} \ar@{}[dr]|{\mbox{\tiny { \bf [U.5(b)]}}}
& A'^\dagger \oa B'^\dagger \ar[dr]^{\lambda_\oa} \ar[u]_{f^\dagger \oa g^\dagger}
& \\ A' \ox B' \ar[rrr]_{\varphi_{A' \ox B'}} \ar[ur]_{\mx}
& & & (A' \ox B')^\dagger \ar[uuu]_{(f \ox g)^\dagger}
}
\]
The inner square commutes because $f$ and $g$ are unitary maps.
Similarly, using {\bf [U.5(b)]}, one can show that if $f$ and $g$ are unitary, then $f \oa g$ is unitary. %unclear

\item
The proof is trivial. \qedhere
\end{enumerate}
\end{proof}

The following lemma will be used to prove that the associator natural isomorphisms are unitary.
\begin{lem}
\label{lemma: auxiliary}
The following diagram commutes:
\[ \xymatrix{
(A \ox B) \ox C \ar[r]^{\mx} \ar[d]_{a_\ox}  & (A \ox B) \oa C  \ar[r]^{\mx \oa 1}  & (A \oa B) \oa C \ar[d]^{a_\oa} \\
A \ox (B \ox C) \ar[r]_{\mx} & A \oa ( B \ox C)  \ar[r]_{1 \oa \mx} & A \oa (B \oa C) }
\]
\end{lem}
\begin{proof}
\[ \xymatrix{
(A \ox B) \ox C \ar[r]^{\mx} \ar[d]_{a_\ox} \ar@{}[dr]|{{\sf \bf mix}~(b)} & (A \ox B) \oa C \ar@{<-}[d]^{\partial^L}  \ar[r]^{\mx \oa 1} \ar@{}[dr]|{{\sf \bf mix}~(a)}  & (A \oa B) \oa C \ar[d]^{a_\oa} \\
A \ox (B \ox C) \ar[r]_{1 \ox \mx}  \ar@/{_1pc}/[dr]_{\mx}  & A \ox (B \oa C) \ar[r]_{\mx} \ar@{}[d]|{nat. \mx} & A \oa (B \oa C) \\
& A \oa ( B \ox C)  \ar@/{_1pc}/[ur]_{1 \oa \mx}  & }  \]
\end{proof}

\begin{lem}
\label{lemma:cohUnitary}
	
Suppose $\X$ is a $\dagger$-isomix category with unitary structure and $A$, $B$, and $C$ are unitary objects then the following are unitary maps:

\begin{multicols}{2}
\begin{enumerate}[(i)]
\item $\lambda_\ox: A^\dagger \ox B^\dagger \rightarrow (A \oa B)^\dagger$
\item $\lambda_\oa:  A^\dagger \oa B^\dagger \rightarrow (A \ox B)^\dagger$
\item $\lambda_\top: \top \rightarrow \bot^\dagger$
\item $\lambda_\bot: \bot \to \top^\dagger$
\item $\varphi_A: A \rightarrow A^\dagger$ 
\item $m: \top \rightarrow \bot$
\item $\mx_{A,B}: A \ox B \rightarrow A \oa B$
\item $\iota : A \rightarrow (A^{\dagger})^\dagger$
\item $a_\ox: (A \ox B) \ox C \rightarrow A \ox (B \ox C)$
\item $a_\oa: (A \oa  B) \oa C \rightarrow A \oa (B \oa C)$
\item $c_\ox: A \ox B \rightarrow B \ox A$
\item $c_\oa: A \oa B \rightarrow B \oa A$
\item $\partial_L: A \ox (B \oa C) \rightarrow (A \ox B) \oa C$
\item $\partial_R: (A \oa B) \ox C \rightarrow A \oa (B \ox C)$
\end{enumerate}
\end{multicols}
\end{lem}

\begin{proof}~
\begin{enumerate}[(i)]
\item $\lambda_\ox: A^\dagger \ox B^\dagger \rightarrow (A \oa B)^\dagger$ is a unitary map because:

\[\xymatrixcolsep{4pc}\xymatrix{
	{}&&&&\\
	A^\dag\ox B^\dag \ar[r]^{\phi_A^{-1}\ox\phi_B^{-1}} \ar[d]^{\lambda_\ox} \ar@{=}@/^3pc/[rr]     \ar@{}[dr]|{\mbox{\tiny {\bf [U.5(a)] }}}
	 &  A\ox B \ar[r]^{\phi_A\ox \phi_B}  \ar[d]^{\mx} 	                     \ar@{}[dr]|{\mbox{\tiny {\bf nat.}}}
	 &  A^\dag \ox B^\dag \ar[r]^{\phi_{A^\dag \ox B^\dag}} \ar[d]^{\mx}									 \ar@{}[dr]|{\mbox{\tiny {\bf [U.5(a)]}}}
	 &  (A^\dag \ox B^\dag)^\dag \ar[d]_{\lambda_\pr^{-1}} \ar@{=}@/^4pc/[ddd]\\
	(A\pr B)^\dag \ar[r]^{\phi_{A\pr B}^{-1}}  \ar[ddr]_{\phi_{(A\pr B)^\dag}} 								  \ar@{}[dr]|{\mbox{\tiny { \bf [U.3]}}}
	 & A\pr B  \ar[r]^{\phi_A\pr\phi_B} \ar@{=}[d]
	 & A^\dag \pr B^\dag \ar[r]^{\phi_{A^\dag}\pr\phi_{B^\dag}}										 \ar@{}[d]|{\mbox{\tiny { \bf [U.3]}}\ \pr\mbox{\tiny { \bf [U.3]}}}
	 & (A^\dag)^\dag \pr (B^\dag)^\dag  \ar@{=}[d]\\
   {}
     & A \pr B  \ar[rr]^{\iota \pr \iota} \ar[d]^{\iota}
     & {}																					\ar@{}[d]|{\mbox{\tiny {\bf [$\dagger$-ldc.5(a)]}}}
     & (A^\dag)^\dag \pr (B^\dag)^\dag  \ar[d]_{\lambda_\pr}\\
   {}
     & ((A \pr B)^\dag)^\dag  \ar[rr]_{\lambda_\ox^\dag}
     & {}
     & (A^\dag \ox B^\dag)^\dag
}\]

\item $\lambda_\oa$ is unitary because:

\[
\xymatrix{
A^\dag\pr B^\dag                            \ar[rrr]^{\phi_{A^\dag\pr B^\dag}} \ar[dr]^{\mx^{-1}}   \ar[ddd]_{\lambda_\pr} \ar@{}[dddr]|{\mbox{\tiny {\bf Lem. \ref{lemma: mixdagger}}}} 
  &
  &
  &
  (A^\dag\pr B^\dag)^\dag               \ar@{}[dddl]|{\mbox{\tiny {\bf (Lem. \ref{lemma: mixdagger})}}^\dag} \\
{}
  & A^\dag\ox B^\dag                       \ar[r]^{\phi_{A^\dag\ox B^\dag}} \ar[d]_{\lambda_\ox}   \ar@{}[ur]|{\mbox{\tiny {\bf Lem. \ref{lemma:cohUnitary} (vi)}}} \ar@{}[dr]|{\mbox{\tiny {\bf Lem. \ref{lemma:cohUnitary} (i)}}}
  & (A^\dag\ox B^\dag)^\dag           \ar[ur]^{(\mx^{-1})^\dag}
  &\\
{}
  & (A\pr B)^\dag                       \ar[r]^{\phi_{(A\pr B)^\dag}}  \ar@{}[dr]|{\mbox{\tiny {\bf Lems. \ref{lemma:cohUnitary} (vi), \ref{lemma:MUCProperties} (i)}}}
  & ((A\pr B)^\dag)^\dag                \ar[u]_{\lambda_\ox^\dag} \ar[dr]^{((\mx^{-1})^\dag)^\dag}
  &\\
(A\ox B)^\dag                            \ar[rrr]^{\phi_{(A\pr B)^\dag}} \ar[ur]^{(\mx^{-1})^\dag} 
  &
  &
  &
  ((A\ox B)^\dag)^\dag              \ar[uuu]_{\lambda_\pr^\dag}
}
\]

\item $\lambda_\bot: \bot \rightarrow \top^\dagger$ is unitary because:
\[ 
\xymatrix{
\bot \ar[d]_{\lambda_\bot} \ar[rr]^{\varphi_\bot} & & \bot^\dagger \ar[d]^{(\lambda_\bot^{-1})^{\dagger}} \\
\top^\dagger \ar[urr]_{m^\dagger} \ar[rr]_{\varphi_{\top^\dagger} = (\varphi_{\top}^{-1})^\dagger} &  &\top^{\dagger \dagger}
}
\]
The left triangle commutes by {\bf [U.4]} and  {\bf [$\dagger$-mix]}.  The right triangle commutes by {\bf [U.4]} and the functoriality of $\dag$.

\item $\lambda_\top: \top \rightarrow \bot^\dagger$ is unitary because:

\[ 
\xymatrix{
\top \ar[d]_{\lambda_\top} \ar[rr]^{\varphi_\top} & & \top^\dagger \ar[d]^{(\lambda_\top^{-1})^{\dagger}} \\
\bot^\dagger \ar[urr]_{(m^{-1})^\dagger} \ar[rr]_{\varphi_{\bot^\dagger} = (\varphi_{\bot}^{-1})^\dagger} &  &\bot^{\dagger \dagger}
}
\]

The left triangle commutes by {\bf [U.4]} and  {\bf [$\dagger$-mix]}.  The right triangle commutes by {\bf [U.4]} and the functoriality of $\dag$.

\item $\varphi_A$ is unitary because the following square commutes by {\bf [U.3]} and {\bf [U.4]}.
\[
\xymatrix{
A \ar[r]^{\varphi_A} \ar[d]_{\varphi_A} & A^\dagger \ar[d]^{(\varphi^{-1})^\dagger} \\
A^\dagger \ar[r]^{\varphi_{A^\dagger}} & A^{\dagger \dagger}
}
\]

\item $m: \bot \rightarrow \top$ is unitary because:
\[
\xymatrix{
\bot \ar[r]^{\varphi_\bot} \ar[d]_{\m} & \bot^\dagger \ar[d]^{(\m^{-1})^\dagger} \\
\top \ar[r]_{\varphi_\top} \ar[ur]^{\lambda_\top} & \top^\dagger
}
\]
The left and right triangles commute by {\bf [U.4]} and {\bf [$\dagger$-mix]} respectively. Hence, the outer squares commutes.

\item $\mx_{A,B}: A \ox B \rightarrow A \oa B$ is unitary as:
\[ \xymatrix{
{}
  &
  &
  &\\
A \ox B \ar[dd]_\mx\ar@/^2.5pc/[rrr]^{\varphi_{A \ox B}} \ar[r]^\mx \ar@/_1.5pc/[drr]_{\varphi_A \ox \varphi_B}  \ar@{}[drr]|{\mbox{\tiny {\bf nat.}}}   \ar@{}[urrr]|{\mbox{\tiny {\bf [U.5(a)]}}}  
  & A \oa B \ar[r]^{\varphi_A \oa \varphi_B} 
  & A^\dag \oa B^\dag \ar[r]^{\lambda_\oa}  
  & (A \ox B)^\dag \\
{}
  &
  & A^\dag \ox B^\dag \ar[u]^\mx \ar[dr]^{\lambda_\ox}   \ar@{}[ur]|{\mbox{\tiny {\bf Lem. \ref{lemma: mixdagger}}}} \\
A \oa B \ar[rrr] _{\varphi_{A\oa B}}    \ar@{}[urr]|{\mbox{\tiny {\bf [U.3]}}} 
  &
  &
  & (A \oa B)^\dag \ar[uu]_{\mx^\dag}
}
\]
                    
\item $\iota: A \rightarrow A^{\dagger \dagger}$ is unitary as in
\[
\xymatrix{
A \ar[d]_{\iota} \ar[r]^{\varphi_A} & A^\dagger \ar[d]^{(\iota^{-1})^\dagger} \ar[ld]_{\varphi_{A^\dagger}} \\ 
A^{\dagger \dagger} \ar[r]_{\varphi_{A^{\dagger \dagger}}} & A^{\dagger \dagger \dagger}
}
\] 
the left triangle commutes by {\bf [U.3]} and the right triangle commutes by: 
\begin{eqnarray*}
(\iota^{-1})^\dagger &= & ((\varphi_{A^\dagger})^{-1} \varphi_A^{-1})^\dagger =  (((\varphi_{A}^{-1})^\dagger)^{-1} \varphi_A^{-1})^\dagger \\
                         & =  & ((\varphi_A^\dagger)(\varphi_A^{-1}) )^\dagger = ( \varphi_A^{-1} )^{\dagger} ( \varphi_A )^{\dagger \dagger} \\
                         & = &  \varphi_{A^\dagger} (\varphi_A)^{\dagger \dagger} = \varphi_{A^\dagger} (\varphi_{A^{\dagger \dagger}})
\end{eqnarray*}

\item $a_\ox$ is unitary as:

\[
\xymatrix{
(A \ox B) \ox C                                           \ar[rrr]^{\varphi_{(A \ox B) \ox C}} \ar[ddddd]_{a_\ox} \ar[dr]_{\mx}   \ar@{}[drrr]|{\mbox{\tiny {\bf [U.5(b)]}}} \ar@{}[dddddr]|{\mbox{\tiny \bf Lem.~ \ref{lemma: auxiliary}}}
 & {} 
 & {}
 & ( (A \ox B) \ox C )^\dagger                     \\ %1
 & {}
(A \ox B) \pr C                                             \ar[r]^{\varphi_{A\ox B} \pr \varphi_{C}} \ar[d]_{\mx\pr 1}   \ar@{}[dr]|{\mbox{\tiny {\bf  [U.5(b)]$\pr$(id)}}} 
 & (A \ox B)^\dag \pr C^\dag                               \ar[d]^{ \lambda_\oa^{-1}\pr 1}   \ar[ur]_{\lambda_\oa}   \ar@{}[dddr]|{\mbox{\tiny {\bf [\dag-ldc.1]}}} 
 & {}  \\ %2
{}
 & (A \pr B) \pr C                                            \ar[r]_{(\varphi_A \oa \varphi_B) \oa \varphi_C} \ar[d]_{a_\oa}  \ar@{}[dr]|{\mbox{\tiny {\bf nat.}}} 
 & (A^\dag \pr B^\dag) \pr C^\dag                         \ar[d]^{a_\oa} 
 & {} \\ %3
{}
 & A \pr (B \pr C)                                           \ar[r]_{\phi_A \pr (\phi_B \pr \phi_C)}   \ar@{}[dr]|{\mbox{\tiny {\bf (id)$\pr$[U.5(b)] }}} 
 & A^\dag \pr (B^\dag \pr C^\dag)       \ar[d]^{1\oa\lambda_\oa}
 & {} \\ %4
{}
 & A \pr (B \ox C)                                          \ar[r]_{\varphi_{A} \pr \varphi_{B\ox C}} \ar[u]^{1\pr \mx}
 & A^\dag \pr (B \ox C)^\dag               \ar[dr]^{\lambda_\oa} 
 & {} \\ %5
A \ox (B \ox C)                                          \ar[rrr]^{\varphi_{A \ox (B \ox C)}} \ar[ur]^{\mx}   \ar@{}[urrr]|{\mbox{\tiny {\bf [U.5(b)]}}} 
 & {}
 & {}
 & ( A \ox (B \ox C)  )^\dagger \ar[uuuuu]_{a_\ox^\dag} %6    
}
\]

\item $a_\oa$ is unitary because:

\[
\xymatrix{
(A \pr B) \pr C                                           \ar[rrr]^{\varphi_{(A \pr B) \pr C}} \ar[ddddd]_{a_\pr} \ar[dr]_{\mx^{-1}}   \ar@{}[drrr]|{\mbox{\tiny {\bf [U.5(a)]}}}  \ar@{}[dddddr]|{\mbox{\tiny \bf Lem.~ \ref{lemma: auxiliary}}} 
 & {} 
 & {}
 & ( (A \pr B) \pr C )^\dagger                     \\ %1
 & {}
(A \pr B) \ox C                                             \ar[r]^{\varphi_{A\pr B} \ox \varphi_{C}} \ar[d]_{\mx^{-1}\ox 1}   \ar@{}[dr]|{\mbox{\tiny {\bf  [U.5(a)]$\ox$(id)}}} 
 & (A \pr B)^\dag \ox C^\dag                               \ar[d]^{ \lambda_\ox^{-1}\ox 1}   \ar[ur]_{\lambda_\pr}   \ar@{}[dddr]|{\mbox{\tiny {\bf [\dag-ldc.1]}}} 
 & {}  \\ %2
{}
 & (A \ox B) \ox C                                            \ar[r]_{(\varphi_A \pr \varphi_B) \pr \varphi_C} \ar[d]_{a_\ox}  \ar@{}[dr]|{\mbox{\tiny {\bf nat.}}} 
 & (A^\dag \ox B^\dag) \ox C^\dag                         \ar[d]^{a_\ox} 
 & {} \\ %3
{}
 & A \ox (B \ox C)                                           \ar[r]_{\phi_A \ox (\phi_B \ox \phi_C)}   \ar@{}[dr]|{\mbox{\tiny {\bf (id)$\ox$[U.5(a)] }}} 
 & A^\dag \ox (B^\dag \ox C^\dag)       \ar[d]^{1\ox\lambda_\ox}
 & {} \\ %4
{}
 & A \ox (B \pr C)                                          \ar[r]_{\varphi_{A} \ox \varphi_{B\pr C}} \ar[u]^{1\ox \mx^{-1}}
 & A^\dag \ox (B \oa C)^\dag               \ar[dr]^{\lambda_\ox} 
 & {} \\ %5
A \oa (B \oa C)                                          \ar[rrr]^{\varphi_{A \oa (B \oa C)}} \ar[ur]^{\mx^{-1}}   \ar@{}[urrr]|{\mbox{\tiny {\bf [U.5(a)]}}} 
 & {}
 & {}
 & ( A \oa (B \oa C)  )^\dagger \ar[uuuuu]_{a_\oa^\dag} %6    
}
\]

\item $c_\ox$ is unitary because:

\[
\xymatrix{
A\ox B                            \ar[rrr]^{\phi_{A\ox B}} \ar[dr]^{\mx} \ar[ddd]_{c_\ox}
 & {}                                  \ar@{}[dr]|{\mbox{\tiny {\bf [U.5(b)]}}}
 & {}
 & (A\ox B)^\dag              \\
{}                                    
 & A \pr B                        \ar[r]^{\phi_A \pr \phi_B} \ar[d]_{c_\pr}   \ar@{}[dr]|{\mbox{\tiny {\bf nat.}}}
 & A^\dag \pr B^\dag       \ar[d]^{c_\pr} \ar[ur]^{\lambda_\pr}           \ar@{}[dr]|{\mbox{\tiny {\bf [$\dagger$-ldc.2(b)]}}}
 & {} \\
{}
 & B \pr A                       \ar[r]^{\phi_B \pr \phi_A}
 & B^\dag \pr A^\dag      \ar[dr]^{\lambda_\pr}
 & {} \\
B\ox A                           \ar[ur]^{\mx} \ar[rrr]^{\phi_{B\ox A}}
& {}                                 \ar@{}[ur]|{\mbox{\tiny {\bf [U.5(b)]}}}
& {}
& (B\ox A)^\dag                 \ar[uuu]_{(c_\ox^{-1})^\dag = c_\ox^\dag}\\
}
\]

where the left square commutes because

$$
\begin{tikzpicture}
	\begin{pgfonlayer}{nodelayer}
		\node [style=circ] (0) at (0.5, -0.25) {};
		\node [style=circ] (1) at (0, -1) {$\top$};
		\node [style=map] (2) at (0, -1.75) {};
		\node [style=circ] (3) at (0, -2.5) {$\bot$};
		\node [style=circ] (4) at (-0.5, -3.25) {};
		\node [style=none] (5) at (0.5, -3.25) {};
		\node [style=none] (6) at (-0.5, -4.25) {};
		\node [style=none] (7) at (0.5, -4.25) {};
		\node [style=none] (8) at (-0.5, 0.5) {};
		\node [style=none] (9) at (0.5, 0.5) {};
	\end{pgfonlayer}
	\begin{pgfonlayer}{edgelayer}
		\draw [densely dotted, in=-90, out=45, looseness=1.00] (4) to (3);
		\draw (3) to (2);
		\draw (2) to (1);
		\draw [densely dotted, in=-135, out=90, looseness=1.00] (1) to (0);
		\draw [style=none] (0) to (9);
		\draw [style=none] (8) to (4);
		\draw [style=none] (0) to (5);
		\draw [style=none, in=90, out=-90, looseness=1.00] (5) to (6);
		\draw [style=none, in=90, out=-90, looseness=1.00] (4) to (7);
	\end{pgfonlayer}
\end{tikzpicture}
=
\begin{tikzpicture}
	\begin{pgfonlayer}{nodelayer}
		\node [style=circ] (0) at (0, -1.25) {$\top$};
		\node [style=map] (1) at (0, -2) {};
		\node [style=circ] (2) at (0.5, -3.5) {};
		\node [style=circ] (3) at (0, -2.75) {$\bot$};
		\node [style=circ] (4) at (-0.5, -0.5) {};
		\node [style=none] (5) at (0.5, -0.5) {};
		\node [style=none] (6) at (-0.5, 0.5) {};
		\node [style=none] (7) at (0.5, 0.5) {};
		\node [style=none] (8) at (0.5, -4.25) {};
		\node [style=none] (9) at (-0.5, -4.25) {};
	\end{pgfonlayer}
	\begin{pgfonlayer}{edgelayer}
		\draw [densely dotted, in=-90, out=150, looseness=1.25] (2) to (3);
		\draw (3) to (1);
		\draw (1) to (0);
		\draw [densely dotted, in=-45, out=90, looseness=1.00] (0) to (4);
		\draw [style=none] (8) to (2);
		\draw [style=none] (9) to (4);
		\draw [style=none, in=-90, out=90, looseness=1.00] (4) to (7);
		\draw [style=none, in=-90, out=90, looseness=1.00] (5) to (6);
		\draw [style=none] (5) to (2);
	\end{pgfonlayer}
\end{tikzpicture}
=
\begin{tikzpicture}
	\begin{pgfonlayer}{nodelayer}
		\node [style=circ] (0) at (0.5, -0.25) {};
		\node [style=circ] (1) at (0, -1) {$\top$};
		\node [style=map] (2) at (0, -1.75) {};
		\node [style=circ] (3) at (0, -2.5) {$\bot$};
		\node [style=circ] (4) at (-0.5, -3.25) {};
		\node [style=none] (5) at (-0.5, -4) {};
		\node [style=none] (6) at (0.5, -4) {};
		\node [style=none] (7) at (-0.5, 0.75) {};
		\node [style=none] (8) at (0.5, 0.75) {};
		\node [style=none] (9) at (-0.5, -0.25) {};
	\end{pgfonlayer}
	\begin{pgfonlayer}{edgelayer}
		\draw [densely dotted, in=-90, out=45, looseness=1.00] (4) to (3);
		\draw (3) to (2);
		\draw (2) to (1);
		\draw [densely dotted, in=-135, out=90, looseness=1.00] (1) to (0);
		\draw [style=none] (6) to (0);
		\draw [style=none] (5) to (4);
		\draw [style=none] (4) to (9);
		\draw [style=none, in=-90, out=90, looseness=1.00] (9) to (8);
		\draw [style=none, in=-90, out=90, looseness=1.00] (0) to (7);
	\end{pgfonlayer}
\end{tikzpicture}
$$

\item $c_\pr$ is unitary because:

\[
\xymatrix{
A\pr B                            \ar[rrr]^{\phi_{A\pr B}}  \ar[dr]^{\mx^{-1}}   \ar[ddd]_{c_\pr}  \ar@{}[drrr]|{\mbox{\tiny {\bf Lem. \ref{lemma:cohUnitary} (vii)}}}
  &
  &
  & (A\pr B)^\dag \\
{}
  & A \ox B                      \ar[r]^{\phi_{A\ox B}} \ar[d]_{c_\ox}    \ar@{}[dr]|{\mbox{\tiny {\bf Lem. \ref{lemma:cohUnitary} (xi)}}}
  & (A\ox B)^\dag     \ar[ur]^{(\mx^{-1})^\dag} &\\
{}
  & B \ox A                      \ar[r]^{\phi_{B\ox A}} 
  & (B \ox A)^\dag    \ar[dr]^{(\mx^{-1})^\dag}  \ar[u]_{c_\ox^\dag} &\\
B\pr A                            \ar[rrr]^{\phi_{B\pr A}}   \ar[ur]^{\mx^{-1}}   \ar@{}[urrr]|{\mbox{\tiny {\bf Lem. \ref{lemma:cohUnitary} (vii)}}}
  &
  &
  & (B\pr A)^\dag              \ar[uuu]_{c_\pr^\dag} 
  }
\]

where the left square commutes for the same reason and the right square is the dagger of the left square.

\item $\partial_L$ is unitary (see Figure \ref{Fig: linear dist. unitary}).

\begin{figure}
\newpage

\begin{sideways}
\scalebox{.74}{
$
\xymatrix{
{}
& {}
& {}
& {}
& {}
& {}
& {}
& {}
& {}\\
{}
& {}
& {}
& {}
& {}
& {}
& {}
& {}
& {}\\
{}
%%%%%%%%%%%%%%%%%%%%%%%%%
%                              0                                     %
%%%%%%%%%%%%%%%%%%%%%%%%%
 & {}
 & A\ox (B \pr C)                      \ar[rr]^{\phi_{A}\ox \phi_{B\pr C}}    \ar[d]_{\mx}       \ar@/^4pc/[rrrrrr]^{\phi_{A\ox (B\pr C)}} \ar@/_6pc/[dddddd]_{\partial_L}
 & {}                                                                                                                                          \ar@{}[d]|{\mbox{\tiny {\bf [U.6(b)]}}}
 & A^\dag \ox (B \pr C)^\dag            \ar@=[r]   \ar@=[d]                                                                \ar@{}[dr]|{\mbox{\tiny {\bf id}}}
 & A^\dag \ox (B \pr C)^\dag            \ar[r]^{\mx}    \ar[d]^{1 \ox \lambda_\ox^{-1} }     \ar@{}[dr]|{\mbox{\tiny {\bf nat.}}}    \ar@{}[u]|{\mbox{\tiny {\bf [U6.(b)]}}}  
 & A^\dag \pr (B \pr C)^\dag            \ar[rr]^{\lambda_\pr}    \ar[d]^{1 \pr \lambda_\ox^{-1}}                                              \ar@{}[ddrr]|{\mbox{\tiny {\bf [$\dagger$-ldc.4(b)]}}}  
 & {}
 & (A\ox (B\pr C))^\dag                         \\
%%%%%%%%%%%%%%%%%%%%%%%%%
%                              1                                     %
%%%%%%%%%%%%%%%%%%%%%%%%%
{}
 & {}
 & A\pr (B \pr C)                      \ar[r]^{\phi_{A\pr (B \pr C)}}     \ar[d]_{a_\pr^{-1}}                                \ar@{}[dr]|{\mbox{\tiny {\bf Lem. \ref{lemma:cohUnitary} (ix)}}}
 & (A\pr (B\pr C))^\dag               \ar[r]^{\lambda_\ox^{-1}}      \ar[d]^{(a_\pr)^\dag}
 & A^\dag \ox (B\pr C)^\dag           \ar[r]^{1\ox \lambda_\ox^{-1}}                                                 \ar@{}[d]|{\mbox{\tiny {\bf [$\dagger$-ldc.1(b)]}}}
 & A^\dag \ox (B^\dag \ox C^\dag)     \ar[r]^{\mx }    \ar[d]^{a_\ox^{-1}}                                                 \ar@{}[dr]|{\mbox{\tiny {\bf mx.~cat.}}}
 & A^\dag \pr (B^\dag \ox C^\dag)             %   \ar[d]^{\partial_R}
 & {}
 & {}\\
%%%%%%%%%%%%%%%%%%%%%%%%%
%                              2                                     %
%%%%%%%%%%%%%%%%%%%%%%%%%
{}
 & {}
 & (A\pr B)\pr C                       \ar[r]^{\phi_{(A\pr B)\pr C}}     \ar@=[d]
 & ((A\pr B)\pr C)^\dag                \ar[r]^{\lambda_\ox^{-1}}                                                            \ar@{}[d]|{\mbox{\tiny {\bf [U.6(a)]}}}
 & (A\pr B)^\dag \ox C^\dag           \ar[r]^{\lambda_\ox^{-1}\ox 1}    \ar[d]_{\mx}                             \ar@{}[dr]|{\mbox{\tiny {\bf nat.}}}
 & (A^\dag \ox B^\dag) \ox C^\dag      \ar[r]^{ \mx\ox 1}     \ar[d]^{\mx}                                            \ar@{}[dr]|{\mbox{\tiny {\bf nat.}}}
 & (A^\dag \pr B^\dag) \ox C^\dag                 \ar[d]^{\mx}\ar@=[r]          \ar[u]_{\partial_R}                                          \ar@{}[ddddr]|{\mbox{\tiny {\bf nat.}}} 
 & (A^\dag \pr B^\dag) \ox C^\dag                 \ar[dddd]^{\lambda_\pr\ox 1}
 & {}\\
%%%%%%%%%%%%%%%%%%%%%%%%%
%                              3                                     %
%%%%%%%%%%%%%%%%%%%%%%%%%
{}
 & {}
 & (A\pr B)\pr C                      \ar[rr]^{\phi_{A\pr B}\pr \phi_{C}}     \ar[d]_{\mx^{-1} \oa 1}                     \ar@{}[ll]|{\mbox{\tiny {\bf mx ~ cat.}}}     \ar@{}[drr]|{\mbox{\tiny{\bf (id)$\ox$(Lem. \ref{lemma:cohUnitary}   (vii))} } }
 & {}                                                                                                                                              
 & (A\pr B)^\dag \pr C^\dag           \ar[r]^{\lambda_\ox^{-1}\pr 1}     \ar[d]_{\mx^\dag\pr 1}         \ar@{}[dr]|{\mbox{\tiny {\bf 1$\ox$(Lem. \ref{lemma: mixdagger})  }}}
 & (A^\dag \ox B^\dag) \pr C^\dag                \ar[d]^{\mx\pr 1}  \ar[r]^{\mx\pr 1}                          \ar@{}[ddr]|{\mbox{\tiny {\bf id}}}
 & (A^\dag \pr B^\dag) \pr C^\dag                \ar@=[dd]
 & {}
 & {}\\
%%%%%%%%%%%%%%%%%%%%%%%%%
%                              4                                     %
%%%%%%%%%%%%%%%%%%%%%%%%%
{}
 & {}
 & (A\ox B)\oa C                       \ar[rr]^{\phi_{A\pr B}\oa \phi_{C}}    \ar@=[d]                                    
 & {}                                                                                                                                              \ar@{}[d]|{\mbox{\tiny {\bf [U.6(a)] }}}
 & (A\ox B)^\dag \oa C^\dag           \ar[r]^{\lambda_\ox^{-1}\pr 1}     \ar[d]_{\lambda_\ox}            \ar@{}[dr]|{\mbox{\tiny {\bf id}}}
 & (A^\dag\oa B^\dag) \oa C^\dag                \ar[d]^{\lambda_\oa \oa 1}
 & {}
 & {}
 & {}\\
%%%%%%%%%%%%%%%%%%%%%%%%%
%                              5                                     %
%%%%%%%%%%%%%%%%%%%%%%%%%
{}
 & {}
 & (A\ox B)\oa C                       \ar[r]^{\mx^{-1}}     \ar@=[d]
 & (A\pr B)\ox C                       \ar[r]^{\phi_{(A\pr B)\pr C}}                                                            \ar@{}[dr]|{\mbox{\tiny {\bf [U.6(a)]}}}
 & ((A\ox B) \ox C)^\dag                \ar[r]^{\lambda_\ox^{-1}}
 & (A\ox B)^\dag \oa C^\dag           \ar[r]^{ \lambda_\oa^{-1} \ox 1}     \ar@=[d]                               \ar@{}[dr]|{\mbox{\tiny {\bf id}}}
 & (A^\dag \oa B^\dag) \oa C^\dag       \ar[d]^{\lambda_\oa \ox 1}
 & {}
 & {}\\
%%%%%%%%%%%%%%%%%%%%%%%%%
%                              6                                     %
%%%%%%%%%%%%%%%%%%%%%%%%%
{}
 & {}
 & (A\ox B) \oa C                      \ar[rrr]^{\phi_{A\ox B} \oa \phi_{C}}                           \ar@/_4pc/[rrrrrr]_{\phi_{(A\ox B) \oa C}}
 & {}
 & {}
 & (A \ox B)^\dag \oa C^\dag          \ar@=[r]                                                                                   \ar@{}[d]|{\mbox{\tiny {\bf [U.6(a)]}}}
 & (A \ox B)^\dag \oa C^\dag          \ar[r]^{\mx^{-1}}
 & (A \ox B)^\dag \ox C^\dag          \ar[r]^{\lambda_\ox}
 &  ((A \ox B) \oa C)^\dag               \ar[uuuuuu]^{\partial_L^\dag}\\
{}
 & {}
 & {}
 & {}
 & {}
 & {}
 & {}
 & {}
 & {}\\
}
$
}

\end{sideways}
\caption{$\partial_L$ is a unitary isomorphism}
\label{Fig: linear dist. unitary}
\newpage
\end{figure}

\item $\partial_R$ is unitary because:

  \scalebox{0.8}{\begin{minipage}{1.2\linewidth}
\[ \hspace{-1.25cm} 
\xymatrix{
{}
  & {}
  & {}
  & {}
  & {}
  & {}\\
(\!A\!\pr\! B\!) \!\ox\! C                           \ar[r]^{\mx} \ar[d]_{\partial_R}    \ar@/^4pc/[rrrrr]^{\phi_{(A\pr B) \ox C }}
  & (\!A\!\pr\! B\!) \!\pr\! C                      \ar[r]^{\mx^{-1}\pr 1}                                                        \ar@{}[d]|{\mbox{\tiny {\bf }}}
  & (\!A\!\ox\! B\!) \!\pr\! C                      \ar[r]^{\phi_{(A\ox B) \pr C}}                                              \ar@{}[dr]|{\mbox{\tiny {\bf Lem. \ref{lemma:cohUnitary}  (xiii)}}}  \ar@{}[ur]|{\mbox{\tiny {\bf Lem. \ref{lemma:cohUnitary} (vii), \ref{lemma:MUCProperties}}}}
  & (\!(\!A\!\ox\! B\!) \!\pr\! C \!)^\dag         \ar[r]^{(\mx^{-1}\pr 1 )^\dag} \ar[d]_{\partial_L^\dag}   
  & (\!(\!A\!\pr\! B\!) \!\pr\! C \!)^\dag         \ar[r]^{\mx^\dag}                                                               \ar@{}[d]|{\mbox{\tiny {\bf }}}
  & (\!(\!A\!\pr\! B\!) \!\ox\! C \!)^\dag \\
%%%%%%%%%
A\!\pr\! (\!B \!\ox\! C\!)                           \ar[r]^{\mx^{-1}}   \ar@/_4pc/[rrrrr]_{\phi_{A\pr (B \ox C) }}
  & A\!\pr\! (\!B \!\pr\! C\!)                      \ar[r]^{1 \pr \mx} 
  & A\!\ox\! (\!B \!\pr\! C\!)                      \ar[r]^{\phi_{A\ox (B \pr C)}} \ar[u]_{\partial_L}  \ar@{}[dr]|{\mbox{\tiny {\bf Lem. \ref{lemma:cohUnitary}  (vii), \ref{lemma:MUCProperties}}}}
  & (\!A\!\ox\! (\!B\!\pr\! C\!) \!)^\dag         \ar[r]^{(\mx \pr 1 )^\dag}
  & (\!A\!\pr\! (\!B \!\pr\! C\!) \!)^\dag         \ar[r]^{(\mx^{-1})^\dag}
  & (\!A\!\pr\! (\!B \!\ox\! C\!) \!)^\dag        \ar[u]_{\partial_R^\dag}\\
{}
  & {}
  & {}
  & {}
  & {}
  & {}
}
\]
\end{minipage}}
\end{enumerate}
\end{proof}

%%%%%%%%%%%%%%%%%%%%%%%%%%%%%%%%%%%%%%%%%%%%%%

\subsection{Unitary categories}

With the notion of unitary objects in place, one can consider $\dagger$-isomix categories in which all the objects are unitary: these are called {\em unitary categories\/}. This section develops the theory of unitary categories.

\begin{defi}
A {\bf unitary category} is a $\dagger$-isomix category with unitary  structure such that every object in the category is a unitary object.
\end{defi}

Clearly, a unitary category must be a compact $\dagger$-LDC, because every object is in the core.

 A $\dagger$-monoidal category is a strict unitary category in which the unitary structure map and the mix map are identity morphisms. Similarily, a $\dagger$-compact closed category is a strict unitary category in which all objects have unitary duals.

 In the rest of this subsection, we show that any unitary category is $\dagger$-linearly equivalent to a conventional dagger monoidal category. A unitary category being a compact LDC is linearly equivalent, using ${\sf Mx}^*_\uparrow: (\X, \ox,\oa) \to (\X,\oa,\oa)$ (see Corollary \ref{compact-mix-functor}) to the underlying monoidal category based on the par (and the tensor). We now show that for a unitary category one can induce a stationary on objects dagger on $(\X,\oa,\oa)$. We denote this dagger by $(\_)^\ddagger$ and define it by $f^\ddagger := \varphi_Bf^\dagger\varphi_A^{-1}$ as illustrated by the left diagram below:
 
 \[ \xymatrix{ B \ar[d]_{\varphi_B}\ar[rr]^{f^\ddagger} \ar@{}[rrd]|{:=}& & A \ar[d]^{\varphi_A} \\
 	B^\dagger \ar[rr]_{f^\dagger} && A^\dagger} ~~~~~~~~~~~~~
 \xymatrix{ A \ar@/_1pc/[dd]_{\iota} \ar[d]^{\varphi_A}\ar[rr]^{f^{\ddagger\ddagger}} & & B \ar[d]_{\varphi_B} \ar@/^1pc/[dd]^{\iota} \\	
 	A^\dagger \ar[d]^{(\varphi_A^{-1})^\dagger} \ar[rr]_{(f^\ddagger)^\dagger} && B^\dagger \ar[d]_{(\varphi_B^{-1})^\dagger} \\ 	
 	A^{\dagger\dagger} \ar[rr]_{f^{\dagger\dagger}} & & B^{\dagger\dagger} } \]
 
 This new dagger clearly preserves composition and is also a stationary on objects involution as proven by the second diagram:
 the lower square of this diagram is the dagger of the inverted definition and the resulting outer square is the naturality of $\iota$ forcing $f^{\ddagger\ddagger} = f$.
 
 Next, we observe that $u: X \to Y$ is a unitary isomorphism in $\X$ if and only if $u^{-1}= u^\ddagger$.  This makes unitary isomorphisms in the traditional sense of categorical quantum mechanics coincide 
 with the notion introduced here.   Thus, $u$ is unitary in the sense here if and only if the diagram below commutes
 \[ \xymatrix{ B \ar[d]_{\varphi_B} \ar[rr]^{u^{-1}} && A  \ar[d]^{\varphi_A} \\ B^\dagger \ar[rr]_{u^\dagger} && A} \]
 but this diagram commutes if and only if $u^{-1} = u^\ddagger$. 
 
 \begin{defi} \label{preserving-unitary-structure}
 	A $\dagger$-Frobenius mix functor, $F: \X \to \Y$, between compact $\dagger$-isomix categories with unitary structure {\bf preserves unitary structure} if
 	\begin{enumerate}[(i)]
 		\item for all unitary objects $A \in \X$, $F(A)$ is a unitary object such that $\varphi_{F(A)} = F(\varphi_A) \rho^F$ 
 		\item Either $n_\bot^F$ or $m_\top^F$ are unitary isomorphisms i.e.,
 		\[
 		\xymatrix{
 			F(\bot) \ar[r]^{F(\varphi_\bot)} \ar[d]_{n_\bot} & F(\bot^\dagger) \ar[r]^{\rho} & F(\bot)^\dagger \\
 			\bot \ar[rr]_{\varphi_\bot} & & \bot^\dagger \ar[u]_{n_\bot^\dagger}
 		} (or)  \xymatrix{
 			\top \ar[rr]^{\varphi_\top} \ar[d]_{m_\top} & & \top^\dagger \\
 			F(\top) \ar[r]_{F(\varphi_\top)} & F(\top^\dagger) \ar[r]_{\rho} & F(\top)^\dagger \ar[u]_{m_\top^\dagger}
 		}
 		\]
 	\end{enumerate}
 \end{defi}
 
Notice that if $F$ preserves unitary structure, it must be an isomix functor by Lemma \ref{Lemma: isomix functor}. Also, when $A \in \X$ is a unitary object,  then $F(A)$ must be a unitary object, and so $F(A)$ is in the core.
 
 We now show that ${\sf Mx}_\uparrow: (\X,\oa,\oa) \to (\X,\ox,\oa)$ provides a unitary structure preserving equivalence of a dagger monoidal category into a unitary category:
 
 \begin{prop}  \label{unitary-2-dagger}
 	Unitary categories are $\dagger$-linearly equivalent via the mix functor ${\sf Mx}_\uparrow: (\X,\oa,\oa) \to (\X,\ox,\oa)$ to the underlying dagger monoidal category on the par.  
 	Furthermore, closed unitary categories under this equivalence become dagger compact closed categories.
 \end{prop}
 
 \begin{proof}
 	We must exhibit a preservator, that is a natural transformation showing that the involution is preserved:
 	\[ \infer={A \to_{\varphi_A} A^\dagger}{{\sf Mx}_\uparrow(A^\ddagger) \to^{\varphi_A} {\sf Mx}_\uparrow(A)^\dagger} \]
 	Note that $\varphi$ is a natural transformation by the definition of $(\_)^\ddagger$ and its coherence requirements 
 	make it a linear natural equivalence.  Making this the preservator immediately means that unitary structure is preserved.
 	
 	Finally, we must show that unitary linear duals under ${\sf Mx}^{*}_\uparrow$ become $\ddagger$-duals.  Given $(\eta,\epsilon): A \dashvv_u B$  we must 
 	show that under ${\sf Mx}^{*}_\uparrow$ this produces a dagger dual.  ${\sf Mx}^{*}_\uparrow(\eta) = {\sf m} ~\eta: \bot \to A \oa B$ and 
 	${\sf Mx}^{*}_\uparrow(\epsilon) = {\sf mx}^{-1} \epsilon: B \oa A \to \bot$
 	We then require that $c_\oa {\sf Mx}^{*}_\uparrow(\epsilon) = {\sf Mx}^{*}_\uparrow(\eta)^\ddagger$. This is provided by:
 	\[ \xymatrix{A \oa B  \ar[d]^{{\sf mx}^{-1}}  \ar@/_2pc/[ddd]_{\varphi_{A \oa B}} \ar[r]^{c_\oa} & B \oa A \ar@/^1pc/[rr]^{{\sf Mx}^{*}_\uparrow(\epsilon)} \ar[r]_{{\sf mx}^{-1}} 
 		& B \ox A \ar[r]_{\epsilon} & \bot \ar[dd]_{{\sf m}} \ar@/^2pc/[ddd]^{\varphi_\bot} \ar@/_/[dddl]_{\lambda_\bot}\\
 		A \ox B \ar[d]^{\varphi_A \ox \varphi_B} \ar@/_/[rru]^{c_\ox} \\
 		A^\dagger \ox B^\dagger \ar[d]^{\lambda_\ox} & & & \top \ar[d]_{\lambda_\top} \\
 		(A \oa B)^\dagger \ar@{}[rrruuu]|{{\rm Defn.} ~\ref{unitary-duals}~(b)} \ar@/_2pc/[rrr]_{{\sf Mx}^{*}_\uparrow(\eta)^\dagger} \ar[rr]_{\eta^\dagger} & & \top^\dagger \ar[r]_{{\sf m}^\dagger} & \bot^\dagger }  \]
 \end{proof}

%%%%%%%%%%%%%%%%%%%%%%%%%%%%%%%%%%
% Unitary construction %
%%%%%%%%%%%%%%%%%%%%%%%%%%%%%%%%%%

\subsection{The unitary construction}

A $\dagger$-isomix category can have many different unitary structures, as we shall describe in this section, 
thus it is {\em structure\/}, and not a property.   The requirements, however, do mean that for a $\dagger$-isomix category, $\X$, there is always the 
smallest unitary structure, referred to as the ``trivial'' unitary structure, that produces a full unitary subcategory in $\X$.  In this subsection, we provide a construction that produces this unitary category from any $\dagger$-isomix category.  This construction, which we call the {\bf unitary construction} provides an important technique for building unitary categories.   The construction is based on identifying objects with pre-unitary structure: the tensor units always have a canonical pre-unitary structure so the construction always produces a non-empty category.   However, to ensure that an application of the construction yields a unitary category in which there are objects which are not isomorphic to the units, one must exhibit concretely such objects.  Fortunately this is often not difficult to do, making the construction quite applicable.

\begin{defi} ~
\label{defn: pre-unitary}
\begin{enumerate}[(i)] 
\item In a $\dagger$-isomix category, a {\bf pre-unitary object} is an object $U \in \Core(\X)$, together with an isomorphism $\alpha: U \to U^\dagger$ such that  $\alpha (\alpha^{-1})^\dagger = \iota$. 

\item Suppose $\X$ is a $\dagger$-isomix category, then define ${\sf Unitary}(\X)$, the {\bf unitary core} of $\X$, as follows:
\begin{description}
\item[Objects] Pre-unitary objects $(U, \alpha)$,
\item[Maps] $(U, \alpha) \to^f (V, \beta)$ where $U \to^f V $ is any map of $\X$.
\end{description}
\end{enumerate}
\end{defi}

We note that any object which is isomorphic to a preunitary object is also pre-unitary:
\begin{lem}
In a $\dagger$-isomix category, if $U$ is a pre-unitary object 
and there exists an isomorphism $f: U \to U'$, then $U'$ is pre-unitary. 
\end{lem}

Our objective is to show that Unitary($\X$) is endowed with all the structure of a unitary category.

\begin{lem}
For any $\dagger$-isomix category, its canonical unitary core is a compact $\dagger$-LDC with tensor and par defined by
\[ (\top,{\sf m}^{-1}\lambda_\bot: \top \to \top^\dagger) ~~~~~(A, \alpha) \ox (B, \beta) := (A \ox B, \mx(\alpha \oa \beta) \lambda_\oa: A \ox B \to (A \ox B)^\dagger)\]
\[ (\bot, {\sf m} ~\lambda_\top: \bot \to \bot^\dagger)  ~~~~~(A, \alpha) \oa (B, \beta) := (A \oa B, \mx^{-1}(\alpha \ox \beta) \lambda_\ox: A \oa B \to (A \oa B)^\dagger) \]
and $(U,\alpha)^\dagger := (U^\dagger, (\alpha^{-1})^\dagger)$.
\end{lem}

\begin{proof}
The proof uses the techniques of Lemma \ref{Lemma: square root tensor unitary}. 

 Note that, as the map and tensor structure is inherited from $\X$, it suffices to show that these objects are all pre-unitary objects.  Starting with $(U \alpha)^\dagger$ 
we have:
\[ (\alpha^{-1})^\dagger (((\alpha^{-1})^\dagger)^{-1})^\dagger = (\alpha^{-1})^\dagger (\alpha^\dagger)^\dagger = (\alpha^\dagger \alpha^{-1})^\dagger = (\iota^{-1})^\dagger = \iota \]
For the tensor and par we have:
\begin{align*}
{\sf m}^{-1}\lambda_\bot (({\sf m}^{-1}\lambda_\bot)^{-1})^\dagger & =  {\sf m}^{-1}\lambda_\bot {\sf m}^\dagger \lambda_\bot^\dagger \\
& \stackrel{\text{\tiny {\bf [$\dagger$-mix]}}}{=}  {\sf m}^{-1} {\sf m} \lambda_\top  \lambda_\bot^\dagger  = \iota \\
\mx^{-1}(\alpha \oa \beta) \lambda_\oa ((\mx^{-1}(\alpha \oa \beta) \lambda_\oa)^{-1})^\dagger 
& =  \mx^{-1}(\alpha \oa \beta) \lambda_\oa (\mx^\dagger) (\alpha^{-1} \oa \beta^{-1})^\dagger (\lambda_\oa^{-1})^\dagger \\
& =  \mx^{-1}(\alpha \oa \beta) \mx \lambda_\ox (\alpha^{-1} \oa \beta^{-1})^\dagger (\lambda_\oa^{-1})^\dagger \\
& =  (\alpha \ox \beta) \lambda_\ox (\alpha^{-1} \oa \beta^{-1})^\dagger (\lambda_\oa^{-1})^\dagger \\
& =  (\alpha \ox \beta) ((\alpha^{-1})^\dagger \ox (\beta^{-1})^\dagger) \lambda_\ox (\lambda_\oa^{-1})^\dagger \\
&\stackrel{\text{\tiny {\bf Defn \ref{defn: pre-unitary}-(i)}}}{=} (\iota \ox \iota) \lambda_\ox (\lambda_\oa^{-1})^\dagger \\
&\stackrel{\text{\tiny {\bf [$\dagger$-ldc.4]}}}{=} \iota 
\end{align*}
\begin{align*}
{\sf m} \lambda_\top (({\sf m} \lambda_\top)^{-1})^\dagger & =  {\sf m} \lambda_\top ({\sf m}^{-1})^\dagger  (\lambda_\top^{-1})^\dagger  \\
& =  {\sf m} ~{\sf m}^{-1} \lambda_\bot (\lambda_\top^{-1})^\dagger = \iota \\
\mx(\alpha \ox \beta) \lambda_\ox ((\mx(\alpha \ox \beta) \lambda_\ox)^{-1})^\dagger 
& =  \mx(\alpha \ox \beta) \lambda_\ox (\mx^{-1})^\dagger (\alpha^{-1} \ox \beta^{-1})^\dagger (\lambda_\ox^{-1})^\dagger \\
& =  (\alpha \oa \beta) \mx~ \mx^{-1} \lambda_\oa (\alpha^{-1} \ox \beta^{-1})^\dagger (\lambda_\ox^{-1})^\dagger \\
&=  (\alpha \oa \beta) ((\alpha^{-1})^\dagger \oa (\beta^{-1})^\dagger) \lambda_\oa (\lambda_\ox^{-1})^\dagger \\
& =  (\iota \oa \iota) \lambda_\oa (\lambda_\ox^{-1})^\dagger  = \iota  \qedhere
\end{align*}
\end{proof}

This makes $\Unitary(\X)$ into a compact $\dagger$-LDC with all the structure inherited directly from $\X$. 
However, more is true: each object now has an obvious unitary structure.  This gives:

\begin{prop}
For any $\dagger$-isomix category, $\X$, $\Unitary(\X)$ is a unitary category with a full and faithful underlying $\dagger$-isomix functor $U: {\sf Unitary}(\X) \to \X$.
\end{prop}

\begin{proof}
The laxors are all identity maps so that the underlying functors is immediately a $\dagger$-mix functor.

It remains to show that every object is unitary:  we set the unitary structure of an object to be $\alpha: (X,\alpha) \to (X,\alpha)^\dagger$.   However, {\bf [U.1]} -- {\bf [U.5]} are immediately satisfied by construction implying this provides unitary structure for every object.
\end{proof}

%%%%%%%%%%%%%%%%%%%%%%%%%%%%%%%%%%%%%%%%%%%%%%%%%%%%
% Couniversality of universal construction %
Next, we prove the couniversal property of the unitary construction. Define ${\sf UCat}$ to be the category of unitary categories and $\dagger$-isomix functors that preserve   unitary structure in the sense of Definition \ref{preserving-unitary-structure}, thus, whenever ${\varphi_A}$ is the unitary structure  then $F'(\varphi_A) \rho^{F'}$ is unitary structure. Define {\sf Kompact} to be the category of compact $\dagger$-LDCs and $\dagger$-isomix functors.

We now show that the unitary construction produces a right adjoint to the underlying functor $U: {\sf UCat} \to {\sf Kompact}$ which is the identity functor. Preliminary to this result we prove that Frobenius functors preserve preunitary objects:

\begin{lem}
\label{Lemma: Frobenius preunitary}
If $F: \X \to \Y$ is a $\dagger$-isomix functor between compact $\dagger$-LDCs and $(A,\varphi)$ is a preunitary object of $\X$, then $(F(A),F(\varphi)\rho)$ is a preunitary object of $\Y$.
\end{lem}
\begin{proof}
To prove that $(F(A),F(\varphi)\rho)$ is a preunitary object, one has the following computation:
\begin{align*}
F(\varphi)\rho ((F(\varphi) \rho)^{-1})^\dagger 
& =  F(\varphi)\rho F(\varphi^{-1})^\dagger (\rho^{-1})^\dagger \\
& =  F(\varphi (\varphi^{-1})^\dagger) \rho (\rho^{-1})^\dagger \\
& =  F(\iota) \rho (\rho^{-1})^\dagger \stackrel{{\bf [\dagger-isomix]}}{=} \iota \qedhere
\end{align*}
\end{proof}

\begin{prop}
\label{Prop: Couniversal}
$U: {\sf UCat} \to {\sf Kompact}$ has a right adjoint ${\sf Unitary}: {\sf Kompact} \to {\sf UCat}; \C \mapsto {\sf Unitary}(\C)$.
\end{prop}
\begin{proof}
The couniversal diagram is as follows:
\[ \xymatrix{ U(\U) \ar[rr]^{F} \ar@{.>}[d]_{U(F^\flat)} && \C \\
                    U({\sf Unitary}(\C)) \ar[urr]_{\epsilon}} \]

Since $F$ is a $\dagger$-isomix functor it preserves preunitary structure (see Lemma \ref{Lemma: Frobenius preunitary}).  This means that each $(U,\varphi_U)$ in $\U$ is carried by $F$ onto a preunitary object in $\C$, $(F(U),F(\varphi)\rho^F)$.  But a preunitary object in $\C$ is an object of ${\sf Unitary}(\C)$ and this determines $F^\flat$.   The functor $F^\flat$ is uniquely determined as it must preserve the unitary structure.
\end{proof}
%%%%%%%%%%%%%%%%%%%%%%%%%%%%%%%%%%%%%%%%%%%%%%%%%%%%%

\subsection{Examples of unitary constructions}

In Section~\ref{daggers-duals-conjugation}, we discussed examples of $\dagger$-isomix categories in which the $\dagger$ is given by composing the conjugation functor and the dualizing functor. In the rest of the section, we apply the unitary construction to each of those examples to construct a unitary category:

\subsubsection{Category of abstract state spaces}

In Section~\ref{Sec: Asp}, we discussed a construction on a $\dagger$-isomix category, $\X$, that produces a category of abstract state spaces, $\Asp(\X)$, which is a $\dagger$-isomix category. In this section, we examine the preunitary objects of $\Asp(\X)$. Since all the basic natural isomorphisms are inherited from $\X$, $\Core(\X)$ determine $\Core(\Asp(\X))$. If $(A ,\alpha)$ is a preunitary object for $\X$, and $(A, e_A, u_A) \in \Asp(\X)$ then, $((A, e_A, u_A), \alpha)$ is a preunitary object for $\Asp(\X)$ if $u_A \alpha = \lambda_\top e_A^\dagger$.

\subsubsection{Category of a group with involution}

We discussed a source of examples of compact $\dagger$-LDCs which are given by groups with conjugation. Applying unitary construction to each of the example categories results in the following unitary categories. It could be noticed that the preunitary objects in each of these categories includes those group elements such that $\overline{g^{-1}} = g$. More explicitly, the preunitary objects are $(g,1)$ such that $\overline{g^{-1}} = g$.

\begin{itemize}
\item In the discrete category of complex numbers, $\D(\C, +, 0)$, \[(a + ib)^\dagger := \overline{(a+ib)^*} = \overline{(-a-ib)} = -a + ib\] The preunitary objects in this category are given by all complex numbers, i.e., $(ib, 1)$. 

\item In the discrete category of non-zero complex numbers, $\D(\C, ., 1)$, the preunitary objects are given by complex numbers on a unit circle.

\item In the discrete category, $\D(P(x), +, 0)$, where $P(x)$ is a polynomial ring, $P(x)^\dagger = -P(-x)$ and the preunitary objects are polynomials $ P(x) = \sum_n a_n x^n$ such that n is odd. 

\item In $\D(\mathbb{M}_2, \cdot, I_2)$ where $\mathbb{M}_2$ is the group of $2 \times 2$ invertible matrices over $\mathbb{C}$. The $\dagger$ structure is as follows:
 \[\left(
  \begin{matrix}
 a+ib & m+in \\
 c+id & p+iq 
 \end{matrix}
 \right)^\dagger := \overline{\left(
 \begin{matrix}
 a+ib & m+in \\
 c+id & p+iq 
 \end{matrix}
 \right)^*} = \left(
 \begin{matrix}
a-ib & c-id \\
m-in & p-iq
 \end{matrix}
 \right)^{-1} \]
 The preunitary objects in this category are given by unitary matrices.
\end{itemize}
%%%%%%%%%%%%%%%%%%%%%%%%%%%%%%%%%%%%%%%%%%%%%%%

\subsubsection{Category of Hopf Modules in a $*$-autonomous category}

In Section~\ref{Sec: HModx}, we described a construction of $\dagger$-isomix categories  from any  symmetric isomix $*$-autonomous category, $\X$, by choosing the Hopf Modules over a  cocommutative $\ox$-Hopf Algebra. We referred to the resulting category as ${\mbox{\bf H-Mod}}_\X$. Now we shall look at the preunitary objects in ${\mbox{\bf H-Mod}}_\X$ in order to apply the unitary construction to this category. We begin by identifying the objects in the core of ${\mbox{\bf H-Mod}}_\X$:

%We know that in any isomix category, $\X$, the core of the category, $\Core({\mbox{\bf H-Mod}}_\X)$ is a compact $\dagger$-isomix category. If the core of ${\mbox{\bf H-Mod}}_\X$ is non-trivial, i.e., the core includes objects other than tensor units too, then, one can apply unitary construction to the core to get a MUC. ${\mbox{\bf H-Mod}}_\X$ being a $\dagger$-isomix category, this technique can be applied to the category to get a MUC. In order to do so, we identify the preunitary objects in $\Core({\mbox{\bf H-Mod}}_\X)$. We begin by identifying the objects in the core of ${\mbox{\bf H-Mod}}_\X$:

\begin{lem}
Suppose $\X$ is a mix $*$-autonomous category and $H$ is a cocommutative Hopf Algebra in $\X$. If $(A, \leftaction{0.4}{white})$ is a H-Module and $A \in \Core(\X)$, then $ (A, \leftaction{0.4}{white}) \in \Core({\mbox{\bf H-Mod}}_\X)$.
\end{lem}
\begin{proof}
The mixor $\mx: A \ox B \to A \oa B$ is inherited directly from $\X$. Hence,  $ (A, \leftaction{0.4}{white}) \in \Core(${\bf H-Mod$_\X)$}.
\end{proof}

Now that we identified the objects in the core, we prove a lemma that will be used later to identify the preunitary objects from the core:

\begin{lem}
\label{Lemma: aux} The following equality holds for a Frobenius algebra:
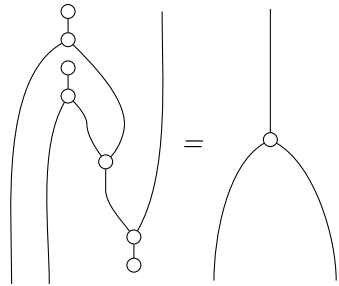
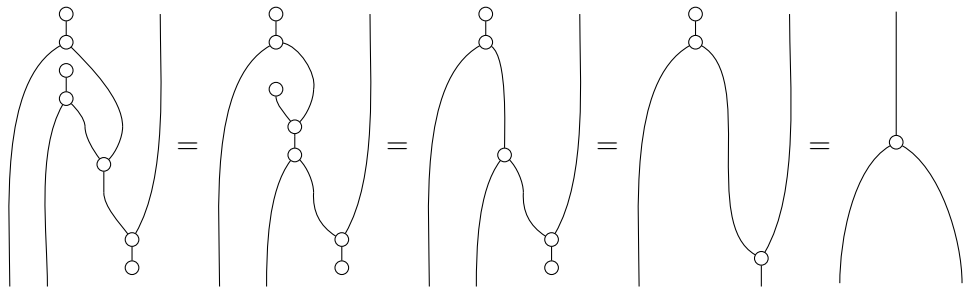
\[ \begin{tikzpicture} %act12
	\begin{pgfonlayer}{nodelayer}
		\node [style=circle] (0) at (-2.5, -0) {};
		\node [style=none] (1) at (-2, 1) {};
		\node [style=none] (2) at (-2.5, -0.75) {};
		\node [style=none] (3) at (-3, 1) {};
		\node [style=none] (4) at (-2.5, -0.75) {};
		\node [style=none] (5) at (-1, 4) {};
		\node [style=circle] (6) at (-1.75, -2.75) {};
		\node [style=circle] (7) at (-1.75, -2) {};
		\node [style=none] (8) at (-3, 1) {};
		\node [style=none] (9) at (-4, -3.25) {};
		\node [style=circle] (10) at (-3.5, 2.5) {};
		\node [style=circle] (11) at (-3.5, 1.75) {};
		\node [style=circle] (12) at (-3.5, 3.25) {};
		\node [style=none] (13) at (-2, 1) {};
		\node [style=none] (14) at (-5, -3.25) {};
		\node [style=circle] (15) at (-3.5, 4) {};
	\end{pgfonlayer}
	\begin{pgfonlayer}{edgelayer}
		\draw [in=-90, out=45, looseness=0.75] (0) to (1.center);
		\draw (0) to (2.center);
		\draw [in=127, out=-90, looseness=0.75] (3.center) to (0);
		\draw [in=-90, out=60, looseness=0.75] (7) to (5.center);
		\draw (7) to (6);
		\draw [in=127, out=-90, looseness=0.75] (4.center) to (7);
		\draw [in=90, out=-45, looseness=0.75] (11) to (8.center);
		\draw (11) to (10);
		\draw [in=-120, out=90, looseness=0.75] (9.center) to (11);
		\draw [in=90, out=-45, looseness=0.75] (12) to (13.center);
		\draw (12) to (15);
		\draw [in=-150, out=90, looseness=0.75] (14.center) to (12);
	\end{pgfonlayer}
\end{tikzpicture}
 = \begin{tikzpicture} %act16
	\begin{pgfonlayer}{nodelayer}
		\node [style=none] (0) at (-1.75, -3.25) {};
		\node [style=circle] (1) at (-3.5, 0.5) {};
		\node [style=none] (2) at (-5, -3.25) {};
		\node [style=none] (3) at (-3.5, 4) {};
	\end{pgfonlayer}
	\begin{pgfonlayer}{edgelayer}
		\draw (1) to (3);
		\draw [in=-150, out=90, looseness=0.75] (2.center) to (1);
		\draw [in=90, out=-30, looseness=0.75] (1) to (0.center);
	\end{pgfonlayer}
\end{tikzpicture}\]
\end{lem}
\begin{proof}
\[\begin{tikzpicture} %act12
	\begin{pgfonlayer}{nodelayer}
		\node [style=circle] (0) at (-2.5, -0) {};
		\node [style=none] (1) at (-2, 1) {};
		\node [style=none] (2) at (-2.5, -0.75) {};
		\node [style=none] (3) at (-3, 1) {};
		\node [style=none] (4) at (-2.5, -0.75) {};
		\node [style=none] (5) at (-1, 4) {};
		\node [style=circle] (6) at (-1.75, -2.75) {};
		\node [style=circle] (7) at (-1.75, -2) {};
		\node [style=none] (8) at (-3, 1) {};
		\node [style=none] (9) at (-4, -3.25) {};
		\node [style=circle] (10) at (-3.5, 2.5) {};
		\node [style=circle] (11) at (-3.5, 1.75) {};
		\node [style=circle] (12) at (-3.5, 3.25) {};
		\node [style=none] (13) at (-2, 1) {};
		\node [style=none] (14) at (-5, -3.25) {};
		\node [style=circle] (15) at (-3.5, 4) {};
	\end{pgfonlayer}
	\begin{pgfonlayer}{edgelayer}
		\draw [in=-90, out=45, looseness=0.75] (0) to (1.center);
		\draw (0) to (2.center);
		\draw [in=127, out=-90, looseness=0.75] (3.center) to (0);
		\draw [in=-90, out=60, looseness=0.75] (7) to (5.center);
		\draw (7) to (6);
		\draw [in=127, out=-90, looseness=0.75] (4.center) to (7);
		\draw [in=90, out=-45, looseness=0.75] (11) to (8.center);
		\draw (11) to (10);
		\draw [in=-120, out=90, looseness=0.75] (9.center) to (11);
		\draw [in=90, out=-45, looseness=0.75] (12) to (13.center);
		\draw (12) to (15);
		\draw [in=-150, out=90, looseness=0.75] (14.center) to (12);
	\end{pgfonlayer}
\end{tikzpicture} = \begin{tikzpicture} %act13
	\begin{pgfonlayer}{nodelayer}
		\node [style=none] (0) at (-1, 4) {};
		\node [style=circle] (1) at (-1.75, -2.75) {};
		\node [style=circle] (2) at (-1.75, -2) {};
		\node [style=circle] (3) at (-3.5, 3.25) {};
		\node [style=none] (4) at (-5, -3.25) {};
		\node [style=circle] (5) at (-3.5, 4) {};
		\node [style=circle] (6) at (-3.5, 2) {};
		\node [style=none] (7) at (-2.5, 2) {};
		\node [style=circle] (8) at (-3, 1) {};
		\node [style=none] (9) at (-2.5, 2) {};
		\node [style=none] (10) at (-3, 0.75) {};
		\node [style=none] (11) at (-3.75, -3.25) {};
		\node [style=none] (12) at (-2.5, -0.75) {};
		\node [style=circle] (13) at (-3, 0.25) {};
		\node [style=none] (14) at (-2.5, -0.75) {};
		\node [style=none] (15) at (-3, 0.75) {};
	\end{pgfonlayer}
	\begin{pgfonlayer}{edgelayer}
		\draw [in=-90, out=60, looseness=0.75] (2) to (0.center);
		\draw (2) to (1);
		\draw (3) to (5);
		\draw [in=-150, out=90, looseness=0.75] (4.center) to (3);
		\draw [in=-90, out=45, looseness=0.75] (8) to (9.center);
		\draw (8) to (10.center);
		\draw [in=127, out=-90, looseness=0.75] (6) to (8);
		\draw [in=90, out=-45, looseness=0.75] (13) to (14.center);
		\draw (13) to (15.center);
		\draw [in=-127, out=90, looseness=0.75] (11.center) to (13);
		\draw [bend right, looseness=1.00] (12.center) to (2);
		\draw [bend left=45, looseness=0.75] (3) to (7.center);
	\end{pgfonlayer}
\end{tikzpicture} =
\begin{tikzpicture} %act14
	\begin{pgfonlayer}{nodelayer}
		\node [style=none] (0) at (-1, 4) {};
		\node [style=circle] (1) at (-1.75, -2.75) {};
		\node [style=circle] (2) at (-1.75, -2) {};
		\node [style=circle] (3) at (-3.5, 3.25) {};
		\node [style=none] (4) at (-5, -3.25) {};
		\node [style=circle] (5) at (-3.5, 4) {};
		\node [style=none] (6) at (-3.75, -3.25) {};
		\node [style=none] (7) at (-2.5, -0.75) {};
		\node [style=circle] (8) at (-3, 0.25) {};
		\node [style=none] (9) at (-2.5, -0.75) {};
	\end{pgfonlayer}
	\begin{pgfonlayer}{edgelayer}
		\draw [in=-90, out=60, looseness=0.75] (2) to (0.center);
		\draw (2) to (1);
		\draw (3) to (5);
		\draw [in=-150, out=90, looseness=0.75] (4.center) to (3);
		\draw [in=90, out=-45, looseness=0.75] (8) to (9.center);
		\draw [in=-127, out=90, looseness=0.75] (6.center) to (8);
		\draw [bend right, looseness=1.00] (7.center) to (2);
		\draw [in=90, out=-36, looseness=0.50] (3) to (8);
	\end{pgfonlayer}
\end{tikzpicture} = 
\begin{tikzpicture} %act15
	\begin{pgfonlayer}{nodelayer}
		\node [style=none] (0) at (-1, 4) {};
		\node [style=none] (1) at (-1.75, -3.25) {};
		\node [style=circle] (2) at (-1.75, -2.5) {};
		\node [style=circle] (3) at (-3.5, 3.25) {};
		\node [style=none] (4) at (-5, -3.25) {};
		\node [style=circle] (5) at (-3.5, 4) {};
	\end{pgfonlayer}
	\begin{pgfonlayer}{edgelayer}
		\draw [in=-90, out=60, looseness=0.75] (2) to (0.center);
		\draw (2) to (1.center);
		\draw (3) to (5);
		\draw [in=-150, out=90, looseness=0.75] (4.center) to (3);
		\draw [in=150, out=-30, looseness=0.75] (3) to (2);
	\end{pgfonlayer}
\end{tikzpicture} = 
\begin{tikzpicture} %act16
	\begin{pgfonlayer}{nodelayer}
		\node [style=none] (0) at (-1.75, -3.25) {};
		\node [style=circle] (1) at (-3.5, 0.5) {};
		\node [style=none] (2) at (-5, -3.25) {};
		\node [style=none] (3) at (-3.5, 4) {};
	\end{pgfonlayer}
	\begin{pgfonlayer}{edgelayer}
		\draw (1) to (3);
		\draw [in=-150, out=90, looseness=0.75] (2.center) to (1);
		\draw [in=90, out=-30, looseness=0.75] (1) to (0.center);
	\end{pgfonlayer}
\end{tikzpicture}\]
\end{proof}

In the following Proposition we identify the preunitary objects in the core:

\begin{prop} 
Suppose $\X$ is a symmetric mix $*$-autonomous category and $H$ is a cocommutative Hopf Algebra in $\X$. If $A \in \Core(\X)$ and $(A, \mulmap{1.2}{white}, \unitmap{1.2}{white}, \comulmap{1.2}{white}, \counitmap{1.2}{white})$ is a cocommutative Frobenius Algebra with an algebra homomorphism $H \to^{h} A$ then, 
\begin{enumerate}[(a)]
\item $(A, \leftaction{0.4}{white})$ is a H-Module where, $\leftaction{0.4}{white}: H \ox A \to A := \begin{tikzpicture}
	\begin{pgfonlayer}{nodelayer}
		\node [style=circle] (0) at (0, -0) {};
		\node [style=none] (1) at (-0.5, 1.75) {};
		\node [style=none] (2) at (0.75, 1.75) {};
		\node [style=none] (3) at (0, -0.75) {};
		\node [style=circle, scale=2] (4) at (-0.5, 1) {};
		\node [style=none] (5) at (-0.5, 1) {$h$};
	\end{pgfonlayer}
	\begin{pgfonlayer}{edgelayer}
		\draw (1.center) to (4);
		\draw [bend left, looseness=1.00] (0) to (4);
		\draw [in=-90, out=15, looseness=1.00] (0) to (2.center);
		\draw (0) to (3.center);
	\end{pgfonlayer}
\end{tikzpicture}$ 
 \item $\overline{(A, \leftaction{0.4}{white})^*} = (A, \leftaction{0.4}{white})$ where $A^*$ is the self-dual Frobenius Algebra with cups and caps defined as
$
\begin{tikzpicture}
	\begin{pgfonlayer}{nodelayer}
		\node [style=circle] (0) at (0, -0) {};
		\node [style=none] (1) at (-0.5, 1) {};
		\node [style=none] (2) at (0.75, 1) {};
		\node [style=circle] (3) at (0, -1) {};
		\node [style=circle] (4) at (0, -1) {};
	\end{pgfonlayer}
	\begin{pgfonlayer}{edgelayer}
		\draw [in=-90, out=15, looseness=1.00] (0) to (2.center);
		\draw [in=150, out=-90, looseness=1.00] (1.center) to (0);
		\draw (0) to (3);
	\end{pgfonlayer}
\end{tikzpicture}  and  \begin{tikzpicture}
	\begin{pgfonlayer}{nodelayer}
		\node [style=circle] (0) at (0, 0) {};
		\node [style=none] (1) at (-0.5, -1) {};
		\node [style=none] (2) at (0.75, -1) {};
		\node [style=circle] (3) at (0, 1) {};
		\node [style=circle] (4) at (0, 1) {};
	\end{pgfonlayer}
	\begin{pgfonlayer}{edgelayer}
		\draw [in=90, out=-15, looseness=1.00] (0) to (2.center);
		\draw [in=-150, out=90, looseness=1.00] (1.center) to (0);
		\draw (0) to (3);
	\end{pgfonlayer}
\end{tikzpicture}
$ respectively. Hence, $A^* = A$ and $(A, \leftaction{0.4}{white})^\dagger =  (A, \leftaction{0.4}{white})$.
\end{enumerate}
\end{prop}
\begin{proof}~
\begin{enumerate}[(a)]
\item $\begin{tikzpicture} %act7
	\begin{pgfonlayer}{nodelayer}
		\node [style=circle] (0) at (0, -0) {};
		\node [style=none] (1) at (0.75, 1.25) {};
		\node [style=none] (2) at (0, -0.75) {};
		\node [style=circle, scale=2] (3) at (-0.5, 0.5) {};
		\node [style=none] (4) at (-0.5, 0.5) {$h$};
		\node [style=none] (5) at (-0.5, 1.25) {};
	\end{pgfonlayer}
	\begin{pgfonlayer}{edgelayer}
		\draw [bend left, looseness=1.00] (0) to (3);
		\draw [in=-90, out=15, looseness=1.00] (0) to (1.center);
		\draw (0) to (2.center);
		\draw (5.center) to (3);
	\end{pgfonlayer}
\end{tikzpicture}: H \ox A \to A$ is a left action because $h: H \to A$ is an algebra homomorphism.
\item $
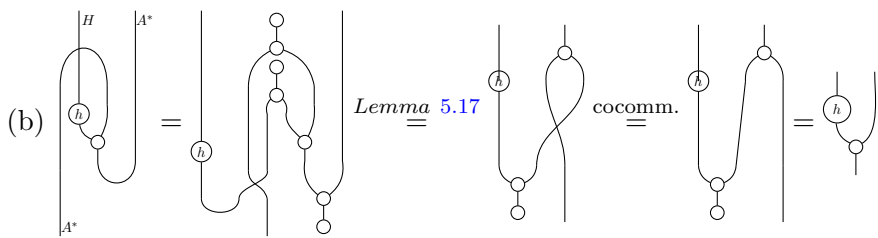
\begin{tikzpicture} %Frob0
	\begin{pgfonlayer}{nodelayer}
		\node [style=none] (0) at (1.25, -1) {};
		\node [style=none] (1) at (2.25, -1) {};
		\node [style=none] (2) at (1.5, 0.5) {};
		\node [style=none] (3) at (1.5, 1) {};
		\node [style=none] (4) at (0.25, 1) {};
		\node [style=none] (5) at (2.25, 3) {};
		\node [style=none] (6) at (0.25, -1.25) {};
		\node [style=none] (7) at (1, 2.75) {$H$};
		\node [style=none] (8) at (2.5, 2.75) {$A^*$};
		\node [style=none] (9) at (0.25, -3) {};
		\node [style=none] (10) at (0.5, -2.75) {$A^*$};
		\node [style=circle, scale=1.5] (11) at (0.75, 0.25) {};
		\node [style=none] (12) at (0.75, 3) {};
		\node [style=circle] (13) at (1.25, -0.5) {};
		\node [style=none] (14) at (0.75, 0.25) {$h$};
	\end{pgfonlayer}
	\begin{pgfonlayer}{edgelayer}
		\draw [in=-90, out=90, looseness=1.00] (2.center) to (3.center);
		\draw [bend left=90, looseness=2.75] (4.center) to (3.center);
		\draw [bend right=90, looseness=2.00] (0.center) to (1.center);
		\draw (5.center) to (1.center);
		\draw (4.center) to (6.center);
		\draw [bend right, looseness=1.00] (11) to (13);
		\draw [bend right=15, looseness=1.00] (13) to (2.center);
		\draw (13) to (0.center);
		\draw (12.center) to (11);
		\draw (6.center) to (9.center);
	\end{pgfonlayer}
\end{tikzpicture} =
\begin{tikzpicture} %Frob2
	\begin{pgfonlayer}{nodelayer}
		\node [style=none] (0) at (1.25, -1) {};
		\node [style=none] (1) at (2.25, -1) {};
		\node [style=none] (2) at (1.5, 0.5) {};
		\node [style=none] (3) at (0.25, 0.25) {};
		\node [style=none] (4) at (1.5, 1) {};
		\node [style=none] (5) at (-0.25, 1) {};
		\node [style=none] (6) at (2.25, 3) {};
		\node [style=none] (7) at (0.25, -1.25) {};
		\node [style=none] (8) at (-0.25, -1.25) {};
		\node [style=none] (9) at (-0.5, -2) {};
		\node [style=none] (10) at (0.25, -2) {};
		\node [style=none] (11) at (-1.5, -2) {};
		\node [style=none] (12) at (-1.5, 3) {};
		\node [style=none] (13) at (0.25, -3) {};
		\node [style=none] (14) at (0.75, 0.25) {};
		\node [style=circle] (15) at (1.25, -0.5) {};
		\node [style=circle, scale=1.5] (16) at (-1.5, -0.75) {};
		\node [style=none] (17) at (-0.25, 1) {};
		\node [style=none] (18) at (1.5, 1) {};
		\node [style=circle] (19) at (0.5, 2.75) {};
		\node [style=circle] (20) at (0.5, 2) {};
		\node [style=none] (21) at (0.25, 0.25) {};
		\node [style=circle] (22) at (0.5, 1.5) {};
		\node [style=none] (23) at (0.75, 0.25) {};
		\node [style=circle] (24) at (0.5, 0.75) {};
		\node [style=none] (25) at (1.25, -1) {};
		\node [style=circle] (26) at (1.75, -2.75) {};
		\node [style=none] (27) at (2.25, -1) {};
		\node [style=circle] (28) at (1.75, -2) {};
		\node [style=none] (29) at (-1.5, -0.75) {$h$};
	\end{pgfonlayer}
	\begin{pgfonlayer}{edgelayer}
		\draw [in=-90, out=90, looseness=1.00] (2.center) to (4.center);
		\draw (6.center) to (1.center);
		\draw (3.center) to (7.center);
		\draw (5.center) to (8.center);
		\draw [in=105, out=-90, looseness=1.25] (8.center) to (10.center);
		\draw [in=90, out=-75, looseness=0.75] (7.center) to (9.center);
		\draw [bend right=90, looseness=1.25] (11.center) to (9.center);
		\draw (13.center) to (10.center);
		\draw [bend right=15, looseness=1.00] (15) to (2.center);
		\draw (15) to (0.center);
		\draw (16) to (11.center);
		\draw [in=-90, out=135, looseness=1.00] (15) to (14.center);
		\draw [bend right, looseness=1.00] (20) to (17.center);
		\draw [bend left, looseness=1.00] (20) to (18.center);
		\draw (19) to (20);
		\draw [bend right, looseness=1.00] (24) to (21.center);
		\draw [bend left, looseness=1.00] (24) to (23.center);
		\draw (22) to (24);
		\draw [bend left, looseness=1.00] (28) to (25.center);
		\draw [bend right, looseness=1.00] (28) to (27.center);
		\draw (26) to (28);
		\draw (12.center) to (16);
	\end{pgfonlayer}
\end{tikzpicture} \stackrel{Lemma ~ \ref{Lemma: aux}}{=} 
\begin{tikzpicture} %Frob3
	\begin{pgfonlayer}{nodelayer}
		\node [style=none] (0) at (0.75, 0.25) {};
		\node [style=none] (1) at (-0.25, 0.25) {};
		\node [style=none] (2) at (0.25, -2) {};
		\node [style=none] (3) at (-1.5, -2) {};
		\node [style=none] (4) at (-1.5, 1.75) {};
		\node [style=none] (5) at (0.25, -3.5) {};
		\node [style=circle, scale=1.5] (6) at (-1.5, 0.25) {};
		\node [style=none] (7) at (-1.5, 0.25) {$h$};
		\node [style=none] (8) at (0.75, 0.25) {};
		\node [style=none] (9) at (0.25, 1.75) {};
		\node [style=circle] (10) at (0.25, 1) {};
		\node [style=none] (11) at (-0.25, 0.25) {};
		\node [style=none] (12) at (-1.5, -2) {};
		\node [style=none] (13) at (-0.5, -2) {};
		\node [style=circle] (14) at (-1, -3.25) {};
		\node [style=circle] (15) at (-1, -2.5) {};
		\node [style=none] (16) at (-1.5, -2) {};
	\end{pgfonlayer}
	\begin{pgfonlayer}{edgelayer}
		\draw (5.center) to (2.center);
		\draw (6) to (3.center);
		\draw (4.center) to (6);
		\draw (7.center) to (6);
		\draw [bend right, looseness=1.00] (10) to (11.center);
		\draw [bend left, looseness=1.00] (10) to (8.center);
		\draw (9.center) to (10);
		\draw [bend left, looseness=1.00] (15) to (16.center);
		\draw [bend right, looseness=1.00] (15) to (13.center);
		\draw (14) to (15);
		\draw [in=90, out=-90, looseness=1.00] (1.center) to (2.center);
		\draw [in=90, out=-90, looseness=1.00] (0.center) to (13.center);
	\end{pgfonlayer}
\end{tikzpicture} \stackrel{\text{cocomm.}}{=} 
\begin{tikzpicture} %Frob4
	\begin{pgfonlayer}{nodelayer}
		\node [style=none] (0) at (0.75, 0.25) {};
		\node [style=none] (1) at (-0.25, 0.25) {};
		\node [style=none] (2) at (0.75, -2) {};
		\node [style=none] (3) at (-1.5, -2) {};
		\node [style=none] (4) at (-1.5, 1.75) {};
		\node [style=none] (5) at (0.75, -3.5) {};
		\node [style=circle, scale=1.5] (6) at (-1.5, 0.25) {};
		\node [style=none] (7) at (-1.5, 0.25) {$h$};
		\node [style=none] (8) at (0.75, 0.25) {};
		\node [style=none] (9) at (0.25, 1.75) {};
		\node [style=circle] (10) at (0.25, 1) {};
		\node [style=none] (11) at (-0.25, 0.25) {};
		\node [style=none] (12) at (-1.5, -2) {};
		\node [style=none] (13) at (-0.5, -2) {};
		\node [style=circle] (14) at (-1, -3.25) {};
		\node [style=circle] (15) at (-1, -2.5) {};
		\node [style=none] (16) at (-1.5, -2) {};
	\end{pgfonlayer}
	\begin{pgfonlayer}{edgelayer}
		\draw (5.center) to (2.center);
		\draw (6) to (3.center);
		\draw (4.center) to (6);
		\draw (7.center) to (6);
		\draw [bend right, looseness=1.00] (10) to (11.center);
		\draw [bend left, looseness=1.00] (10) to (8.center);
		\draw (9.center) to (10);
		\draw [bend left, looseness=1.00] (15) to (16.center);
		\draw [bend right, looseness=1.00] (15) to (13.center);
		\draw (14) to (15);
		\draw (1.center) to (13.center);
		\draw (0.center) to (2.center);
	\end{pgfonlayer}
\end{tikzpicture} =
\begin{tikzpicture} %act11
	\begin{pgfonlayer}{nodelayer}
		\node [style=circle] (0) at (0, -0) {};
		\node [style=none] (1) at (0.5, 2) {};
		\node [style=none] (2) at (0, -0.75) {};
		\node [style=circle, scale=2] (3) at (-0.5, 1) {};
		\node [style=none] (4) at (-0.5, 1) {$h$};
		\node [style=none] (5) at (-0.5, 2) {};
	\end{pgfonlayer}
	\begin{pgfonlayer}{edgelayer}
		\draw [bend left, looseness=1.00] (0) to (3);
		\draw [in=-90, out=30, looseness=0.75] (0) to (1.center);
		\draw (0) to (2.center);
		\draw (3) to (5.center);
	\end{pgfonlayer}
\end{tikzpicture}$
\end{enumerate}
\end{proof}

\begin{cor}
$(((A, \mulmap{1.2}{white}, \unitmap{1.2}{white}, \comulmap{1.2}{white}, \counitmap{1.2}{white}), \leftaction{0.4}{white}), 1)$ is a preunitary object.
\end{cor}

Thus, we have a source of non-trivial preunitary objects so that we can form a non-trivial unitary category.

\subsection{Mixed unitary categories}

We are now ready for the definition of mixed unitary categories, which is the key structure developed in this paper.

\begin{defi}
A {\bf mixed unitary category} (MUC) is a $\dagger$-isomix category, $\C$, equipped with a strong $\dagger$-isomix functor 
$M: \U \to \C$ from a unitary category $\U$ to $\C$ such that there exists the following natural transformations:
\[ \mx': M(U) \oa X \to M(U) \ox X  \text{ with } \mx  ~\mx' = 1 \text{ and }\mx' ~ \mx = 1 \]
\[ \mx'': X \oa M(U) \to X \ox M(U) \text{ with } \mx  ~\mx'' = 1 \text{ and }\mx'' ~ \mx = 1 \]
A mixed unitary category, $M: \U \to \C$ is {\bf symmetric} if the functor $M$, the 
unitary category $\U$, and the $\dagger$-isomix category $\C$ are symmetric.
\end{defi}

In the definition of a MUC, the requirement of a transformation $\mx'$ which is inverse to $\mx$ ensures that the functor $M: \U \to \C$ factors through the $\Core(\C)$. We discuss examples of MUCs in the next section.  First we show that the unitary construction on a $\dagger$-isomix category produces a mixed unitary category (MUC) which is couniversal.

Mix unitary categories organize themselves into a 2-category ${\sf MUC}$ (although we shall not discuss the 2-cell structure):
\begin{description}
\item[0-cells]  Are mix unitary categories $M: \U \to \X$;
\item[1-cells]  Are MUC morphisms: these are squares of $\dagger$-isomix functors $(F',F,\gamma): M \to N$ commuting up to a $\dagger$-linear natural isomorphism $\gamma$:
 \[ \xymatrix{ \U \ar[d]_{F'} \ar@{}[drr]|{\Downarrow~\gamma} \ar[rr]^M & & \X \ar[d]^F \\ \V \ar[rr]_{N} & & \Y} \]
 The functor $F': \U \to \V$ is between unitary categories and we demand of it that it preserves unitary structure in the 
 sense of Definition \ref{preserving-unitary-structure}, thus, whenever ${\varphi_A}$ is the unitary structure  then $F'(\varphi_A) \rho^F$ is unitary structure.
 \item[2-cells] These are ``pillows''  of natural transformations. $(\beta, \beta') : (F, F', \gamma_F) \Rightarrow (G, G', \gamma_G)$ 
is a 2-cell if and only if it satisfies the following equality:
\[ \xymatrix{ \U \ar@/_1pc/[dd]_{G'} \ar@/^1pc/[dd]^{F'}   \ar@{}[ddrr]|
{\Downarrow~\gamma_F} \ar[rr]^M & & \X \ar@/^1pc/[dd]^F \\ 
{\xLeftarrow{\beta'}} & & \\ 
\V \ar[rr]_{N} & & \Y} ~ \xymatrix{ \\ = \\ } ~ \xymatrix{ \U \ar@/_1pc/[dd]_{G'} \ar@{}[ddrr]|
{\Downarrow~\gamma_G} \ar[rr]^M & & \X \ar@/_1pc/[dd]_G \ar@/^1pc/[dd]^F \\ 
 & & {\xLeftarrow{\beta}} \\ 
\V \ar[rr]_{N} & & \Y} \]
 \end{description}
 
 We remark that we have observed that any MUC can be ``simplified'' to a dagger monoidal category with a strong $\dagger$-mix Frobenius functor into a $\dagger$--isomix category: this is achieved by precomposing with ${\sf Mx}_\downarrow$.   This may seem a worthwhile simplification, but it should be recognized that it simply transfers complexity from the unitary category itself onto the preservator which must now ``create'' unitary structure:
 \[ \xymatrix{\U  \ar[d]_{{\sf Mx}^{*}_\downarrow} \ar[rr]^M & & \C \ar@{=}[d] \\
                     \U_\downarrow \ar[rr]_{{\sf Mx}_\downarrow;M} & & \C} \]
 Here $\U_\downarrow = (\U,\oa,\oa)$ is viewed as a dagger monoidal category and ${\sf Mx}_\downarrow^{*}$ is the inverse of ${\sf Mx}_\downarrow$.  
 The point is that the preservator of the lower arrow ${\sf Mx}_\downarrow;M$ is non-trivial as it must encode the unitary structure of $\U$.
 
Our objective is now to show that the unitary construction of the previous section gives rise to a right adjoint to the underlying 2-functor $U: {\sf MUC} \to {\sf MCC}$  where 
the 2-category ${\sf MCC}$ is defined as:

\begin{description}
\item[0-cells] Its objects are  {\bf mixed $\dagger$-compact categories} (MCC), that is strong $\dagger$-Frobenius functors $V: \C \to \Y$ where $\C$ is a compact $\dagger$-LDC, $\Y$ is a $\dagger$-isomix category, 
and $V$ factors through the core of $\Y$ i.e, for all $\forall$ objects $C \in \C$, $Y \in \Y$,   $\exists$ $\mx': V(C) \oa Y \to V(C) \ox Y$ such that $\mx~\mx' = 1$ and $\mx' ~ \mx = 1$.
\item[1-cells]   The 1-cells are squares of mix Frobenius functors which commute up to a linear natural isomorphism;
\item[2-cells]    Are pillows of natural transformations (which we shall ignore).
\end{description}

An example of a mix $\dagger$-compact category is, of course, the inclusion of the core into a $\dagger$-isomix category $C:{\sf Core}(\X) \hookrightarrow \X$;

\begin{prop}
$U: {\sf MUC} \to {\sf MCC}$ has a right adjoint ${\sf Unitary}: {\sf MCC} \to {\sf MUC}; (\C \to^V \X) \mapsto ({\sf Unitary}(\C) \to^{U;V} \X)$.
\end{prop}

\begin{proof}

The couniversal diagram is as follows:
\[ \xymatrix{ {\U \to^M \X} \ar[rr]^{(F,G,\gamma)} \ar[d]_{(F^\flat,G,\gamma^\flat)} && {\C \to^V \Y} \\
                    {{\sf Unitary}(\C) \to_{U;V} \Y} \ar[urr]_{\epsilon} } \]
where $\epsilon$ is the square on the left and $(F^\flat,G,\gamma^\flat)$ is the square on the right:
 \[ \xymatrix{{\sf Unitary}(\C) \ar[d]_U \ar[r]^{~~~~U} & \C \ar[r]^V & \Y \ar@{=}[d] \\
                         \C \ar[rr]_V && \Y} 
     ~~~~~~~~
     \xymatrix{\U \ar[dr]_F \ar[d]^{F^\flat} \ar@{}[drr]|{~~~~~~~~~\uparrow~\gamma} \ar[rr]^M & & \X \ar[d]^{G} \\
                      {\sf Unitary}(\C) \ar[r]_{~~~~U} & \C \ar[r]_V & \Y} \]
                      
It follows from Proposition \ref{Prop: Couniversal} that the couniversal diagram commute.
\end{proof}

This proposition means that in building a non-trivial MUC from a mixed $\dagger$-compact category it suffices to show that the compact $\dagger$-LDC contains non-trivial pre-unitary objects. 

%%%%%%%%%%%%%%%%%%%%%%%%%%%%%%%%%%%%%%%%%%%%%%%%%%%%%%%

\subsection{Examples of mixed unitary categories}
\label{Sec: MUC examples}

%%%%%%%%%%%%%%%%%%%%%%%%%%%%%%%%%%%%%%%%%%%%%%%%%%%%%%%

In this section we present a number of examples of MUCs.  We have already noted that dagger monoidal categories are automatically unitary categories in which the unitary 
structure is given by identity maps.  The identity functors then give a rather trivial MUC.  More excitingly one can take the bicompletion of the $\dagger$-monoidal category: this is a non-trivial $\dagger$-isomix $*$-autonomous category 
extension of the original $\dagger$-monoidal category and provides, thus, an interesting example of how MUCs arise.  

Our purpose in this section is to exhibit some non-trivial manifestations of the various structural components of a MUC.  To this end we discuss in some detail three basic examples.

\subsubsection{Finite dimensional framed vector spaces}

In this section we show that the example ${\sf FFVec}_K$, the category of finite dimensional framed vector spaces defined in 
Section~\ref{subsection:fdfv} is a unitary category (hence is immediately a mixed unitary category). The unitary structure map 
on each object $(V, {\cal V})$ is defined as follows:
\[ \varphi_{(V,{\cal  V})}: (V,{\cal  V}) \to (V,{\cal  V})^\dag; v_i \mapsto \widetilde{v_i} \]
and it remains to check the coherences {\bf [U.3]}--{\bf [U.6]}.  First note that {\bf [U.4]} holds immediately by the observation above that 
$\iota(v_i) = \widetilde{\widetilde{v_i}}$.  For {\bf [U.3]} we require that $\varphi_{A^\dag}(\widetilde{a_i}) = (\varphi_A^{-1})^\dag (\widetilde{a_i})$ 
the result is a higher-order term, so we may check that the evaluations are the same on basis elements:
\begin{eqnarray*}
	(\varphi_{A^\dag}(\widetilde{a_i}) ) (\widetilde{a_j}) & = & \widetilde{\widetilde{a_i}}(\widetilde{a_j}) = \partial_{i,j} \\
	((\varphi_{A}^{-1})^\dag(\widetilde{a_i}))(\widetilde{a_j}) & = & \widetilde{a_i} (\varphi_{A}^{-1}(\widetilde{a_j})) = \widetilde{a_i}(a_j) =  \partial_{i,j}
\end{eqnarray*}
Note that {\bf [U.5]}(a) and {\bf [U.5]}(b), in this example, require $\lambda_\top = \varphi_\top$ which can easily be verified as each reduces to conjugation.
{\bf [U.6]}(a) and {\bf [U.6]}(b), in this example, are the same requirement which is verified by:
\[ \lambda_\ox(\varphi_A \ox \varphi_B(a_i \ox b_j) ) = \lambda_\ox (\widetilde{a_i} \ox \widetilde{b_j}) = \widetilde{a_i \ox b_j} = \varphi_{A \ox B} (a_i \ox b_j) \]

This gives:

\begin{prop}
	${\sf FFVec}_K$ with the unitary structure above is a MUC.
\end{prop}

This raises the question of what precisely the unitary maps of this example are.  To elucidate this we note that  a functor can easily be constructed
$U:{\sf FFVec}_K \to {\sf Mat}(K)$ where, for each object in ${\sf FFVec}_K$ we choose a total order on the elements of the basis and note that 
any map is then given by a matrix acting on the bases: thus a matrix in ${\sf Mat}(K)$ with the appropriate dimensions.  We now observe:

%typed wrong. A!=B.  Do you mean that they have the same underlying vector space associated to them.  <-- No as they have the same dimension in Mat(\X) they are the same objects!
\begin{lem} 
	An isomorphism $u: (A,{\cal A}) \to (B,{\cal B})$ in ${\sf FFVec}_K$ is unitary if and only if $U(f)$ is unitary in ${\sf Mat}(\X)$.
\end{lem}

\proof While $U$ does not preserve $(\_)^\dag$ on the nose it does so up to the natural equivalence determined by $U(\varphi_A)$ which being a basis 
permutation is a unitary equivalence.   Thus, it is not hard to see that the following diagram commutes:
\[ \xymatrix{ U(B,{\cal B})  \ar[d]_{U(f)^\dag} \ar[rr]^{U(\varphi_B)} & &U((B,{\cal B})^\dag) \ar[d]^{U(f^\dag)} \\
	U(A,{\cal A}) \ar[rr]_{U(\varphi_A)} & & U((A,{\cal A})^\dag) } \]
Recall that in the category of matrices, the dagger is stationary on objects so $U(B,{\cal B}) = U(B,{\cal B})^\dag$.  

Now suppose $u$ is unitary in ${\sf FFVec}_K$  then $u^{-1} = \varphi_B u^\dagger \varphi_A^{-1}$ so that 
\[ U(u)^{-1} = U(u^{-1}) = U(\varphi_B u^\dagger \varphi_A^{-1}) = U(\varphi_B) U(u^\dagger) U(\varphi_A^{-1}) = U(u)^\dagger \]
so that its underlying map is unitary. Conversely, if $U(u)$ is unitary then 
\[ U(u^{-1}) = U(u)^{-1} = U(u)^\dag = U(\varphi_B u^\dagger \varphi_A^{-1}) \]
which immediately implies, as $U$ is faithful, that $u$ is unitary in ${\sf FFVec}_K$.
\endproof

One might reasonably regard this as a rather roundabout way to describe the standard notion of a unitary map.  However, two things of importance have been 
achieved.  First an example of a unitary category with a non-stationary dagger and, thus, a non-identity unitary structure, has been exhibited.  Second we have 
shown how the standard unitary structure may be re-expressed in this formalism using non-stationary constructs.

%%%%%%%%%%%%%%%%%%%%%%%%%%%%%%%%%%%%%%%%%%%%%%%%%%%%%%%%%%%%%%%%%%%%
\subsubsection{Finiteness matrices}

In Section~\ref{Sec: Finiteness matrices}, we described the category of finiteness matrices, ${\sf FMat}(\C)$. The core of ${\sf FMat}(\C)$ is the subcategory determined by objects whose webs are finite sets, that is the objects are $X = (|X|, P(X))$ where $|X|$ is a finite set. 
Clearly, $\Core({\sf FMat}(\C))$ is then equivalent to the category of finite dimensional matrices, ${\sf Mat}(\C)$.  This is a well-known $\dagger$-compact closed category, which is a unitary category with unitary structure given by identity maps (as $(\_)^\dagger$ is stationary on objects).  

The inclusion ${\cal I}: {\sf Mat}(\C) \to {\sf FMat}(\C)$ provides an important example of a MUC.   

%%%%%%%%%%%%%%%%%%%%%%%%%%%%%%%%%%%%%%%%%%%%%%%%%%%%%%%%%%%%%%%%%%%%%

\subsubsection{The embedding of finite-dimensional Hilbert Spaces into Chu spaces}

In Section~\ref{Section: Chu}, we showed that the Chu construction applied to a symmetric conjugative closed monoidal category, $\X$, with pullbacks gives a $\dagger$-isomix category. Recall that the dagger in the resulting category of Chu spaces is given by composing the conjugation with the dualizing functor.  In this section, we start by discussing, in  general, the construction of a mixed unitary category from a Chu category ${\sf Chus}_\X(I)$.  A crucial step in this is to identify objects which are in the core of this category.

Recall that a  symmetric monoidal closed category, $\X$, is (degenerately) a compact linearly distributive category and, thus, there may be objects which have linear adjoints: these are called {\bf nuclear} objects \cite{higgs&rowe}.  Explicitly a nuclear object $A$ in a symmetric monoidal closed  category is an object with $A \multimap B \cong A^{*} \ox B$, where 
$A^* := A \multimap I$.   The nuclear objects form a compact closed subcategory of $\X$ which is conjugative when $\X$ is conjugative.  In ${\sf Vec}_\mathbb{C}$ the nucleus consist precisely of the finite dimensional vector spaces.  If $(\eta, \epsilon): A \dashv\!\!\!\dashv B$ is witness that $A$ (and $B$) are nuclear in $\X$ then the object $(A,B,\epsilon,c_\otimes\epsilon)$ is in the core of ${\sf Chus}_\X(I)$ because in the second component of the tensor product with any other object $(X,Y,\nu,c_\ox \nu)$ one has the degenerate pullback:
\[ \xymatrix{
& Y \ox B \ar[rd]^{\simeq} \ar[ld] & \\ 
X \multimap B \ar[rd]^{\simeq}  & & Y \ox A^{*} \ar[ld]  \\
& X \multimap A^* \to^{\simeq} (X \multimap I) \ox A &
} \]
where we use the isomorphism $B \to^{\simeq} A^{*}$.

In this manner the nuclear objects of ${\sf Nuclear}(\X)$, which form a compact closed category with a dagger, may be embedded into the core of ${\sf Chus}_\X(I)$.  To obtain a unitary category it suffices then to use the unitary construction for which, to obtain a non-trivial result, we need to show that there are non-trivial examples of pre-unitary objects. To achieve this we consider an object $H$ for which $(e, n): H \dashv\!\!\!\dashv \overline{H}$ and such that $e$ satisfies:
\[ \xymatrix{\overline{\overline{H}} \otimes \overline{H} \ar[d]_{\varepsilon \otimes 1} \ar[r]^{\chi} & \overline{\overline{H} \otimes \overline{H}} \ar[dd]^{\overline{e}}\\
                   H \otimes \overline{H} \ar[d]_{e} \\
                   I \ar[r]_{\chi^{\!\!\!\circ}} & \overline{I} } \]
               
 For such an object we note:
 \begin{align*}
{ \overline{(H,\overline{H},e,c_\otimes e)}}^* & =  (\overline{H},\overline{\overline{H}},\chi\overline{c_\otimes e} (\chi^{\!\!\!\circ})^{-1},\chi\overline{e} (\chi^{\!\!\!\circ})^{-1})^* \\
 & =  (\overline{\overline{H}},\overline{H},\chi\overline{e} (\chi^{\!\!\!\circ})^{-1},\chi\overline{c_\otimes e} (\chi^{\!\!\!\circ})^{-1})
 \end{align*}
 
This makes
 $$(\varepsilon^{-1},1) : (H,\overline{H},e,c_\otimes e) \to (\overline{\overline{H}},\overline{H},\chi\overline{e} (\chi^{\!\!\!\circ})^{-1},\chi\overline{c_\otimes e} (\chi^{\!\!\!\circ})^{-1})$$
a preunitay map.  Note that it is a Chu map by the commuting diagram above and as $\overline{\varepsilon} =\varepsilon$ we have 
$$(\varepsilon^{-1},1) (1,\overline{\varepsilon}) = (\varepsilon^{-1},1) (1,\varepsilon) = (\varepsilon^{-1},\varepsilon)$$
where $(\varepsilon^{-1},\varepsilon)$ is the involutor.

In ${\sf Vec}_{\mathbb{C}}$ a map $e: H \otimes \overline{H} \to \mathbb{C}$ is a ``sesquilinear form'' and the diagram above asserts that it is in addition a symmetric form.  Any  Hilbert space with its inner product, thus, satisfies the above conditions.   Thus,  it is clear that the embedding of the category of finite dimensional Hilbert Spaces into Chu spaces, ${\sf FHilb} \hookrightarrow {\sf Chus}_{{\sf Vec}_\C} (\C)$ is a mixed unitary category.  The embedding is in fact a full and faithful embedding which extends to {\em all\/} Hilbert spaces (although only the finite dimensional ones land in the core).  

Explicitly the embedding is defined as follows: suppose $H$ is a (finite dimensional)) Hilbert Space, then the corresponding Chu Space is given by $(H, \overline{H}, \langle - | - \rangle_H)$, where $\langle - | - \rangle_H: H \ox \overline{H} \to \C$ is the inner product. For any linear map $H \to^{f} K$ between Hilbert Spaces, the corresponding Chu map is given by $(f, f^\dagger): (H, \overline{H}, \langle - | - \rangle_H) \to (K, \overline{K}, \langle - | - \rangle_K)$, where $f^\dagger$ is the Hermitian adjoint of $f$ so, $\langle f(a) | b \rangle = \langle a | f^\dagger(b) \rangle$.

Furthermore, observe that $(H, \overline{H}, \langle - | - \rangle_H)^\dagger :=  \overline{(H, \overline{H}, \langle - | - \rangle_H)^*} =  (H, \overline{H}, \langle - | - \rangle_H)$. Hence, this embedding preserves the (stationary) dagger for all Hilbert spaces.  However, the par of two infinite dimensional Hillbert spaces in this Chu category is not a Hilbert space so that the duality cannot be seen within the category of Hilbert spaces.

%%%%%%%%%%%%%%%%%%%%%%%%%%%%%%%%%%%%%%%%%%%%%%%%%%%%%%%%%%%%%%%%%%%%%%
\subsubsection{Constructing MUCs using the unitary construction}
One can construct a MUC from any $\dagger$-isomix category using the unitary construction: for any $\dagger$-isomix category, $\X$, ${\sf Unitary}(\Core(\X)) \to^{U} \Core(\X) \hookrightarrow \X$ is a MUC. In this manner we have already many examples of MUCs:

\begin{itemize}
	\item The inclusion $\C \hookrightarrow \D( \C, +, 0)$
	\item ${\sf Unitary}(\Core({\mbox{\bf H-Mod}}_\X)) \to^{U} \Core(\mbox{\bf H-Mod}) \hookrightarrow {\mbox{\bf H-Mod}}_\X$
	\item ${\sf Unitary}(\Core({\sf Chus_\X}(I)) \to^{U} \Core({\sf Chus_\X}(I)) \hookrightarrow {\sf Chus_\X}(I)$
\end{itemize}

\section{Conclusion}

In this paper, we have extended the theory of $\dagger$-monoidal categories and $\dagger$-compact closed categories to linearly distributive 
and *-autonomous settings to obtain the semantics of (multiplicative) $\dagger$-linear logic.  In these linear settings, the 
two different tensor products (tensor and par) must be flipped by the dagger.  Thus, one cannot have a stationary (identity on objects) 
dagger, and hence one is forced to replace the conventional dagger by a contravariant structure-preserving involution.  This has coherence  
consequences: almost two thirds of this paper is dedicated to understanding the details of these coherences.

If multiplicative  $\dagger$-linear logic is to provide a semantics for a generalized categorical quantum mechanics (CQM), then notions 
such as isometry and unitary isomorphism, which are central to CQM, should have an expression in this logic.  Here we showed that 
with additional ``unitary structure'' one can recapture classical CQM as a ``unitary core''  of multiplicative  $\dagger$-linear logic.
Furthermore, we showed how, from any $\dagger$-isomix category, it is always possible to extract  
a ``unitary core'' which is, up to equivalence, a $\dagger$-monoidal category (i.e a classical semantic setting for CQM).   

This led to the notion of a mixed unitary category (MUC) given by a $\dagger$-isomix category with a chosen unitary 
core as our proposal for an extension of CQM.   A MUC can be viewed as an extension much as a $K$-algebra 
extends a field $K$ and permits the expression of properties which are difficult to express within $K$ itself.  
In the extended setting of a MUC -- finiteness matrices with its core for example -- provides an  
extension of the classical CQM setting in which infinite dimensional types, such as those given by the exponential modalities,  
 are present.   Furthermore, in the extended setting one can bend, and yank wires without the category being compact.

The fact that a unitary category is a component of a MUC allows one to mimic the construction of completely positive 
maps, ${\sf CP}^\infty$, see \cite{CoH16}, in a way which displays 
some interesting features.  To start with the ancillary objects (which are to be traced out) must now, necessarily, be chosen from the 
unitary core and these
 it can be supposed are an essentially small class even though the overall category may be large.  This keeps the number completely 
 positive maps between any two objects small.  The resulting category is under reasonable assumptions a MUC (see \cite{CS19}) 
 which has an appropriate analogue of environment structure \cite{Coecke10}.  Furthermore, in the presence of duals the 
 whole construction is functorial.  
 
An important observation of CQM is that an orthogonal basis, for a Hilbert space, correspond to a special commutative Frobenius algebra 
\cite{CPV12}.  This allows one to replace the notion of a basis by algebraic structure.  A significant consequence of this has been the 
algebraic expression of ``uncertainty" using complementary Frobenius algebras \cite{Coecke&Duncan}, which, in turn, led to the 
formulation of the ZX-calculus.   As was mentioned in the introduction, linear settings allow the expression of structures which parallel Frobenius 
algebras and this can be exploited to allow an expression complementarity in MUCs and, furthermore, to link this to the exponential 
modalities in $\dagger$-linear logic \cite{CoS20}.  

 An instructive source of examples of MUCs, which was mentioned in the introduction and left for future work, uses Joyal's bicompletion procedure \cite{Joy95}:  
 here, starting with a $\dagger$-monoidal category, or a compact $\dagger$-isomix category, $\C$, one can form a MUC $\iota: \C \to \Lambda(\C)$ by 
 simply bicompleting.  Furthermore, the bicompletion is a (non-compact) $\dagger$-isomix category which, when the starting point, $\C$, is $\dagger$-compact 
 closed, is a $\dagger$-isomix $*$-autonomous category.

\bibliographystyle{alpha}
\bibliography{dagger-frob}
\end{document}